\documentclass[leqno, 11pt]{amsart}
\usepackage{amsmath, amsthm, amssymb, latexsym}
  \newlength{\auxwidth}
  \newlength{\auxheight}
\newcounter{theorem}
\newcounter{proposition}
\newcounter{lemma}
\newcounter{corollary}
\newcounter{conjecture}

\newtheorem{definition}{Definition}[section]
\newtheorem{theorem}[definition]{Theorem}
\newtheorem{proposition}[definition]{Proposition}
\newtheorem{lemma}[definition]{Lemma}

\theoremstyle{remark}
\newtheorem{example}[definition]{Example}
\newtheorem{remark}[definition]{Remark}

  \newtheorem*{acknowledgements}{Acknowledgements}

\numberwithin{equation}{section}

\newcommand{\op}[1]{\operatorname{#1}}

\newcommand{\acou}[2]{\ensuremath{\langle #1 , #2 \rangle}}

\newcommand{\brak}[1]{\ensuremath{\langle\! #1\!\rangle}}

\newcommand{\Tr}{\ensuremath{\op{Tr}}}
\newcommand{\tr}{\op{tr}}
\newcommand{\Tra}{\ensuremath{\op{Trace}}}

\newcommand{\Hol}{\op{Hol}}

\newcommand{\C}{\ensuremath{\mathbb{C}}} 
\newcommand{\bH}{\ensuremath{\mathbb{H}}} 
\newcommand{\N}{\ensuremath{\mathbb{N}}} 
\newcommand{\R}{\ensuremath{\mathbb{R}}} 
\newcommand{\Z}{\ensuremath{\mathbb{Z}}}

\newcommand{\Rd}{\ensuremath{\R^{d+1}}}

\newcommand{\Rdo}{\R^{d+1}\!\setminus\! 0}

\newcommand{\URd}{U\times\R^{d+1}}

\newcommand{\URdo}{U\times(\R^{d+1}\!\setminus\! 0)}

\newcommand{\Ca}[1]{\ensuremath{\mathcal{#1}}}

\newcommand{\cD}{\ensuremath{\mathcal{D}}}
\newcommand{\cE}{\Ca{E}}
\newcommand{\cF}{\Ca{F}}
\newcommand{\cG}{\ensuremath{\mathcal{G}}}
\newcommand{\cH}{\ensuremath{\mathcal{H}}}

\newcommand{\cK}{\ensuremath{\mathcal{K}}}
\newcommand{\cL}{\ensuremath{\mathcal{L}}}

\newcommand{\cS}{\ensuremath{\mathcal{S}}}

\newcommand{\cU}{\ensuremath{\mathcal{U}}}

\newcommand{\fg}{\ensuremath{\mathfrak{g}}}
\newcommand{\fh}{\ensuremath{\mathfrak{h}}}

\newcommand{\psivdo}{$\Psi_{H}$DO}
\newcommand{\psivdos}{$\Psi_{H}$DO's}

\newcommand{\pvdo}{\ensuremath{\Psi_{H}}} 
\newcommand{\pvhdo}{\ensuremath{\Psi_{H,\op{v}}}}

\newcommand{\psido}{$\Psi$DO} 
\newcommand{\psidos}{$\Psi$DO's} 
\newcommand{\psinf}{\ensuremath{\Psi^{-\infty}}}

\newcommand{\Svb}{S_{\scriptscriptstyle{| |}}}
\newcommand{\SvbU}[1]{S_{\scriptscriptstyle{| |}}^{#1}(U\times\Rd)}

\newcommand{\ah}{\text{ah}}
\newcommand{\reg}{{\text{reg}}}
\newcommand{\ord}{{\op{ord}}}

\newcommand{\xiy}{{\xi\rightarrow y}}
\newcommand{\xitauyt}{{(\xi,\tau)\rightarrow (y,t)}}
\newcommand{\yxi}{{y\rightarrow\xi}}

\newcommand{\supp}{\op{supp}}

\newcommand{\rk}{\op{rk}}
\newcommand{\im}{\op{im}}
 
\newcommand{\dom}{\op{dom}}

\newcommand{\End}{\ensuremath{\op{End}}}

\newcommand{\hotimes}{\hat\otimes}

\renewcommand{\Box}{\square}
\newcommand{\subsubset}{\subset\!\subset}

\newcommand{\Sp}{\op{Sp}}

\newcommand{\dbarb}{\bar\partial_{b}}

\begin{document}
\title{Functional calculus and spectral asymptotics for hypoelliptic operators on Heisenberg Manifolds. I.} 

\author{Rapha\"el Ponge}

\address{Department of Mathematics, Ohio State University, Columbus, USA.}
\email{ponge@math.ohio-state.edu}
 \keywords{Heisenberg calculus, complex powers, heat equation, spectral asymptotics, analysis on CR and contact manifolds, hypoelliptic operators}
 \subjclass[2000]{Primary 58J40, 58J50; Secondary 58J35, 32V10, 35H10, 53D10}


 \begin{abstract}
This paper is part of a series papers devoted to geometric and spectral theoretic applications of the hypoelliptic calculus on Heisenberg manifolds.
More specifically, in this paper we make use of the Heisenberg calculus of Beals-Greiner and Taylor to analyze the spectral theory of hypoelliptic operators on 
 Heisenberg manifolds. The main results of this paper include: (i) Obtaining complex powers of hypoelliptic operators as holomorphic families of \psivdos, 
 which can be used to define a scale of weighted Sobolev spaces interpolating the weighted Sobolev spaces of Folland-Stein and providing us with sharp 
 regularity   estimates for hypoelliptic operators on Heisenberg manifolds; (ii) Criterions on the principal symbol of $P$ to invert the heat operator 
 $P+\partial_{t}$ and to derive the small time heat kernel 
 asymptotics for $P$; (iii) Weyl asymptotics for hypoelliptic operators which can be reformulated geometrically for the main geometric operators on CR and contact 
 manifolds, that is, the Kohn Laplacian, the horizontal sublaplacian and its conformal powers, as well as the contact Laplacian. For dealing we 
 cannot make use of the standard approach of Seeley, so we rely on a new approach based on the pseudodifferential approach representation of the heat 
 kernel. This is especially suitable for dealing with positive hypoelliptic operators. We will deal with more general operator in a forthcoming paper using 
 another new approach. The results of this paper will be used in another forthcoming paper dealing with an analogue for the Heisenberg calculus of 
 the noncommutative geometry which, in particular, will allow us to make use in the Heisenberg setting of Connes' noncommutative 
 geometry, including the operator theoretic framework for the local index formula of Connes-Moscovici.
\end{abstract}

\maketitle 

\section{Introduction} 
This paper is part of a series papers devoted to geometric and spectral theoretic applications of the hypoelliptic calculus on Heisenberg manifolds. 
Recall that a Heisenberg manifold $(M,H)$ consists of a manifold $M$ together with a distinguished hyperplane bundle $H\subset TM$. This definition 
 covers many examples: Heisenberg group and its quotients by cocompact lattices, (codimension 1) foliations, CR and contact manifolds and the 
 confolations of~\cite{ET:C}. 
 
 The name Heisenberg manifold comes from the fact that  
 the relevant tangent structure for a Heisenberg manifold $(M,H)$ is rather that of a bundle $GM$ of two-step nilpotent Lie 
 groups, whose Lie group structure is encoded by an intrinsic Levi form $\cL:H\times H\rightarrow TM/H$  
 (see~\cite{BG:CHM}, \cite{Be:TSSRG}, \cite{EMM:HAITH}, \cite{FS:EDdbarbCAHG}, \cite{Gr:CCSSW}, \cite{Po:Pacific1}, \cite{Ro:INA}).  
 
 In this context the most natural operators include H{\"o}rmander's sum of squares, the Kohn Laplacian, the horizontal sublaplacian and its conformal 
 powers,  as well as the contact Laplacian (see Section~\ref{sec:Operators} and the references therein for an overview of these operators). 
 Although these operators may be 
 hypoelliptic, they are definitely not elliptic. 
 Therefore, the classical pseudodifferential calculus cannot be used efficiently to study these operators. 
 
The relevant substitute to the standard pseudodifferential calculus is provided by the Heisenberg calculus of 
 Beals-Greiner~\cite{BG:CHM} and Taylor~\cite{Ta:NCMA} (see also \cite{BdM:HODCRPDO},  \cite{CGGP:POGD}, \cite{Dy:POHG}, \cite{Dy:APOHSC}, 
 \cite{EMM:HAITH}, \cite{FS:EDdbarbCAHG}, \cite{RS:HDONG}).  
The idea in the Heisenberg calculus, which goes back to Elias Stein, is the following. Since the relevant notion of tangent structure for a Heisenberg 
manifold $(M,H)$ is that of a bundle $GM$ of 2-step nilpotent graded Lie groups, it stands for reason to construct a pseudodifferential calculus which at every point 
$x \in M$ is well modelled by the calculus of 
 convolution operators on the nilpotent tangent group $G_{x}M$. 

The result is a class of pseudodifferential operators, 
the \psivdos, which are locally \psidos\ of type 
$(\frac{1}{2},\frac{1}{2})$, but unlike the latter possess a full symbolic calculus and makes sense on a general Heisenberg manifold. In particular, a \psivdo\ 
admits a parametrix in the Heisenberg calculus if, and only if, its principal symbol is invertible, and then the \psivdo\ is hypoelliptic with a loss of 
derivatives controlled by its order (see Section~\ref{sec:PsiHDO} for a more detailed overview of the Heisenberg calculus). 

 In~\cite{Po:Pacific1} a tangent groupoid has been associated to any Heisenberg manifold $(M,H)$ as the differentiable groupoid encoding a smooth 
 deformation of $M\times M$ (see also~\cite{VE:PhD}). This is the analogue for Heisenberg manifolds of Connes' tangent groupoid of a manifold which 
 plays a pivotal role in his proof in~\cite{Co:NCG} of the index theorem of Atiyah-Singer~\cite{AS:IEO1}. The results of~\cite{Po:Pacific1} have been used 
 subsequently in~\cite{Po:BSM1} to give an intrinsic definition of the principal symbol for the Heisenberg calculus, which was not done in~\cite{BG:CHM} or 
 in~\cite{Ta:NCMA}. 
  
 In this paper, and its sequel~\cite{Po:JFA3}, we make use of the Heisenberg calculus to analyze the spectral theory of hypoelliptic operators on 
 Heisenberg manifolds. The main results of this paper include:\smallskip
 
 - Obtaining complex powers of hypoelliptic operators as holomorphic families of \psivdos, which can be used to define a scale of $W_{H}^{s}$, $s\in 
 \R$, of weighted Sobolev spaces interpolating the weighted Sobolev spaces of Folland-Stein and providing us with sharp estimates for 
 \psivdos.\smallskip 
 
 - Criteria on the principal symbol of $P$ to invert the heat operator $P+\partial_{t}$ and to derive the small time heat kernel 
 asymptotics for $P$.\smallskip 
 
 - Weyl asymptotics for hypoelliptic operators which can be reformulated geometrically for the main geometric operators on CR and contact 
 manifolds, that is, the Kohn Laplacian, the horizontal sublaplacian and its conformal powers, as well as the contact Laplacian.\smallskip  
 
 These results will be important ingredients in~\cite{Po:GAFA1} to construct an analogue for the Heisenberg calculus of the noncommutative residue 
trace of  Wodzicki~(\cite{Wo:LISA}, \cite{Wo:NCRF}) and Guillemin~\cite{Gu:NPWF} and to study the zeta and eta functions of hypoelliptic operators. 
In turn this has several geometric consequences. In particular, this allows us to make use of the framework of Connes' noncommutative geometry, 
including the operator theoretic framework for the local index formula of Connes-Moscovici~\cite{CM:LIFNCG}. 

Let us also mention that the lack of microlocality of the Heisenberg calculus does not allow us to carry out in the Heisenberg setting the standard approach 
of Seeley~\cite{Se:CPEO} (see also~\cite{Sh:POST}, \cite{Gr:FCPDBP}) to complex powers of elliptic operators. 
Instead, we rely on a new approach, based on the pseudodifferential representation of the heat in~\cite{BGS:HECRM}, 
which is quite suitable for dealing with positive differential operators.  
A similar approach has also been used  by Mathai-Melrose-Singer~\cite{MMS:FAI} and Melrose~\cite{Me:SPLLB} in the context of 
projective pseudodifferential operators on Azamaya bundles. 
In~\cite{Po:JFA3} we will use another approach to deal with \psivdos\ 
which are not positive differential operators.

\subsection{Holomorphic families of \psivdos} 
Prior to dealing with complex powers of hypoelliptic operators we define holomorphic families of \psivdos\ and check their main properties. 
 
In a local Heisenberg chart $U\subset \Rd$ the definition of a holomorphic family of \psivdos\ parametrized by an open $\Omega\subset \C$ 
is similar to that of the standard definition of a holomorphic family 
of \psidos\ in~\cite[7.14]{Wo:LISA} and \cite[p.~189]{Gu:GLD} (see also~\cite{KV:GDEO}). In particular, we allow the order of the family of \psidos\ to vary 
analytically. Several properties of \psivdos\ on $U$ extend \emph{mutatis standis} to the setting of holomorphic families of \psivdos. 
In particular, concerning the product of \psivdos\ we have: 

\begin{proposition}
    For $j=1,2$ let $(P_{j,z})_{z \in \Omega}\subset \pvdo^{*}(M,\cE)$ be a holomorphic family of \psivdos\ and supposed that $(P_{1,z})_{z \in 
    \Omega}$ or  $(P_{2,z})_{z \in \Omega}$ is uniformly properly supported with respect to $z$. Then the family of products $(P_{1,z}P_{2,z})_{z \in 
    \Omega}$ is a holomorphic family of \psivdos. 
\end{proposition}

There is, however, a difficulty when trying to extend the definition to arbitrary Heisenberg manifolds. More precisely, 
the proof of the invariance of the Heisenberg calculus by Heisenberg diffeomorphisms relies on a 
characterization of the distribution kernels of \psivdos\ by means of a suitable class of kernels 
$\cK^{*}(\URd)= \sqcup_{m\in \C}\cK^{m}(\URd)\subset \cD'(\URd)$. Each distribution $K \in \cK^{m}(\URd)$ 
admits an asymptotic expansion in the sense of distributions, 
\begin{equation}
    K\sim \sum_{j\geq 0} K_{m+j}, \qquad K_{l}\in \cK_{l}(\URd),
\end{equation}
where  $\cK_{l}(\URd)$ consists of kernels which are smooth for $y\neq 0$ and homogeneous of degree $l$ if $l \not \in\Z$ 
and homogeneous of degree $l$ up to logarithmic terms otherwise. In particular, the definition of $\cK_{l}(\URd)$ depends upon whether $l$ is an 
integer or in not, which causes trouble for defining holomorphic families of kernels in $\cK^{*}(\URd)$ whose order may take integer values. 

We resolve this issue by giving an alternative description of the class $\cK^{*}(\URd)$ in terms of what we call \emph{almost homogeneous} kernels, 
as they are homogenous modulo smooth terms and, under the Fourier 
transform, they correspond to the almost homogeneous symbols considered in~\cite{BG:CHM}. 

Since the definition of an almost homogeneous symbol of degree $l$ does not depend on whether $l$ is an integer or not, there is no trouble anymore 
to define holomorphic families of almost homogeneous kernels. Therefore, we can make use of the characterization of $\cK^{*}(\URd)$ in terms of almost 
homogenous kernels to define holomorphic families with values in $\cK^{*}(\URd)$ (see~Definition~\ref{def:HolPHDO.kernels-families} for the precise 
definition). 

We show that the distribution of kernel holomorphic families of \psivdos\ can be characterized in terms of holomorphic families with values in 
$\cK^{*}(\URd)$. As a consequence we can extend to the setting of holomorphic families of \psivdos\ the arguments in the proof 
in~\cite{BG:CHM} and~\cite{Po:BSM1}  of the invariance by Heisenberg diffeomorphisms of the Heisenberg calculus, so that we get: 

\begin{proposition}
    Let $\phi:U\rightarrow \tilde{U}$ be a change of Heisenberg chart and let $(\tilde{P}_{z\in \Omega})_{z \in \Omega}$ be a holomorphic family of 
    \psivdos\ on $\tilde{U}$. Then the family $(P_{z})_{z\in\Omega}:=(\phi^{*}P_{z})_{z\in\Omega}$ is a holomorphic family of \psivdos\ on $U$.
\end{proposition}
 
This allows us to define holomorphic families of \psivdos\ on a general Heisenberg manifold $(M,H)$ and acting on the sections of a vector 
bundle on $\cE$. In this setting the main properties of holomorphic families of \psivdos\ are summarized below. 

\begin{proposition}
     Let $(P_{z})_{z \in \Omega}\subset \pvdo^{*}(M,\cE)$ be a holomorphic family of \psivdos. Then:\smallskip
     
     1) The family of principal symbols $(\sigma_{*}(P_{z}))_{z \in \Omega}$ belongs to  $\Hol(\Omega,C^{\infty}(\fg^{*}M\setminus 0, 
     \End \cE))$.\smallskip
     
     2) The family of transpose operators $(P_{z}^{t})_{z \in \Omega}\subset \pvdo^{*}(M,\cE^{*})$ is a holomorphic family of \psivdos.\smallskip 
     
     3) The family of adjoints $(P_{z}^{*})_{z \in \Omega}$ is an anti-holomorphic family of \psivdos, that is, $(P_{\bar 
     z}^{*})_{z \in \Omega}$ is a holomorphic family of \psivdos.
\end{proposition}

Finally, let us mention that the almost homogeneous approach to the Heisenberg calculus will have further applications in~\cite{Po:JFA3} 
for constructing a class of \psivdos\ with parameter containing the resolvents of hypoelliptic \psivdos\ (see also~\cite{Po:IJM1}).  

\subsection{Complex powers of hypoelliptic operators}
One the main aims of this paper is to realize complex powers of hypoelliptic operators as holomorphic families of \psivdos. As alluded to above we cannot carry out the 
standard approach of Seeley~\cite{Se:CPEO} in the Heisenberg setting. Therefore, we have to rely on another approach based on the pseudodifferential 
representation of the heat kernel of~\cite{BGS:HECRM}. This is quite suitable for dealing differential operators, which is enough to cover 
may examples. The general case will be dealt  with in~\cite{Po:JFA3} using a different approach using a new kind of pseudodifferential representation 
of the resolvent. 

%
%

Let $P:C^{\infty}(M,\cE)\rightarrow C^{\infty}(M,\cE)$ be a hypoelliptic selfadjoint differential operator and assume that $P$ is bounded from below. 
It is well known 
the heat semigroup $e^{-tP}$, $t\geq 0$, allows us to invert the heat operator $P+\partial_{t}$. Conversely, 
constructing a suitable pseudodifferential calculus nesting parametrices for $P+\partial_{t}$ allows us to give a pseudodifferential representation for the 
heat kernel of $P$. 

In the elliptic setting the relevant pseudodifferential calculus is the Volterra calculus of Greiner~\cite{Gr:AEHE} and 
Piriou~\cite{Pi:COPDTV}. The latter is only a simple modification of the classical pseudodifferential calculus in order to take into account the parabolicity and 
the Volterra property with respect to the time variable of the heat equation. Moreover, this approach holds in fairly greater generality 
and has many applications~(\cite{BGS:HECRM}, \cite{BS:HEPNP1}, \cite{BS:HEPNP2}, \cite{Gr:AEHE}, \cite{Kr:PhD}, \cite{Kr:VFPDO}, \cite{KS:IPSPDE}, 
\cite{Me:APSIT}, \cite{Pi:COPDTV}, \cite{Po:CMP1}, \cite{Po:PAMS2}).

The Greiner's approach has been extended to the Heisenberg calculus in~\cite{BGS:HECRM}, mostly with the purpose of deriving the small time heat 
kernel asymptotics for the Kohn Laplacian on CR manifolds. In particular, a class of Volterra \psivdos\ is obtained that contains parametrices for the heat
operator $P+\partial_{t}$. As a consequence, once the principal symbol of $P+\partial_{t}$  is an invertible 
Volterra-Heisenberg symbol, the inverse of $P+\partial$ is a Volterra-\psivdo\ which, in turn, yields a pseudodifferential representation of the heat 
kernel of $P$. 

Assume now that $P$ is positive. Thanks to the spectral theorem we can define the complex power of $P^{s}$, $s\in \C$,  as an unbounded operator on $L^{2}(M,\cE)$ 
which is bounded for $\Re s\leq 0$. Moreover, for $\Re s<0$ the Mellin formula holds, 
\begin{equation}
    P^{s}=\Gamma(s)^{-1}\int_{0}^{\infty}t^{s}(1-\Pi_{0}(P))e^{-tP}\frac{dt}{t},
    \label{eq:Intro1.Mellin-formula}
\end{equation}
where $\Pi_{0}(P)$ denotes the orthogonal projection onto the kernel of $P$. 
Combining this  formula with the pseudodifferential representation of the heat kernel of $P$ allows us to prove: 

\begin{theorem}\label{thm:Intro1.complex-powers}
    Let $P:C^{\infty}(M,\cE)\rightarrow C^{\infty}(M,\cE)$ be a positive differential operator of Heisenberg  order $m$  such that 
    the principal symbol of the heat operator $P+\partial_{t}$ is an invertible Volterra-Heisenberg symbol. 
    Then the family $(P^{s})_{s\in \C}$ of the complex powers of $P$, defined by functional calculus on $L^{2}(M,\cE)$, is a holomorphic 1-parameter group 
    of \psivdos\ such that $\ord P^{s}=ms$ for any $s\in \C$.  
\end{theorem}
 
In~\cite{BGS:HECRM} it was shown that for a selfadjoint sublaplacian $\Delta:C^{\infty}(M,\cE)\rightarrow C^{\infty}(M,\cE)$
the principal symbol of $\Delta+\partial_{t}$ is invertible under a condition closely related to the Rockland condition and the invertibility of the 
principal symbol of $\Delta$ (see~\cite[5.23]{BGS:HECRM} and Section~\ref{sec.powers1}). Therefore, in the case of a sublaplacian we get: 

\begin{theorem}\label{thm:Intro1.complex-powers-sublaplacian}
    Let $\Delta:C^{\infty}(M,\cE)\rightarrow C^{\infty}(M,\cE)$ be a positive sublaplacian satisfying the condition of~\cite[5.23]{BGS:HECRM}. 
    Then  the family $(\Delta^{s})_{s\in \C}$ of the complex powers of $\Delta$ is a holomorphic 1-parameter group 
    of \psivdos\ such that $\ord \Delta^{s}=2s$ for any $s\in \C$.  
\end{theorem}

This theorem holds for the following kinds of sublaplacians:\smallskip 

(a) A selfadjoint sum of squares $\nabla_{X_{1}}^{*}\nabla_{X_{1}}+\ldots+\nabla_{X_{m}}^{*}\nabla_{X_{m}}$, 
where the vector fields $X_{1},\ldots,X_{m}$ linearly span $H$ and $\nabla$ is a connection on $\cE$, 
provided that the Levi form of $(M,H)$  is nonvanishing;\smallskip 

(b) The Kohn Laplacian on a CR manifold acting on $(p,q)$-forms under condition $Y(q)$.\smallskip 

(c) The horizontal Laplacian on a Heisenberg manifold acting on horizontal forms of degree $k$  under condition $X(k)$.\smallskip 

For more general operators we show that, provided that the Levi form of $(M,H)$ has constant rank, the Rockland condition for $P$ is enough to 
deal with invertibility of the principal symbol of $P+\partial_{t}$ (see Theorem~\ref{thm:Intro1.heat} below). This allows us to get:

\begin{theorem}\label{thm:Intro1.complex-powers-Rockland}
    Suppose that the Levi form of $(M,H)$ has constant rank and 
    let $P:C^{\infty}(M,\cE)\rightarrow C^{\infty}(M,\cE)$ be a positive differential operator of Heisenberg  order $m$ such that $P$ satisfies the 
    Rockland condition at every point. Then the family $(P^{s})_{s\in \C}$ of the complex powers of $P$ is a holomorphic 1-parameter group of \psivdos\ 
    such that $\ord P^{s}=ms$ for any $s\in \C$.  
\end{theorem}

In particular, Theorem~\ref{thm:Intro1.complex-powers-Rockland} is valid for the contact Laplacian on a contact manifold. In this context 
this allows us to fill a technical gap in the proof of by Julg-Kasparov~\cite{JK:OKTGSU} of the Baum-Connes conjecture for $SU(n,1)$ (see~\cite{Po:Crelle1}). 

Finally, let us mention that a similar approaches to complex powers of (hypo)elliptic operators have been used  by 
Mathai-Melrose-Singer~\cite{MMS:FAI} and Melrose~\cite{Me:SPLLB} in the context of projective pseudodifferential operators on Azamaya bundles. 

\subsection{Weighted Sobolev spaces}
Assuming that the Levi form~of $(M,H)$ is nowhere vanishing, we can  
construct a scale $W_{H}^{s}(M)$, $s \in \R$, of 
weighted Sobolev spaces obtained as follows. 

Let $X_{1},\ldots,X_{m}$ be vectors fields 
spanning $H$ and let $\Delta_{X}=X_{1}^{*}X_{1}+\ldots+X_{m}^{*}X_{m}$. Then 
by Theorem~\ref{thm:Intro1.complex-powers-sublaplacian} the complex powers $(1+\Delta_{X})^{s}$ gives rise to a 1-parameter group of \psivdos. We 
then define $W_{H}^{s}(M)$ as the Hilbert space consisting of distributions $u \in 
\cD'(M)$ such that $(1+\Delta_{X})^{\frac{s}{2}}u\in L^{2}(M)$ together with the Hilbertian norm, 
\begin{equation}
    \|u\|_{W_{H}^{s}}=\|(1+\Delta_{X})^{\frac{s}{2}}u\|_{L^{2}}, \qquad u \in W_{H}^{s}(M). 
\end{equation}

Neither the underlying space of $W_{H}^{s}(M)$, $s\in \R$, nor its topology depend on the choice of the vector fields $X_{1},\ldots,X_{m}$ 
(Proposition~\ref{prop:Sobolev.embeddings}). Moreover, these weighted Sobolev spaces are nicely related to the standard Sobolev spaces $L^{2}_{s}(M)$ and 
to the weighted Sobolev spaces of Folland-Stein~\cite{FS:EDdbarbCAHG}, as we have:

\begin{proposition}
1) For $k=1,2,\ldots$ the weighted Sobolev spaces $W_{H}^{k}(M)$ and $S_{k}^{2}(M)$ agree as spaces and bear the same topology.\smallskip
    
2) For $s \in \R$ the following continuous embedding hold: 
   \begin{equation}
   \begin{array}{rcl}
       L^2_{s}(M)  & \hookrightarrow W_{H}^{s}(M) \hookrightarrow L^{2}_{s/2}(M) & \qquad \text{if $s\geq 0$},\\
       L^{2}_{s/2}(M)  & \hookrightarrow W_{H}^{s}(M) \hookrightarrow L^2_{s}(M) & \qquad \text{if $s< 0$}.
   \end{array}
         \label{eq:Intro1.embeddings}
    \end{equation}
\end{proposition}

We can similarly define weighted Sobolev spaces $W_{H}^{s}(M,\cE)$ of distributional sections of $\cE$ and the embedding~(\ref{eq:Intro1.embeddings})
hold \emph{verbatim} in this setting. In addition, we can localize the Sobolev spaces $W_{H}^{s}(M,\cE)$, that is, we can speak about distributional 
sections which are $W_{H}^{s}$ near a point (see Definition~\ref{def:Sobolev.loaclization-WHs} and Lemma~\ref{lem:Sobolev.localization-WHs}).

The main interest of the weighted Sobolev spaces $W_{H}^{s}(M,\cE)$, $s \in \R$, with respect to the standard Sobolev spaces is that they allow us 
to get sharper regularity estimates for \psivdos, for we have:

\begin{proposition}\label{prop:Intro1.regularity-PsiHDOs}
    Let $P\in \pvdo^{m}(M,\cE)$, set $k=\Re m$ and let $s\in \R$.\smallskip 
    
    1) The operator $P$ extends to a continuous linear mapping from $W_{H}^{s+k}(M,\cE)$ to $W_{H}^{s}(M,\cE)$. 
    \smallskip 
    
    2) Assume that the principal symbol of $P$ is invertible. Then, for any $u \in \cD'(M,\cE)$ and any $x_{0}\in M$, we have
    \begin{gather}
        Pu \in W_{H}^{s}(M,\cE) \Longrightarrow u \in W_{H}^{s+k}(M,\cE),\\
        \text{$Pu$ is $W_{H}^{s}$ near $x_{0}$} \Longrightarrow \text{$u$ is $W_{H}^{s+k}$ near $x_{0}$}.
         \label{eq:Intro1.hypoellipricity-WHs}
    \end{gather}
     In fact, for any $s'\in \R$ we have the estimate,
    \begin{equation}
        \|u\|_{W^{s+k}_{H}} \leq C_{ss'}(\|Pu\|_{W^{s}_{H}}+\|u\|_{W^{s'}_{H}}), \qquad u \in W^{s+k}_{H}(M,\cE).
    \end{equation}
\end{proposition}

We can also obtain a version of Proposition~\ref{prop:Intro1.regularity-PsiHDOs} for holomorphic families and complex powers of hypoelliptic operators as follows.

\begin{proposition}
1)  Let $\Omega \subset \C$ be open and let $(P_{z})_{z\in \Omega}\in \Hol(\Omega, \Psi^{*}_{H}(M,\cE))$. Assume there exists a real $m$ 
    such that $\Re\ord P_{z}\leq m<\infty$. Then for any $s\in \R$ the family $(P_{z})_{z\in \Omega}$ defines 
    a holomorphic family with values in $\cL(W_{H}^{s+m}(M,\cE),W_{H}^{s}(M,\cE))$.\smallskip
 
2) Let $P:C^{\infty}(M,\cE)\rightarrow C^{\infty}(M,\cE)$ be a positive differential operator of Heisenberg  order $m$ such that the principal symbol 
of $P+\partial_{t}$ is an invertible Volterra-Heisenberg symbol. Then, for any reals $a$ and $s$, the complex powers of $P$ satisfy
\begin{equation}
    (P^{z})_{\Re z<a} \in \Hol(\Re z<a, \cL(W_{H}^{s+ma}(M,\cE),W_{H}^{s}(M,\cE))).
\end{equation}
\end{proposition}

\subsection{Rockland condition and the heat equation} Since Theorem~\ref{thm:Intro1.complex-powers} and the main results of~\cite{BGS:HECRM} 
hold under the condition that the principal symbol of the heat operator $P+\partial_{t}$ is an invertible Volterra-Heisenberg symbol, it stands to reason 
to find criterions for the invertibility of the latter. We show here that, to a large extent, the Rockland condition and the positivity of the principal symbol 
of $P$ are enough to imply the invertibility of the principal symbol of $P+\partial_{t}$. In particular, the results results of~\cite{BGS:HECRM} and 
Theorem~\ref{thm:Intro1.complex-powers-sublaplacian} are valid for a very wide class of hypoelliptic operators.  

Here we say that a homogeneous symbol 
$p_{m}$ of degree $m$ is positive if it can be put into the form $\overline{q_{\frac{m}{2}}}*q_{\frac{m}{2}}$ for some symbol 
$q_{\frac{m}{2}}$ homogeneous of degree $\frac{m}{2}$. It may be difficult in practice to check whether an operator has a positive principal symbol only 
by looking at this principal symbol. Nevertheless, for hypoelliptic operators we have: 

\begin{proposition}\label{prop:Intro1.positivity-criterion}
Assume that the Levi form of $(M,H)$ has constant rank  and let $P:C^{\infty}(M,\cE)\rightarrow C^{\infty}(M,\cE)$ 
be a selfadjoint differential operator of even Heisenberg order $m$ with an invertible principal symbol. 
Then the following are equivalent:\smallskip
   
   (i) The operator $P$ is bounded from below.\smallskip
   
   (ii) The principal symbol $\sigma_{m}(P)$ of $P$ is positive.\smallskip
   
\noindent In particular, for any lower order selfadjoint perturbation $R$ the operator $P+R$ remains bounded from below.   
\end{proposition}

Next, we prove: 

\begin{theorem}\label{thm:Intro1.heat}
Assume that the Levi form of $(M,H)$ has constant rank and 
let $P:C^{\infty}(M,\cE)\rightarrow C^{\infty}(M,\cE)$ be a selfadjoint differential operator of even Heisenberg order $m$ which is bounded from below 
and satisfies the Rockland condition at every point. Then:\smallskip
   
   1) The principal symbol of $P+\partial_{t}$ is an invertible Volterra-Heisenberg symbol.\smallskip 
  
  2)  The heat operator  $P+\partial_{t}$ admits a Volterra-\psivdo\ inverse;\smallskip 
   
   3) The heat kernel $k_{t}(x,y)$ of $P$ has an asymptotics in  $C^{\infty}(M,(\End \cE)\otimes|\Lambda|(M))$ of the form
    \begin{equation}
     k_{t}(x,x) \sim_{t\rightarrow 0^{+}} t^{-\frac{d+2}{m}} \sum t^{\frac{2j}{m}} a_{j}(P)(x),
    \label{eq:Intro1.Rockland-Heat.heat-kernel-asymptotics}
    \end{equation}
where the density $a_{j}(P)(x)$ is locally computable in terms of the symbol $q_{-m-2j}(x,\xi,\tau)$ of degree $-m-2j$ of any Volterra-\psivdo\ parametrix 
for $P+\partial_{t}$.
  \end{theorem}

  Notice that we only have to is to prove the first statement since the results of~\cite{BGS:HECRM} allows us to deduce the last two statements from the 
 first one. This is done by making use of  Theorem~\ref{thm:Intro1.complex-powers} and of results of Christ-Geller-G\l owacki-Polin~\cite{CGGP:POGD} and 
 Folland-Stein~\cite{FS:HSHG} for Rockland operators on nilpotent graded Lie groups. 

 The above theorem is valid for the following operators:\smallskip  
  
(a) The conformal powers $\boxdot_{\theta}^{(k)}$ of the (scalar) Tanaka Laplacian acting on functions on a strictly pseudoconvex CR 
manifold;\smallskip 

(b) The contact Laplacian on a contact manifold.\smallskip 

These examples are not covered by the results of~\cite{BGS:HECRM}. 

On the other hand, if $\Delta$ is a sublaplacian then, when the Levi form of $(M,H)$ is nondegenerate, it is shown in~\cite{BG:CHM} that the 
invertibility of the principal symbol of $\Delta$ is equivalent to the Rockland condition, but is weaker than the 
condition considered in~\cite[5.23]{BGS:HECRM} (see~\cite[18.4]{BG:CHM} and Section~\ref{sec:PsiHDO} for the precise statement of the condition). 
Anyway, using Theorem~\ref{thm:Intro1.heat} we get: 

 \begin{theorem}\label{thm:Intro1.heat-sublaplacian}
     Assume that the Levi form of $(M,H)$ is nondegenerate and let $\Delta:C^{\infty}(M,\cE)\rightarrow C^{\infty}(M,\cE)$ 
     be a selfadjoint sublaplacian which is bounded from below and satisfies the weaker condition of~\cite[18.4]{BG:CHM}. Then the conclusions of 
     Theorem~\ref{thm:Intro1.heat} are valid for $\Delta+\partial_{t}$ and the heat kernel of $\Delta$.
 \end{theorem}

\subsection{Spectral asymptotics for hypoelliptic operators} 
One interesting application of the heat kernel asymptotics in~\cite{BGS:HECRM} and in~(\ref{eq:Intro1.Rockland-Heat.heat-kernel-asymptotics}) 
is to allows us to derive spectral asymptotics for hypoelliptic operators as follows. 

Let $P:C^{\infty}(M,\cE)\rightarrow C^{\infty}(M,\cE)$ be a selfadjoint differential operator of Heisenberg order $m$ which is bounded from below and 
such that the principal symbol of $P+\partial_{t}$ is an invertible Volterra-Heisenberg symbol. 

Let $\lambda_{0}(P)\leq \lambda_{1}(P)\leq \ldots$ denote the eigenvalues of $P$ counted with multiplicity and  
let $N(P;\lambda)$ denote its counting function, that is, 
\begin{equation}
     N(P;\lambda)=\#\{k\in \N;\ \lambda_{k}(P)\leq \lambda \}, \qquad \lambda\geq 0.
\end{equation}
Then we prove: 
 \begin{theorem}\label{thm:Intro1.spectral-asymptotics}
Under the above assumptions the following hold.\smallskip 
     
     1) As $t\rightarrow 0^{+}$ we have
     \begin{equation}
         \Tr e^{-tP} \sim t^{-\frac{d+2}{m}} \sum t^{\frac{2j}{m}} A_{j}(P), \qquad  
   A_{j}(P)=\int_{M}\tr_{\cE} a_{j}(P)(x).
          \label{eq:Intro1.heat-trace-asymptotics}
     \end{equation}
 where the density $a_{j}(P)(x)$ is the coefficient of $t^{\frac{j-d+2}{m}}$ in the 
 asymptotics~(\ref{eq:Rockland-Heat.heat-kernel-asymptotics}).\smallskip 
      
    2) We have $A_{0}(P)>0$.\smallskip 
     
2) As $\lambda\rightarrow \infty$ we have 
\begin{equation}
    N(P;\lambda) \sim \nu_{0}(P)\lambda^{\frac{d+2}{m}}, \qquad 
    \nu_{0}(P)=\Gamma(1+\frac{d+2}{m})^{-1}A_{0}(P).
    \label{eq:Intro1.counting-function-asymptotics}
\end{equation}

3) As $k\rightarrow \infty$ we have 
     \begin{equation}
         \lambda_{k}(P)\sim \left(\frac{k}{\nu_{0}(P)}\right)^{\frac{m}{d+2}}. 
         \label{eq:Intro1.eigenvalue-asymptotics}
     \end{equation}     
 \end{theorem}

 The first asymptotics is an immediate consequence of~(\ref{eq:Intro1.Rockland-Heat.heat-kernel-asymptotics}). Moreover, once proved that we have 
 $A_{0}(P)>0$ the asymptotics~(\ref{eq:Intro1.counting-function-asymptotics}) and~(\ref{eq:Intro1.eigenvalue-asymptotics}) follows 
 from~(\ref{eq:Intro1.heat-trace-asymptotics}) and Karamata's Tauberian theorem. The bulk of the proof is thus to prove the second assertion, which is
 obtained by making use of spectral theoretic arguments.
 
By relying on other pseudodifferential calculi several authors have obtained Weyl asymptotics closely related to~(\ref{eq:Intro1.counting-function-asymptotics}) 
in the more general setting of hypoelliptic with multicharacteristics~(see \cite{II:PDPEAABSFSP}, \cite{Me:HOCVC2WE}, \cite{MS:ECH2}, 
\cite{Mo:ESOHCM1}, \cite{Mo:ESOHCM2}). 

As far as the Heisenberg setting is concerned, the approach using the Volterra-Heisenberg calculus has two main advantages. 
 
First, the pseudodifferential analysis is significantly simpler. In particular, the Volterra-Heisenberg calculus yields for free 
 the heat kernel asymptotics once the principal symbol of the heat operator is shown to be invertible, for which it is enough to use the Rockland condition 
 in many cases.  
 
 Second, as the Volterra-Heisenberg calculus take fully into account the underlying Heisenberg geometry of the manifold and is invariant by change of 
 Heisenberg coordinates, we can very effectively deal with operators admitting normal forms (see Proposition~\ref{prop:Spectral.normal-form-nu0P} on this point). 
 In particular, as explained below, we can reformulate in geometric terms the Weyl asymptotics~(\ref{eq:Intro1.counting-function-asymptotics}) for the main 
 geometric operators on CR and contact manifolds. 

\subsection{Weyl asymptotics and CR geometry} 
Let $M^{2n+1}$ be a compact $\kappa$-strictly pseudoconvex CR manifold, together with a pseudohermitian contact form $\theta$, that is, a contact form 
such that the associated Levi form has signature $(n-\kappa,\kappa,0)$. We also endow $M$ with a Levi metric compatible with $\theta$.  
Then the volume of $M$ with 
respect to this  Levi metric is independent of the choice of the Levi form and is equal to 
\begin{equation}
        \op{vol}_{\theta}M=\frac{(-1)^{\kappa}}{n!2^{n}} \int_{M} \theta \wedge d\theta^{n}.
\end{equation}
We call $\op{vol}_{\theta}M$ the pseudohermitian volume of $(M,\theta)$. We can relate the Weyl asymptotics~(\ref{eq:Intro1.counting-function-asymptotics}) 
for the Kohn Laplacian and the 
horizontal sublaplacian to $\op{vol}_{\theta}M$ as follows. 

First, for $\mu\in (-n,n)$ we let 
\begin{equation}
    \nu(\mu)= (2\pi)^{-(n+1)} \int_{-\infty}^{\infty}e^{-\mu\xi_{0}}(\frac{\xi_{0}}{\sinh \xi_{0}})^{n}d\xi_{0}.
\end{equation}

Second, in the case of the Kohn Laplacian  we obtain: 
\begin{theorem}
Let $\Box_{b}:C^{\infty}(M,\Lambda^{*,*})\rightarrow C^{\infty}(M,\Lambda^{*,*})$ denote the Kohn Laplacian associated to the Levi metric on $M$. 
Then, for $p,q=0,\ldots,n$ with $q \neq \kappa$ and $q\neq n-\kappa$, as $\lambda \rightarrow \infty$ we have 
\begin{gather}
    N(\Box_{b|_{\Lambda^{p,q}}};\lambda) \sim \alpha_{n\kappa pq}(\op{vol}_{\theta}M)\lambda^{n+1},\\  
    \alpha_{n\kappa pq}= \frac{1}{2^{n+1}} \binom{n}{p} \sum_{\max(0,q-\kappa)\leq  k\leq \min(q,n-\kappa)} \binom{n-\kappa}{k}\binom{\kappa}{q-k} 
    \nu(n-2(\kappa-q+2k)).
     \label{eq:Intro1.alphapq}
\end{gather}
\end{theorem}
In the strictly pseudoconvex case, i.e.~$\kappa=0$, the result is an easy consequence of~\cite[Thm.~8.31]{BGS:HECRM}, but in the case $\kappa \geq 1$ 
this seems to be new. 

Next, in the CR setting the horizontal sublaplacian preserves the bidegree and, in the same way as with the Kohn Laplacian, we get: 

\begin{theorem}\label{thm:Intro1.Deltab}
 Let $\Delta_{b}$ be the horizontal sublaplacian associated to the Levi metric on $M$. Then for $p,q=0,\ldots,n$, with 
  $(p,q)\neq (\kappa,n-\kappa)$ and $(p,q)\neq(n-\kappa,\kappa)$, as $\lambda  \rightarrow \infty$ we have 
    \begin{gather}
  N(\Delta_{b|_{\Lambda^{p,q}}};\lambda) \sim \beta_{npq}(\op{vol}_{\theta}M) \lambda^{n+1}, \label{eq:Intro1.Weyl-Deltab-CR}\\  
    \beta_{n\kappa pq}= \!  \!  \!  \!\sum_{\substack{\max(0,q-\kappa)\leq  k\leq \min(q,n-\kappa)\\ \max(0,p-\kappa)\leq l\leq \min(p,n-\kappa)}}  \!  
    \!  \!  \!  \binom{n-\kappa}{l}\binom{\kappa}{p-l} \binom{n-\kappa}{k}\binom{\kappa}{q-k} \nu(2(q-p)+4(l-k)).
     \label{eq:Intro1.betapq}
\end{gather}
\end{theorem}

Finally, as the Weyl asymptotics~(\ref{eq:Intro1.counting-function-asymptotics}) depends only on the principal symbol of the operator $P$, in the 
strictly pseudoconvex case, i.e.~$\kappa=0$, we can also deal with the conformal 
powers of the horizontal Laplacian.
\begin{theorem}
 Assume that $(M^{2n+1},\theta)$ is a strictly pseudoconvex pseudohermitian manifold and  for 
$k=1,\ldots,n+1$ let $\boxdot_{\theta}^{(k)}:C^{\infty}(M)\rightarrow C^{\infty}(M)$ be a $k$'th conformal power of $\Delta_{b}$ associated to $\theta$. 
Then  as $\lambda 
    \rightarrow \infty$ we have 
    \begin{equation}
         N(\boxdot_{\theta}^{(k)};\lambda) \sim \nu(0)(\op{vol}_{\theta}M) \lambda^{\frac{n+1}{k}}.    
     \end{equation}
\end{theorem}
 
 \subsection{Weyl asymptotics and contact geometry}
 Let $(M^{2n+1},\theta)$ be a compact orientable contact manifold. We assume that the Heisenberg bundle $H=\ker \theta$ has a calibrated almost complex
 structure $J$ so that $d\theta(X,JX)=-d\theta(JX,X)>0$ for any section $X$ of $H\setminus 0$. We then endow $M$ with the Riemannian metric 
 $g_{\theta}=d\theta(.,J.)+\theta^{2}$. The volume of $M$ with respect to $g_{\theta}$ depends only on $\theta$ and is equal to:
 \begin{equation}
     \op{vol}_{\theta}M=\frac{1}{n!}\int_{M}d\theta^{n}\wedge \theta.
 \end{equation}
 We call $\op{vol}_{\theta}M$ the contact volume of $M$. 
  
 We can relate the Weyl asymptotics for the horizontal sublaplacian to the contact volume to get:
 
 \begin{theorem}
   Let $\Delta_{b}:C^{\infty}(M,\Lambda^{*}_{\C}H^{*})\rightarrow  C^{\infty}(M,\Lambda^{*}_{\C}H^{*})$ be the horizontal sublaplacian associated to 
   the metric $g_{\theta}$. Then, for $k=0,\ldots,2n$ with $k\neq n$, as $\lambda \rightarrow \infty$ we have 
    \begin{gather}
         N(\Delta_{b|_{\Lambda^{k}_{\C}H^{*}}};\lambda) \sim \gamma_{nk} (\op{vol}_{\theta}M) \lambda^{n+1},\qquad 
          \label{eq:Intro1.Weyl-Deltab-contact}
 \gamma_{nk}=2^{-n}\sum_{p+q=k}\binom{n}{p} \binom{n}{q}\nu(p-q).         
    \end{gather}
\end{theorem}
 
Note that  when $(M,\theta)$ is a strictly pseudoconvex pseudohermitian manifold the asymptotics~(\ref{eq:Intro1.Weyl-Deltab-contact}) 
is compatible with~(\ref{eq:Intro1.Weyl-Deltab-CR}) because the contact volume differs from the pseudohermitian volume by a factor of $2^{-n}$.

Finally, we can also deal with the contact Laplacian as follows. 
\begin{theorem}
    Let $\Delta_{R}:C^{\infty}(M,\Lambda^{*}\oplus \Lambda^{n}_{*})\rightarrow C^{\infty}(M,\Lambda^{*}\oplus \Lambda^{n}_{*})$ 
    be the contact Laplacian on $M$.\smallskip 
    
    1)  For $k=0,\ldots,2n$ with $k\neq n$ there exists a  universal constant $\nu_{nk}>0$ depending only on $n$ and $k$ 
    such that as $\lambda \rightarrow \infty$ we have
     \begin{equation}
           N(\Delta_{R|_{\Lambda^{k}}})\sim \nu_{nk} (\op{vol}_{\theta}M)\lambda^{n+1}. 
                \label{eq:Intro1.Weyl-contact-Laplacian1}
      \end{equation}
    
    2) For $j=1,2$ there exists a  universal constant $\nu_{n}^{(j)}>0$ depending only on $n$ and $j$  such that as $\lambda \rightarrow \infty$ we have 
    \begin{equation}
      N(\Delta_{R|_{\Lambda^{n}_{j}}})\sim \nu_{n}^{(j)} (\op{vol}_{\theta}M)\lambda^{\frac{n+1}{2}}. 
         \label{eq:Intro1.Weyl-contact-Laplacian2}
    \end{equation}
\end{theorem}

\subsection{Organization of the paper} The paper is organized as follows. We start by reviewing the background needed in this paper. In 
Section~\ref{sec.Heisenberg} we recall the main definitions and examples concerning Heisenberg manifolds and their tangent Lie group bundles following 
the point of view of~\cite{Po:Pacific1}. In Section~\ref{sec:Operators} we recall the definitions of main operators on Heisenberg manifolds: sum of squares, Kohn 
Laplacian, horizontal sublaplacian and its conformal powers and the contact Laplacian. In Section~\ref{sec:PsiHDO} we give a detailed overview of the Heisenberg 
calculus of~\cite{BG:CHM} and~\cite{Ta:NCMA}, following closely the expositions of~\cite{BG:CHM} and~\cite{Po:BSM1}. 

In Section~\ref{sec.HolPHDO} we define holomorphic families of \psivdos\ and check their main properties. In Section~\ref{sec.powers1}, after 
having recalled the pseudodifferential representation of the heat kernel of~\cite{BGS:HECRM}, we make use of the latter to 
construct complex powers of 
positive hypoelliptic operators as holomorphic families of \psivdos.  
In Section~\ref{sec.Sobolev} we define the weighted Sobolev spaces $W_{H}(M,\cE)$ and show that they give 
sharp regularity results for \psivdos. In Section~\ref{sec:Rockland-heat} we show that if $P$ is selfadjoint differential operator bounded from below, then 
the fact that $P$ satisfies the Rockland condition at every point is 
enough to insure us the invertibility of the principal symbol of $P+\partial_{t}$. This show that the results of~\cite{BGS:HECRM} and Section~\ref{sec.powers1} 
are valid for a wide class of hypoelliptic operators. 

In the last three section we deal with spectral asymptotics for hypoelliptic operators on Heisenberg manifolds. In Section~\ref{sec:Spectral} we derive 
general spectral asymptotics for such operators on a general Heisenberg manifold. We specialize them in Section~\ref{sec:Spectral-CR} to the 
Kohn Laplacian and the horizontal sublaplacian and its conformal powers on a CR manifold. Finally, in Section~\ref{sec:Spectral-contact} look at these 
asymptotics in the special cases of the horizontal sublaplacian and of the contact Laplacian on a contact manifold. 

\begin{acknowledgements} 
The author is grateful  to Alain Connes, Charles Epstein, Bernard Helffer, Henri Moscovici and Michel Rumin for helpful and 
 stimulating discussions. He wish also to thank for its hospitality the Mathematics Department of Princeton 
 University, where this paper was finally completed. 
 
 In addition, the research of the paper was partially supported by a graduate fellowships from the French Ministry of National Education and Research, by a 
 fellowship from the RT Network geometric analysis HPCRN-CT-1999-00118, and by the NSF grant DMS 0409005. 
 
 Finally, some of the results of this paper were announced in~\cite{Po:CRAS1} and presented as part of 
 the author's PhD thesis at University of Paris-Sud (Orsay, France), made under the supervision of Alain Connes.  
 \end{acknowledgements}

\section{Heisenberg manifolds and their tangent Lie group bundles} 
\label{sec.Heisenberg}
In this section we recall the main facts about Heisenberg manifolds and their tangent Lie group bundle. The exposition here follows closely that of~\cite{Po:Pacific1}.

\begin{definition}
   1) A Heisenberg manifold is a smooth manifold $M$ equipped with a distinguished hyperplane bundle $H \subset TM$. \smallskip 
   
   2) A Heisenberg diffeomorphism $\phi$ from a Heisenberg manifold $(M,H)$ onto another Heisenberg manifold 
   $(M,H')$ is a diffeomorphism $\phi:M\rightarrow M'$ such that $\phi^{*}H = H'$. 
\end{definition}

\begin{definition}
   Let $(M^{d+1},H)$ be a Heisenberg manifold. Then:\smallskip 
   
   1) A (local) $H$-frame for $TM$ is  a (local)  frame $X_{0}, X_{1}, \ldots, X_{d}$ of $TM$ so that $X_{1}, \ldots, 
   X_{d}$ span~$H$.\smallskip  
   
   2) A local Heisenberg chart is a  local chart with a local $H$-frame of $TM$ over its domain.
\end{definition}

The main examples of Heisenberg manifolds are the following.\smallskip 

\emph{a) Heisenberg group}. The $(2n+1)$-dimensional Heisenberg group
$\bH^{2n+1}$ is $\R^{2n+1}=\R \times \R^{n}$ equipped with the 
group law, 
\begin{equation}
    x.y=(x_{0}+y_{0}+\sum_{1\leq j\leq n}(x_{n+j}y_{j}-x_{j}y_{n+j}),x_{1}+y_{1},\ldots,x_{2n}+y_{2n}).  
\end{equation}
A left-invariant basis for its Lie algebra $\fh^{2n+1}$ is then
provided by the vector-fields, 
\begin{equation}
    X_{0}=\frac{\partial}{\partial x_{0}}, \quad X_{j}=\frac{\partial}{\partial x_{j}}+x_{n+j}\frac{\partial}{\partial 
    x_{0}}, \quad X_{n+j}=\frac{\partial}{\partial x_{n+j}}-x_{j}\frac{\partial}{\partial 
    x_{0}}, \quad 1\leq j\leq n,
     \label{eq:Examples.Heisenberg-left-invariant-basis}
\end{equation}
which  for $j,k=1,\ldots,n$ and $k\neq j$ satisfy the relations,
\begin{equation}
    [X_{j},X_{n+k}]=-2\delta_{jk}X_{0}¥, \qquad [X_{0},X_{j}]=[X_{j},X_{k}]=[X_{n+j},X_{n+k}]=0.
     \label{eq:Examples.Heisenberg-relations}
\end{equation}
In particular, the subbundle spanned by the vector fields 
$X_{1},\ldots,X_{2n}$ yields a left-invariant Heisenberg structure on 
$\bH^{2n+1}$.\smallskip

 \emph{b) Foliations.} Recall that a (smooth) foliation is a manifold $M$ together with a subbundle $\cF \subset TM$ 
which is integrable in the Froebenius' sense, i.e.~so that
$[\cF,\cF]\subset \cF$. Therefore, any codimension 1 foliation is a Heisenberg manifold.\smallskip  

 \emph{c) Contact manifolds}. 
Opposite to foliations are contact manifolds: a \emph{contact
structure} on a manifold $M^{2n+1}$ is given by a global non-vanishing $1$-form $\theta$ on $M$ such that
$d\theta$ is non-degenerate on $H=\ker \theta$. In particular, $(M,H)$ is a Heisenberg manifold. In fact, by
Darboux's theorem any contact manifold $(M^{2n+1},\theta)$ is locally
contact-diffeomorphic to the Heisenberg group $\bH^{2n+1}$ equipped with its standard contact
form $\theta^{0}= dx_{0}+\sum_{j=1}^{n}(x_{j}dx_{n+j}-x_{n+j}dx_{j})$.\smallskip

 \emph{d) Confoliations}. According to Elyashberg-Thurston~\cite{ET:C} a \emph{confoliation structure} on an oriented manifold
$M^{2n+1}$ is given by a global non-vanishing $1$-form $\theta$ on $M$ such that
$(d\theta)^{n}\wedge \theta\geq 0$. In particular, when $d\theta
\wedge \theta=0$ (resp.~$(d\theta)^{n}\wedge \theta>0$) we are
in presence of a foliation (resp.~a contact structure). In any case the hyperplane bundle $H=\ker \theta$ defines a 
Heisenberg structure on $M$.\smallskip

 \emph{e) CR manifolds.} A CR
structure on an orientable manifold $M^{2n+1}$ is given by a rank $n$
complex subbundle $T_{1,0}\subset T_{\C}M$ which is integrable in  Froebenius' sense and such that 
$T_{1,0}\cap T_{0,1}=\{0\}$, where $T_{0,1}=\overline{T_{1,0}}$. 
Equivalently, the subbundle $H=\Re (T_{1,0}\oplus T_{0,1})$ has the 
structure of a complex bundle of (real) dimension $2n$. In
particular, the pair $(M,H)$ forms a Heisenberg manifold. 

Moreover, since $M$ and $H$ is orientable by means of its complex structure the normal bundle $TM/H$ is orientable, hence admits a global 
nonvanishing section $T$. Let $\theta$ be a global section of $T^{*}M/H^{*}$ such that $\theta(T)=1$ and $\theta$ annihilates $H$. Then Kohn's Levi 
form is the form $L_{\theta}$ on $T_{1,0}$ such that, for for sections $Z$ and $W$ of $T_{1,0}$, we have
\begin{equation}
    L_{\theta}(Z,W)=-id\theta(Z,\bar{W})=i\theta([Z,\bar{W}]). 
     \label{eq:Heisenberg.Kohn-Levi-form}
\end{equation}

We say that $M$ is strictly pseudoconvex (resp.~nondegenerate, $\kappa$-strictly pseudoconvex) when for some choice of $\theta$  the Levi form 
$L_{\theta}$  is everywhere positive definite (resp.~is everywhere non-degenerate, has everywhere signature $(n-\kappa,\kappa,0)$). In particular, 
when $M$ is nondegenerate the 1-form $\theta$ is a contact form on $M$. 

The main example of a CR manifold is that of the (smooth) boundary $M=\partial D$ of a complex domain $D \subset \C^{n}$. In particular, 
when $D$ is strongly pseudoconvex (or strongly pseudoconcave) then $M$ is strictly pseudoconvex.

\subsection{The tangent Lie group bundle}
A simple description of the tangent Lie group bundle of a Heisenberg manifold $(M^{d+1},H)$ is given as follows.

\begin{lemma}[\cite{Po:Pacific1}]
The Lie bracket of vector fields induces on $H$ a 2-form with values in $TM/H$, 
\begin{equation}
    \cL: H\times H \longrightarrow TM/H,
     \label{eq:Bundle.Levi-form1}
\end{equation}
so that for any sections $X$ and $Y$ of $H$ near a point $a\in M$ we have
\begin{equation}
    \cL_{a}(X(a),Y(a)) = [X,Y](a) \quad \bmod H_{a}.
     \label{eq:Bundle.Levi-form2}
\end{equation}
\end{lemma}

\begin{definition}
 The $2$-form  $\cL$ is called the Levi form of $(M,H)$.
\end{definition}

The Levi form $\cL$ allows us to define a bundle $\fg M$ of graded Lie algebras  by endowing $(TM/H)\oplus H$ 
with the smooth fields of Lie Brackets and gradings such that
\begin{equation}
    [X_{0}+X',Y_{0}+Y']_{a}=\cL_{a}(X',Y') \qquad \text{and} \qquad t.(X_{0}+X')=t^{2}X_{0}+tX' \quad t \in \R,
    \label{eq:Heisenberg.intrinsic-Lie-algebra-structure}
\end{equation}
for $a\in M$ and $X_{0}$, $Y_{0}$ in $T_{a}M/H_{a}$ and $X'$, $Y'$ in $H_{a}$. 

\begin{definition}
    The bundle $\fg M$ is called the tangent Lie algebra bundle of $M$.
\end{definition}

As we can easily check $\fg M$ is a bundle of $2$-step nilpotent Lie algebras which contains the normal bundle $TM/H$ in its center.
Therefore, its associated 
graded Lie group bundle $GM$ can be described as follows. As a bundle $GM$ is $(TM/H)\oplus H$ and the exponential 
map is merely the identity. In particular, the grading of $GM$ is as in~(\ref{eq:Heisenberg.intrinsic-Lie-algebra-structure}). 
Moreover, since  $\fg M$ is 
2-step nilpotent the Campbell-Hausdorff formula gives 
\begin{equation}
    (\exp X)(\exp Y)= \exp(X+Y+\frac{1}{2}[X,Y]) \qquad \text{for  sections $X$, $Y$ of $\fg M$}.
\end{equation}
From this we deduce that the product on $GM$ is such that  
\begin{equation}
    (X_{0}+X').(Y_{0}+X')=X_{0}+Y_{0}+\frac{1}{2}\cL(X',Y')+X'+Y',    
    \label{eq:Bundle.Lie-group-law}
\end{equation}
for  sections $X_{0}$, $Y_{0}$ of $TM/H$  and sections $X'$, $Y'$ of $H$.

\begin{definition}
    The bundle  $GM$ is called the tangent Lie group bundle of $M$. 
\end{definition}

In fact, the fibers of $GM$ as classified by the Levi form $\cL$ as follows.

\begin{proposition}[\cite{Po:Pacific1}]\label{prop:Bundle.intrinsic.fiber-structure}
  1) Let $a\in M$. Then $\cL_{a}$ has rank $2n$ if, and only if, as a 
  graded Lie group $G_{a}M$ is isomorphic to $\bH^{2n+1}\times \R^{d-2n}$.\smallskip 
  
  2) The Levi form $\cL$ has constant rank $2n$ if, and only if, $GM$ is  a fiber bundle with typical fiber 
  $\bH^{2n+1}\times \R^{d-2n}$.
\end{proposition}
 
Now, let $\phi:(M,H)\rightarrow (M',H')$ be a Heisenberg diffeomorphism from $(M, H)$ onto another Heisenberg manifold 
$(M',H')$. Since $\phi_{*}H=H'$ we see that $\phi'$ induces a smooth vector bundle isomorphism 
$\overline{\phi}:TM/H\rightarrow TM'/H'$. 

\begin{definition}
We let  $\phi_{H}':(TM/H)\oplus 
  H \rightarrow (TM'/H')\oplus H'$ denote the vector bundle isomorphism such that
    \begin{equation}
    \phi'_{H}(a)(X_{0}+X')=\overline{\phi}'(a)X_{0}+\phi'(a)X',
     \label{eq:Bundle.Intrinsic.Phi'H}
\end{equation}
for any $a\in M$ and any $X_{0}\in T_{a}/H_{a}$ and $X'\in H_{a}$.
\end{definition}

\begin{proposition}[\cite{Po:Pacific1}]\label{prop:Bundle.Intrinsic.Isomorphism}
The vector bundle isomorphism  $\phi'_{H}$ is an isomorphism of graded Lie group bundles from $GM$ onto $GM'$. In particular, the Lie group bundle isomorphism 
class of $GM$ depends only on the Heisenberg diffeomorphism class of $(M,H)$.  
\end{proposition}

\subsection{Heisenberg coordinates and nilpotent approximation of vector fields}
It is interesting to relate the intrinsic description of $GM$ above with the more extrinsic description of~\cite{BG:CHM} (see also~\cite{Be:TSSRG}, 
\cite{EMM:HAITH}, \cite{FS:EDdbarbCAHG},  \cite{Gr:CCSSW}, \cite{Ro:INA}) in terms of the Lie group 
associated to a nilpotent Lie algebra of model vector fields. 

First, let $a\in M$ and let us describe $\fg_{a}M$ as the graded Lie algebra of left-invariant vector fields on $G_{a}M$  
by identifying any $X \in \fg_{a}M$ with the left-invariant vector fields $L_{X}$ on $G_{a}M$ given by 
\begin{equation}
    L_{X}f(x)= \frac{d}{dt}f[(t\exp X).x]_{|_{t=0}}= \frac{d}{dt}f[(tX).x]_{|_{t=0}}, \qquad f \in C^{\infty}(G_{a}M).
\end{equation}
This allows us to associate to any vector fields $X$ near $a$ a unique left-invariant vector fields $X^{a}$ on $G_{a}M$ 
such that 
\begin{equation}
    X^{a}= \left\{ 
    \begin{array}{ll}
        L_{X_{0}(a)} & \text{if $X(a)\not \in H_{a}$},  \\
        L_{X(a)} & \text{otherwise,} 
    \end{array}\right.
     \label{eq:Bundle.intrinsic.model-vector-fields}
\end{equation}
where $X_{0}(a)$ denotes the class of $X(a)$ modulo $H_{a}$. 

\begin{definition}
    The left-invariant vector fields $X^{a}$ is called the model vector fields of $X$ at $a$.
\end{definition}

Let us look at the above construction in terms of a $H$-frame $X_{0},\ldots,X_{d}$ near 
$a$, i.e.~of a local trivialization of the vector bundle $(TM/H)\oplus H$. For $j,k=1,\ldots,d$ we let 
\begin{equation}
    \cL(X_{j},X_{k})=[X_{j},X_{k}]X_{0}=L_{jk}X_{0} \quad \bmod H.
\end{equation}
With respect to the coordinate system $(x_{0},\ldots,x_{d})\rightarrow x_{0}X_{0}(a)+\ldots+x_{d}X_{d}(a)$ we can 
write the product law of $G_{a}M$ as 
\begin{equation}
    x.y=(x_{0}+\frac{1}{2}\sum_{j,k=1}^{d}L_{jk}x_{j}x_{k},x_{1},\ldots,x_{d}).
     \label{eq:Heisenberg.productGmM-coordinates}
\end{equation}
Then the vector fields $X_{j}^{a}$, $j=1,\ldots,d$, in~(\ref{eq:Bundle.intrinsic.model-vector-fields}) 
is just the left-invariant vector fields corresponding to the vector $e_{j}$ of the canonical basis 
of $\Rd$, that is, we have
\begin{equation}
    X_{0}^{a}=\frac{\partial}{\partial x_{0}} \quad \text{and}  \quad X_{j}^{a}=\frac{\partial}{\partial x_{j}} 
    -\frac{1}{2}\sum_{k=1}^{d}L_{jk}x_{k}\frac{\partial}{\partial x_{0}}, \quad 1\leq j\leq d.
     \label{eq:Heisenberg.Xjm.coordinates}
\end{equation}
In particular, for $j,k=1,\ldots,d$ we have the relations, 
\begin{equation}
    [X_{j}^{a},X_{k}^{a}]=L_{jk}(a)X_{0}^{a}, \qquad [X_{j}^{a},X_{0}^{a}]=0.
     \label{eq:Heisenberg.constant-structures.Gm}
\end{equation}

Now, let $\kappa:\dom \kappa \rightarrow U$ be a Heisenberg chart near $a=\kappa^{-1}(u)$ and let 
$X_{0},\ldots,X_{d}$ be the associated $H$-frame of $TU$.  
Then there exists a unique affine coordinate change $x \rightarrow \psi_{u}(x)$ such that 
$\psi_{u}(u)=0$ and $\psi_{u*}X_{j}(0)=\frac{\partial}{\partial x_{j}}$ for 
$j=0,1,\ldots,d$. Indeed, if for $j=1,\ldots,d$ we set $X_{j}(x)=\sum_{k=0}^{d}B_{jk}(x)\frac{\partial}{\partial x_{k}}$ then 
we have
\begin{equation}
    \psi_{u}(x)=A(u)(x-u), \qquad A(u)=(B(u)^{t})^{-1}.
\end{equation}

\begin{definition}\label{def:Heisenberg.extrinsic.u-coordinates}
1) The coordinates provided by $\psi_{u}$ are called the privileged coordinates at $u$ 
with respect to the $H$-frame $X_{0},\ldots,X_{d}$. 

2) The map $\psi_{u}$ is called the privileged-coordinate map with respect to the $H$-frame $X_{0},\ldots,X_{d}$.
\end{definition}
\begin{remark}
    The privileged coordinates at $u$ are called $u$-coordinates in~\cite{BG:CHM}, but they correspond to the privileged coordinates 
    of~\cite{Be:TSSRG} and \cite{Gr:CCSSW} in the special case of a Heisenberg manifold. 
\end{remark}

Next, on $\Rd$ we consider the dilations 
\begin{equation}
    \delta_{t}(x)=t.x=(t^{2}x_{0},tx_{1}, \ldots, tx_{d}), \qquad t \in \R,
    \label{eq:Heisenberg.dilations}
\end{equation}
with respect to which $\frac{\partial}{\partial_{x_{0}}}$ is 
homogeneous of degree $-2$ and $\frac{\partial}{\partial_{x_{1}}},\ldots,\frac{\partial}{\partial_{x_{d}}}$ 
is homogeneous of degree~$-1$. 

Since in the privileged coordinates at $u$ we have $X_{j}(0)=\frac{\partial}{\partial x_{j}}$ we can write
\begin{equation}
    X_{j}= \frac{\partial}{\partial_{x_{j}}}+ \sum_{k=0}^{d} a_{jk}(x) \frac{\partial}{\partial_{x_{k}}}, 
    \qquad j=0,1,\ldots d,
\end{equation}
where the $a_{jk}$'s are smooth functions such that $a_{jk}(0)=0$. Therefore, we may define
\begin{gather}
    X_{0}^{(u)}= \lim_{t\rightarrow 0} t^{2}\delta_{t}^{*}X_{0}= \frac{\partial}{\partial_{x_{0}}},
    \label{eq:Heisenberg.X0u}\\
     X_{j}^{(u)}= \lim_{t\rightarrow 0} t^{-1}\delta_{t}^{*}X_{j}= 
     \frac{\partial}{\partial_{x_{j}}}+\sum_{k=1}^{d}b_{jk}x_{k} \frac{\partial}{\partial_{x_{0}}}, \quad 
     j=1,\ldots,d, \label{eq:Heisenberg.Xju}     
\end{gather}
where for $j,k=1,\ldots,d$ we have set $b_{jk}= \partial_{x_{k}}a_{j0}(0)$. 

Observe that $X_{0}^{(u)}$ is homogeneous of degree $-2$ and $X_{1}^{(u)},\ldots,X_{d}^{(u)}$ are homogeneous of degree $-1$.
 Moreover, for $j,k=1,\ldots,d$ we have 
\begin{equation}
    [X_{j}^{(u)},X_{0}^{(u)}]=0 \quad \text{and} \quad [X_{j}^{(u)},X_{0}^{(u)}]=(b_{kj}-b_{jk})X_{0}^{(u)}.
     \label{eq:Heisenberg.constant-structures.Gu1}
\end{equation}
Thus, the linear space spanned by $X_{0}^{(u)},X_{1}^{(u)}, \ldots, X_{d}^{(u)}$ is a graded 2-step nilpotent 
Lie algebra $\fg^{(u)}$. In particular, $\fg^{(u)}$ is the Lie algebra of left-invariant vector fields over the graded Lie group $G^{(u)}$ 
consisting of $\Rd$ equipped with the grading~(\ref{eq:Heisenberg.dilations}) and the group law,
\begin{equation}
    x.y=(x_{0}+\sum_{j,k=1}^{d}b_{kj}x_{j}x_{k},x_{1},\ldots,x_{d}).
\end{equation}

Now, if near $a$ we let $\cL(X_{j},X_{k})=[X_{j},X_{k}]=L_{jk}(x)X_{0}\bmod H$ then we have 
\begin{equation}
    [X_{j}^{(u)},X_{k}^{(u)}]=\lim_{t\rightarrow 0}[t\delta_{t}^{*}X_{j},t\delta_{t}^{*}X_{k}] = 
    \lim_{t\rightarrow 0} t^{2}\delta_{t}^{*}(L_{jk}(\circ \kappa^{-1}(x))X_{0})=L_{jk}(a)X_{0}^{(u)}.
     \label{eq:Heisenberg.constant-structures.Gu2}
\end{equation}
Comparing this with~(\ref{eq:Heisenberg.constant-structures.Gm}) and~(\ref{eq:Heisenberg.constant-structures.Gu1}) 
then shows that $\fg^{(u)}$ has the same the constant structures as those of 
$\fg_{a}M$, hence is isomorphic to $\fg_{a}M$. Consequently, the Lie groups $G^{(u)}$ and $G_{a}M$ are isomorphic. 
In fact, as it follows from~\cite{BG:CHM} and~\cite{Po:Pacific1} an explicit isomorphism is given by 
\begin{equation}
     \phi_{u}(x_{0},\ldots,x_{d})= (x_{0}-\frac{1}{4}\sum_{j,k=1}^{d}(b_{jk}+b_{kj})x_{j}x_{k},x_{1},\ldots,x_{d}).
     \label{eq:Bundle.Extrinsic.Phiu}
\end{equation}

\begin{definition}\label{def:Bundle.extrinsic.normal-coordinates}
Let $\varepsilon_{u}=\phi_{u}\circ \psi_{u}$. Then:\smallskip 

1) The new coordinates provided by $\varepsilon_{u}$  are called Heisenberg 
coordinates at $u$ with respect to the $H$-frame $X_{0},\ldots,X_{d}$.\smallskip  

2) The map $\varepsilon_{u}$ is called the $u$-Heisenberg coordinate map.
\end{definition}

\begin{remark}
       The Heisenberg coordinates at $u$ have been also considered in~\cite{BG:CHM} as a technical tool 
       for inverting the principal symbol of a hypoelliptic sublaplacian.
\end{remark}

  Next, as it follows from~\cite[Lem.~1.17]{Po:Pacific1} we also have
\begin{equation}
    \phi_{*}X_{0}^{(u)}=\frac{\partial}{\partial x_{0}}=X_{0}^{a} \quad \text{and} \quad 
    \phi_{*}X_{j}^{(u)}=\frac{\partial}{\partial x_{j}}-\frac{1}{2}\sum_{k=1}^{d}L_{jk}x_{k}\frac{\partial}{\partial x_{0}}=X_{j}^{a}, \quad 
    j=1,\ldots,d. 
    \label{eq:Heisenberg.Xu-Xm}
\end{equation}
Since $\phi_{u}$ commutes with the Heisenberg dilations~(\ref{eq:Heisenberg.dilations}) 
using~(\ref{eq:Heisenberg.X0u})--(\ref{eq:Heisenberg.Xju}) we get
\begin{equation}
    \lim_{t\rightarrow 0} t^{2}\delta_{t}^{*}\phi_{u*}X_{0}^{(u)}=X^{a}_{0} \quad \text{and} \quad 
    \lim_{t\rightarrow 0} t\delta_{t}^{*}\phi_{u*}X_{j}^{(u)}=X^{a}_{j}, \quad j=1,\ldots,d.
\end{equation}

In fact, as shown in~\cite{Po:Pacific1} 
for any vector fields $X$ near $a$, as $t\rightarrow 0$ and in Heisenberg coordinates at $a$, we have
\begin{equation}
   \delta_{t}^{*}X=  \left\{ 
   \begin{array}{ll}
       t^{-2}X^{a} +\op{O}(t^{-1})& \text{if $X(a)\in H_{a}$},\\
       t^{-1}X^{a} +\op{O}(1) & \text{otherwise}. 
   \end{array}\right. 
   \label{eq:Bundle.Extrinsic.approximation-normal}
\end{equation}
Therefore, we obtain:

\begin{proposition}[\cite{Po:Pacific1}]\label{prop:Bundle.equivalent-descriptions}
    In the Heisenberg coordinates centered at $m=\kappa^{-1}(u)$ the tangent Lie group $G_{a}M$ coincides with $G^{(u)}$ and for any vector fields 
    $X$ the model vector fields $X^{a}$ approximates $X$ near $a$ in the sense of~(\ref{eq:Bundle.Extrinsic.approximation-normal}).
\end{proposition}

One consequence of the equivalence between the two approaches to $GM$ is a tangent approximation for Heisenberg diffeomorphisms as follows. 

Let  $\phi:(M,H)\rightarrow (M',H')$ be a Heisenberg diffeomorphism
from $(M,H)$ to another Heisenberg manifold $(M',H')$. We also endow $\Rd$ with the pseudo-norm,
\begin{equation}
    \|x\|= (x_{0}^{2}+(x_{1}^{2}+\ldots+x_{d}^{2})^{2})^{1/4}, \qquad x\in \Rd,
\end{equation}
so that for any $x \in \Rd$ and any $t \in \R$ we have 
\begin{equation}
    \|t.x\|=|t|\, \|x\| . 
     \label{eq:Bundle.homogeneity-pseudonorm}
\end{equation}

\begin{proposition}[{\cite[Prop.~2.21]{Po:Pacific1}}]\label{prop:Heisenberg.diffeo}
   Let $a\in M$ and set $m'=\phi(a)$. Then, in Heisenberg coordinates at $a$ and at $a'$   
   the diffeomorphism $\phi(x)$ has a behavior near 
   $x=0$ of the form 
   \begin{equation}
       \phi(x)= \phi_{H}'(0)x+(\op{O}(\|x\|^{3}), \op{O}(\|x\|^{2}),\ldots,\op{O}(\|x\|^{2})).
        \label{eq:Bundle.Approximation-diffeo}
   \end{equation}
   In particular, there is no term of the form $x_{j}x_{k}$, $1\leq j,k\leq d$, in the Taylor expansion of $\phi_{0}(x)$ 
   at $x=0$.
\end{proposition}
   \begin{remark}
  An asymptotics similar to~(\ref{eq:Bundle.Approximation-diffeo}) is given  
  in~\cite[Prop.~5.20]{Be:TSSRG} in privileged coordinates at $u$ and 
   $u'=\kappa_{1}(a')$, 
  but the leading term there is only a Lie algebra isomorphism from $\fg^{(u)}$ onto $\fg^{(u')}$. This is only 
  in Heisenberg coordinates that we recover the Lie group isomorphism $\phi'_{H}(a)$ as the leading term of the asymptotics.
\end{remark}
   
\begin{remark}
    An interesting application  of Proposition~\ref{prop:Heisenberg.diffeo} in~\cite{Po:Pacific1} is the construction of the tangent groupoid 
    $\cG_{H}M$ of $(M,H)$ as the differentiable groupoid encoding the smooth deformation of  $M\times M$ to $GM$. 
   This groupoid is the analogue in the Heisenberg setting of Connes' tangent groupoid of a manifold~(\cite[II.5]{Co:NCG}, \cite{HS:MKOEFFTK}) and its shows 
   that $GM$ is tangent to $a$ in a differentiable  fashion (compare~\cite{Be:TSSRG}, \cite{Gr:CCSSW}). 
\end{remark}

\section{Most common operators on Heisenberg manifolds} 
\label{sec:Operators}
In the his section we recall the constructions of the most common operators on a Heisenberg manifold. At the exception of the contact Laplacian, 
all these operators are sublaplacians or coincides with an integer power of a sublaplacian up to a lower order term. 

Recall that a sublaplacian on a Heisenberg manifold $(M^{d+1},H)$ acting on the sections of a vector bundle $\cE$ over $M$ is 
differential operator $\Delta:C^{\infty}(M,\cE)\rightarrow C^{\infty}(M,\cE)$ such that, near any $a\in M$, there exists a 
$H$-frame $X_{0},X_{1},\ldots,X_{d}$ of $TM$ so that $\Delta$ takes the form
 \begin{equation}
    \Delta=-\sum_{j=1}^{d} X_{j}^{2} + \text{lower order terms}.
    \label{eq:Heisenberg.sublaplacian}
\end{equation}
%
%

\subsection{H\"ormander's sum of squares.}
 Let $X_{1},\ldots,X_{m}$ be (real) vector fields on a manifold $M^{d+1}$ and consider the sum of squares,
\begin{equation}
    \Delta=-(X_{1}^{2}+\ldots+X_{m}^{2}).
\end{equation}
By a celebrated theorem of H\"ormander~\cite{Ho:HSODE} the operator $\Delta$ is hypoelliptic provided that the following bracket condition is 
satisfied:the vector fields $X_{0},X_{1},\ldots,X_{m}$ together with their successive Lie brackets $[X_{j_{1}}, [X_{j_{2}},\ldots ,
X_{j_{1}}]\ldots]]$ span the tangent bundle $TM$ at every point. 

When $X_{1},\ldots,X_{m}$ span a hyperplane bundle $H$ then $\Delta$ is a sublaplacian with \emph{real} coefficients and the bracket condition 
reduces to $H+[H,H]=TM$ or, equivalently, the Levi form of $(M,H)$ is nonvanishing. 

In fact, given a vector bundle $\cE$, the theorem of H\"ormander holds for more general sublaplacians $\Delta:C^{\infty}(M,\cE)\rightarrow 
C^{\infty}(M,\cE)$ of the form
\begin{equation}
    \Delta=-(\nabla_{X_{1}}^{2}+\ldots+\nabla_{X_{m}}^{2})+L,
     \label{eq:Operators.generalized-sum-of-squares}
\end{equation}
where $\nabla$ is a connection on $\cE$ and $L=L(X_{1},\ldots,X_{m})$ is a first order polynomial with \emph{real} coefficients in the vector fields 
$X_{1},\ldots,X_{m}$. In particular, if $M$ is endowed with a positive density and $\cE$ with a Hermitian metric, this includes the selfadjoint sum of squares,
\begin{equation}
    \Delta=\nabla_{X_{1}}^{*}\nabla_{X_{1}}+\ldots+\nabla_{X_{m}}^{*}\nabla_{X_{m}}^{2}.
\end{equation}

\subsection{The Kohn Laplacian on a CR manifold}  
In~\cite{KR:EHFBCM} Kohn-Rossi showed that the Dolbeault complex on a bounded complex domain  induces on 
its boundary a horizontal complex of differential forms. This was later extended by Kohn~\cite{Ko:BCM} to the general setting of a CR 
manifold $M^{2n+1}$ as follows.

Let $M^{2n+1}$ be a CR manifold with CR bundle $T_{1,0}\subset T_{\C}M$ and set $T_{0,1}=\overline{T_{1,0}}$. Then the subbundle   
$H=\Re (T_{1,0}\oplus T_{0,1})\subset TM$ admits an integrable complex structure and the splitting $H\otimes \C=T_{1,0}\oplus T_{0,1}$ 
gives rise to a decomposition  $\Lambda H^{*}\otimes \C = \oplus_{0\leq p,q\leq n } \Lambda^{p,q}$ into forms of bidegree $(p,q)$.

Assume that $T_{\C}M$ is endowed with a Hermitian metric such 
that $T_{1,0}$ and $T_{0,1}$ are orthogonal subspaces and complex conjugation is an (antilinear) isometry. This Hermitian metric gives rise to a 
Hermitian metric on $\Lambda^{*}T^{*}_{\C}M$ with respect to which the decomposition  $\Lambda H^{*}\otimes \C = \oplus_{0\leq p,q\leq n } 
\Lambda^{p,q}$ becomes orthogonal. Let $\Pi_{p,q}:\Lambda^{*}T_{\C}^{*}M \rightarrow \Lambda^{p,q}$ be the orthogonal projection onto $\Lambda^{p,q}$. 
Then the Kohn-Rossi operator $\bar \partial_{b}:  C^{\infty}(M,\Lambda^{p,q}) \rightarrow C^{\infty}(M,\Lambda^{p,q+1})$ is given by 
\begin{equation}
    \bar\partial_{b} \eta = \Pi_{p,q+1}(d\eta),\qquad  \eta \in C^{\infty}(M,\Lambda^{p,q}). 
\end{equation}
Since the integrability of $T_{1,0}$ implies that $\bar\partial_{b}^{2}=0$, this yields chain complexes 
$\bar\partial_{b}:C^{\infty}(M,\Lambda^{p,*}) \rightarrow C^{\infty}(M,\Lambda^{p,*+1})$. 

Endowing $M$ with a smooth density $\rho>0$ we let  $\bar\partial_{b} ^{*}$ denote 
the formal adjoint of $\bar\partial_{b}$. Then the Kohn Laplacian is   
\begin{equation}
    \Box_{b}=( \bar\partial_{b}+\bar\partial^{*}_{b})^{2}=\bar\partial^{*}_{b}\bar\partial_{b} +  \bar\partial_{b}\bar\partial^{*}_{b}.
\end{equation}

The Kohn Laplacian is a sublaplacian (see~\cite[Sect.~13]{FS:EDdbarbCAHG}, \cite[Sect.~20]{BG:CHM}), so is not elliptic. Nevertheless,  
Kohn~\cite{Ko:BCM} proved that 
under a geometric condition on the Levi form~(\ref{eq:Heisenberg.Kohn-Levi-form}), the so-called condition $Y(q)$, 
the operator $\Box_{b}$ acting on $(p,q)$-forms is hypoelliptic with loss of one derivative, 
i.e.~for any compact $K\subset M$ we have 
\begin{equation}
    \|u \|_{s+1}  \leq C_{Ks}( \|\Box_{b}u \|_{s}  + \|u\|_{0}) \qquad \forall C_{K}^{\infty}(M,\Lambda^{p,q}). 
     \label{eq:Operators.hypoellipticity}
\end{equation} 
where $\|.\|_{s}$ denotes the norm of the Sobolev space $L^{2}_{s}(M,\Lambda^{p,q})$. 

The condition $Y(q)$ at point $x\in M$ means that if we let $(r(x)-\kappa(x),\kappa(x),n-r(x))$ be the signature of $L_{\theta}$ at $x$, so that $r(x)$ is 
the rank of $L_{\theta}$ and $\kappa(x)$ the number of its negative eigenvalues, then we must have
\begin{equation}
    q\not \in \{\kappa(x),\kappa(x)+1,\ldots,\kappa(x)+n-r(x)\}\cup \{r(x)-\kappa(x),r(x)-\kappa(x)+1,\ldots,n-\kappa(x)\}.
     \label{eq:Operators.Y(q)-condition}
\end{equation}
Kohn's theorem then tells us that when the $Y(q)$-condition holds everywhere $\Box_{b}$ is hypoelliptic on $(p,q)$-forms. 

For instance,  when $M$ is $\kappa$-strictly pseudoconvex, the $Y(q)$-condition exactly means that we must have $q\neq \kappa$ and $q\neq n-\kappa$. In general 
this condition is equivalent to the existence of a parametrix within the Heisenberg calculus (\emph{cf.}~\cite{BG:CHM}, \cite{Po:BSM1}; 
see also~\cite{BdM:HODCRPDO}, \cite{FS:EDdbarbCAHG}), from which we recover the hypoellipticity of $\Box_{b}$.  

Finally, let us mention that the condition $Y(q)$ is a sufficient, but not a necessary condition for the hypoellipticity of the Kohn Laplacian 
(see, e.g.,~\cite{Ko:SEPDCRM}, \cite{Ko:OMSHETCROBL}, \cite{Ni:PhD}).

\subsection{The horizontal Laplacian on a Heisenberg manifold}
Let $(M^{d+1},H)$ be a Heisenberg manifold endowed with a Riemannian metric. Then the horizontal sublaplacian is the differential operator 
$\Delta_{b}:C^{\infty}(M,\Lambda^{*}_{\C}H^{*})\rightarrow C^{\infty}(M,\Lambda^{*}_{\C}H^{*})$ is given by 
\begin{equation}
    \Delta_{b}=d_{b}^{*}d_{b}+d_{b}d_{b}^{*}, \qquad d_{b}\alpha=\pi_{b}(d\alpha), 
\end{equation}
where $\pi_{b}:\Lambda^{*}_{\C}T^{*}M\rightarrow \Lambda^{*}_{\C}H^{*}$ denotes the orthogonal projection onto 
$\Lambda^{*}_{\C}H^{*}$. 

This operator was first introduced by Tanaka~\cite{Ta:DGSSPCM} in the CR setting, but versions of this operator acting on functions 
were independently defined by Greenleaf~\cite{Gr:FESPM} and Lee~\cite{Le:FMPHI}. Moreover, it can be shown that $d_{b}^{2}=0$ if, and only if, 
the subbundle $H$ is integrable, so in general $\Delta_{b}$ is not the Laplacian of a chain complex. 

 On functions $\Delta_{b}$ is a sum of squares modulo a lower order term, hence is hypoelliptic by H\"ormander's theorem. 
When acting on horizontal forms of higher degree, that is, on sections of $\Lambda^{k}_{\C}H^{*}$ with $k\geq 1$, it is was shown, 
in~\cite{Ta:DGSSPCM} and~\cite{Ru:FDVC} in the contact case and in~\cite{Po:BSM1} in the general case, that the operator $\Delta_{b}$ is hypoelliptic 
when the condition $X(k)$ holds everywhere. 

In the terminology of~\cite{Po:BSM1} the condition $X(k)$ at a point $x \in M$ means that, 
if we let $2r(x)$ denote the rank of the Levi form $\cL$ at $x$, then we must have 
\begin{equation}
    k\not \in\{r(x),r(x)+1,\ldots,d-r(x)\}.
    \label{eq:Operators.X(k)-condition}
\end{equation}

For instance, if $M^{2n+1}$ is a contact manifold or a nondegenerate CR manifold then the Levi form is everywhere nondegenerate, so $r(x)=2n$ and 
the $X(k)$-condition becomes $k\neq n$. 

Assume now that $M$ is a CR manifold with Heisenberg structure $H=\Re(T_{1,0}\oplus T_{0,1})$ and assume that 
$T_{\C}M$ is endowed with a Hermitian metric with respect to which $T_{1,0}$ and $T_{0,1}$ are orthogonal subspaces and complex conjugation is an isometry. 
Then  we have  $d_{b}=\dbarb +\partial_{b}$, where $\partial_{b}$ denotes the conjugate of $\dbarb$, that is, $\partial_{b}\alpha =\overline{\dbarb \bar \omega}$ 
for any $\omega\in C^{\infty}(M,\Lambda_{\C}^{*}H^{*})$. Moreover, as $\dbarb \partial_{b}^{*}+\partial_{b}^{*} \dbarb= \dbarb^{*} \partial_{b}+\partial_{b} 
 \dbarb^{*}=0$ (see~\cite{Ta:DGSSPCM}) we get 
 \begin{equation}
     \Delta_{b}=\Box_{b}+ \overline{\Box}_{b},
      \label{eq:Operators.Tanaka-Kohn}
 \end{equation}
 where $\overline{\Box}_{b}$ is the conjugate of $\Box_{b}$. In particular, we see that the horizontal sublaplacian preserves the bidegree, i.e., it 
 acts on $(p,q)$-forms. 
 
  In fact, as explained in~\cite[p.~34]{Po:BSM1}, one can show that the operator $\Delta_{b}$ acting on $(p,q)$-forms  is hypoelliptic 
  when at any point $x \in M$ the condition~$Y(p,q)$ 
 is satisfied. This means that at any point $x\in M$ we have  
  \begin{equation}
      (p,q)\not\in \{(\kappa(x)+j,r(x)-\kappa(x)+k);\ \max(j,k)\leq n-r(x)\} \cup  \{(r(x)-\kappa(x)+j,\kappa(x)+k);\ \max(j,k)\leq n-r(x)\},
  \end{equation}
  where $(r(x)-\kappa(x),\kappa(x),n-r(x))$ is the signature at $x$ of the Levi form $L_{\theta}$ associated to some non-vanishing real 1-form $\theta$ 
  anihilating $T_{1,0}\oplus T_{0,1}$. In particular, when $M$ is $\kappa$-strictly pseudoconvex the $Y(p,q)$ reduces to $(p,q)\neq (\kappa,n-\kappa)$ 
  and $(p,q)\neq (n-\kappa,\kappa)$.
 
 
\subsection{The conformal powers of the horizontal sublaplacian on a strictly pseudoconvex CR manifold}
Let $(M^{2n+1},\theta)$ be a strictly pseudoconvex pseudohermitian manifold. Let $T_{1,0}\subset T_{\C}M$ be the CR bundle of $M$  and define 
$T_{0,1}=\overline{T_{1,0}}$ and $H=\Re(T_{1,0}\oplus T_{0,1})$. 

According to Webster~\cite{We:PHSRH} a pseudohermitian structure on $M$ is given a choice 
of a 1-form $\theta$ annihilating $H$ such that the associated Levi form $L_{\theta}$ is positive definite on $T_{1,0}$, so that $\theta$ is a contact 
form with respect to which the complex structure of $H$ is calibrated, i.e.~we have $d\theta(X,JX)>0$ for any non-zero section of $H$. 

We extend the $L_{\theta}$ into a Hermitian metric $h_{\theta}$ on $T_{\C}M$ that $T_{1,0}$ and $T_{0,1}$, are orthogonal subspaces, 
complex conjugation is an (antilinear) isometry and $h_{\theta|_{H^{\perp}}}=\theta^{2}$. Then as shown by Tanaka~\cite{Ta:DGSSPCM} 
and Webster~\cite{We:PHSRH} there is a unique unitary connection on $T_{\C}M$ preserving the pseudohermitian structure. 
Note that the  contact form $\theta$ is unique up to a conformal change $f \rightarrow e^{2f} \theta$, $f\in C^{\infty}(M,\R)$. In order to study 
the analogue in CR geometry of the Yamabe problem  Jerison-Lee~\cite{JL:YPCRM} (see also~\cite{JL:ESIHGCRYP}, \cite{JL:ICRNCCRYP}) 
introduced a conformal version of the scalar horizontal sublaplacian acting on functions,
\begin{equation}
   \boxdot_{\theta} :C^{\infty}(M) \longrightarrow C^{\infty}(M), \qquad \boxdot_{\theta}=\Delta_{b}+ \frac{n}{n+2} R_{n},
\end{equation}
where $R_{n}$ denotes the scalar curvature of the Tanaka-Webster connection. This is a conformal operator in the sense that 
\begin{equation}
     \boxdot_{e^{2f}\theta} = e^{-(n+2)f} 
     \boxdot_{\theta} e^{nf}\qquad \forall f\in C^{\infty}(M, \R). 
\end{equation}

Assume now that $M$ is strictly pseudoconvex. Then the construction of Jerison-Lee 
has been recently generalized by Gover-Graham~\cite{GG:CRIPSL}. Assuming the existence of a $(n+2)$'th root for the canonical bundle 
$\Lambda^{n,1}$ their constructions yield for $k=1,2,\ldots,n+1$ a selfadjoint differential operator,  
 \begin{equation}
    \boxdot_{\theta}^{(k)}: C^{\infty}(M)\longrightarrow C^{\infty}(M), 
\end{equation}
such that the leading part (in the Heisenberg sense) of $\boxdot_{\theta}^{(k)}$ agrees with that of $\Delta_{b}^{k}$ and we have 
\begin{equation}
     \boxdot_{e^{2f}\theta}^{(k)}= e^{-(n+1+k)f}\boxdot_{\theta}^{(k)} 
    e^{(n+1-k)f}, \qquad \forall f\in C^{\infty}(M, \R).
\end{equation}
This is essentially the operator $P_{w,w'}$ of Theorem~1.1 of~\cite{GG:CRIPSL} with $w=w'=\frac{k-n-1}{2}$, noticing that when $w'=w$ the density bundles 
$\cE(w,w)$, on which the Gover-Graham's operators act, are in fact trivializable. 
For $k=1$ we recover the operator $\boxdot_{\theta}$ of Jerison-Lee. 

\subsection{Contact complex and contact Laplacian}
Given an orientable contact manifold $(M^{2n+1},\theta)$ the contact complex of Rumin~\cite{Ru:FDVC} can be seen as an attempt to get on $M$ 
a complex of horizontal differential forms by forcing out the equalities $d_{b}^{2}=0$ and $(d_{b}^{*})^{2}$ as follows.

Let $H=\ker \theta$ and assume that $H$ is endowed a calibrated almost complex structure $J\in \End_{\R}H$, $J^{2}=-1$, so that 
$d\theta(X,JX)=-d\theta(JX,X)>0$ for any section $X$ of $H$. We then get a Riemannian metric on $M$ by letting 
\begin{equation}
    g_{\theta}=d\theta(.,J.)+\theta^{2}.
\end{equation}

Let $T$ be the Reeb fields of $\theta$. Then we have 
\begin{equation}
    d_{b}^{2}=-\cL_{T}\varepsilon(d\theta)=\varepsilon(d\theta)\cL_{T},
\end{equation}
where $\varepsilon(d\theta)$ denotes the exterior multiplication by $d\theta$. 

There are two ways of modifying the space $\Lambda^{*}_{\C}H^{*}$ of horizontal forms to get a complex. The first one is to force the equality $d_{b}^{2}=0$ by
restricting the operator $d_{b}$ to $\Lambda^{*}_{2}:=\ker \varepsilon(d\theta) \cap \Lambda^{*}_{\C}H^{*}$ since this 
bundle is stable under $d_{b}$ and on there $d_{b}^{2}$ vanishes. 

The second way by similarly forcing the equality $(d_{b}^{*})^{2}$ via the restriction of $d_{b}^{*}$ to  
$\Lambda^{*}_{1}:=\ker 
\iota(d\theta)\cap \Lambda^{*}_{\C}H^{*}=(\im \varepsilon(d\theta))^{\perp}\cap \Lambda^{*}_{\C}H^{*}$, where $\iota(d\theta)$ denotes the interior product 
with $d\theta$. This means that we replace $d_{b}$ by the operator $\pi_{1}\circ d_{b}$, where $\pi_{1}$ is the orthogonal projection onto $\Lambda^{*}_{1}$.

On the other hand, since $d\theta$ is nondegenerate on $H$ the operator $\varepsilon(d\theta):\Lambda^{k}_{\C}H^{*}\rightarrow \Lambda^{k+2}_{\C}H^{*}$  is 
injective for $k\leq n-1$ and surjective for $k\geq n+1$. This implies that $\Lambda_{2}^{k}=0$ for $k\leq n$ and $\Lambda_{1}^{k}=0$ for $k\geq n+1$. 
Therefore, we only have two halves of complexes. As observed by Rumin we get a full complex by connecting 
the two halves by means of the (second order) differential operator, 
\begin{equation}
    D_{R}:\Lambda_{1}^{n} \longrightarrow\Lambda_{2}^{n+1}, \qquad 
    D_{R}=\cL_{T}-d_{b}\varepsilon(d\theta)^{-1}d_{b},
\end{equation}
where $\varepsilon(d\theta)^{-1}$ is the inverse of $\varepsilon(d\theta):\Lambda^{n-1}_{\C}H^{*}\rightarrow \Lambda^{n+1}_{\C}H^{*}$. In other words, letting 
$\Lambda^{k}=\Lambda^{k}_{1}$ for $k=0,\ldots,n-1$ and $\Lambda^{k}=\Lambda^{k}_{2}$ for $k=n+1,\ldots,2n$, we have the  chain complex:  
\begin{equation}
    C^{\infty}(M)\stackrel{d_{R}}{\rightarrow}
    \ldots \stackrel{d_{R}}{\rightarrow} C^{\infty}(M,\Lambda^{n}_{1})\stackrel{D_{R}}{\rightarrow} C^{\infty}(M,\Lambda^{n}_{2}) 
    \stackrel{d_{R}}{\rightarrow}
    \ldots \stackrel{d_{R}}{\rightarrow} C^{\infty}(M,\Lambda^{2n}),
     \label{eq:Operators.contact-complex}
\end{equation}
where $d_{R}:C^{\infty}(M,\Lambda^{k})\rightarrow C^{\infty}(M,\Lambda^{k+1})$ is equal to $\pi_{1}\circ d_{b}$ for $k=0,\ldots,n-1$ and to $d_{b}$ for 
$k=n+1,\ldots,2n$. This complex is called the contact complex of $M$. 

The contact Laplacian $\Delta_{R}:C^{\infty}(M,\Lambda^{*}\oplus \Lambda^{n}_{*})\rightarrow C^{\infty}(M,\Lambda^{*}\oplus \Lambda^{n}_{*})$ 
is given by the formulas,
\begin{equation}
    \Delta_{R}=\left\{
    \begin{array}{ll}
        (n-k)d_{R}d_{R}^{*}+(n-k+1) d_{R}^{*}d_{R}& \text{on $\Lambda^{k}$, $k=0,\ldots,n-1$},\\
        (d_{R}^{*}d_{R})^{2}+D_{R}^{*}D_{R} &  \text{on $\Lambda^{n}_{1}$},\\
        D_{R}D_{R}^{*}+  (d_{R}d_{R}^{*})& \text{on $\Lambda^{n}_{2}$},\\
         (n-k+1)d_{R}d_{R}^{*}+(n-k) d_{R}^{*}d_{R} & \text{on $\Lambda^{k}$, $k=n+1,\ldots,2n$}.
    \end{array}\right.
\end{equation}

Next, the almost complex structure $J$ of $H$ defines a bigrading on $\Lambda^{*}_{\C}H^{*}$. More precisely, we have an 
   orthogonal splitting $H\otimes \C=T_{1,0}\oplus T_{0,1}$ with $T_{1,0}=\ker(J+i)$ and $T_{0,1}=\overline{T_{1,0}}=\ker(J-i)$.
   Therefore, if we consider the subbundles $\Lambda^{1,0}=T_{1,0}^{*}$ and $\Lambda^{0,1}=T_{0,1}^{*}$ of $H^{*}\otimes \C\subset T^{*}_{\C}M$ 
   then we have the orthogonal decomposition $  \Lambda^{*}_{\C}H^{*}=\bigoplus_{0\leq p,q\leq n}\Lambda^{p,q}$ with 
   $\Lambda^{p,q}=(\Lambda^{1,0})^{p}\wedge (\Lambda^{0,1})^{q}$. We then get a bigrading on $\Lambda_{j}^{n}$, $j=1,2$, by letting 
 \begin{equation}
     \Lambda_{j}^{n}=\bigoplus_{p+q=n}\Lambda^{p,q}_{j}, \qquad \Lambda^{p,q}_{j}= \Lambda^{n}_{j}\cap \Lambda^{p,q}.
 \end{equation}

For $k=0,\ldots,2n$ with $k\neq n$ there exists a first order differential operator $P_{k}$ such that 
\begin{gather}
    (n-k+2)\Delta_{R}= (n-k)(n-k+1)\Delta_{b} +P_{k}^{*}P_{k}\qquad \text{on $\Lambda^{k}$, $k=0,..,n-1$},
    \label{eq:Operators.equalities-DeltaR-Deltab.1}\\
  (l-n+2)\Delta_{R}= (k-n)(k-n+1)\Delta_{b} +P_{l}^{*}P_{l} \qquad \text{on $\Lambda^{l}$,  $k=n+1,\ldots,2n$}. 
\end{gather}
For $j=1,2$ and $p+q=n$ there exist second order differential operators $P_{p,q}^{(j)}$ so that we have 
\begin{gather}
    4\Delta_{R}=\Delta_{b}^{2}+(P_{p,q}^{(j)})^{*}P_{p,q}^{(j)} \qquad \text{on $\Lambda^{p,q}_{\pm}$ with $\sup(p,q)\geq 1$},\\
    \Delta_{R}=(\Delta_{b}+iT)^{2}  \qquad \text{on $\Lambda^{n,0}_{j}$},\\
    \Delta_{R}=(\Delta_{b}-iT)^{2} \qquad \text{on $\Lambda^{0,n}_{j}$}.
    \label{eq:Operators.equalities-DeltaR-Deltab.5}
\end{gather}

These formulas enabled Rumin~\cite{Ru:FDVC} to show that $\Delta_{R}$ satisfied at every point a condition, the Rockland condition~\cite{Ro:HHGRTC}, 
which then allowed him to 
apply results of Helffer-Nourrigat~(\cite{HN:HMOPCV}) to show that $\Delta_{R}$ was maximal hypoelliptic. Alternatively, the facts that 
the Rockland condition is satisfied at every point and that the manifold is contact insure us that $\Delta_{R}$ admits a parametrix within the 
Heisenberg calculus, hence is hypoelliptic (see~\cite{Po:BSM1}).

\section{Heisenberg calculus}
\label{sec:PsiHDO}
In this section we gather the main facts about the Heisenberg calculus, following closely the expositions in~\cite{BG:CHM} and~\cite{Po:BSM1} (see 
also~\cite{EMM:HAITH}, \cite{Ta:NCMA}). 

 \subsection{Left-invariant pseudodifferential operators}
 Let $(M^{d+1},H)$ be a Heisenberg manifold and let $G=G_{a}M$ be the tangent Lie group of $M$ at a point $a \in M$. We review here the main results about 
left-invariant pseudodifferential operators on $G$ (see also~\cite{BG:CHM}, \cite{CGGP:POGD}, \cite{Ta:NCMA}). 
 
 Recall that for any finite dimensional  vector space $E$ the Schwartz class $\cS(E)$ is a Fr\'echet space and the Fourier transform is the continuous 
 isomorphism of $\cS(E)$ onto $\cS(E^{*})$ given by
 \begin{equation}
     \hat{f}(\xi)=\int_{E}e^{i\acou{\xi}{x}}f(x)dx, \qquad f \in \cS(E), \quad \xi \in E^{*},
 \end{equation}
 where $dx$ denotes the Lebesgue measure of $E$. 
 \begin{definition}
     $\cS_{0}(E)$ is the closed subspace of $\cS(E)$ consisting of $f \in \cS(E)$ such that for any differential operator $P$ on $E^{*}$ we have $(P\hat{f})(0)=0$. 
 \end{definition}
 
 Since $G$ has the same underlying set as that of its Lie algebra $\fg=\fg_{x}M$ we can let $\cS(G)$ and $\cS_{0}(G)$ denote the Fr\'echet spaces 
 $\cS(E)$ and $\cS_{0}(E)$ associated to the underlying linear space $E$ of $\fg$ 
 (notice that the Lebesgue measure of $E$ coincides with the Haar measure of $G$ since $G$ is nilpotent). 
 
 Next, for  $\lambda \in \R$ and $\xi=\xi_{0}+\xi'$ in $\fg^{*}= (T_{a}^{*}M/H^{*}_{a})\oplus H_{a}$ we let 
  \begin{equation}
     \lambda.\xi=  \lambda.(\xi_{0}+\xi')=\lambda^{2}\xi_{0}+\lambda \xi'.
      \label{eq:PsiHDO.Heisenberg-dilation-fg*}
 \end{equation}

 \begin{definition}
     $S_{m}(\fg^{*})$, $m \in \C$, is the space of functions $p \in C^{\infty}(\fg^{*}\setminus 0)$ which are homogeneous of degree $m$, in the sense 
     that, for any $\lambda>0$ we have 
\begin{equation}
         p(\lambda.\xi)=\lambda^{m}p(\xi) \qquad \forall \xi\in \fg^{*}\setminus 0. 
\end{equation}
    In addition $S_{m}(\fg^{*})$ is endowed with the Fr\'echet space topology induced from that of $C^{\infty}(\fg^{*}\setminus 0)$. 
 \end{definition}

 Note that the image $\hat{\cS}_{0}(G)$ of $\cS(G)$ under the Fourier transform consists of functions $v\in \cS(\fg^{*})$ such that, given any norm 
 $|.|$ on $G$,  near $\xi=0$ we have $|g(\xi)|=\op{O}(|\xi|^{N})$ for any integer $N\geq 0$. Thus, any $p \in S_{m}(\fg^{*})$ defines an element of
 $\hat{\cS}_{0}(\fg^{*})'$ by letting 
 \begin{equation}
     \acou{p}{g}= \int_{\fg^{*}}p(\xi)g(\xi)d\xi, \qquad g \in \hat{\cS}_{0}(\fg^{*}).
 \end{equation}
 This allows us to define the inverse Fourier transform of $p$ as the element $\check{p}\in \cS_{0}(G)'$ such that 
 \begin{equation}
     \acou{\check{p}}{f}=\acou{p}{\check{f}} \qquad \forall f \in \cS_{0}(G).
       \label{eq:PsiHDO.inverse-Fourier-transform-symbol}
 \end{equation}

\begin{proposition}[\cite{BG:CHM}, \cite{CGGP:POGD}]\label{prop:PsiHDO.convolution-symbols-group}
    1) For any $p \in S_{m}(\fg^{*})$ the left-convolution by $\check{p}$, 
\begin{equation}
   \check{p}*f(x):=\acou{\check{p}(y)}{f(x.y^{-1})}, \qquad f \in \cS_{0}(G), 
     \label{eq:PsiHDO.convolution-operator}
\end{equation}
defines a continuous endomorphism of $\cS_{0}(G)$.\smallskip 

  2) There is a continuous bilinear product,
 \begin{equation}
     *:S_{m_{1}}(\fg^{*})\times S_{m_{2}}(\fg^{*}) \longrightarrow S_{m_{1}+m_{2}}(\fg^{*}),
 \end{equation} 

such that, for any  $p_{1}\in S_{m_{1}}(\fg^{*})$ and $p_{2}\in S_{m_{2}}(\fg^{*})$, the composition of the left-convolution operators by 
$\check{p_{1}}$ and $\check{p_{2}}$ is the left-convolution operator by $(p_{1}*p_{2})^{\vee}$, that is, we have
\begin{equation}
    \check{p_{1}}*(\check{p_{2}}*f)=(p_{1}*p_{2})^{\vee}*f  \qquad \forall f \in \cS_{0}(G) .
\end{equation}
\end{proposition}

Let us also mention that if $p \in S_{m}(\fg^{*})$ then the convolution operator $Pu=\check{p}*f$ is a pseudodifferential operator. 
Indeed, let $X_{0}(a),\ldots,X_{d}(a)$ be a (linear) basis of $\fg$ so that $X_{0}(a)$ is in $T_{a}M/H_{a}$ and 
$X_{1}(a),\ldots,X_{d}(a)$ span $H_{a}$. For $j=0,\ldots,d$ let $X_{j}^{a}$ be the left-invariant vector fields on $G$ such that 
$X^{w}_{j|_{x=0}}=X_{j}(a)$. The basis 
$X_{0}(a),\ldots,X_{d}(a)$ yields a linear isomorphism $\fg\simeq \Rd$, hence a global chart of $G$. In this chart $p$ is a 
homogeneous symbol on $\Rd\setminus 0$ with respect to the dilations 
\begin{equation}
    \lambda.x=(\lambda^{2}x_{0},\lambda x_{1},\ldots,\lambda x_{d}), \qquad x\in \Rd, \quad \lambda>0.
     \label{eq:PsiHDO.Heisenberg-dilations-Rd}
\end{equation}

Similarly, each vector fields $\frac{1}{i}X_{j}^{a}$, $j=0,\ldots,d$, corresponds to a vector fields on $\Rd$ whose symbol is denoted 
$\sigma_{j}^{a}(x,\xi)$. Then, setting $\sigma=(\sigma_{0},\ldots,\sigma_{d})$, it can be shown that in the above chart $P$ is 
the operator 
\begin{equation}
    Pf(x)=\int_{\Rd} e^{ix.\xi}p(\sigma^{a}(x,\xi))\hat{f}(\xi), \qquad f \in \cS_{0}(\Rd).
     \label{eq:PsiHDO.PsiDO-convolution}
\end{equation}
In other words $P$ is the pseudodifferential operator $p(-iX^{a}):=p(\sigma^{a}(x,D))$ acting on $\cS_{0}(\Rd)$.

\subsection{$\mathbf{\Psi_{H}}$DO's on an open subset of $\Rd$} 

Let $U$ be an open subset of $\Rd$ together with a hyperplane bundle $H \subset TU$ and a $H$-frame $X_{0},X_{1},\ldots,X_{d}$ of $TU$. 
Then the class of \psivdos\ on $U$  is a class of pseudodifferential operators modelled on that of homogeneous convolution operators on the fibers of $GU$.

\begin{definition} $S_{m}(\URd)$, $m\in\C$, is the space of symbols 
    $p(x,\xi)\in C^{\infty}(U\times\Rdo)$ that are homogeneous of degree $m$ with respect to the $\xi$-variable, that is, 
    \begin{equation}
        p(x,\lambda.\xi)=\lambda^m p(x,\xi) \qquad \text{for any $\lambda>0$},
    \end{equation}
    where $\xi \rightarrow \lambda.\xi$ denotes the Heisenberg dilation~(\ref{eq:PsiHDO.Heisenberg-dilations-Rd}).
\end{definition}

Observe that the homogeneity of $p\in S_{m}(\URd)$ implies that, for any compact $K \subset U$, it satisfies estimates
\begin{equation}
     | \partial^\alpha_{x}\partial^\beta_{\xi}p(x,\xi)|\leq C_{K\alpha\beta}\|\xi\|^{\Re m-\brak\beta}, \qquad x\in K, \quad \xi \neq 0,
     \label{eq:PsiVDO.estimates-homogeneous-symbols}
\end{equation}
where $ \|\xi\|=(|\xi_{0}|^{2}+|\xi_{1}|^{4}+\ldots +|\xi_{d}|^{4})^{1/4}$ and $\brak\alpha = 2\alpha_{0}+\alpha_{1}+\ldots+ 
\alpha_{d}$. 

\begin{definition}$S^m(\URd)$,  $m\in\C$, consists of symbols  $p\in C^{\infty}(\URd)$ with
an asymptotic expansion $ p \sim \sum_{j\geq 0} p_{m-j}$, $p_{k}\in S_{k}(\URd)$,
in the sense that, for any integer $N$ and for any compact $K \subset U$, we have
\begin{equation}
    | \partial^\alpha_{x}\partial^\beta_{\xi}(p-\sum_{j<N}p_{m-j})(x,\xi)| \leq 
    C_{\alpha\beta NK}\|\xi\|^{\Re m-\brak\beta -N}, \qquad  x\in K, \quad \|\xi \| \geq 1.
    \label{eq:PsiVDO.asymptotic-expansion-symbols}
\end{equation}
\end{definition}

Next, for $j=0,\ldots,d$ let  $\sigma_{j}(x,\xi)$ denote the symbol of $\frac{1}{i}X_{j}$ (in the 
classical sense) and set  $\sigma=(\sigma_{0},\ldots,\sigma_{d})$. For any $p \in S^{m}(\URd)$ it can be shown that the symbol 
$p_{\sigma}(x,\xi):=p(x,\sigma(x,\xi))$ is in the H\"ormander class of symbols of type $(\frac{1}{2}, \frac{1}{2})$ (see~\cite[Prop.~10.22]{BG:CHM}). 
Therefore, we define a continuous linear operator from $C^{\infty}_{c}(U)$ to $C^{\infty}(U)$ by letting 
\begin{equation}
          p(x,-iX)f(x)= (2\pi)^{-(d+1)} \int e^{ix.\xi} p(x,\sigma(x,\xi))\hat{f}(\xi)d\xi,
    \qquad f\in C^{\infty}_{c}(U).
\end{equation}

In the sequel we let $\Psi^{-\infty}(U)$ denotes the class of smoothing operators, i.e.~of operators given by smooth kernels. 
    
\begin{definition}
   $\pvdo^{m}(U)$, $m\in \C$, consists of operators $P:C^{\infty}_{c}(U)\rightarrow C^{\infty}(U)$ of the form
\begin{equation}
         P= p(x,-iX)+R,
\end{equation}
 with $p$ in $S^{m}(\URd)$, called the symbol of $P$, and $R$ is smoothing.
\end{definition}

    The above definition of the symbol of $P$ differs from that of~\cite{BG:CHM}, since 
    there the authors defined it to be $p_{\sigma}(x,\xi)=p(x,\sigma(x,\xi))$. Note also that $p$ is unique 
    modulo $S^{-\infty}(\URd)$. 

\begin{lemma}\label{lem:PsiHDO.asymptotic-completeness}
   For $j=0,1,\ldots$ let $p_{m-j}\in S_{m-j}(\URd)$. Then there exists $P\in \pvdo^{m}(U)$ with symbol $p\sim \sum_{j\geq 0}p_{m-j}$. Moreover, the 
   operator $P$ is unique modulo smoothing operators. 
\end{lemma}
   
    The class $\pvdo^{m}(U)$ does not depend on the choice of the $H$-frame $X_{0}, \ldots, X_{d}$ (see~\cite[Prop.~10.46]{BG:CHM}). Moreover, 
since it is contained in the class 
of \psidos\ of type $(\frac{1}{2},\frac{1}{2})$ we get:

\begin{proposition}\label{prop:PsiHDO.Sobolev-regularity}
    Let $P$ be a \psivdo\ of order $m$ on $U$.\smallskip 
    
    1)  $P$ extends to a continuous mapping from 
    $\cE'(U)$ to $\cD'(U)$ and has a distribution kernel which is smooth off the diagonal. \smallskip 
    
    3) Let $k=\Re m$ if $\Re m\geq 0$ and $k=\frac{1}{2}\Re m$ otherwise. Then for any $s \in \R$ the operator $P\in \pvdo^{m}(U)$ extends
    to a continuous mapping from $L^{2}_{s,\op{comp}}(U)$ to $L^{2}_{s-k,\op{loc}}(U)$. 
\end{proposition}

\subsection{Composition of $\mathbf{\Psi_{H}}$DO's}
Recall that there is no symbolic calculus for \psidos\ of type $(\frac{1}{2},\frac{1}{2})$ since the product of two such \psidos\ needs not be 
again a \psido\ of type $(\frac{1}{2},\frac{1}{2})$. However, using  the fact that the \psivdos\ are modelled on left-invariant pseudodifferential operators 
allows us  to construct a symbolic calculus for \psivdos.

First, for $j=0,\ldots,d$ let $X_{j}^{(x)}$ be the leading homogeneous part of $X_{j}$ in privileged coordinates centered at $x$ 
defined according to~(\ref{eq:Heisenberg.X0u}) and~(\ref{eq:Heisenberg.Xju}). These vectors span a nilpotent Lie algebra of left-invariant vector fields
on a nilpotent graded Lie group $G^{x}$ which corresponds 
to $G_{x}U$ by pulling back the latter from the Heisenberg coordinates at $x$ to the privileged coordinates at $x$. 

As alluded to above the product law of $G^{(x)}$ defines a convolution product for symbols, 
\begin{equation}
    *^{(x)}: S_{m_{1}}(\Rd) \times S_{m_{2}}(\Rd) \longrightarrow S_{m_{1}+m_{2}}(\Rd).
     \label{eq:PsiHDO.convolution-symbol-pointwise}
\end{equation}
such that, with the notations of~(\ref{eq:PsiHDO.PsiDO-convolution}), on $\cL(\cS_{0}(\Rd))$ we have  
\begin{equation}
    p_{1}(-iX^{(x)})p_{2}(-iX^{(x)})=(p_{1}*^{(x)}p_{2})(-iX^{(x)}) \qquad \forall p_{j}\in S_{m_{j}}(\Rd).
\end{equation}

As it turns out the product $*^{(x)}$ depends smoothly on $x$ (see~\cite[Prop.~13.33]{BG:CHM}). Therefore, we get a continuous bilinear product,
\begin{gather}
    *: S_{m_{1}}(\URd) \times S_{m_{2}}(\URd) 
        \rightarrow S_{m_{1}+m_{2}} (\URd),\\   
        p_{1}*p_{2}(x,\xi)=(p_{1}(x,.)*^{(x)}p_{2}(x,.))(\xi), \qquad p_{j}\in S_{m_{j}}(\URd).
 \label{eq:PsiHDO.convolution-symbols-URd}
\end{gather}

\begin{proposition}[{\cite[Thm.~14.7]{BG:CHM}}] \label{prop:PsiHDO.composition}
    For $j=1,2$ let $P_{j}\in \pvdo^{m_{j}}(U)$ have  symbol $p_{j}\sim \sum_{k\geq 0} p_{j,m_{j}-k}$ and assume  
    that one of these operators is properly supported. Then the operator  
$P=P_{1}P_{2}$ is a \psivdo\ of order $m_{1}+m_{2}$ and has symbol  $p\sim \sum_{k\geq 0} p_{m_{1}+m_{2}-k}$, with  

\begin{equation}
     p_{m_{1}+m_{2}-k}(x,\xi) = \sum_{k_{1}+k_{2}\leq k} \sum_{\alpha,\beta,\gamma,\delta}^{(k-k_{1}-k_{2})}
            h_{\alpha\beta\gamma\delta} (x)  (D_{\xi}^\delta p_{1,m_{1}-k_{1}})* (\xi^\gamma 
            \partial_{x}^\alpha \partial_{\xi}^\beta p_{2,m_{2}-k_{2}})(x,\xi),    
\end{equation}
where $\underset{\alpha,\beta,\gamma,\delta}{\overset{\scriptstyle{(l)}}{\sum}}$ denotes the sum over the indices such that 
$|\beta|=|\gamma|$ and $|\alpha|+|\beta| \leq \brak\beta -\brak\gamma+\brak\delta = l$, and the functions 
$h_{\alpha\beta\gamma\delta}(x)$'s are  polynomials in the derivatives of the coefficients of 
the vector fields $X_{0},\ldots,X_{d}$.
\end{proposition}

\begin{remark}
It follows from~(\ref{eq:PsiHDO.convolution-symbols-URd}) that for any $x \in U$ the $x$-symbol $p_{1,m_{1}}*p_{2,m_{2}}(x,.)$ 
depends only on $p_{m_{1}}(x,.)$ and $p_{m_{2}}(x,.)$. However, the value of  $(p_{1,m_{1}}*p_{2,m_{2}})(x,\xi)$ at $(x,\xi)\in \URd$ depends on all 
the $\eta$-values $p_{m_{1}}(x,\eta)$ and $p_{m_{2}}(x,\eta)$ as $\eta$ ranges over $\Rd$. Thus we may localize the product of Heisenberg symbols with 
respect to $x$, but not respect to $(x,\xi)$, that is, the product of \psivdos\ is local, but is not microlocal. 
\end{remark}

\subsection{The distribution kernels of $\mathbf{\Psi_{H}}$DO's}
An important fact about \psidos\ is their characterization in terms of their distribution kernels.

  First, we extend the notion of homogeneity of functions to distributions. For $K\in \cS'(\Rd)$ and $\lambda >0$ we let $K_{\lambda}$ denote 
  the element of $\cS'(\Rd)$ such that

     \begin{equation}
           \acou{K_{\lambda}}{f}=\lambda^{-(d+2)} \acou{K(x)}{f(\lambda^{-1}.x)} \qquad \forall f\in\cS(\Rd). 
            \label{eq:PsiHDO.homogeneity-K-m}
      \end{equation}
In the sequel we will also use the notation $K(\lambda.x)$ for denoting $K_{\lambda}(x)$. We then say that $K$ is homogeneous of degree $m$, $m\in\C$, 
when $K_{\lambda}=\lambda^m K$ for any $\lambda>0$.

\begin{definition}\label{def:PsiHDO.regular-distributions}
  $\cS'_{\reg}(\Rd)$  consists of tempered distributions on $\Rd$ which are smooth outside the 
  origin. We equip it with the weakest  topology such that  the inclusions of $\cS'_{\reg}(\Rd)$  into $\cS'(\Rd)$ and $C^{\infty}(\Rdo)$ are continuous. 
\end{definition}

\begin{definition}
  $\cK_{m}(\URd)$, $m\in\C$, consists of distributions $K(x,y)$ in $C^\infty(U)\hotimes\cS'_{\reg}(\Rd)$ such that  for some functions $c_{\alpha}(x) \in 
   C^{\infty}(U)$,  $\brak\alpha=m$, we have  
    \begin{equation}
        K(x,\lambda.y)= \lambda^m K(x,y) + \lambda^m\log\lambda
                \sum_{\brak\alpha=m}c_{\alpha}(x)y^\alpha \qquad \text{for any $\lambda>0$}.
    \end{equation}
\end{definition}
 
The interest of considering the distribution class $\cK_{m}(\URd)$ stems from:

\begin{lemma}[{\cite[Prop.~15.24]{BG:CHM}},~{\cite[Lem.~I.4]{CM:LIFNCG}}]\label{prop:PsiHDO.Sm-Km}
   1) Any  $p\in S_{m}(\URd)$ agrees on $\URdo$ with  a distribution $\tau(x,\xi)$  in $C^{\infty}(U)\hotimes \cS'(\Rd)$ 
    such that $\check \tau_{\xiy}$ is in $\cK_{\hat m}(\URd)$, $\hat m=-(m+d+2)$. 
    
    2) If $K(x,y)$ is  in $\cK_{\hat m}(\URd)$ then the restriction of 
    $\hat{K}_{\yxi }(x,\xi)$ to $\URdo$ belongs to $S_{m}(\URd)$. 
\end{lemma}

This result is a consequence of the solution to the problem of extending a homogeneous function $p\in C^{\infty}(\Rd\setminus 0)$ into a 
homogeneous distribution on $\Rd$ and of the fact that for $\tau \in \cS'(\Rd)$ we have 
\begin{equation}
    (\hat{\tau})_{\lambda} = |\lambda|^{-(d+2)}(\tau_{\lambda^{-1}})^{\wedge}\qquad \forall \lambda \in \R\setminus 0. 
     \label{eq:PsiHDO.dilations-Fourier-transform}
\end{equation}
In particular, if $\tau$ is homogeneous of degree $m$ then $\hat{\tau}$ is homogeneous of $-(m+d+2)$.

The relevant class of kernels for the Heisenberg calculus is the following.
\begin{definition}\label{def:PsiHDO.K*}
$\cK^{m}(\URd)$, $m\in \C$, consists of distributions $K\in \cD'(\URd)$ with an asymptotic expansion 
     $K\sim \sum_{j\geq0}K_{m+j}$,  $K_{l}\in \cK_{l}(\URd)$,

 in the sense that,  for any integer $N$, as soon as $J$ is large enough we have 
   \begin{equation}
K-\sum_{j\leq J}K_{m+j}\in  C^{N}(\URd). 
     \label{eq:PsiHDO.asymptotics-kernel}
 \end{equation} 
\end{definition}

Since under the Fourier transform the asymptotic expansion~(\ref{eq:PsiVDO.asymptotic-expansion-symbols}) 
for symbols corresponds to that for distributions in~(\ref{eq:PsiHDO.asymptotics-kernel}),  
using Lemma~\ref{prop:PsiHDO.Sm-Km} we get: 

\begin{lemma}[{\cite[pp.~133--134]{BG:CHM}}]\label{lem:PsiHDO.characterization.Km}
    Let $K \in \cD'(\URd)$. Then the following are equivalent:\smallskip
    
    (i) The distribution $K$ belongs to $\cK^{m}(\URd)$;\smallskip
   
   (ii) We can put $K$ into the form
  \begin{equation}
      K(x,y)=\check{p}_{\xiy}(x,y)+R(x,y),
  \end{equation}
 for some $p\in S^{\hat{m}}(\URd)$, $\hat{m}=-(m+d+2)$, and some $R\in 
    C^{\infty}(\URd)$.\smallskip 
    
    Moreover, if (i) and (ii) holds and we expand $K \sim \sum_{j\geq 0} K_{m+j}$, $K_{l}\in \cK_{l}(\URd)$, 
    then we have $p\sim \sum_{j\geq 0}p_{\hat{m}-j}$ where $p_{\hat{m}-j} \in S_{\hat{m}-j}(\URd)$ is the restriction to $\URdo$ of 
    $(K_{m+j})^{\wedge}_{\yxi}$. 
\end{lemma}

Next, for $x \in U$ let $\psi_{x}$ denote the affine change to the privileged coordinates at $x$ and let us write  
$(A_{x}^{t})^{-1}\xi=\sigma(x,\xi)$ with  $A_{x}\in \op{GL}_{d+1}(\R)$. 
Since $\psi_{x}(x)=0$ and $\psi_{x*}X_{j}=\frac{\partial}{\partial 
y_{j}}$ at $y=0$ for $j=0,\ldots,d$,  one checks that $\psi_{x}(y)=A_{x}(y-x)$. 

Let $p\in S^{m}(\URd)$. As $p(x,-iX)=p^{\sigma}(x,D)$ with $p^{\sigma}(x,\xi)=p(x,\sigma(x,\xi))=p(x,(A_{x}^{t})^{-1}\xi)$ 
the distribution kernel 
$k_{p(x,-iX)}(x,y)$ of $p(x,-iX)$ is represented by the oscillating integrals
\begin{equation}
   (2\pi)^{-(d+1)} \int e^{i(x-y).\xi} p(x,(A_{x}^{t})^{-1}\xi)d\xi =   (2\pi)^{-(d+1)}|A_{x}|\int e^{iA_{x}(x-y).\xi} p(x,\xi)d\xi.
\end{equation}
Since $\psi_{x}(y)=A_{x}(y-x)$ we deduce that 
\begin{equation}
    k_{p(x,-iX)}(x,y)=|\psi_{x}'| \check{p}_{\xiy}(x,-\psi_{x}(y)).
    \label{eq:PsiHDO.kernel-quantization-symbol-psiy}
\end{equation}
Combining this with Lemma~\ref{lem:PsiHDO.characterization.Km} then gives: 

\begin{proposition}[{\cite[Thms.~15.39, 15.49]{BG:CHM}}]\label{prop:PsiVDO.characterisation-kernel1}
 Let $P:C_{c}^\infty(U)\rightarrow C^\infty(U)$ be a continuous linear operator with distribution kernel $k_{P}(x,y)$. Then the following are 
 equivalent:\smallskip  
 
 (i) $P$ is a \psivdo\ of order $m$, $m\in \C$.\smallskip 
 
 (ii) There exist $K\in \cK^{\hat{m}}(\URd)$, $\hat{m}=-(m+d+2)$, and $R \in C^{\infty}(U\times U)$ such that 
 \begin{equation}
     k_{P}(x,y)=|\psi_{x}'|K(x,-\psi_{x}(y)) +R(x,y) .
      \label{eq:PsiHDO.characterization-kernel.privilegedx}
 \end{equation}
 Furthermore, if (i) and (ii) hold and $K\sim \sum_{j\geq 0}K_{\hat{m}+j}$, 
$K_{l}\in \cK_{l}(\URd)$, then $P$ has symbol $p\sim \sum_{j\geq 0} p_{m-j}$, $S_{l}\in S_{l}(\URd)$, where 
$p_{\hat{m}-j}$ is the restriction to $\URdo$ of 
    $(K_{m+j})^{\wedge}_{\yxi}$. 
\end{proposition}

In the sequel we will need a version of Proposition~\ref{prop:PsiVDO.characterisation-kernel1} in Heisenberg coordinates.  
To this end let $\varepsilon_{x}$ denote the coordinate change to the Heisenberg coordinates at $x$ and set $\phi_{x}=\varepsilon_{x}\circ \psi_{x}^{-1}$. 
Recall that $\phi_{x}$ is  a Lie group isomorphism from $G^{(x)}$ to 
$G_{x}U$ such that $\phi_{x}(\lambda.y)=\lambda.\phi_{x}(y)$ for any $\lambda \in \R$. Moreover, using~(\ref{eq:Bundle.Extrinsic.Phiu}) 
one can check that $|\phi_{x}'|=1$ and $\phi_{x}^{-1}(y)=-\phi_{x}(-y)$. Therefore,  
from~(\ref{eq:PsiHDO.kernel-quantization-symbol-psiy}) we see that we can put $ k_{p(x,-iX)}(x,y)$ into the form
\begin{equation}
    k_{p(x,-iX)}(x,y)=|\varepsilon_{x}'| K_{P}(x,-\varepsilon_{x}(y)), \quad K_{P}(x,y)=  
    \check{p}_{\xiy}(x,-\phi_{x}(-y))=\check{p}_{\xiy}(x,\phi_{x}^{-1}(y)).
    \label{eq:PsiHDO.kernel-quantization-symbol}
\end{equation}

In fact, the coordinate changes $\phi_{x}$ give rise to an action on distributions on $\URd$ given by
   \begin{equation}
     K(x,y) \longrightarrow \phi_{x}^{*}K(x,y), \qquad \phi_{x}^{*}K(x,y)=K(x,\phi_{x}^{-1}(y)).
       \label{eq:PsiHD.isomorphism-cK}
   \end{equation}
  Since $\phi_{x}$ depends smoothly on $x$, this action induces a continuous linear isomorphisms of $C^{N}(\URd)$, $N \geq 0$, and $C^{\infty}(\URd)$ 
  onto themselves. As $\phi_{x}(y)$ is 
  polynomial in $y$ in such way that $\phi_{x}(0)=0$  and $\phi_{x}(\lambda.y)=\lambda.\phi_{x}(y)$ for every $\lambda \in \R$, we deduce that the above 
  action also yields a continuous linear isomorphism of $C^{\infty}(U)\hotimes \cS_{\reg}'(\Rd)$ onto itself and, for every $\lambda>0$,  we have 
   \begin{equation}
       (\phi_{x}^{*}K)(x,\lambda.y)= \phi_{x}^{*}[K(x,\lambda.y)], \qquad K \in \cD'(\URd).
        \label{eq:PsiHDO.homogeneity.phix*}
   \end{equation}

   Furthermore, as $\phi_{x}(y)$ is polynomial in $y$ we see that for every $\alpha \in \N^{d+1}$ we can write 
$\phi_{x}(y)^{\alpha}$ in the form $\phi_{x}(y)^{\alpha}=\sum_{\brak\beta=\brak \alpha}d_{\alpha\beta}(x)y^{\beta}$ with $d_{\alpha\beta}\in C^{\infty}(\URd)$.  
   It then follows that, for every $m \in \C$,  the map $K(x,y) \rightarrow \phi_{x}^{*}K(x,y)$ induces a linear isomorphisms of $\cK_{m}(\URd)$ and 
   $\cK^{m}(\URd)$ onto themselves. Combining this with~(\ref{eq:PsiHDO.kernel-quantization-symbol}) and 
   Proposition~\ref{prop:PsiVDO.characterisation-kernel1} then gives:

\begin{proposition}\label{prop:PsiVDO.characterisation-kernel2}
 Let $P:C_{c}^\infty(U)\rightarrow C^\infty(U)$ be a continuous linear operator with distribution kernel $k_{P}(x,y)$. Then the following are 
 equivalent:\smallskip  
 
 (i) $P$ is a \psivdo\ of order $m$, $m\in \C$.\smallskip 
 
 (ii) There exist $K_{P}\in \cK^{\hat{m}}(\URd)$, $\hat{m}=-(m+d+2)$, and $R \in C^{\infty}(U\times U)$ such that 
 \begin{equation}
     k_{P}(x,y)=|\varepsilon_{x}'|K_{P}(x,-\varepsilon_{x}(y)) +R(x,y) .
      \label{eq:PsiHDO.characterization-kernel.Heisenberg}
 \end{equation}
Furthermore, if (i) and (ii) hold and $K_{P}\sim \sum_{j\geq 0}K_{P,\hat{m}+j}$,  
$K_{l}\in \cK_{l}(\URd)$, then $P$ has symbol $p\sim \sum_{j\geq 0} p_{m-j}$, $S_{l}\in S_{l}(\URd)$, where 
$p_{\hat{m}-j}$ is the restriction to $\URdo$ of 
$[K_{P, \hat{m}+j}(x,\phi_{x}^{-1}(y))]^{\wedge}_{\yxi }$.
\end{proposition}

\begin{remark}\label{rem:PsiHDO.principal-symbol-at-x}
 Let $a\in U$. Then from~(\ref{eq:PsiHDO.characterization-kernel.Heisenberg}) 
 we see that the distribution kernel of $\tilde{P}=(\varepsilon_{a})_{*}P$ at $x=0$ is 
 \begin{equation}
     k_{\tilde{P}}(0,y)=|\varepsilon_{a}'|^{-1} k_{P}(\varepsilon_{a}^{-1}(0),\varepsilon_{a}^{-1}(y))=K_{P}(a,-y).
      \label{eq:PsiHDO.characterization-kernel.Heisenberg-origin}
 \end{equation}

 On the other hand, as we are in Heisenberg coordinates already, we have $\psi_{0}=\varepsilon_{0}=\phi_{0}=\op{id}$. Thus, in the 
 form~(\ref{eq:PsiHDO.characterization-kernel.Heisenberg}) for $\tilde{P}$ we have $K_{\tilde{P}}(0,y)=K_{P}(a,y)$. 
Therefore, if we let 
 $p_{m}(x,\xi)$ denote the principal symbol of $P$ and let  $K_{P,\hat{m}}\in \cK_{\hat{m}}(\URd)$ denote the leading kernel of $K_{P}$ 
 then by Proposition~\ref{prop:PsiVDO.characterisation-kernel2} we have 
 \begin{equation}
     p_{m}(0,\xi)=[K_{P,\hat{m}}]^{\wedge}_{\yxi }(a,\xi),
 \end{equation}
This shows that $[K_{P,\hat{m}}]^{\wedge}_{\yxi }(a,\xi)$ is the principal symbol of $P$ at $x=0$ in Heisenberg coordinates centered at $a$.
\end{remark}

\subsection{$\mathbf{\Psi_{H}}$DO's on a general Heisenberg manifold}
Let $(M^{d+1},H)$ be a Heisenberg manifold. As alluded to before the \psivdos\ on an subset of $\Rd$ are \psidos\ of type 
$(\frac{1}{2},\frac{1}{2})$. However, the latter don't make sense on a general manifold, for their class is not preserved by an arbitrary change of chart. Nevertheless, 
when dealing with \psivdos\ this issue is resolved if we restrict ourselves to changes of Heisenberg charts. Indeed, we have:

 \begin{proposition}[{\cite{BG:CHM}}, \cite{Po:BSM1}]\label{prop:PsiHDO.invariance}
     Let $U$ (resp.~$\tilde{U}$) be an open subset of $\Rd$ together with a hyperplane bundle $H\subset TU$ (resp.~$\tilde{H}\subset T\tilde{U}$) and a 
    $H$-frame of $TU$ 
    (resp.~a $\tilde{H}$-frame of $T\tilde{U}$). Let $\phi:(U,H)\rightarrow (\tilde{U},\tilde{H})$ be a Heisenberg diffeormorphism and let $\tilde{P}\in 
    \Psi_{\tilde{H}}^{m}(\tilde{U})$.\smallskip 
    
   1) The operator $P=\phi^{*}\tilde{P}$ is a \psivdo\ of order $m$ on $U$.\smallskip 
    
   2) If the distribution kernel of $\tilde{P}$ is of the form~(\ref{eq:PsiHDO.characterization-kernel.Heisenberg}) 
   with $K_{\tilde{P}}(\tilde{x},\tilde{y})\in \cK^{\hat{m}}(\tilde{U}\times \Rd)$ then the distribution kernel of 
   $P$ can be written in the form~(\ref{eq:PsiHDO.characterization-kernel.Heisenberg}) with $K_{P}(x,y)\in \cK^{\hat{m}}(\URd)$ such that 
   \begin{equation}
       K_{P}(x,y) \sim \sum_{\brak\beta\geq \frac{3}{2}\brak\alpha} \frac{1}{\alpha!\beta!} 
       a_{\alpha\beta}(x)y^{\beta}(\partial_{\tilde{y}}^{\alpha}K_{\tilde{P}})(\phi(x),\phi_{H}'(x)y),
        \label{eq:PsiHDO.asymptotic-expansion-KP}
   \end{equation}
   where we have let $a_{\alpha\beta}(x)=\partial^{\beta}_{y}[|\partial_{y}(\tilde{\varepsilon}_{\phi(x)}\circ \phi\circ \varepsilon_{x}^{-1})(y)| 
   (\tilde{\varepsilon}_{\phi(x)}\circ \phi\circ \varepsilon_{x}^{-1}(y)-\phi_{H}'(x)y)^{\alpha}]_{|_{y=0}}$ and $\tilde{\varepsilon}_{\tilde{x}}$ 
   denote the change to the Heisenberg coordinates at $\tilde{x}\in \tilde{f}$. In particular, we have 
   \begin{equation}
       K_{P}(x,y)=|\phi_{H}'(x)| K_{\tilde{P}}(\phi(x),\phi_{H}'(x)y) \qquad \bmod \cK^{\hat{m}+1}(\URd).
       \label{eq:PsiHDO.asymptotic-expansion-KP-principal}
   \end{equation}
\end{proposition}


 As a consequence of Proposition~\ref{prop:PsiHDO.invariance} we can define \psivdos\ on $M$  acting on the sections of 
a vector bundle $\cE$ over $M$.

\begin{definition}
    $\pvdo^{m}(M,\cE)$,  $m\in \C$,  consists of continuous operators $P:C^{\infty}_{c}(M,\cE) \rightarrow 
    C^{\infty}(M,\cE)$ such that:\smallskip 
    
    (i) The distribution kernel of $P$ is smooth off the diagonal;\smallskip 
    
    (ii) For any trivialization $\tau:\cE_{|_{U}}\rightarrow U\times \C^{r}$ over a 
     local Heisenberg chart $\kappa:U \rightarrow V\subset \Rd$ the operator $\kappa_{*}\tau_{*}(P_{|_{f}})$ belongs to
     $\pvdo^{m}(V, \C^{r}):=\pvdo^{m}(V)\otimes \End \C^{r}$. 
\end{definition}

All the aforementioned properties of \psidos\  on an open subset of $\Rd$ hold \emph{mutatis standis} in this setting. 

\subsection{Transposes and adjoints of $\mathbf{\Psi_{H}}$DO's}
Let us now look at the transpose and adjoints of \psivdos. First, given a Heisenberg chart $U \subset \Rd$ we have:

\begin{proposition}[\cite{BG:CHM}, \cite{Po:BSM1}]\label{prop:PsiHDO.transpose-chart}
    Let $P\in \pvdo^{m}(U)$. Then:\smallskip 
    
    1) The transpose operator $P^{t}$ is a \psivdo\ of order $m$ on $U$.\smallskip 
    
    2) If we write the distribution kernel of 
    $P$ in the form~(\ref{eq:PsiHDO.characterization-kernel.Heisenberg}) 
   with $K_{P}\in \cK^{\hat{m}}(\URd)$ then $P^{t}$ can be written in the form~(\ref{eq:PsiHDO.characterization-kernel.Heisenberg}) with 
   $K_{P^{t}}\in  \cK^{\hat{m}}(\URd)$ such that 
   \begin{equation}
       K_{P^{t}}(x,y) \sim  \sum_{\frac{3}{2}\brak\alpha \leq \brak \beta} \sum_{|\gamma|\leq |\delta| \leq 2|\gamma|} 
       a_{\alpha\beta\gamma\delta}(x) y^{\beta+\delta}  
       (\partial^{\gamma}_{x}\partial_{y}^{\alpha}K_{P})(x,-y), 
        \label{eq:PsiHDO.transpose-expansion-kernel}
   \end{equation}
   where 
   $a_{\alpha\beta\gamma\delta}(x)=\frac{|\varepsilon_{x}^{-1}|}{\alpha!\beta!\gamma!\delta!}
   [\partial_{y}^{\beta}(|\varepsilon_{\varepsilon_{x}^{-1}(-y)}'|(y-\varepsilon_{\varepsilon_{x}^{-1}(y)}(x))^{\alpha})
   \partial_{y}^{\delta}(\varepsilon_{x}^{-1}(-y)-x)^{\gamma}](x,0)$. In particular, 
   \begin{equation}
       K_{P^{t}}(x,y)=K_{P}(x,-y) \qquad \bmod \cK^{\hat m+1}(\URd). 
        \label{eq:PsiHDO.transpose-principal-kernels}
   \end{equation}
\end{proposition}


Thanks to this result we can prove: 

\begin{proposition}[\cite{BG:CHM}, \cite{Po:BSM1}]\label{prop:PsiHDO.transpose-adjoint}
  Let $P:C^{\infty}(M,\cE)\rightarrow C^{\infty}(M,\cE)$ be a \psivdo\ of order $m$. Then:
 
  1) The transpose operator $P^{t}:\cE'(M,\cE^{*})\rightarrow \cD'(M,\cE^{*})$ is a \psivdo\ of order $m$;\smallskip 
  
 2) If $M$ is endowed with a smooth positive density and $\cE$ with a Hermitian metric then the adjoint $P^{*}:C^{\infty}(M,\cE)\rightarrow 
 C^{\infty}(M,\cE)$ is a \psivdo\ of order $\bar{m}$. 
\end{proposition}

 \subsection{Principal symbol and model operators.} Let us now give an intrinsic definition of the principal symbol of  a \psivdo\ and of its model 
operator at a point. 

Let $\pi:\fg^{*}M \rightarrow M$ be the canonical projection of the bundle $\fg^{*}M$ onto $M$. 

\begin{definition}
    $S_{m}(\fg^{*}M,\cE)$, $m\in \C$, is the space of sections $p\in C^{\infty}(\fg^{*}M\setminus 0, \End \pi^{*}\cE)$ which are homogeneous of 
    degree $m$ in the sense that, for any $\lambda>0$, we have 
    \begin{equation}
        p(x,\lambda.\xi)=\lambda^{m}p(x,\xi) \qquad \forall (x,\xi) \in \fg^{*}M\setminus 0,
    \end{equation}
    where $\xi \rightarrow \lambda.\xi$ denotes the dilation~(\ref{eq:PsiHDO.Heisenberg-dilation-fg*}). 
\end{definition}

Let $P\in \pvdo^{m}(M)$ and  for $j=1,2$ let $\kappa_{j}$ be a Heisenberg chart with domain 
$V_{j}\subset M$. We let $\phi:U_{1}\rightarrow U_{2}$ be the corresponding transition map, with $U_{j}=\kappa_{j}(V_{1}\cap V_{2})\subset \Rd$, and 
for $j=1,2$ we define $P_{j}:=\kappa_{j*}(P_{|_{V_{1}\cap V_{2}}})$,
so that $P_{1}=\phi^{*} P_{2}$. 
Since $P_{j}$ belongs to $\pvdo^{m}(U_{j})$ its distribution kernel is of 
the form~(\ref{eq:PsiHDO.characterization-kernel.Heisenberg}) with $K_{P_{j}}\in \cK^{\hat m}(U_{j}\times \Rd)$.  
Let $K_{P_{j},\hat{m}}\in \cK_{\hat m}(U_{j}\times\Rd)$ be the leading kernel of $K_{P_{j}}$ and define 
\begin{equation}
    p_{j,m}(x,\xi)=[K_{P_{j},\hat{m}}]^{\wedge}_{\yxi }(x,\xi), \qquad (x,\xi)\in U_{j}\times \Rdo. 
     \label{eq:PsiHDO.global-principal-symbol}
\end{equation}

Recall that by  Remark~\ref{rem:PsiHDO.principal-symbol-at-x} for any $a \in U_{j}$ the symbol $p_{j}(a,.)$ 
yields in Heisenberg coordinates centered at $a$ the principal symbol of $P_{j}$ at $x=0$.  
Moreover, as shown in~\cite{Po:BSM1}, it follows from Proposition~\ref{prop:PsiHDO.invariance} that we have
\begin{equation}
    p_{1,m}(x,\xi)= p_{2,m}(\phi(x),[\phi_{H}'(x)]^{-1t}\xi).
\end{equation}
This shows that $p_{m}\!\!:=\kappa_{1}^{*}p_{1,m}$ is an element of $S_{m}(\fg^{*}(V_{1}\cap V_{1}))$ which is independent of the choice of the chart 
$\kappa_{1}$. Since $S_{m}(\fg^{*}M)$ is a sheaf this gives rise to a uniquely defined symbol $\sigma_{m}(P)\in S_{m}(\fg^{*}M)$. 
As we can similarly deal with \psivdos\ acting on $\cE$ we obtain:

%

\begin{proposition}[\cite{Po:BSM1}]\label{prop:PsiHDO.principal-symbol}
    For any $P \in \pvdo^{m}(M,\cE)$ there is a unique symbol $\sigma_{m}(P)\in S_{m}(\fg^{*}M, \cE)$ such that if in a local trivializing Heisenberg 
    chart $U\subset \Rd$ we let $K_{P,\hat{m}}(x,y)\in \cK_{\hat{m}}(\URd)$ be the leading kernel for the kernel $K_{P}(x,y)$  in the 
    form~(\ref{eq:PsiHDO.characterization-kernel.Heisenberg})  
    for $P$,  then we have
    \begin{equation}
       \sigma_{m}(P)(x,\xi)=[K_{P,\hat{m}}]^{\wedge}_{\yxi }(x,\xi), \qquad (x,\xi)\in U\times \Rdo. 
    \end{equation}
    
    Equivalently, for any $x_{0} \in M$ the symbol $\sigma_{m}(P)(x_{0},.)$ agrees in trivializing Heisenberg coordinates centered at $x_{0}$ with 
  the principal symbol of $P$ at $x=0$. 
\end{proposition}

\begin{definition}\label{def:PsiHDO.principal-symbol}
   For $P\in \pvdo^{m}(M,\cE)$ the symbol $\sigma_{m}(P)\in S_{m}(\fg^{*}M,\cE)$ provided by Proposition~\ref{prop:PsiHDO.principal-symbol}
     is called the principal symbol of $P$. 
\end{definition}

    Since we have two notions of principal symbol we shall distinguish between them by saying that  $\sigma_{m}(P)$ is the global 
    principal symbol of $P$ and that in a local trivializing chart the principal symbol $p_{m}$  of $P$ in the sense of~(\ref{eq:PsiVDO.asymptotic-expansion-symbols}) 
    is the local principal symbol of $P$ in this chart. 
    
    In fact, as shown in~\cite{Po:BSM1} the symbols $p_{m}(x,\xi)$ and $\sigma_{m}(P)(x,\xi)$ are related by
    \begin{gather}
         p_{m}(x,\xi)= (\hat{\phi}_{x}^{*}\sigma_{m}(P))(x,\xi), 
           \label{eq:PsiHDO.local-global-principal-symbol} \\
       (\hat{\phi}_{x}^{*}\sigma_{m}(P))= 
       [[\sigma_{m}(P)]^{\vee}_{\xiy}(x,\phi_{x}^{-1}(y))]^{\wedge}_{\yxi }=  [\phi_{x}^{*}[\sigma_{m}(P)]^{\vee}_{\xiy}]^{\wedge}_{\yxi },
       \label{eq:PsiHDO.isomorphism-symbols}
     \end{gather}
    where $\phi_{x}^{*}$ is the isomorphism map~(\ref{eq:PsiHD.isomorphism-cK}).
    In particular, since the latter is a linear isomorphism of $\cK_{m}(\URd)$ onto itself 
    the map  $p\rightarrow \hat{\phi}_{x}^{*}p$ is a linear isomorphism of $S_{m}(\URd)$ onto itself. 

 On the other hand, the principal map is surjective. More precisely, we have:
    \begin{proposition}[\cite{Po:BSM1}]\label{prop:PsiHDO.surjectivity-principal-symbol-map}
     For every $m\in \C$ the principal symbol map $\sigma_{m}: \pvdo^{m}(M,\cE) \rightarrow S_{m}(\fg^{*}M,\cE)$ gives rise to a linear isomorphism
    $\pvdo^{m}(M,\cE)/\pvdo^{m-1}(M,\cE)\stackrel{\sim}{\longrightarrow}  S_{m}(\fg^{*}M,\cE)$. 
 \end{proposition}

  Now, granted the above definition of the principal symbol, we can define the model operator at a point as follows. 
  
  \begin{definition}\label{def:PsiHDO.model-operator}
    Let $P\in \pvdo^{m}(M,\cE)$ have (global) principal symbol $\sigma_{m}(P)$.  Then the model operator of $P$ at $a\in M$ is the left-invariant 
    \psivdo-operator $P^{a}:\cS_{0}(G_{a}M,\cE_{a})\rightarrow 
    \cS_{0}(G_{a}M,\cE_{x})$ with symbol $\sigma_{m}(P)^{\vee}_{\xiy}(a,.)$, so that we have  
    \begin{equation}
        P^{a}f(x)=\acou{\sigma_{m}(P)^{\vee}_{\xiy}(a,y)}{f(x.y^{-1})}, \qquad f \in 
        \cS_{0}(G_{a}M,\cE_{a}).   
    \end{equation}
\end{definition}


For $a \in M$ we let $*^{a}:S_{m_{1}}(G_{a}M)\times S_{m_{2}}(G_{a}M)\rightarrow S_{m_{1}+m_{2}}(G_{a}M)$ be the 
convolution product for symbols defined by the product law of $G_{a}M$. Then it is proved in~\cite{Po:BSM1} that in a local trivializing chart we 
have  
%
%
%
%
%
%
\begin{equation}
    p_{m_{1}}*^{a}p_{m_{2}}=(\hat{\phi}_{a})_{*}[(\hat{\phi}_{a}^{*}p_{m_{1}})*^{(a)}(\hat{\phi}_{a}^{*}p_{m_{2}})] \qquad \forall p_{m_{j}}\in 
    S_{m_{j}}(\Rd),
    \label{eq:PsiHDO.global-local-convolution-symbols}
\end{equation}
where $(\hat{\phi}_{a})_{*}$ denotes the inverse of $\hat{\phi}_{a}^{*}$. Since $\hat{\phi}_{a}^{*}$, its inverse and $*^{(a)}$ depend smoothly on 
$a$, we deduce that that so does $*^{a}$.  This allows us to get:

\begin{proposition}[\cite{Po:BSM1}]
1)  The group laws on the fibers of $GM$ give rise to a convolution product,
    \begin{gather}
        *:S_{m_{1}}(\fg^{*}M,\cE)\times S_{m_{2}}(\fg^{*}M,\cE) \longrightarrow S_{m_{1}+m_{2}}(\fg^{*}M,\cE),
    \end{gather}
such that for symbols $p_{m_{j}}\in S_{m_{j}}(\fg^{*}M,\cE)$, $j=1,2$, we have
    \begin{gather}
        p_{m_{1}}*p_{m_{2}}(x,\xi)=[p_{m_{1}}(x,.)*^{x}p_{m_{2}}(x,.)](\xi), \qquad (x,\xi)\in \fg^{*}M\setminus 0,
    \end{gather}
 where $*^{x}$ denote the convolution product for symbols on $G_{x}M$.\smallskip
\end{proposition}

\begin{proposition}[\cite{Po:BSM1}]\label{prop:PsiHDO.composition2}
 For $j=1,2$ let $P_{j}\in \pvdo^{m_{j}}(M,\cE)$ and suppose that $P_{1}$ or $P_{2}$ 
    is properly supported. Then:\smallskip 
   
    1) We have $\sigma_{m_{1}+m_{2}}(P_{1}P_{2})=\sigma_{m_{1}}(P)*\sigma_{m_{2}}(P)$.\smallskip
    
    2) At every point $a\in M$ the model operator of $P_{1}P_{2}$ is $(P_{1}P_{2})^{a}=P^{a}_{1}P_{2}^{a}$.
\end{proposition}

Using Proposition~\ref{prop:PsiHDO.transpose-chart} we can also deal with the principal symbols and model operators of the adjoints and transposes of \psivdos. 

\begin{proposition}[\cite{BG:CHM}, \cite{Po:BSM1}]\label{prop:PsiHDO.transpose-global}
  Let $P \in \pvdo^{m}(M,\cE)$ have principal symbol $\sigma_{m}(P)$. Then:
 
 1) The principal symbol of the transpose $P^{t}$ is $\sigma_{m}(P^{t})(x,\xi)= \sigma_{m}(x,-\xi)^{t}$;\smallskip 
  
  2) If $P^{a}$ is the model operator of $P$ at $a$, then the model operator of $P^{t}$ at $a$ is the transpose operator 
  $(P^{a})^{t}: \cS_{0}(G_{x}M,\cE_{x}^{*})\rightarrow \cS_{0}(G_{x}M,\cE_{x}^{*})$.
\end{proposition}

Assume now that $M$ and $\cE$ are endowed with a positive density and a Hermitian metric respectively and let $L^{2}(M,\cE)$ be the associated 
$L^{2}$-Hilbert space. 

\begin{proposition}[\cite{Po:BSM1}]\label{prop:PsiHDO.adjoint-manifold}
    Let $P \in \pvdo^{m}(M,\cE)$ have principal symbol $\sigma_{m}(P)$. Then:\smallskip
   
 1) The principal symbol of the adjoint $P^{*}$ is $\sigma_{\bar{m}}(P^{*})(x,\xi)=\sigma_{m}(P)(x,\xi)^{*}$.\smallskip 
  
  2) If $P^{a}$ denotes the model operator of $P$ at $a\in M$ then the model operator of $P^{*}$ at $a$ is  
  the adjoint  $(P^{a})^{*}$ of $P^{a}$. 
 \end{proposition}

\subsection{Hypoellipticity, parametrices and Rockland condition}
Let $P:C^{\infty}_{c}(M,\cE) \rightarrow C^{\infty}(M,\cE)$ be a \psivdo\ of order $m$ such that $k:=\Re m>0$. First, we have:

\begin{proposition}[\cite{BG:CHM}, \cite{Po:BSM1}]\label{thm:PsiHDO.hypoellipticity}
   The following are equivalent:\smallskip 

   1) The principal symbol $\sigma_{m}(P)$ of $P$ is invertible with respect to the convolution product for homogeneous 
    symbols;\smallskip 
    
    2) The operator $P$ admits a parametrix $Q$ in $\pvdo^{-m}(M,\cE)$, i.e.~$PQ=QP=1  \bmod \psinf(M,\cE)$.\smallskip 
 
 \noindent Furthermore, if 1) and 2) hold then $P$ is hypoelliptic with loss of $\frac{k}{2}$-derivatives, i.e.,~for any $s\in \R$ and any compact 
 $K\subset M$  we have estimates
    \begin{equation}
        \| f\|_{L^{2}_{s+k/2}} \leq C_{Ks}(\|P f\|_{H^{s}}+\|f\|_{L^{2}_{s}}) \qquad \forall f \in C^{\infty}_{K}(M,\cE).  
         \label{eq:PsiHDO.subellipticity.subelliptic-estimates}
    \end{equation}
 \end{proposition}
    
When $M$ is compact, combining this with 
the compactness of the embedding of  $H^{k/2}(M,\cE)$ into $L^{2}(M,\cE)$ we get: 
\begin{proposition}\label{prop:PsiHDO.spectrum}
    Suppose $M$ compact and assume that $P$ has an invertible principal symbol and a 
    spectrum different from $\C$. Then: \smallskip 

    1) The spectrum of $P$ consists of isolated eigenvalues with finite multiplicities.\smallskip 

    2) For any $\lambda\in \op{Sp}P$ the eigenspace $\ker (P-\lambda)$ is a finite dimensional subspace of $C^{\infty}(M,\cE)$.  
\end{proposition}

Now, assume that $M$ is endowed with a positive density and $\cE$ with a Hermitian metric.  Let $P^{a}$ be the model operator of $P$ at 
a point $a\in M$ 
and let  $\pi: G\rightarrow \cH_{\pi}$  be a (nontrivial) unitary representation of $G=G_{a}M$. We define the symbol 
$\pi_{P^{a}}$ as follows (see also~\cite{Ro:HHGRTC}, \cite{Gl:SSGMCAINGHG}, \cite{CGGP:POGD}).

Let $\cH_{\pi}^{0}(\cE_{a})$ be the subspace of $\cH_{\pi}(\cE_{a}):=\cH_{\pi}\otimes \cE_{a}$ spanned by the vectors of the form
\begin{equation}
     \pi_{f}\xi=\int_{G}(\pi_{x}\otimes 1_{\cE_{a}})(\xi \otimes f(x)) dx, \qquad  f \in 
\cS_{0}(G), \quad  \eta \in \cH_{\pi}(\cE_{a}).
\end{equation}
with $f$ in $\cS_{0}(G,\cE_{a})=\cS_{0}(G)(\cE_{a})$ and $\xi \in \cH_{\pi}$. 
Then we let $\pi_{P^{a}}$ denote the (unbounded) operator of $\cH_{\pi}(\cE_{a})$ with domain 
$\cH_{\pi}^{0}(\cE_{a})$ such that
\begin{equation}
    \pi_{P^{a}}(\pi_{f}\xi)=\pi_{P^{a}f}\xi \qquad \forall f\in \cS_{0}(G,\cE_{a})\quad \forall \xi \in \cH_{\pi}.
\end{equation}

One can check that $\pi_{P^{a*}}$ is the adjoint of $\pi_{P^{a}}$ on $\cH_{\pi}^{0}$, hence 
is densely defined. Thus $\pi_{P^{a}}$ is closable and we can let $\overline{\pi_{P^{a}}}$ denotes its closure. 

In the sequel we let $C^{\infty}_{\pi}(\cE_{a})=C^{\infty}_{\pi}\otimes \cE_{a}$, where $C^{\infty}_{\pi}\subset \cH_{\pi}$ denotes the space of 
smooth vectors of $\pi$ (i.e.~the subspace of vectors $\xi \in \cH_{\pi}$ so that $x \rightarrow \pi(x)\xi$ is smooth from $G$ to $\cH_{\pi}$).

\begin{proposition}[\cite{CGGP:POGD}]\label{PsiHDO.properties-symbol-representation}
    1) The domain of $\overline{\pi_{P^{a}}}$ always contains $C^{\infty}_{\pi}(\cE_{a})$.\smallskip 
   
   2)  If $\Re m \leq 0$ then the operator $\overline{\pi_{P^{a}}}$ is bounded.\smallskip 
   
   3) We have $\overline{(\pi_{P^{a}})^{*}}=(\overline{\pi_{P^{a}}})^{*}$.\smallskip 
   
   4) If $P_{1}$ and $P_{2}$ are \psidos\ on $M$ then $\overline{\pi_{(P_{1}P_{2})^{a}}}=\overline{\pi_{P^{a}_{1}}}\, \overline{\pi_{P^{a}_{2}}}$.  
\end{proposition}

\begin{remark}
   If $\cE_{a}=\C$ and $P^{a}$ is a differentiable operator then, as it is left-invariant, $P^{a}$ belongs to the enveloping algebra 
$\cU(\fg)$ of the Lie algebra $\fg=\fg_{a}M$ of $G$. In this case $\overline{\pi_{P^{a}}}$ coincides on $C^{\infty}_{\pi}$ with the operator 
$d\pi(P^{a})$, where $d\pi$ is the representation of $\cU(\fg)$ induced by $\pi$.  
\end{remark}

\begin{definition}
    We say that $P$ satisfies the Rockland condition at $a$ if for any nontrivial unitary irreducible representation 
$\pi$ of $G_{a}M$ the operator $\overline{\pi_{P^{a}}}$ is injective on $C^{\infty}_{\pi}(\cE_{a})$.
\end{definition}

\begin{proposition}[\cite{Po:BSM1}]\label{thm:PsiHDO.Rockland-Parametrix}
    Suppose that the Levi form  of $(M,H)$ has constant rank. Then the following are equivalent:\smallskip 
    
    (i) $P$ and $P^{t}$ satisfy the Rockland condition at every point of $M$;\smallskip 
    
   (ii) $P$ and $P^{*}$ satisfy the Rockland condition at every point of $M$;\smallskip 
    
    (iii) The principal symbol of $P$ is invertible.\smallskip 
    
   \noindent  In particular,  if $P$ is selfadjoint then principal symbol of $P$ is invertible if, and only if, $P$ satisfies the Rockland condition 
   at every point of $M$. In any case, if the conditions (i), (ii), (iii) hold $P$ admits a parametrix in $\pvdo^{-m}(M,\cE)$ and is hypoelliptic in the sense 
   of~(\ref{eq:PsiHDO.subellipticity.subelliptic-estimates}).   
\end{proposition}
%
%
%
%

Finally, for a sublaplacian $\Delta:C^{\infty}(M,\cE) \rightarrow C^{\infty}(M,\cE)$ the Rockland condition can be reformulated as follows. Near a 
point $a \in M$ let $X_{0},X_{1},\ldots,X_{d}$ be a local $H$-frame of $TM$ with respect to which $\Delta$ takes the form,
 \begin{equation}
    \Delta=-(X_{1}^{2}+\ldots+X_{d}^{2})- i\mu(x) X_{0}+ \op{O}_{H}(1),
    \label{eq:Rockland.sublaplacian.preferred-frame}
\end{equation}
where  $\mu(x)$  is a local section of $\End \cE$ 
and the notation $\op{O}_{H}(1)$  means a differential operator of Heisenberg order~$\leq 1$.

Let $L(x)=(L_{jk}(x))$ be the matrix of 
$\cL$ with respect to the $H$-frame $X_{0},X_{1},\ldots,X_{d}$, so that for $j,k=1,\ldots,d$ we have 
\begin{equation}
    \cL(X_{j},X_{k})=[X_{j},X_{k}]=L_{jk}(x)X_{0} \quad \bmod H.
\end{equation}

Let $2n$ be the rank of $\cL_{a}$ and $L(a)$, let  $\lambda_{1},\ldots,\lambda_{d}$ denote the eigenvalues of $L(a)$ and 
consider the condition, 
 \begin{equation}
    \Sp \mu(a) \cap \Lambda_{a}=\emptyset,
     \label{eq:Rockland.sublaplacian}
 \end{equation}
 where the singular set $\Lambda_{a}$ is defined as follows, 
 \begin{gather} 
 \Lambda_{a}=     (-\infty, -\frac12 \Tra |L(a)|]\cup [\frac12 \Tra 
   |L(a)|,\infty) \qquad \text{if $2n<d$},\\
    \Lambda_{a}=\{\pm(\frac12 \Tra |L(a)|+\sum_{j=1}^{d}\alpha_{j}|\lambda_{j}|); \alpha_{j}\in \N^{d}\}\qquad \text{if $2n=d$}.
 \end{gather}

As is turns out this condition makes sense independently of the choice of the $H$-frame and is the relevant condition to look at in the case of 
 a sublaplacian. More precisely, we have: 

\begin{proposition}[{\cite[Thm.~18.4]{BG:CHM}}, \cite{Po:BSM1}]\label{thm:PsiDO.Rockland-sublaplacian}
    1) The condition~(\ref{eq:Rockland.sublaplacian}) makes sense intrinsically for any $a \in M$.\smallskip
    
    2) At every point $a\in M$ the Rockland conditions for $\Delta$ and $\Delta^{t}$ are equivalent to~(\ref{eq:Rockland.sublaplacian}).\smallskip
    
    3) The principal symbol of $\Delta$ is invertible if, and only if, the condition~(\ref{eq:Rockland.sublaplacian}) holds at every point of 
    $M$, so when the latter occurs $\Delta$ admits a parametrix in $\pvdo^{-2}(M,\cE)$ and is hypoelliptic with loss of one derivative.
\end{proposition}

\begin{remark}
 In the sequel we will also make use of the weaker condition, 
  \begin{gather}
    \Sp \mu(a) \cap \Lambda_{a}'=\emptyset, 
    \label{eq:Rockland.sublaplacian'}\\
     \Lambda_{a}'= (-\infty, -\frac12 \Tra |L(a)|]\cup [\frac12 \Tra  |L(a)|,\infty).
 \end{gather}
This condition was considered in~\cite{BGS:HECRM} and implies the existence of a parametrix for the heat operator associated to a sublaplacian in the 
Volterra-Heisenberg calculus of~\cite{BGS:HECRM} (see~Section~\ref{sec.powers1}).
\end{remark}

\subsection{Examples}
We now briefly explained how Proposition~\ref{thm:PsiHDO.Rockland-Parametrix} and Proposition~\ref{thm:PsiDO.Rockland-sublaplacian} 
can be used in the case of the examples described in Section~\ref{sec:Operators}. 

\subsubsection*{(a) Generalized sums of squares} If a sublaplacian $\Delta$ is a generalized sum of squares of the 
form~(\ref{eq:Operators.generalized-sum-of-squares}) 
then in~(\ref{eq:Rockland.sublaplacian.preferred-frame})
the matrix $\mu(x)$ is is diagonal with purely imaginary values, so 
that the condition~(\ref{eq:Rockland.sublaplacian}) is satisfied if, and only if, the Levi form does not vanish at 
$x$. Thus $\Delta$ has an invertible principal symbol 
if, and only if, the Levi form of $(M,H)$ is nonvanishing or, equivalently, the  
bracket condition  $H+[H,H]=TM$ is satisfied. In this case $\Delta$ admits a \psivdo\ parametrix and is hypoelliptic, so  
in the special case of a Heisenberg manifold we recover the hypoellipticity result of~\cite{Ho:HSODE}.

\subsubsection*{(b) Kohn Laplacian.} Let $\Box_{b}:C^{\infty}(M,\Lambda^{p,q})\rightarrow C^{\infty}(M,\Lambda^{p,q})$ 
be the Kohn Laplacian on a CR manifold $M^{2n+1}$ acting on $(p,q)$-forms, $0\leq p,q\leq n$. Then at point $x$ the 
condition~(\ref{eq:Rockland.sublaplacian}) for $\Box_{b}$ is equivalent to the $Y(q)$-condition
(see~\cite{BG:CHM},  \cite{Po:BSM1}). 
Thus $\Box_{b}$ has an 
invertible principal symbol  if, and only if, the condition $Y(q)$ holds everywhere. 
In this case $\Box_{b}$ admits a parametrix in $\pvdo^{-2}(M,\Lambda^{p,q})$ and we recover the theorem of Kohn~\cite{Ko:BCM} on the hypoellipticity of 
$\Box_{b}$.   

\subsubsection*{(c) Horizontal sublaplacian.} Let $\Delta_{b}:C^{\infty}(M,\Lambda^{k}_{\C}H^{*})\rightarrow C^{\infty}(M,\Lambda^{k}_{\C}H^{*})$  
be the horizontal sublaplacian on a Heisenberg manifold $(M^{d+1},H)$ acting on horizontal forms of degree $k$, $0\leq k\leq d$. Then at a point 
$x \in M$ the condition~(\ref{eq:Rockland.sublaplacian}) for $\Delta_{b}$ 
are equivalent to the condition~$X(k)$ (see~\cite{Po:BSM1}). Therefore, $\Delta_{b}$ has an invertible principal symbol if, and only 
if, the condition~$X(k)$  holds at every point. In particular, in case of a non-vanishing Levi form $\Delta_{b}$ has always an 
invertible principal when acting on functions. In any case, when the condition $X(k)$ holds everywhere $\Delta_{b}$ admits a  parametrix in 
$\pvdo^{-2}(M,\Lambda^{k}_{\C}H^{*})$ and is hypoelliptic with loss of one derivative.   

\subsubsection*{(d) Conformal powers of the horizontal sublaplacian.} Let $(M^{2n+1},\theta)$ be a strictly pseudoconvex CR manifold equipped with a 
pseudohermitian structure and for $k=1,\ldots,n$ let $\boxdot_{\theta}^{(k)}$ be a conformal power of $\Delta_{b}$ acting on functions as in~\cite{JL:YPCRM}  
and~\cite{GG:CRIPSL}. As on functions $\Delta_{b}$ satisfies the Rockland condition  at every point and has invertible principal symbol, the same is true for all 
its integer powers. As $\boxdot_{\theta}^{(k)}$  and $\Delta_{b}^{k}$ has same principal symbol, it follows that $ \boxdot_{\theta}^{(k)}$ satisfies the 
Rockland condition at every point and has an invertible principal symbol. In particular, the operator $\boxdot_{\theta}^{(k)}$ admits a parametrix 
in $\pvdo^{-2k}(M)$ and is hypoelliptic with loss of $k$ derivatives.

\subsubsection*{(e) Contact Laplacian} Let $(M^{2n+1},\theta)$ be an orientable contact manifold. For $k=0,1,\ldots,2n$ with $k\neq n$ 
let $\Delta_{R}:C^{\infty}(M,\Lambda^{k})\rightarrow C^{\infty}(M,\Lambda^{k})$  be the contact Laplacian on $M$ 
acting on contact forms of degree $k$. Then $\Delta_{b}$ satisfies the Rockland condition at every point  (see~\cite[p.~100]{Ru:FDVC}). Since in the contact 
case the Levi form of $(M,H)$ has constant rank $2n$, we then may apply Proposition~\ref{thm:PsiHDO.Rockland-Parametrix} to deduce that 
$\Delta_{R}$ has an invertible principal symbol, hence admits a parametrix in   $\pvdo^{-2}(M,\Lambda^{k})$ and is hypoelliptic with loss of one 
derivative. 

Likewise, the contact Laplacian $\Delta_{R}:C^{\infty}(M,\Lambda^{n}_{*})\rightarrow C^{\infty}(M,\Lambda^{n}_{*})$ acting on contact forms of degree $n$ 
satisfies the Rockland condition at every point and has an invertible principal symbol, so it admits a parametrix in $\pvdo^{-4}(M,\Lambda^{n}_{*})$ 
and is hypoelliptic with loss of two derivatives.

\section{Holomorphic families of $\mathbf{\Psi}$DOs}
\label{sec.HolPHDO}

In this section we define  holomorphic families of \psivdos\ and check their main properties. To this end we make use of 
an "almost homogeneous'' approach to the Heisenberg calculus, described in the first subsection.

\subsection{Almost homogeneous approach to the Heisenberg calculus} 
In this subsection we explain how the \psidos\ can be described in terms of symbols and 
kernels which are almost homogeneous, in the sense that there are homogeneous modulo infinite order terms.

\begin{definition}\label{def:HolPHDO.almost-homgenous-symbols}
A symbol $p(x,\xi)\in C^{\infty}(\URd)$ is almost homogeneous of degree $m$, $m \in \C$, when we have
\begin{equation}
     p(x,\lambda.\xi)-\lambda^m p(x,\xi)\in S^{-\infty}(\URd)\qquad \text{for any $\lambda>0$}.
\end{equation}
We let $S_{\ah}^m(\URd)$ denote the space of almost homogeneous symbols of degree $m$.
\end{definition}\begin{lemma}[{\cite[Prop.~12.72]{BG:CHM}}]\label{lem:HolPHDO.almost-homogeneity}
Let  $q(x,\xi) \in C^{\infty}(\URd)$. Then  we have equivalence: \smallskip 

(i) The symbol $q(x,\xi)$ is almost homogeneous of degree $m$;\smallskip 

(ii) The symbol $q(x,\xi)$ is in $S^{m}(\URd)$ and we have $q\sim 
p_{m}$ with $p_{m}\in S_{m}(\URd)$ (i.e.~the only nonzero homogeneous term in the asymptotic 
expansion~(\ref{eq:PsiVDO.asymptotic-expansion-symbols}) for $q$ is $p_{m}$). 
\end{lemma}
 
Granted this we shall now prove: 
\begin{lemma}\label{lem:HolPHDO.almost-homogeneous-characterization}
    Let $p \in C^{\infty}(\URd)$. Then we have equivalence:\smallskip
    
    (i) The function $p$ belongs to $S^{m}(\URd)$.\smallskip 
    
    (ii) For $j=0,1,..$ there exists $q_{m-j}\in S_{\ah}^{m-j}(\URd)$ such that $p\sim \sum_{j\geq 0} q_{m-j}$.
\end{lemma}
\begin{proof}
    Suppose that for $j=0,1,..$ there exists $q_{m-j}\in S_{\ah}^{m-j}(\URd)$ such that $p\sim \sum_{j\geq 0} q_{m-j}$.
    By Lemma~\ref{lem:HolPHDO.almost-homogeneity} there exists $p_{m-j}\in S_{m-j}(\URd)$ such that $q_{m-j}\sim p_{m-j}$. Then 
    we have $p\sim \sum_{j\geq 0}p_{m-j}$ and so  $p$ belongs to $S^{m}(\URd)$. Thus (ii) implies (i).
        
    Conversely, assume that $p$ belongs to $S^{m}(\URd)$ and let us write $p\sim \sum_{j\geq 0} 
    p_{m-j}$, $p_{l}\in S_{l}(\URd)$. Let $\varphi \in C^{\infty}_{c}(\Rd)$ be such that $\varphi(\xi)=1$ near $\xi=1$ and $\varphi(\xi)=0$ for 
    $\|\xi\|\leq 1$. For $j=0,1,..$ set 
    $q_{m-j}(x,\xi)=(1-\varphi(\xi))p_{m-j}(x,\xi)$. As for any $t>0$ the function $q_{m-j}(x,t.\xi)-t^{m-j}q_{m-j}(x,\xi)$ 
    is equal to $(\varphi(t.\xi)-\varphi(\xi))p_{m-j}(x,\xi)$, so belongs to $S^{-\infty}(\URd)$, we see that $q_{m-j}$ lies in 
    $S^{m-j}_{\ah}(\URd)$. Moreover, as we have $q_{m-j}(x,\xi)=p_{m-j}(x,\xi)$ for $\|\xi\|>1$ we see that $p\sim \sum_{j\geq 0} 
    \tilde{p}_{m-j}$. Hence (i) implies (ii). 
\end{proof}

The almost homogeneous symbols have been considered in~\cite[Sect.~12]{BG:CHM} already. In the sequel it will be important to have a 
"dual'' notion of almost homogeneity for distributions as follows. \begin{definition}
    The space $\cD'_{\reg}(\Rd)$ consists of the distributions on $\Rd$ that are smooth outside the 
    origin.  It is  endowed with
    the weakest topology that makes continuous the inclusions of  $\cD'_{\reg}(\Rd)$ into 
    $\cD'(\Rd)$ and  $C^{\infty}(\Rdo)$.
\end{definition}
\begin{definition}\label{def:HolPHDO.almost-homogeneous-kernels}
 A distribution $K(x,y)\in C^\infty(U)\hotimes \cD'_{\reg}(\Rd)$  is said to be almost homogeneous of degree $m$, $m \in \C$, when
    \begin{equation}
         K(x, \lambda.y)-\lambda^{m}K(x,y)\in C^{\infty}(\URd) \quad \text{for any  $\lambda>0$}. 
     \end{equation}
We let $\cK_{\ah}^m(\URd)$ denote the space of almost homogeneous distributions of degree $m$.
 \end{definition}

\begin{proposition}[Compare~{\cite[pp.~18-21]{Ta:NCMA}}]\label{prop:HolPHDO.characterization-almost-homogeneous-kernels}
 Let $K(x,y)\in C^\infty(U)\hotimes \cD'_{\reg}(\Rd)$ and set $\hat{m}=-(m+d+2)$. Then we have equivalence:\smallskip 
 
(i) The distribution $K$ belongs to $\cK_{\ah}^m(\URd)$.\smallskip

(ii) We can put $K(x,y)$ into the form,
\begin{equation}
    K(x,y)=K_{m}(x,y)+R(x,y),
\end{equation}
for some $K_{m}\in \cK_{m}(\URd)$ and $R\in C^{\infty}(\URd)$.\smallskip

(iii) We can put $K(x,y)$ into the form,
\begin{equation}
    K(x,y)=\check{p}_{\xiy}(x,y)+R(x,y),
\end{equation}
for some $p\in S_{\ah}^{\hat{m}}(\URd)$  and $R\in C^{\infty}(\URd)$. 
\end{proposition}
\begin{proof}
    First, if $K_{m}\in \cK_{m}(\URd)$ then~(\ref{eq:PsiHDO.homogeneity-K-m}) 
    implies that, for any  $\lambda>0$, the distribution $K(x, \lambda.y)-\lambda^{m}K(x,y)$ is in $C^{\infty}(\URd)$. 
    Thus (ii) implies~(i). 
    
    Second, let $p\in S_{\ah}^{\hat{m}}(\URd)$. By Lemma~\ref{lem:HolPHDO.almost-homogeneity} there is $p_{\hat{m}}\in S_{m}(\URd)$ such that 
    $p\sim p_{\hat{m}}$.  
    Thanks to Lemma~\ref{prop:PsiHDO.Sm-Km} we extend $p_{\hat{m}}$ into a distribution $\tau(x,\xi)$  in  
 $C^{\infty}(U) \hotimes \cS'(\Rd)$ such that $\check{\tau}_{\xiy}(x,y)$ is in $\cK_{m}(\URd)$. Let  
 $\varphi \in C^{\infty}_{c}(\Rd)$ be such that $\varphi =1$ near the origin. Then we can write 
 \begin{equation}
     p=\tau+\varphi(p-\tau) +(1-\varphi)(p-p_{\hat{m}}).
 \end{equation}
 Here $\varphi(\xi)(p(x,\xi)-\tau(x,\xi))$ belongs to $C^{\infty}(U)\otimes \cD'(\Rd)$ and is supported on a fixed compact set with respect to $\xi$, 
 so $[ \varphi(p-\tau)]^{\vee}_{\xiy}$ is smooth. 
 Moreover, as $p\sim p_{\hat{m}}$ both $(1-\varphi)(p-p_{m})$ and  $[(1-\varphi)(p-p_{m})]^{\vee}_{\xiy}$ are in $S^{-\infty}(\URd)$. It then follows  
 that $\check{p}_{\xiy}$ coincides with $\check{\tau}_{\xiy}$ up to an element of $C^{\infty}(\URd)$.  Since $\check{\tau}_{\xiy}(x,y)$ is in 
 $\cK_{m}(\URd)$ we deduce from this that (iii) implies~(ii).
 
 To complete the proof it remains to show that (i) implies (iii). Assume that $K(x,y)$ belongs to $\cK_{\ah}^m(\URd)$. 
 Let  $\varphi(y) \in C^{\infty}_{c}(\Rd)$ be  
 so that $\varphi(y) =1$ near $y=0$ and set $p=(\varphi K)^{\wedge}_{\yxi}(x,y)$. Then $p$ is a smooth and has slow growth with respect to $\xi$. 
 Moreover $\check{p}_{\xiy}(x,y)$ differs from 
 $K(x,y)$ by the smooth function $(1-\varphi(y))K(x,y)$. 
 
 Next, using~(\ref{eq:PsiHDO.dilations-Fourier-transform}) we see that, for any $\lambda>0$, the function $ p(x,\lambda.\xi)-\lambda^{\hat m}p(x,\xi)$ 
 is equal to
\begin{multline}
    \lambda^{-(d+2)}[\varphi(\lambda^{-1}.y)K(x,\lambda^{-1}.y)-\lambda^{-m}\varphi(y) K(x,y)]^{\wedge}_{\yxi}\\ 
    =  \lambda^{-(d+2)}[(\varphi(\lambda^{-1}.y)-\varphi(y))K(x,\lambda^{-1}.y)+
    \varphi(y) (K(x,\lambda^{-1}.y)-\lambda^{-m}K(x,y))]^{\wedge}_{\yxi}.
 \end{multline}
Note that $(\varphi(\lambda^{-1}.y)-\varphi(y))K(x,\lambda^{-1}.y)+
    \varphi(y) (K(x,\lambda^{-1}.y)-\lambda^{-m}K(x,y))$ belongs $C^{\infty}(\URd)$ and is  compactly supported with respect to 
$y$, so it belongs to $S^{-\infty}(\URd)$. Since the latter is also true for Fourier transform with respect to $y$ we see that 
 $ p(x,\lambda.\xi)-\lambda^{\hat 
m}p(x,\xi) $ is in $S^{-\infty}(\URd)$ for any $\lambda>0$, that is, the symbol $p$ is almost homogeneous of degree $\hat{m}$. It then follows that 
the distribution $K$ satisfies~(iii). This proves that (i) implies (iii). The proof is thus achieved.
\end{proof}
\subsection{Holomorphic families of $\mathbf{\Psi}_{H}$DO's} From now on we let $\Omega$ denote an open subset of $\C$. 
\begin{definition}\label{def:complex.symbols}
A family $(p_{z})_{z\in\Omega}\subset S^*(\URd)$ is said to be holomorphic when:  \smallskip  
     
     (i) The order $m(z)$ of $p_{z}$ depends analytically on $z$; \smallskip 

     (ii) For any $(x,\xi)\in \URd$ the function $z\rightarrow p_{z}(x,\xi)$ is holomorphic on $\Omega$; \smallskip 
     
     (iii) The bounds of the asymptotic expansion $p_{z} \sim \sum_{j\geq 0} p_{z, m(z)-j}$, 
         $p_{z,l}\in S_{l}(\URd)$, are locally uniform with respect to $z$, that is,~for 
         any integer $N$ and any compacts $K\subset U$ and $L\subset \Omega$ we have 
\begin{equation}
   | \partial_{x}^\alpha\partial_{\xi}^\beta (p_{z}-\sum_{j<N}  
            p_{z,m(z)-j})(x,\xi)| \leq C_{NKL\alpha\beta} \|\xi\|^{\Re m(z)-N-\brak\beta}, \quad x\in K, \ \|\xi\|\geq  1, \ z\in 
            L . 
    \label{eq:complex.symbols.asymptotic-expansion}
\end{equation}
We let $\Hol(\Omega,S^*(\URd))$ denote the set of the families $(p_{z})_{z\in\Omega}\subset S^*(\URd)$ that are holomorphic. 
\end{definition}
\begin{remark}\label{rem:HolPHDO.analyticity-homogeneous-symbols}
    If $(p_{z})_{z\in \Omega}$ is a holomorphic family of symbols then the homogeneous symbols 
    $p_{z,m(z)-j}$ depend   
    analytically on $z$. Indeed, for $\xi\neq 0$ we have  
    \begin{equation}
        p_{z,m(z)}(x,\xi)=\lim_{\lambda \rightarrow \infty} \lambda^{-m(z)}  
            p_{z}(x,\lambda.\xi).
    \end{equation}
     Since  the above axioms imply that the family $(\lambda^{-m(z)}  p_{z}(x,\lambda.\xi) )_{\lambda \geq 1}$ is bounded in the Fr\'echet-Montel space $\Hol(\Omega, 
C^\infty(\URdo)$ the convergence actually holds in $\Hol(\Omega, C^\infty(\URdo)$. Hence $p_{z,m(z)}$ depends 
analytically on $z$. Moreover, as for $\xi\neq 0$ we also have 
    \begin{equation}
        p_{j,z}(x,\xi)=\lim_{\lambda \rightarrow \infty} \lambda^{j-m(z)} 
                     (p_{z}(x,\lambda.\xi) -\sum_{l< j}\lambda^{m(z)-l}p_{z,m(z)-l}(x,\xi)),\quad j=1,2,\ldots,
    \end{equation}
an easy induction shows that all the symbols $p_{z,m_{z}-j}$ depend analytically on $z$.
\end{remark}Recall that $\psinf(U)=\cL(\cE'(U),C^{\infty}(U))$ is naturally a Fr\'echet space which is isomorphic to $C^{\infty}(U\times U)$ by Schwartz's kernel 
theorem. Therefore holomorphic families of smoothing operators makes well sense and 
we may define holomorphic families of \psivdos\ as follows.

\begin{definition}\label{def:HolPHDO.holomorphic-families}
 A family $(P_{z})_{z\in \Omega}\subset \pvdo^{m}(U)$ is holomorphic when it can be put into the form 
\begin{equation}
     P_{z} = p_{z}(x,-iX) + R_{z} \qquad z \in \Omega, 
\end{equation}
  for some family $(p_{z})_{z\in \Omega}\in \Hol(\Omega, S^{*}(\URd))$ and some family $(R_{z})_{z\in \Omega}\in \Hol(\Omega, \Psi^{-\infty}(U))$. 
  We let $\Hol(\Omega,\pvdo^{*}(U))$ denote the set of holomorphic families of \psivdos.
\end{definition}
For technical sake it will be useful to consider the symbol class below.

\begin{definition}
   $\Svb^k(\URd)$, $k \in \R$, consists of symbols $p(x,\xi)\in C^{\infty}(\URd)$ such that, for any 
   compact $K\subset U$,  we have 
   \begin{equation}
     |\partial_{x}^\alpha\partial^\beta_{\xi}p(x,\xi)| \leq 
    C_{K\alpha\beta}(1+\|\xi\|)^{k-\brak\beta}, \qquad (x,\xi)\in K\times\Rd.  
    \label{eq:HolPHDO.SvbU.estimates}
   \end{equation} 
 Its topology  is given by the sharpest constants $C_{K\alpha\beta}$'s 
 in the estimates~(\ref{eq:HolPHDO.SvbU.estimates}). 
\end{definition}

Note that the estimates~(\ref{eq:PsiVDO.asymptotic-expansion-symbols})  imply that $S^{m}(\URd)$, $m \in \C$, 
is contained in $\Svb^k(\URd)$ for any $k \geq \Re m$.  

\begin{proposition}\label{prop:Complexpowers.operators.properties}
    Let $(P_{z})_{z\in \Omega}$ be a holomorphic family of \psivdo's. Then: \smallskip 
    
    1) $(P_{z})_{z\in \Omega}$ gives rise to families in 
        $\Hol(\Omega,\cL(C_{c} ^\infty(U),C^\infty(U)))$ and $\Hol(\Omega,\cL(\cE'(U)),\cD'(U)))$. \smallskip 
	
        2) Off the diagonal of $U\times U$ the distribution kernel of $P_{z}$ is represented by a holomorphic family of smooth  functions. 
    \end{proposition}
\begin{proof}
    Without any loss of generality we may suppose that $P_{z}=p_{z}(x,-iX)$, with $(p_{z})_{z\in \Omega}$ in 
    $\Hol(\Omega, S^{*}(\URd))$. Moreover, shrinking $\Omega$ if necessary,  we may also assume that the degree $m_{z}$ of 
    $p_{z}$ stays bounded, as much so   $(p_{z})_{z\in \Omega}$ is contained in $\Svb^k(\URd)$ for some real $k\geq 0$.  
Let $\sigma_{j}(x,\xi)$ denote the 
classical symbol of $-iX_{j}$ and set $\sigma=(\sigma_{0},\ldots,\sigma_{d})$. 
Then the proof of~\cite[Prop.~10.22]{BG:CHM} shows that the map $p(x,\xi) \rightarrow p^{\sigma}(x,\xi):=p(x,\sigma(x,\xi))$ is 
continuous from 
$\Svb^{k}(\URd)$ to $S^{k}_{\frac{1}{2},\frac{1}{2}}(\URd)$. Thus, the family $(p^{\sigma}_{z})_{z \in \Omega}$ 
belongs to $\Hol(\Omega,S^{k}_{\frac{1}{2},\frac{1}{2}}(\URd))$.

Next, it follows from the proof of~\cite[Thm.~2.2]{Ho:PDOHE} that:\smallskip   

(i) The  quantization map $q \rightarrow q(x,D)$ induces continuous $\C$-linear maps from $S^{k}_{\frac{1}{2},\frac{1}{2}}(\URd)$ 
 to $\cL(C_{c} ^\infty(U),C^\infty(U))$ and to $\cL(\cE'(U)),\cD'(U))$; 
\smallskip 

(ii) The linear map $q(x,\xi) \rightarrow \check{q}_{\xiy}(x,y)$ is continuous from $S^{k}_{\frac{1}{2},\frac{1}{2}}(\URd)$ to $C^{\infty}(U)\otimes 
\cD'_{\reg}(\Rd)$, in such way that for any $q\in S^{k}_{\frac{1}{2},\frac{1}{2}}(\URd)$ the distribution kernel $\check{q}_{\xiy}(x,x-y)$ 
of $q(x,D)$ is represented off the diagonal by a smooth function depending continuously on $q$.\smallskip 

As  a continuous $\C$-linear map is  analytic it follows that, on the 
one hand,  $(p_{z}(x,-iX))_{z\in \Omega}$ gives rise to  elements of $\Hol(\Omega,\cL(C_{c} ^\infty(U),C^\infty(U)))$ and 
$\Hol(\Omega,\cL(\cE'(U)),\cD'(U)))$ and, 
 on the other hand,  the distribution kernel of $P_{z}$ is represented outside the diagonal of $U\times U$ by a holomorphic family of smooth  functions.
\end{proof}

\begin{definition}
    Let $(P_{z})_{z \in \Omega}\subset \cL(C^{\infty}_{c}(U),C^{\infty}(U))$ and for $z \in \Omega$ let $k_{P_{z}}(x,y)$ denote the distribution kernel of 
    $P_{z}$. Then the family $(P_{z})_{z \in \Omega}$ is said to be uniformly properly supported when, for any 
compact $K\subset U$ there exist compacts $L_{1}\subset U$ and $L_{2}\subset K$ such that for any $z\in \Omega$ we have
\begin{equation}
    \supp k_{P_{z}}(x,y)\cap (U\times K) \subset L_{1} \qquad \text{and} \qquad \supp k_{P_{z}}(x,y)\cap (K\times U) \subset L_{2}.
     \label{eq:HolPHDO.uniformly-properly-supported}
\end{equation}
\end{definition}

\begin{proposition}\label{prop:Complexpowers.operators.properties2}
Let $(P_{z})_{z\in \Omega}$ be a holomorphic family of \psivdos.  \smallskip 
                   
 1) We can write $P_{z}$ in the form $P_{z}=Q_{z}+R_{z}$ with  $(Q_{z})_{z\in \Omega}\in \Hol(\Omega, \pvdo^{*}(U)$ uniformly properly supported 
 and $(R_{z})_{z\in \Omega}\in \Hol(\Omega, \Psi^{-\infty}(U)$.\smallskip 

2) If the family $(P_{z})_{z \in \Omega}$ is uniformly properly supported then it gives rise to 
 holomorphic families of continuous endomorphisms of $C_{c}^\infty(U)$ and $C^\infty(U)$ and of $\cE'(U)$ and $\cD'(U)$.
\end{proposition}
\begin{proof}
    Let $(\varphi_{i})_{i \geq 0} \subset C^{\infty}_{c}(U)$ be a partition of unity which is subordinated to a locally finite covering 
    $(U_{i})_{i \geq 0}$ of $U$ by 
    relatively compact open subsets. For each $i\geq 0$ let $\psi_{i}\subset C^{\infty}_{c}(U)_{i}$ be such that 
    $\psi_{i}=1$ near $\supp \varphi_{i}$ and set $\chi(x,y)=\sum \varphi_{i}(x)\psi_{i}(y)$. Then $\chi$ is a smooth function on $U\times U$ which is 
    properly supported and such that $\chi(x,y)=1$ near the diagonal of $U\times U$. 
    
    For $z \in \Omega$ let $k_{P_{z}}(x,y)$ denote the distribution kernel of $P_{z}$ and let $Q_{z}$ and $R_{z}$ be the elements of 
    $\cL(C^{\infty}_{c}(U), C^{\infty}(U))$ with respective distribution kernels
     \begin{equation}
        k_{Q_{z}}(x,y)=\chi(x,y)k_{P_{z}}(x,y) \quad \text{and} \quad k_{R_{z}}(x,y)=(1-\chi(x,y))k_{P_{z}}(x,y).
    \end{equation}
    Notice that since $\chi$ is properly supported the family $(Q_{z})_{z\in \Omega}$ is uniformly properly supported.
    As  by Proposition~\ref{prop:Complexpowers.operators.properties} the distribution $k_{P_{z}}(x,y)$ is represented outside the 
        diagonal of $U\times U$ by a holomorphic family of smooth  functions, we see that $(k_{R_{z}}(x,y))$ is a holomorphic family of smooth kernels, 
        i.e.~$(R_{z})_{z\in \Omega}$ is a holomorphic family of smoothing operators. Since $Q_{z}=P_{z}-R_{z}$ it follows that $(Q_{z})_{z\in \Omega}$ 
        is a holomorphic family of \psivdos. Hence the first assertion.
        
       Assume now that $(P_{z})_{z\in \Omega}$ is uniformly properly supported. Thanks to Proposition~\ref{prop:Complexpowers.operators.properties}
        we already know that  $(P_{z})_{z \in \Omega}$ gives rise to holomorphic families with 
        values in $\cL(C^{\infty}_{c}(U),C^{\infty}(U))$ and $\cL(\cE'(U),\cD'(U))$. Let $K$ be a compact subset of $U$. 
        Then~(\ref{eq:HolPHDO.uniformly-properly-supported})  implies that  
        there exists a compact $L\subset U$ such that for every $z \in \Omega$ the operator $P_{z}$ maps 
        $C_{K}^{\infty}(U)$ to $C^{\infty}_{L}(U)$ and $\cE'_{K}(U)$ to $\cE'_{L}(U)$, in such way that  
        $(P_{z})_{z\in \Omega}$ gives rise to holomorphic families with values in $\cL(C_{K}^{\infty}(U),C^{\infty}_{L}(U))$ and $\cL(\cE'_{K}(U),\cE'_{L}(U))$. 
        In view of the definitions of the topologies of $C^{\infty}_{c}(U)$ and $\cE'(U)$ as the inductive limit topologies of $C^{\infty}_{K}(U)$ and 
        $\cE_{K}(U)$ as $K$ ranges over compacts of $U$, this shows that the family 
        $(P_{z})_{z\in \Omega}$ gives rise to elements of $\Hol(\Omega,C_{c}^\infty(U))$ and $\Hol(\Omega, \cE'(U))$.
       
        Next, let $(\varphi_{i})_{i\geq 0}\subset C^{\infty}_{c}(U)$ be a partition of unity. For each index $ i$ let $K_{i}$ be a compact 
        neighborhood of $\supp \varphi_{i}$. Then~(\ref{eq:HolPHDO.uniformly-properly-supported}) 
        implies that there exists a compact $L_{i}\subset U$ such that $\supp 
        k_{P_{z}}(x,y)\cap(K_{i}\times U)\subset L_{i}$ for every $z \in \Omega$. Let $\psi_{i}\in C^{\infty}(U)$ be such that $\psi_{i}=1$ near 
        $K_{i}$. Then we have 
        \begin{equation}
            P_{z}=\sum_{i\geq 0}\varphi_{i}P_{z}=\sum_{i\geq 0}\varphi_{i}P_{z}\psi_{i}.
             \label{eq:HolPHDO.uniformly-properly-supported-locally-finite-sum}
        \end{equation}
        Since each family $(\varphi_{i}P_{z}\psi_{i})_{z\in \Omega}$ is holomorphic with values in $\cL(C^{\infty}(U))$ and $\cL(\cD'(U))$ and the 
        sums are locally finite this shows that $(P_{z})_{z\in \Omega}$ gives rise to elements of $\Hol(\Omega,C^\infty(U))$ and $\Hol(\Omega,\cD'(U))$. 
\end{proof}

\subsection{Composition of holomorphic families of \psivdos}
Let us now look at the analyticity of the composition of \psivdos. To this 
end we need to deal with holomorphic families of almost homogeneous symbols as follows. 

\begin{definition}
    A family $(q_{z})_{z\in \Omega}\subset S_{\ah}^{*}(\URd)$ is holomorphic when:\smallskip
    
    (i) The degree $m(z)$ of $q_{z}$ is a holomorphic function on $\Omega$;\smallskip
    
    (ii) The family $(q_{z})_{z\in \Omega}$ belongs to $\Hol(\Omega,  C^{\infty}(\URd))$;\smallskip 
    
    (iii) For any $t>0$ the family $(q_{z}(x,t.\xi) -t^{m(z)}q_{z}(x,\xi))_{z\in \Omega}$ is in $\Hol(\Omega, S^{-\infty}(\URd))$.\smallskip
    
   \noindent We let $\Hol(\Omega,  S_{\ah}^{*}(\URd))$ denote the set of holomorphic families of almost homogeneous symbols.
\end{definition}\begin{lemma}\label{lem:HolPHDO.characterization-almost-homogeneous-symbols}
 Let $(q_{z})_{z\in\Omega}\in \Hol(\Omega, C^{\infty}(\URd))$. Then  we have equivalence:\smallskip 
 
 (i) The family $(q_{z})_{z\in\Omega}$ belongs to $ \Hol(\Omega, S_{\ah}^{*}(\URd))$ and has degree 
  degree $m(z)$;\smallskip

 (ii) The family $(q_{z})_{z\in \Omega}$ lies in $\Hol(\Omega, S^{*}(\URd))$ and, in the sense of~(\ref{eq:complex.symbols.asymptotic-expansion}), 
 we have $q_{z} \sim p_{z}$ where, for 
 every $z \in \Omega$, the symbol $p_{z}$ belongs to $S_{m(z)}(\URd)$.
\end{lemma}
\begin{proof}
    Assume that  $(q_{z})_{z\in \Omega}$ is in $\Hol(\Omega, S^{*}(\URd))$ and, in the sense of~(\ref{eq:complex.symbols.asymptotic-expansion}), 
      we have $q_{z} \sim p_{z}$ where, for 
 every $z \in \Omega$, the symbol $p_{z}$ belongs to $S_{m(z)}(\URd)$. Then the order $m(z)$ of $q_{z}$ is a holomorphic function on $\Omega$ and,
for any compact subset $K \subset U$, any integer $N$ and any open $\Omega' \subsubset \Omega$, we have
     \begin{equation}
          | \partial_{x}^\alpha\partial^\beta_{\xi}(q_{z}-p_{z})(x,\xi)\leq C_{NK\Omega'\alpha\beta}\|\xi\|^{-N}, 
         \label{eq:HolPHDO.asymptotics-symbol-single-term}
     \end{equation}
for $x\in K$, $\|\xi\|\geq 1$ and $z \in \Omega'$. It follows that, for any $t>0$, the family 
$\{q_{z}(x,t.\xi)-t^{m(z)}q_{z}(x,\xi)\}_{z\in \Omega}=\{(q_{z}(x,t.\xi)-q_{z}(x,t.\xi))-t^{m(z)}(q_{z}(x,\xi)-q_{z}(x,\xi))\}_{z\in \Omega}$  is 
contained in $\Hol(\Omega, S^{-\infty}(\URd))$. Thus $(q_{z})_{z\in \Omega}$ belongs to $ \Hol(\Omega, S_{\ah}^{*}(\URd))$.
         
  Conversely, suppose that $(q_{z})_{z\in\Omega}$ is contained in $ \Hol(\Omega, S_{\ah}^{*}(\URd))$ and has 
  degree $m(z)$. Then, for any $t>0$, any compact 
  $K \subset U$, any integer $N$ and any open $\Omega'\subsubset \Omega$, we have   
  \begin{equation}   
      |\partial_{x}^\alpha\partial^\beta_{\xi}(q_{z}(x,t.\xi) -  t^{m(z)} q_{z}(x,\xi))| \leq 
     C_{tNK\Omega'\alpha\beta}  (1+\|\xi\|)^{-{N}},  
      \label{eq:HolPHDO.almost-homogeneity}
  \end{equation}
  for $(x,\xi)\in K \times\Rd$ and $z\in \Omega'$.
  Then replacing $\xi$ by $s.\xi$, 
 $s>0$, in~(\ref{eq:HolPHDO.almost-homogeneity}) shows that when $N\geq \sup_{z\in\Omega'}\Re \hat{m}(z)$ we have 
 \begin{multline}
     | \partial_{x}^\alpha\partial^\beta_{\xi}(s^{m(z)}q_{z}(x,st.\xi) - (st)^{m(z)}q_{z}(x,s.\xi))|  \leq 
        C_{tNK\Omega'\alpha\beta}   s^{\Re m(z)-N}\|\xi\|^{-{N}}\\ \leq  C_{tNK\Omega'\alpha\beta}   s^{-1}\|\xi\|^{-{N}}.
     \label{eq:HolPHDO.twisted-almost-homogeneity}
 \end{multline}
 for $(x,\xi)\in K \times\Rdo$ and $z\in \Omega'$. 
 
 Next, for $ k\in\N$ let $q_{z,k}(x,\xi)=(2^{k})^{-m}q_{z}(x,2^{k}.\xi)$. Then, for any compact 
$K \subset U$, any open $\Omega'\subsubset \Omega$ and any integer $N \geq \sup_{z\in\Omega'}\Re \hat{m}(z)$, we have
 \begin{equation}
     | \partial_{x}^\alpha\partial^\beta_{\xi}(q_{z,k+1}(x,\xi) - q_{z,k}(x,\xi))|\leq C_{2NK\Omega'\alpha\beta} 2^{-k}\|\xi\|^{-{N}},
 \end{equation}
for $(x,\xi)\in K \times\Rdo$ and $z\in \Omega'$.
This shows that the series 
 $\sum_{k \geq 0} (q_{z,k+1}-q_{z,k})$ is convergent in $\Hol(\Omega,C^\infty(\URdo))$. Hence the sequence $(q_{z,k})_{k\geq 0}$ 
 converges in $\Hol(\Omega,C^\infty(\URdo))$  to some  family $(p_{z})_{z\in \Omega}$. 
 In fact, taking $s=2^{k}$ in~(\ref{eq:HolPHDO.twisted-almost-homogeneity})  
  and letting $k\rightarrow \infty$ with $t$ fixed 
 shows that $q_{z}$ is homogeneous of degree $m(z)$ with respect to the $\xi$-variable. Moreover, for any compact 
$K \subset U$, any open $\Omega'\subsubset \Omega$ and any integer $N \geq \sup_{z\in\Omega'}\Re \hat{m}(z)$, we have 
\begin{equation}
     | \partial_{x}^\alpha\partial^\beta_{\xi}(q_{z}- p_{z})(x,\xi)| \leq \sum  |\partial_{x}^\alpha\partial^\beta_{\xi}(q_{z,k+1}-q_{z,k})(x,\xi)|  
     \leq C_{2NK\Omega'\alpha\beta} \|\xi\|^{-{N}}, 
\end{equation}
for $(x,\xi)\in K \times(\Rdo)$ and $z\in \Omega'$, i.e.~we have $q_{z}\sim p_{z}$ in the sense of~(\ref{eq:complex.symbols.asymptotic-expansion}). 
\end{proof}

Using Lemma~\ref{lem:HolPHDO.characterization-almost-homogeneous-symbols} and arguing as in the proof of 
Lemma~\ref{lem:HolPHDO.almost-homogeneous-characterization} we get the following 
characterization of holomorphic families of symbols. \begin{lemma}\label{lem:complex.symbols.criterion-almost-homogeneity}
    Let $(p_{z})_{z\in \Omega} \in \Hol(\Omega,C^{\infty}(\URd))$. Then 
    we have equivalence:\smallskip 
    
    (i) The family $(p_{z})_{z\in \Omega}$ belongs to $\Hol(\Omega, S^{*}(\URd))$ and has holomorphic order $m(z)$.\smallskip 
    
    (ii) For $j=0,1,\ldots$ there exists $(q_{j,z})_{z\in \Omega}\in \Hol(\Omega, S_{\ah}(\URd))$ almost homogeneous of degree 
    $m(z)-j$ so that we have $p_{z}\sim \sum_{j \geq 0} q_{j,z}$ in the sense of~(\ref{eq:complex.symbols.asymptotic-expansion}).
\end{lemma}

Next, it is shown in~\cite[Sect.~12]{BG:CHM} that, as for homogeneous symbols 
in~(\ref{eq:PsiHDO.convolution-symbol-pointwise})--(\ref{eq:PsiHDO.convolution-symbols-URd}), there is a continuous bilinear product, 
\begin{equation}
    *:\Svb^{k_{1}}(\URd) \times \Svb^{k_{2}}(\URd) \longrightarrow \Svb^{k_{1}+k_{2}}(\URd),  
\end{equation}
which is homogeneous in the sense that, for any $\lambda \in \R$, we have 
\begin{equation}
    (p_{1}*p_{2})_{\lambda}=p_{1,\lambda}*p_{2,\lambda}, \qquad p_{j}\in \Svb^{k_{j}}(\URd).
     \label{eq:HolPHDO.homogeneity-*product}
\end{equation}
This product is related to the product of homogeneous symbols as follows.

\begin{lemma}[{\cite[Sect.~13]{BG:CHM}}]
   For $j=1,2$ let $p_{j}\in S^{m_{j}}(\URd)$ have principal symbol 
   $p_{m_{j}}\in S_{m_{j}}(\URd)$. Then $p_{1}*p_{1}$ lies in $S^{m_{1}+m_{2}}(\URd)$ 
   and has principal symbol $p_{m_{1}}*p_{m_{2}}$. 
\end{lemma}

Furthermore, this product is holomorphic, for we have: 

\begin{lemma}\label{lem:HolPHDO.*convolution} 
   For $j=1,2$ let $(p_{j,z})_{z \in \Omega} \subset S^{*}(\URd)$ be a holomorphic family of symbols. 
   Then $(p_{1,z}*p_{2,z})_{z \in \Omega}$ is a holomorphic family of symbols as well.  \end{lemma}
\begin{proof}
    For $j=1,2$ let $m_{j}(z)$ be the order of $p_{j,z}$. Since $m_{1}(z)$ and $m_{2}(z)$ are holomorphic functions, possibly by 
    shrinking $\Omega$, we may assume that $\sup_{z\in \Omega}m_{j}(z)\leq k<\infty$. Then
   each family $(p_{j,z})_{z\in \Omega}$ belongs to $\Hol(\Omega, S_{\|}^{k}(\URd))$. 
    Since $*$ is  a continuous $\C$-bilinear map from  
$S_{\|}^{k}(\URd)\times S_{\|}^{k}(\URd)$ to $S_{\|}^{2k}(\URd)$ we 
see that $p_{1,z}*p_{2,z}$ is in $\Hol(\Omega, S_{\|}^{2k}(\URd))$, hence in $\Hol(\Omega, C^{\infty}(\URd))$.

Now, assume that $p_{j,z}$, $j=1,2$, is almost homogeneous of degree $m_{j}(z)$. Then 
using~(\ref{eq:HolPHDO.homogeneity-*product}) 
we see that, for any $\lambda>0$, the symbol $(p_{1,z}*p_{2,z})_{\lambda} -\lambda^{m_{1}(z)+m_{2}(z)} p_{1,z}*p_{2,z}$ is equal to
\begin{multline}
  p_{1,z,\lambda}*p_{1,z,\lambda} -\lambda^{m_{1}(z)+m_{2}(z)} p_{1,z}*p_{2,z}\\ 
 = (p_{1,z,\lambda}-\lambda^{m_{1}(z)}p_{1,z})*p_{2,z}+\lambda^{m_{1}(z)} p_{1,z,\lambda}*(p_{2,z,\lambda}-\lambda^{m_{2}(z)}p_{2,z}).
\end{multline}
Since $(p_{1,z,\lambda}-\lambda^{m_{1}(z)}p_{1,z})_{z\in \Omega}$ and $(p_{2,z,\lambda}-\lambda^{m_{2}(z)}p_{2,z})_{z\in \Omega}$ both belong to 
$\Hol(\Omega,S^{-\infty}(\URd))$ 
combining this with the analyticity of $*$ on $\Svb^{*}(\URd)$ shows that, for any $\lambda>0$, the family   
$(p_{1,z}*p_{2,z})_{\lambda} -\lambda^{m_{1}(z)+m_{2}(z)} p_{1,z}*p_{2,z}$ belongs to 
$\Hol(\Omega,S^{-\infty}(\URd))$. Then Lemma~\ref{lem:complex.symbols.criterion-almost-homogeneity}  
 implies that $p_{1,z}*p_{2,z}$ is a holomorphic family of almost homogeneous 
of symbols of degree $m_{1}(z)+m_{2}(z)$. 

In general, by Lemma~\ref{lem:complex.symbols.criterion-almost-homogeneity} we have
$p_{j,z} \sim \sum_{l\geq 0} p_{j,z,l}$, with $(p_{j,z,l})_{z\in 
\Omega}\in \Hol(\Omega, S^{*}_{\ah}(\URd))$ of degree 
$m_{j}(z)-l$ and $\sim$ taken in the sense of~(\ref{eq:complex.symbols.asymptotic-expansion}). 
In particular, for any integer $N$  we have 
$ p_{j,z}=\sum_{l<N} p_{j,z,l}\bmod \Hol(\Omega, S_{\|}^{k-N}(\URd))$. 
Thus, 
\begin{equation}
    p_{1,z}*p_{2,z}=\sum_{l+p<N} p_{1,z,l}*p_{2,z,p} \quad \bmod \Hol(\Omega, S_{\|}^{2k-N}(\URd)).
     \label{eq:HolPHDO.analyticity-*convolution-symbols}
\end{equation}
As  explained above $(p_{1,z,l}*p_{2,z,p})_{z\in \Omega}$ is a holomorphic family of almost homogeneous 
of symbols of degree $m_{1}(z)+m_{2}(z)-l-p$. It then follows from 
Lemma~\ref{lem:complex.symbols.criterion-almost-homogeneity} that $(p_{1,z}*p_{2,z})_{z\in \Omega}$ 
is a holomorphic family of symbols.
\end{proof}

We are now ready to prove: 
\begin{proposition}\label{prop:HolPHDO.composition}
    For $j=1,2$ let $(P_{j,z})_{z\in \Omega}$ be in 
    $\Hol(\Omega,\pvdo^{*}(U))$ and suppose that at least one the families $(P_{1,z})_{z 
    \in \Omega}$ or $(P_{2,z})_{z \in \Omega}$ is uniformly properly supported. 
    Then the family $(P_{1,z}P_{2,z})_{z  \in \Omega}$ is a holomorphic family of \psivdos. 
\end{proposition}
\begin{proof}
By assumption $(P_{1,z})_{z\in \Omega}$ or $(P_{2,z})_{z\in \Omega}$ 
is a uniformly properly supported holomorphic family of \psivdos, 
hence gives rise to elements of $\Hol(\Omega,\cL(C^{\infty}(U)))$ and $\Hol(\Omega,\cL(\cE'(U)))$ 
by Proposition~\ref{prop:Complexpowers.operators.properties2}. 
Moreover, Proposition~\ref{prop:Complexpowers.operators.properties2} 
tells us that the other family at least coincides with a uniformly properly supported holomorphic family of \psivdos\ 
 up to a holomorphic family of smoothing operators. 
 It thus follows that $(P_{1,z}P_{2,z})_{z\in \Omega}$ is the product of two uniformly properly supported holomorphic families of 
 \psivdos\ up to a holomorphic families of smoothing operators. 

As a consequence we may assume that the families 
$(P_{1,z})_{z\in \Omega}$ 
and $(P_{2,z})_{z\in \Omega}$ are both uniformly properly supported. Thanks 
to~(\ref{eq:HolPHDO.uniformly-properly-supported-locally-finite-sum})   
this allows us to write 
\begin{equation}
    P_{1,z}P_{2,z}=\sum_{i \geq 0} \varphi_{i} P_{1,z}\psi_{i}P_{2,z},
    \label{eq:PsiHDO.composition-locally-finite-sum}
\end{equation}
where $(\varphi_{i})_{i \geq 0}\subset C^{\infty}_{c}(U)$ and $(\psi_{i})_{i\geq 0} \subset  
   C^{\infty}_{c}(U)$ are locally finite families 
  such that $(\varphi_{i})$ is a partition of the unity and $\psi_{i}=1$ near $\supp\varphi_{i}$. 
  
  Next, for $j=1,2$ let us write $P_{j,z}=p_{j,z}(x,-iX)+R_{j,z}$, with  $(p_{j,z})_{z\in \Omega}$ in 
  $\Hol(\Omega, S^{*}(\URd))$ and 
    $(R_{j,z})_{z\in \Omega}$ in $\Hol(\Omega, \psinf(U))$. Since by Proposition~\ref{prop:Complexpowers.operators.properties}
    each family $(p_{j,z}(x,-iX))_{z\in \Omega}$ is holomorphic with values in 
    $\cL(C_{c}^{\infty}(U), C^{\infty}(U))$ and $\cL(\cE'(U),\cD'(U))$ 
    using~(\ref{eq:PsiHDO.composition-locally-finite-sum}) we see that 
    \begin{equation}
        P_{1,z}P_{2,z}=\sum \varphi_{i} 
        p_{1,z}(x,-iX)\psi_{i}p_{2,z}(x,-iX) \quad \bmod 
        \Hol(\Omega, \Psi^{-\infty}(U)).
    \end{equation}
    
    At this stage we make appeal to: 
    
    \begin{lemma}[{\cite[Prop.~14.45]{BG:CHM}}] \label{lem:HolPHDO.proof-composition}
 For $j=1,2$ let $p_{j}\in \SvbU{k_{j}}$ and let $\psi \in 
 C_{c}^\infty(U)$. Then:
   \begin{equation}
     p_{1}(x,-iX) \psi p_{2}(x,-iX)= p_{1}\#_{\psi} p_{2}(x,-iX), 
 \end{equation}
 where $ \#_{\psi}$ is a continuous bilinear map from 
 $\SvbU{k_{1}}\times\SvbU{k_{2}}$ to $\SvbU{k_{1}+k_{2}}$. Moreover, for any integer $N\geq 1$ we have  
   \begin{equation}
      p_{1}\#_{\psi}p_{2} = \sum_{j<N}
        \sum_{\alpha\beta\gamma\delta}^{(j)} h_{\alpha\beta\gamma\delta} \psi (D_{\xi}^\delta p_{1})* (\xi^\gamma 
            \partial_{x}^\alpha \partial_{\xi}^\beta p_{2}) + R_{N, \psi}(p_{1},p_{2}),
  \end{equation}
  where the notations are the same as in 
  Proposition~\ref{prop:PsiHDO.composition} and $R_{N, \psi}$ is a 
  continuous bilinear map from  $\SvbU{k_{1}}\times\SvbU{k_{2}}$ to $\SvbU{k_{1}+k_{2}-N}$.
\end{lemma}
\begin{remark}
 The continuity contents of Lemma~\ref{lem:HolPHDO.proof-composition} is explicitly stated in Proposition~14.45 of \cite{BG:CHM}, but they follow 
 from its proof or from a standard use of the closed graph theorem.   
\end{remark}

Now, thanks to Lemma~\ref{lem:HolPHDO.proof-composition} we have
\begin{equation}
    P_{1,z}P_{2,z}=p_{z}(x,-iX)+R_{z}, \qquad p_{z}= \sum \varphi_{i}p_{1,z}\#_{\psi_{i}}p_{2,z}. 
     \label{eq:HolPHDO.composition-PsiHDO}
\end{equation}
Furthermore, possibly by shrinking $\Omega$, we may assume that there is a real $k$ such that 
for $j=1,2$ we have $\Re \ord p_{j,z}\leq k$ for any $z \in \Omega$. Then 
the continuity contents of Lemma~\ref{lem:HolPHDO.proof-composition} 
imply that $(p_{z})_{z\in \Omega}$ belongs to $\Hol(\Omega, \Svb^{2k}(\URd))$ 
and for any integer $N\geq 1$ we can write
\begin{equation}
       p_{z} =  \sum_{r<N}  q_{r,z} +R_{N,z}, \quad 
        q_{j,z}= \overset{(r-s-t)}{\underset{\alpha\beta\gamma\delta}{\sum}}
h_{\alpha\beta\gamma\delta} (D_{\xi}^\delta p_{1,z})* (\xi^\gamma 
            \partial_{x}^\alpha \partial_{\xi}^\beta p_{2,z}),
\end{equation}
where $(R_{N,z})_{z\in \Omega}:=(\sum_{i\geq 0}\varphi_{i} R_{N, 
\psi_{i}}(p_{1,z},p_{2,z}))_{z\in \Omega}$ is in $\Hol(\Omega,S_{\|}^{2k-N}(\URd))$. 
Thanks to Lemma~\ref{lem:HolPHDO.*convolution} 
            the family $(q_{r,z})_{z\in \Omega}$ is in $\Hol(\Omega, S^{*}(\URd))$  
            and has order $m_{1}(z)+m_{2}(z)-r$. Therefore, we have $p_{z} \sim \sum_{j\geq 0} q_{j,z}$ 
            in the sense of~(\ref{eq:complex.symbols.asymptotic-expansion}), 
            which by Lemma~\ref{lem:HolPHDO.*convolution} implies that $(p_{z})_{z\in \Omega}$ 
            belongs to $\Hol(\Omega, S^{*}(\URd))$.
           Combining this with~(\ref{eq:HolPHDO.composition-PsiHDO}) 
           then shows that $(P_{1,z}P_{2,z})_{z \in \Omega}$ is a holomorphic family of \psivdos. \end{proof}

\subsection{Kernel characterization of holomorphic families of $\mathbf{\Psi}$DO's}
We shall now give a characterization of holomorphic families of \psivdos\ in terms of holomorphic families with values in 
$\cK^{*}(\URd)$. Since the 
latter is defined in terms of asymptotic expansions of kernels in 
$\cK_{*}(\URd)$ (\emph{cf.} Definition~\ref{def:PsiHDO.K*}) a 
technical difficulty occurs, because for a family $(K_{z})_{z\in 
\Omega}\subset \cK_{*}(\URd)$ logarithmic singularities may appear as 
the order of $K_{z}$ crosses non-negative integer values. 

This issue is resolved by making use of holomorphic families of almost homogeneous distributions as follows.
\begin{definition}
   A family $(K_{z})_{z \in \Omega}\subset \cK^{*}_{\ah}(\URd)$ is holomorphic when:\smallskip 
   
   (i) The degree $m(z)$ of $K_{z}$ is a holomorphic function on $\Omega$;\smallskip 
   
   (ii) The family $(K_{z})_{z\in \Omega}$ belongs to $\Hol(\Omega,C^{\infty}(U)\otimes \cD'_{\reg}(\Rd))$;\smallskip 
  
   (iii) For any $\lambda >0$ the family $\{K_{z}(x,\lambda.y)-\lambda^{m(z)}K_{z}(x,y)\}_{z \in \Omega}$ is in 
   $\Hol(\Omega, C^{\infty}(\URd))$.\smallskip 
   
\noindent We let $\Hol(\Omega,\cK^{*}_{\ah}(\URd))$ denote the set of holomorphic 
families of almost homogeneous distributions.
\end{definition}

\begin{lemma}\label{lem:HolPHDO.almost-homogeneous-kernel-symbol}
    Let $(K_{z})_{z\in \Omega}\in \Hol(\Omega, C^{\infty}(U)\otimes 
    \cD'_{\reg}(\Rd))$.
    Then we have equivalence:\smallskip 
    
    (i) The family $(K_{z})_{z \in \Omega}$ belongs to $\Hol(\Omega, \cK_{\ah}^{*}(\URd))$ and has degree $m(z)$.\smallskip 
    
    (ii) We can put $(K_{z})_{z \in \Omega}$ into the form,
    \begin{equation}
        K_{z}(x,y)=\check{p}_{z,\xiy}(x,y) +R_{z}(x,y), \qquad z\in \Omega,
    \end{equation}
    for some family $(p_{z})_{z \in \Omega}\in 
    \Hol(\Omega,S^{*}_{\ah}(\URd))$ of degree $\hat{m}(z):=-(m(z)+d+2)$ and some family $(R_{z})_{z\in \Omega}\in \Hol(\Omega, 
    C^{\infty}(\URd))$.   
\end{lemma}
\begin{proof}
 Assume that the family $(K_{z})_{z \in \Omega}$ belongs to $\Hol(\Omega, \cK_{\ah}^{*}(\URd))$ and has degree 
 $m(z)$.  Let $\varphi \in 
 C^{\infty}_{c}(\URd)$ be such that $\varphi(y)=1$ near $y=0$ and set $p_{z}=(\varphi(y)K_{z}(x,y))_{\yxi}$. Thus,
 \begin{equation}
     K_{z}(x,y)= \check{p}_{z,\xiy}(x,y) 
 +(1-\varphi(y))K_{z}(x,y) = \check{p}_{z,\xiy}(x,y) \ \bmod \Hol(\Omega, C^{\infty}(\URd)).
 \end{equation}
 As $(\varphi(y)K_{z}(x,y))_{z\in \Omega}$ is in 
 $\Hol(\Omega, C^{\infty}(U)\otimes\cE'_{L}(\Rd)$ with $L=\supp \varphi$, we see that $(p_{z})_{z\in \Omega}$ belongs to 
 $\Hol(\Omega,C^{\infty}(\URd))$. 
 Furthermore, thanks to~(\ref{eq:PsiHDO.dilations-Fourier-transform}) for any $\lambda>0$ we have: 
 \begin{multline}
     p_{z}(x,\lambda.\xi)-\lambda^{\hat{m}(z)}p_{z}(x,\xi)=\\ \lambda^{-(d+2)}
     [(\varphi(\lambda^{-1}.y)-\varphi(y))K_{z}(x,\lambda^{-1}.y)+
    \varphi(y) (K_{z}(x,\lambda^{-1}.y)-\lambda^{-m}K_{z}(x,y))]_{\yxi}^{\wedge}(x,\xi).
 \end{multline}
Since the r.h.s.~above is the Fourier transform with respect to $y$ of an element of 
$\Hol(\Omega, C_{c}^{\infty}(\URd))\subset \Hol(\Omega, 
    S^{-\infty}(\URd))$ we see that $(p_{z}(x,\lambda.\xi)-\lambda^{\hat{m}(z)}p_{z}(x,\xi))_{z \in \Omega}$ is in 
    $\Hol(\Omega, S^{-\infty}(\URd))$ for any $\lambda>0$. Combining this with Lemma~\ref{lem:HolPHDO.characterization-almost-homogeneous-symbols}
    then shows that 
    $(p_{z})_{z\in \Omega}$ is a family of almost homogeneous symbols of degree $\hat{m}(z)$.     
    Conversely, let $(p_{z})_{z \in \Omega}\in  \Hol(\Omega, S^{*}_{\ah}(\URd))$ be almost homogeneous of degree 
    $\hat{m}(z)$. Possibly by shrinking $\Omega$ we may assume that we have $\sup_{z \in \Omega}\Re \hat{m}(z)\leq k<\infty$. 
    Then the family $(p_{z})_{z \in \Omega}$ belongs to $\Hol(\Omega, \Svb^{k}(\URd))$, 
    hence to $\Hol(\Omega, S_{\frac{1}{2}\frac{1}{2}}^{k}(\URd))$. As mentioned in the course of the proof of 
    Proposition~\ref{prop:Complexpowers.operators.properties}, the map $q(x,\xi) \rightarrow \check{q}_{\xiy}(x,y)$ is analytic from 
    $S_{\frac{1}{2}\frac{1}{2}}^{k}(\URd)$ to $C^{\infty}(U)\otimes \cD'_{\reg}(\Rd)$, so the family 
    $(\check{p}_{z,\xiy})_{z\in \Omega}$ lies in $\Hol(\Omega,C^{\infty}(U)\otimes \cD'_{\reg}(\Rd))$. 
    
    Next, using~(\ref{eq:PsiHDO.dilations-Fourier-transform}) we see that, for any $\lambda>0$, we have
    \begin{equation}
        \check{p}_{z,\xiy}(x,\lambda.y)-\lambda^{m(z)}\check{p}_{z,\xiy}(x,y)= 
        [p_{z}(x,\lambda^{-1}.\xi)-\lambda^{-\hat{m}(z)}p_{z}(x, \xi)]_{\xiy}^{\wedge}(x,y).
         \label{eq:HolPHDO.almost-homogeneity-checkp}
    \end{equation}
    Since by Lemma~\ref{lem:HolPHDO.characterization-almost-homogeneous-symbols}
    the r.h.s.~of~(\ref{eq:HolPHDO.almost-homogeneity-checkp}) is the inverse Fourier transform with respect to $\xi$ of an element 
    of $\Hol(\Omega,S^{-\infty}(\URd))$ we see that the family
    $\{\check{p}_{z,\xiy}(x,\lambda.y)-\lambda^{m(z)}\check{p}_{z,\xiy}(x,y)\}_{z \in \Omega}$ is contained in $\Hol(\Omega, C^{\infty}(\URd))$. 
    It then follows that $(\check{p}_{z,\xiy})_{z \in \Omega}$ is a holomorphic family of almost homogeneous distributions of degree $m(z)$. 
\end{proof}

\begin{definition}\label{def:HolPHDO.kernels-families}
   A family $(K_{z})_{z\in \Omega}\subset \cK^{*}(\URd)$ is holomorphic when:\smallskip 
    
    (i) The order $m_{z}$ of $K_{z}$ is a 
    holomorphic function of $z$;\smallskip 
    
   (ii) For $j=0,1,..$ there exists $(K_{j,z})\in \Hol (\Omega, \cK_{\ah}^{*}(\URd))$ of degree 
   $m(z)+j$ such that $K_{z}\sim \sum_{j \geq 0} K_{j,z}$  in the sense that, for any open 
    $\Omega' \subsubset \Omega$ and any integer $N$, as soon as $J$ is large enough we have 
 \begin{equation}
         K_{z}- \sum_{j\leq J} K_{z,m_{z}+j}\in \Hol(\Omega',C^{N}(\URd)).
                      \label{eq:HolPHDO.kernel-asymptotic}
 \end{equation}
\end{definition}\begin{proposition}\label{prop:HolPHDO.characterization-Hol-cK*}
   For a family  $(K_{z})_{z\in \Omega}\subset \cK^{*}(\URd)$ the following are equivalent:\smallskip
   
   (i) The family   $(K_{z})_{z\in \Omega}$ is holomorphic and has order $m(z)$.\smallskip
   
   (ii) We can put $(K_{z})_{z\in \Omega}$ into the form, 
   \begin{equation}
       K_{z}(x,y)=(p_{z})^{\vee}_{\xiy}(x,y)+R_{z}(x,y),
       \label{eq:HolPHDO.kernel-symbol}
   \end{equation}
  for some family $(p_{z})_{z \in \Omega}\in \Hol(\Omega,S^{*}(\URd))$ of order $\hat{m}(z):=-(m(z)+d+2)$ and 
  some family $(R_{z})_{z\in \Omega}\in \Hol(\Omega, C^{\infty}(\URd))$.   
   \end{proposition}
 \begin{proof}
     Assume that $(K_{z})_{z\in \Omega}$ belongs to $\Hol(\Omega, \cK^{*}(\URd))$. Let $\varphi\in C^{\infty}(\Rd)$ be such that 
    $\varphi(y)=1$ near $y=0$ and for $z\in \Omega$ let $p_{z}=(\varphi(y)K_{z}(x,y))^{\wedge}_{\yxi}$. Since $(K_{z})_{z\in \Omega}$ lies in 
    $\Hol(\Omega, C^{\infty}(U)\hotimes \cD_{\reg}'(\Rd))$ we see that $(p_{z})_{z\in \Omega}$ belongs to $\Hol(\Omega, C^{\infty}(\URd))$ and we 
    have 
      \begin{equation}
       K_{z}(x,y)=(p_{z})^{\vee}_{\xiy}(x,y)+(1-\varphi(y))K_{z}(x,y)=(p_{z})^{\vee}_{\xiy}(x,y) \quad \bmod \Hol(\Omega,C^{\infty}(U\times U)).
       \label{eq:HolPHDO.kernel-symbol2}
   \end{equation}
    
   Let us write $K_{z}\sim \sum_{j\geq 0}K_{j,z}$ with 
  $(K_{j,z})_{z \in \Omega}\in \Hol(\Omega, \cK_{\ah}^{*}(\URd))$ of 
  degree $\hat{m}(z)+j$ and $\sim$ taken in the sense of~(\ref{eq:HolPHDO.kernel-asymptotic}). For $j=0,1,\ldots$ we let 
  $p_{j,z}=(\varphi(y)K_{j,z})_{z \in \Omega}$. Then arguing as in the proof of Lemma~\ref{lem:HolPHDO.almost-homogeneous-kernel-symbol} shows 
  that $(p_{j,z})_{z\in \Omega}$ is a family in $\Hol(\Omega, 
  S^{*}_{\ah}(\URd))$ of degree $\hat{m}(z)-j$. 
  
  Next, in the sense of~(\ref{eq:HolPHDO.kernel-asymptotic}) we have  $(p_{z})^{\vee}_{\xiy}(x,y)\sim \sum_{j \geq 0}(p_{j,z})^{\vee}_{\xiy}(x,y)$.
 Under the Fourier transform with respect to $y$ this shows  that, 
 for any compact $L\subset U$, any integer $N$ and any open $\Omega' \subsubset \Omega$, as soon as $J$ is large enough 
 we have estimates,
     \begin{equation}
         |\partial_{x}^{\alpha}\partial_{\xi}^{\beta}(p_{z}- \sum_{j\leq J} (\varphi(y) K_{j,z}(x,y))_{\yxi}^{\wedge})(x,\xi)| \leq 
         C_{\Omega'NJL\alpha\beta}(1+|\xi|^{2})^{-[N/2]},
          \label{eq:HolPHDO.kernel-characterization.asymptotic-expansion-symbol}
     \end{equation}
 for $(x,\xi)\in L\times\Rd$ and $z \in \Omega'$. Hence $p_{z}\sim \sum_{j\geq 0}p_{j,z}$ in the sense 
 of~(\ref{eq:complex.symbols.asymptotic-expansion}).  It thus follows that
 $(p_{z})_{z\in\Omega}$ is in $\Hol(\Omega, S^{*}(\URd))$ and has order $\hat{m}(z)$, so using~(\ref{eq:HolPHDO.kernel-symbol2}) we see that the family $(K_{z})_{z\in 
 \Omega}$ is of the form~(\ref{eq:HolPHDO.kernel-symbol}). 
 
      Conversely, assume that $(K_{z})_{z \in \Omega}$ is of the form  $K_{z}(x,y)=(p_{z})^{\vee}_{\xiy}(x,y)+R_{z}(x,y)$ for some family
      $(p_{z})_{z \in \Omega}$ in $\Hol(\Omega,S^{*}(\URd))$ of order $\hat{m}(z):=-(m(z)+d+2)$ and some family 
  $(R_{z})_{z\in \Omega}$ in $\Hol(\Omega, C^{\infty}(\URd))$. 
      By Lemma~\ref{lem:complex.symbols.criterion-almost-homogeneity} 
    we have $p_{z}\sim \sum_{j\geq 0}p_{j,z}$ with $(p_{j,z})_{z\in \Omega}\in 
    \Hol(\Omega,S^{*}_{\ah}(\URd))$ of degree $m(z)-j$ and $\sim$ taken in the sense 
    of~(\ref{eq:complex.symbols.asymptotic-expansion}). Thus, under 
    the inverse Fourier transform with respect to $\xi$, we get an asymptotic expansion $(p_{z})^{\vee}_{\xiy}\sim \sum_{j\geq 
    0}(p_{j,z})^{\vee}_{\xiy}$ in the sense of~(\ref{eq:HolPHDO.kernel-asymptotic}). 
    As Lemma~\ref{lem:HolPHDO.almost-homogeneous-kernel-symbol} tells us that  
    $((p_{j,z})^{\vee}_{\xiy})_{z \in \Omega}$ is a holomorphic family of almost 
    homogeneous distributions of degree $\hat{m}(z)+j$, it follows that 
    $((p_{z})^{\vee}_{\xiy})_{z\in \Omega}$ is a family in $\Hol(\Omega,\cK^{*}(\URd))$ of order $\hat{m}(z)$. Since $(K_{z})_{z \in \Omega}$ agrees with 
    $((p_{z})^{\vee}_{\xiy})_{z\in \Omega}$ up to an element of $\Hol(\Omega,C^{\infty}(U\times U))$, the same is true for $(K_{z})_{z \in \Omega}$. 
\end{proof}
   We are now ready to prove the kernel characterization of holomorphic families of \psivdos. 

\begin{proposition}\label{prop:HolPHDO.kernel.characterization}
 Let $(P_{z})_{z\in \Omega} \in \Hol(\Omega,\cL(C_{c}^{\infty}(U),C^{\infty}(U)))$. Then the following are 
 equivalent:\smallskip 
 
 (i) The family $(P_{z})_{z\in \Omega}$ is a holomorphic family of \psivdos\ of order $m(z)$.\smallskip 
 
 (ii) There exist $(K_{z})_{z\in \Omega}\in \Hol(\Omega,\cK^{*}(\URd))$ of order $\hat{m}(z):=-(m(z)+d+2)$ 
  and $(R_{z})_{z\in \Omega}\in \Hol(\Omega,C^{\infty}(U\times U))$ 
 such that the distribution kernel $k_{P_{z}}(x,y)$ of $P_{z}$ is of the form
  \begin{equation}
     k_{P_{z}}(x,y) = |\psi_{x}'| P_{z}(x,-\varepsilon_{x}(y)) + R_{z}(x,y),
     \label{eq:HolPHDO.kernel.characterization-psix}
 \end{equation}
 where  $\psi_{x}$ denotes the coordinate change to the privileged coordinates at $x$.\smallskip   

 (iii) There exist $(K_{P_{z}})_{z\in \Omega}\in \Hol(\Omega,\cK^{*}(\URd))$ of order $\hat{m}(z):=-(m(z)+d+2)$ 
  and $(R_{z})_{z\in \Omega}\in \Hol(\Omega,C^{\infty}(U\times U))$ 
 such that the distribution kernel $k_{P_{z}}(x,y)$ of $P_{z}$ is of the form
  \begin{equation}
     k_{P_{z}}(x,y) = |\varepsilon_{x}'| P_{z}(x,-\varepsilon_{x}(y)) + R_{z}(x,y).
     \label{eq:HolPHDO.kernel.characterization}
 \end{equation}
 where  $\varepsilon_{x}$ denotes the coordinate change to the Heisenberg coordinates at $x$.   
 \end{proposition}
\begin{proof}
   First, it follows from~(\ref{eq:PsiHDO.kernel-quantization-symbol-psiy}) and 
   Proposition~\ref{prop:HolPHDO.characterization-Hol-cK*} that (i) and (ii) are equivalent. 
   
   Next, for $x\in U$ let $\phi_{x}$ denote the transition map from the privileged coordinates at $x$ to the Heisenberg coordinates at $x$. 
Recall that the coordinate changes $\phi_{x}$, $x \in U$, give rise to an action on distributions on $\URd$ given by
   \begin{equation}
     K(x,y) \longrightarrow \phi_{x}^{*}K(x,y), \qquad \phi_{x}^{*}K(x,y)=K(x,\phi_{x}^{-1}(y)).
        \label{eq:HolPHDO.action-phix*}
   \end{equation}

   Since $\phi_{x}$ depends smoothly on $x$, this action gives rise to continuous linear isomorphisms of $C^{N}(\URd)$, $N \geq 0$, and $C^{\infty}(\URd)$ 
  onto themselves, hence to analytic isomorphisms. Moreover, since $\phi_{x}(0)=0$  this also yields an analytic isomorphism of $C^{\infty}(U)\hotimes 
  \cD_{\reg}'(\Rd)$ onto itself. Combining this with the homogeneity property~(\ref{eq:PsiHDO.homogeneity.phix*}) 
  we then deduce that~(\ref{eq:HolPHDO.action-phix*}) induces linear isomorphisms of 
  $\Hol(\Omega,\cK^{*}_{\ah}(\URd))$ and  $\Hol(\Omega,\cK^{*}(\URd))$ onto themselves. Together 
  with~(\ref{eq:PsiHDO.kernel-quantization-symbol}) this shows that the statements 
  (ii) and (iii) are equivalent.  
\end{proof}
\subsection{Holomorphic families of $\mathbf{\Psi}$DO's on a general Heisenberg manifold}
 Let us now define holomorphic families of \psivdos\ a general Heisenberg manifold. 
 
 First, the following holds. 

\begin{lemma}\label{lem:HolPHDO.invariance}
    Let $(K_{z})_{z\in \Omega}\in \Hol(\Omega, \cK*(\URd))$ and assume there exists integer $N$ such that $\inf_{z\in \Omega}\Re \ord K_{z}\geq 
    2N$.   Then the family $(K_{z})_{z\in \Omega}$ is contained in $\Hol(\Omega, C^{N}(\URd))$. 
\end{lemma}
\begin{proof}
    Thanks to Proposition~\ref{prop:HolPHDO.kernel.characterization} we may assume that $K_{z}$ is of the form 
    $K_{z}(x,y)=\check{p}_{z,\xiy}(x,y)$ with $(p_{z})_{z\in \Omega}$ in $\Hol(\Omega, 
    S^{*}(\URd))$.  As we have $-(\Re \ord p_{z}+d+2)= \Re \ord K_{z}\geq 2N$ we see that $(p_{z})_{z\in \Omega}$ is contained in $\Hol(\Omega, 
    \Svb^{-(2N+d+2)}(\URd))$. Since the map 
         $p\rightarrow \check{p}_{\xiy}$ is continuous from $S_{||}^{k}(\URd)$ to $C^{\infty}(U)\hotimes C^{N}(\URd)$ (see~\cite{Po:BSM1}), 
    it follows that $(K_{z})_{z \in \Omega}$ lies in $\Hol(\Omega, C^{N}(\URd))$. 
\end{proof}

\begin{lemma}\label{lem:Appendix.Heisenberg.invariance}
  1)  Assume $k<-(d+2)$ and set $N=[-\frac{k+d+2}2]$. Then the map 
         $p\rightarrow \check{p}_{\xiy}$ is continuous from $S_{||}^{k}(\URd)$ to $C^{\infty}(U)\hotimes C^{N}(\URd)$.  \smallskip 
         
         2) For $\Re m>0$ we have $\cK^{m}(\URd)\subset C^{\infty}(U)\hotimes C^{[\frac{\Re m}{2}]}(\Rd)$.  
\end{lemma}
%

\begin{proposition}\label{prop:HolPHDO.invariance}
Let $\tilde{U}$ be an open subset of $\Rd$ together with a hyperplane bundle $\tilde{H}\subset T\tilde{U}$ and a $\tilde{H}$-frame of $T\tilde{U}$ 
and let $\phi:(U,H)\rightarrow (\tilde{U},\tilde{H})$ be a Heisenberg diffeomorphism. Then for any family 
$(\tilde{P}_{z})_{z \in \Omega}\in \Hol(\Omega, \Psi_{\tilde{H}}^{*}(\tilde{U}))$ the family 
$(P_{z})_{z \in \Omega}:=(\phi^{*}\tilde{P}_{z})_{z \in \Omega}$ is contained in $\Hol(\Omega, \Psi_{H}^{*}(U))$.
\end{proposition}
\begin{proof}
For $x \in U$ and $\tilde{x}\in \tilde{U}$ let $\varepsilon_{x}$ and $\tilde{\varepsilon}_{\tilde{x}}$ 
denote the coordinate changes to the Heisenberg coordinates at $x$ and $\tilde{x}$ respectively. Then
by Proposition~\ref{prop:HolPHDO.kernel.characterization} the distribution kernel $k_{\tilde{P}_{z}}(\tilde{x},\tilde{y})$ of $\tilde{P}_{z}$ 
is of the form
  \begin{equation}
      k_{\tilde{P}_{z}}(\tilde{x},\tilde{y}) = |\tilde{\varepsilon}_{\tilde{x}}'|K_{\tilde{P}_{z}}(\tilde{x},-\tilde{\varepsilon}_{\tilde{x}}(\tilde{y})) + 
      \tilde{R}_{z}(\tilde{x},\tilde{y}),
  \end{equation}
with $(K_{\tilde{P}_{z}})_{z\in \Omega}$ in  $\Hol(\Omega, \cK^{*}(\tilde{U}\times \tilde{U}))$ and $(\tilde{R}_{z})_{z\in \Omega}$ in $\Hol(\Omega, 
C^{\infty}(\tilde{U}\times\Rd))$. Then the proof of Proposition~\ref{prop:PsiHDO.invariance} in~\cite{Po:BSM1}
shows that the distribution kernel 
$k_{P_{z}}(x,y)$ of $P_{z}$ takes the form, 
\begin{multline}
    k_{P_{z}}(x,y)= |\varepsilon_{x}'| K_{P_{z}}(x,-\varepsilon_{x}(y)) \\
    + (1-\chi(x,-\varepsilon_{x}(y))) |\tilde{\varepsilon}_{\phi(x)}'| 
    K_{\tilde{P}_{z}}(\phi(x),-\tilde{\varepsilon}_{\phi(x)}(\phi(y))) +\tilde{R}_{z}(\phi(x),\phi(y)),
\end{multline}
where we have let
\begin{equation}
    K_{P_{z}}(x,y)= \chi(x,y) |\partial_{y}\Phi(x,y)| K_{\tilde{P}_{z}}(\phi(x),\Phi(x,y)), \qquad \Phi(x,y)=-\tilde{\varepsilon}_{\phi(x)}\circ \phi\circ 
    \varepsilon_{x}^{-1}(-y),
    \label{eq:HolPHDO.kernel-characterization.KPz}
\end{equation}
and $\chi(x,y)\in C^{\infty}(\URd)$ is supported on the open subset $\cU=\{ (x,y)\in \URd; \  \varepsilon_{x}^{-1}(-y)\in U\}$, is properly 
supported with respect to $x$ and satisfies $\chi(x,y)=1$ near $U\times \{0\}$. In particular, we have

\begin{equation}
    k_{P_{z}}(x,y)= |\varepsilon_{x}'| K_{P_{z}}(x,-\varepsilon_{x}(y)) \quad \bmod \Hol(\Omega,C^{\infty}(U\times U)).
      \label{eq:HolPHDO.invariance.kPz}
\end{equation}

Let us now prove that $(K_{P_{z}})_{z \in \Omega}$ is an element of $\Hol(\Omega, \cK^{*}(\URd))$. To this end, possibly by shrinking 
$\Omega$, we may assume that $\inf_{z \in \Omega} \Re \hat{m}(z)\geq \mu >-\infty$.
Moreover, the proof of Proposition~\ref{prop:PsiHDO.invariance} in~\cite{Po:BSM1} shows that for any integer $N$  we have 
\begin{gather}
    K_{P_{z}}(x,y)= \sum_{\brak \alpha<N} \sum_{\frac32\brak\alpha \leq \brak\beta < \frac{3}{2}N} 
    K_{\alpha\beta, z}(x,y)+ \sum_{j=1}^{3} R_{N, z}^{(j)}(x,y),
    \label{eq:HolPHDO.invariance-KPz}\\
    K_{\alpha\beta, z}(x,y)=a_{\alpha\beta}(x) y^\beta 
    (\partial^\alpha_{\tilde{y}}K_{\tilde{P}_{z}})(\phi(x),\phi_{H}'(x)y),
\end{gather}
where the smooth functions $a_{\alpha\beta}(x)$ are as in Proposition~\ref{prop:PsiHDO.invariance} and the remainder terms 
   $R_{N,z}^{(j)}(x,y)$, $j=1,2,3$, take the forms:\smallskip 
   
     - $R_{N, z}^{(2)}(x,y)= \sum_{\brak \alpha<N} \sum_{\brak\beta \dot{=} \frac{3}{2}N} r_{M\alpha}(x,y) y^{\beta} 
    (\partial^\alpha_{\tilde{y}}K_{\tilde{P}_{z}})(\phi(x),\phi_{H}'(x)y)$ for some functions $r_{\alpha\beta}(x,y)$ in 
    $C^{\infty}(\URd)$;\smallskip
    
    - $R_{N, z}^{(3)}(x,y)= \sum_{\brak\alpha=N}\sum_{\brak\beta\dot{=}\frac{3}{2}N} \int_{0}^{1} 
    r_{\alpha\beta}(t,x,y)(\tilde{y}^{\beta}\partial^\alpha_{\tilde{y}}K_{\tilde{P}_{z}})(\phi(x),t\Phi(x,y)+(1-t)\phi_{H}'(x)y)dt$, 
    for some functions $r_{\alpha\beta}(t,x,y)$ in $C^{\infty}([0,1]\times \URd)$;\smallskip
    
  - $R_{N, z}^{(3)}(x,y)= \sum_{\brak \alpha<N} \sum_{\frac32\brak\alpha \leq \brak\beta < \frac{3}{2}N} (1-\chi(x,y))K_{\alpha\beta, z}(x,y)$.\smallskip 
    
Observe that the map $\Phi(x,y)=(\phi(x),\phi'_{H}(x)y)$ is a smooth diffeomorphism from $\URd$ onto $\tilde{U}\times \Rd$ such that $\Phi(x,0)=(\phi(x),0)$ 
and $\Phi(x,\lambda.y)=(\phi(x),\lambda.\phi'_{H}(x)y)$, so along similar lines as that of the proof of Proposition~\ref{prop:HolPHDO.kernel.characterization} 
we can prove that the map  
\begin{equation}
 \cD'(\URd) \ni K(x,y)\longrightarrow K(\phi(x),\phi'_{H}(x)y) \in \cD'(\tilde{U}\times \Rd)
 \end{equation}
gives rise to a linear map from $\Hol(\Omega, \cK^{*}(\URd))$ to $\Hol(\Omega, \cK^{*}(\tilde{U}\times \Rd))$ preserving the order. Therefore, 
the family $(K_{\alpha\beta, z}(x,y))_{z \in \Omega}$ is contained in $\Hol(\Omega, \cK^{*}(\tilde{U}\times \Rd))$ and has order 
$\hat{m}(z)+\brak\beta-\brak\alpha$.  Incidentally, the term $(R_{N, z}^{(3)})_{z \in \Omega}$ 
belongs to $\Hol(\Omega, C^{\infty}(\URd))$.

On the other hand, if $\frac32\brak\alpha \leq \brak\beta  \dot{=} \frac{3}{2}N$ then the order 
$\hat{m}_{\alpha\beta}(z)=\hat{m}(z)+\brak\beta-\brak\alpha$ of the family $(\tilde{y}^{\beta}\partial^\alpha_{\tilde{y}}K_{\tilde{P}_{z}})_{z \in 
\Omega}\in \Hol(\Omega, \cK^{*}(\URd))$ satisfies  
$ \Re\hat{m}_{\alpha\beta}(z)\geq   \Re \hat{m}(z)+\frac{1}{3}\brak\beta\geq \mu +\frac{N}{2}$. Therefore, it follows from 
Lemma~\ref{lem:HolPHDO.invariance} that, for any integer $J$, as soon as $N$ is  large enough 
$(\tilde{y}^{\beta}\partial^\alpha_{\tilde{y}}K_{\tilde{P}_{z}})_{z \in \Omega}$ is in $\Hol(\Omega, C^{J}(\tilde{U}\times \Rd))$, as much so the 
remainder terms $(R_{N, z}^{(2)})_{z \in \Omega}$ and $(R_{N, z}^{(3)})_{z \in \Omega}$ are in $\Hol(\Omega, C^{J}(\URd))$.
   
All this shows that in the sense of~(\ref{eq:HolPHDO.kernel-asymptotic}) we have 
$K_{P_{z}}(x,y)\sim \sum_{\frac32\brak\alpha \leq \brak\beta} K_{\alpha\beta, z}(x,y)$, 
which implies that $(K_{P_{z}})_{z \in \Omega}$ belongs to $\Hol(\Omega, \cK^{*}(\URd))$. Combining this 
    with~(\ref{eq:HolPHDO.invariance.kPz}) and Proposition~\ref{prop:HolPHDO.kernel.characterization}
    then shows that $(P_{z})_{z\in \Omega}$ is a holomorphic family of \psivdos.
\end{proof}

Now, let $(M^{d+1},H)$ be a Heisenberg manifold and let $\cE$ be a smooth vector bundle over $M$. Then Proposition~\ref{prop:HolPHDO.invariance}
allows us to define 
holomorphic families with values in $\Psi^{*}_{H}(M,\cE)$ as follows.
\begin{definition}
 A  family $(P_{z})_{z\in \Omega}\subset \pvdo^{*}(M,\cE)$ is holomorphic when:\smallskip 
 
     (i) The order $m(z)$ of $P_{z}$ is a holomorphic function of $z$; \smallskip 
 
     (ii)  For $\varphi$ and  $\psi$ in $C^\infty_{c}(M)$ with disjoint supports $(\varphi P_{z}\psi)_{z\in \Omega}$ is 
     a holomorphic family of smoothing operators (i.e.~is given by a holomorphic family of smooth distribution kernels); \smallskip 
 
     (iii) For any trivialization $\tau:\cE_{|_{U}}\rightarrow U\times \C^{r}$ over a 
     local Heisenberg chart $\kappa:U \rightarrow V\subset \Rd$ the family $(\kappa_{*}\tau_{*}(P_{z|_{U}}))_{z\in\Omega}$ belongs to 
     $\Hol(\Omega, \pvdo^{*}(V,\C^{r})):=\Hol(\Omega, \pvdo^{*}(V))\otimes \End \C^{r}$.
\end{definition}
All the preceding properties of holomorphic families of \psivdos\ on an open subset of $\Rd$ 
hold \emph{verbatim} for holomorphic families with values in $\Psi_{H}^{*}(M,\cE)$. Moreover, we have: 

\begin{proposition}
    The principal symbol map $\sigma_{*}:\pvdo^{*}(M,\cE) \rightarrow S_{*}(\fg^{*}M,\cE)$ is analytic, in the sense that for any holomorphic family 
    $(P_{z})_{z\in \Omega}\subset \pvdo^{*}(M,\cE)$ the family of symbols $(\sigma_{*}(P_{z}))_{z\in\Omega}$ is in 
    $\Hol(\Omega, C^{\infty}(\fg^{*}M\setminus 0,\End \cE))$. 
\end{proposition}
\begin{proof}
   Let $(P_{z})_{z\in \Omega}\subset \pvdo^{*}(M,\cE)$ be a holomorphic family of \psivdos\ of order $m(z)$ and let us show that 
   the family of symbols $(\sigma_{*}(P_{z}))_{z\in\Omega}$ belongs to $\Hol(\Omega, C^{\infty}(\fg^{*}M\setminus 0,\End \cE))$. 
   Since this a purely local issue we may as well assume that  $(P_{z})_{z\in \Omega}$ is a holomorphic family of scalar \psivdos\ on a local trivializing 
   Heisenberg chart $U\subset \Rd$. 
   
   By Proposition~\ref{prop:HolPHDO.kernel.characterization} we can put the distribution kernel of $P_{z}$ into the form, 
   \begin{equation}
         k_{P_{z}}(x,y) = |\varepsilon_{x}'| K_{P_{z}}(x,-\varepsilon_{x}(y)) + R_{z}(x,y),
   \end{equation}
   with $(K_{P_{z}})_{z\in \Omega}\in \Hol(\Omega,\cK^{*}(\URd))$ of order $\hat{m}(z)=-(m(z)+d+2)$ and 
   $(R_{z})_{z\in \Omega}\in \Hol(\Omega,C^{\infty}(U\times U))$. Let $\varphi \in C^{\infty}(\Rd)$ be such that $\varphi(y)=1$ near $y=0$ and let 
   $p_{z}=(\varphi(y)K_{P_{z}}(x,y))^{\wedge}_{\xiy}$. Then the proof of Proposition~\ref{prop:HolPHDO.characterization-Hol-cK*}
   shows that $(p_{z})_{z\in \Omega}$ is a holomorphic family 
   of symbols. Moreover, we have 
   \begin{equation}
       K_{P_{z}}(x,y)= (p_{z})^{\vee}_{\xiy}(x,y)+(1-\varphi(y))K_{P_{z}}(x,y)=  (p_{z})^{\vee}_{\xiy}(x,y) \ \bmod \Hol(\Omega,C^{\infty}(U\times U)).
        \label{eq:HolPHDO.KPz-local-symbol}
   \end{equation}
   
   Let $z\in \Omega$ and let $K_{\hat{m}(z)}  \in \cK_{\hat{m}(z)}(\URd)$ be the principal kernel of $K_{P_{z}}$. 
   Then~(\ref{eq:HolPHDO.KPz-local-symbol}) and  
   Proposition~\ref{lem:PsiHDO.characterization.Km} show 
   that the leading symbol of $p_{z}$ is the restriction to $\URdo$ of $(K_{\hat{m}(z)})^{\wedge}_{\yxi}$. Since the latter is equal to $\sigma_{m(z)}(P_{z})$, 
   we see that the leading symbol of $p_{z}$ is just $\sigma_{m(z)}(P_{z})$.
   Since $(p_{z})_{z\in \Omega}$ is a holomorphic family of symbols it then follows from Remark~\ref{rem:HolPHDO.analyticity-homogeneous-symbols} 
   that the family $(\sigma_{m(z)}(P_{z}))_{z \in \Omega}$ 
   belongs to $\Hol(\Omega, C^{\infty}(\URdo))$. The proof is thus achieved.
\end{proof}

\subsection{Transposition and holomorphic families of \psivdos}
Let us now look at the analyticity and anti-analyticity of taking transposes and adjoints of \psivdos. 

\begin{proposition}\label{prop:HolPHDO.transpose}
    Let $(P_{z})_{z\in \Omega} \subset \pvdo^{*}(M,\cE)$ be a holomorphic family of \psivdos. Then the transpose family $(P_{z}^{t})\subset 
    \pvdo^{*}(M,\cE^{*})$ is a holomorphic family of \psivdos.  
\end{proposition}
\begin{proof}
    For $z \in \Omega$ let $k_{P_{z}}(x,y)$ denote the distribution kernel of $P_{z}$. The distribution kernel of $P_{z}^{t}$ is 
    $k_{P_{z}}(x,y)=k_{P_{z}}(y,x)^{t}$, hence is represented outside the diagonal by a holomorphic family of smooth kernels. Therefore, we need only to prove 
    the statement for a holomorphic 
    family of scalar \psivdos\ on a Heisenberg chart $U\subset \Rd$, as we shall now suppose that the family $(P_{z})_{z\in \Omega}$ is. In addition, there 
    is no loss of generality in assuming that the order $m(z)$ of $P_{z}$ is such that there exists $\mu\in \R$ so that $\Re \hat{m}(z)\geq \mu$ for any 
    $z \in \Omega$. 
    
    Next, thanks to Proposition~\ref{prop:HolPHDO.kernel.characterization}
    the kernel of $P_{z}$ is of the form
    \begin{equation}
        k_{P_{z}}=|\varepsilon_{x}'|K_{P_{z}}(x,-\varepsilon_{x}(y))+R_{z}(x,y),
    \end{equation}
    with $(K_{P_{z}})_{z \in \Omega}$ in $\Hol(\Omega,\cK^{*}(\URd))$ and $(R_{z})_{z \in \Omega}$ in $\Hol(\Omega, C^{\infty}(U\times U))$. 
    Then the proof of 
    Proposition~\ref{prop:PsiHDO.transpose-chart} in~\cite{Po:BSM1} shows that we can write 
    \begin{equation}
        k_{P_{z}}(x,y)= |\varepsilon_{x}'|K_{P_{z}^{t}}(x,-\varepsilon_{x}(y)) +  
        (1-\chi(x,-\varepsilon_{x}(y)))|\varepsilon_{y}'|K_{P_{z}}(y,-\varepsilon_{y}(x))+R_{z}(y,x),
    \end{equation}
where we have let 
    \begin{equation}
        K_{P_{z}^{t}}(x,y)= \chi(x,y)|\varepsilon_{x}'|^{-1}|\varepsilon_{y}'|K_{P_{z}}(\varepsilon_{x}^{-1}(-y), 
        -\varepsilon_{\varepsilon_{x}^{-1}(-y)}(x)),  
    \end{equation}
and the function  $\chi(x,y)\in C^{\infty}(\URd)$ is supported in the open subset $\cU=\{(x,y)\in \URd; \ \varepsilon_{x}^{-1}(-y)\in U\}$ and is in such 
way to be properly supported with respect $x$ and to be equal to $1$ near $U\times \{0\}$. 
    In particular, we have 
    \begin{equation}
        k_{P_{z}^{t}}(x,y)= |\varepsilon_{x}'|K_{P_{z}^{t}}(x,-\varepsilon_{x}(y)) \qquad \bmod \Hol(\Omega,C^{\infty}(U\times U)).
         \label{eq:HolPHDO.kernel-transpose}
    \end{equation}
    Moreover, as shown in~\cite{Po:BSM1}, for any integer $N$ we have
    \begin{gather}
        K_{P_{z}^{t}}(x,y)= \sum_{\alpha,\beta,\gamma,\delta}^{(N)} K_{\alpha\beta\gamma\delta,z}+ \sum_{j=1}^{4}R_{N,z}(x,y),\\ 
     K_{\alpha\beta\gamma\delta,z}=  a_{\alpha\beta\gamma\delta}(x) y^{\beta+\delta}  
       (\partial^{\gamma}_{x}\partial_{y}^{\alpha}K_{P_{z}})(x,-y)
\end{gather}
where the smooth functions $a_{\alpha\beta\gamma\delta}(x)$ are as in Proposition~\ref{prop:PsiHDO.transpose-chart}, the summation goes over all the multi-orders 
$\alpha$, $\beta$, $\gamma$ and $\delta$ such that $\brak\alpha<N$, $\frac{3}{2}\brak\alpha \leq 
\brak \beta < \frac{3}{2}N$ and $|\gamma|\leq |\delta| \leq 2|\gamma|<2N$ and the remainder terms $R_{N,z}^{(j)}(x,y)$ take the forms:\smallskip 

- $R_{N,z}^{(1)}=\sum_{\brak\alpha=N}\sum_{\brak \beta\dot{=}\frac{3}{2}\brak \alpha} 
    |\varepsilon_{x}'|^{-1}|\varepsilon_{y}'| \int_{0}^{1}r_{N\alpha\beta}(t,x,y)
    (y^{\beta}\partial_{y}^{\alpha}K_{P_{z}})(\varepsilon_{x}^{-1}(-y), \Phi_{t}(x,y))$, where 
the functions $r_{N\alpha\beta}(t,x,y)$ are in $C^{\infty}([0,1]\times \URd)$, the equality $k\dot{=}\frac{3}{2}l$ means that $k$ is equal to $\frac{3}{2}l$ if 
$\frac{3}{2}l$ is integer and  to $\frac{3}{2}l+\frac{1}{2}$ otherwise, and we have let 
$\Phi_{t}(x,y)=-y+t(y-\varepsilon_{\varepsilon_{x}^{-1}(-y)}(x))$;\smallskip

- $R_{N,z}^{(2)}(x,y)= \sum_{\brak \beta \dot{=}\frac{3}{2}N}r_{N\alpha}(x,y)y^{\beta} (\partial_{y}^{\alpha}K_{P_{z}})(\varepsilon_{x}^{-1}(-y),-y)$ 
with $r_{N\alpha}(x,y)$ in $C^{\infty}(\URd)$;\smallskip

- $R_{N,z}^{(3)}(x,y)=\sum_{|\gamma|=N}\sum_{N\leq |\delta|\leq 2N} a_{\alpha\beta\gamma\delta}(x)y^{\beta+\delta} \int_{0}^{1}(1-t)^{N-1} 
  (\partial^{\gamma}_{x}\partial_{y}^{\alpha}K_{P_{z}})(\varepsilon_{t}(x,y),-y)$;\smallskip

  - $R_{N,z}^{(4)}(x,y)=\sum_{\alpha,\beta,\gamma,\delta}^{(N)} (1-\chi(x,y))K_{\alpha\beta\gamma\delta,z}$.\smallskip

Each family $(K_{\alpha\beta\gamma\delta,z})_{z \in \Omega}$ belongs to $\Hol(\Omega,\cK^{*}(\URd))$. Moreover, the remainder term
 $R_{N,z}^{(4)}$ belongs to $\Hol(\Omega, C^{\infty}(\URd))$ and, along similar lines as that of the proof of 
Proposition~\ref{prop:HolPHDO.invariance}, we can show that for any integer $J$ the other remainder terms $(R_{N,z}^{(j)})_{z\in \Omega}$ are 
in $\Hol(\Omega, C^{J}(\URd))$ as soon $N$ is large enough. Therefore, we have 
$K_{P_{z}^{t}} \sim  \sum_{\frac{3}{2}\brak\alpha \leq \brak \beta} \sum_{|\gamma|\leq |\delta| \leq 2|\gamma|}  K_{\alpha\beta\gamma\delta,z}$
in the sense of~(\ref{eq:HolPHDO.kernel-asymptotic}),  which means that $K_{P_{t}^{z}}$ belongs to $\Hol(\Omega,\cK^{*}(\URd))$. Combining this 
with~(\ref{eq:HolPHDO.kernel-transpose}) and Proposition~\ref{prop:HolPHDO.kernel.characterization} 
shows that $(P_{z}^{t})_{z\in \Omega}$ is a holomorphic family of \psivdos. 
\end{proof}

Assume now that $M$ is endowed with a density~$>0$ and $\cE$ with a Hermitian metric. Then Proposition~\ref{prop:HolPHDO.transpose} 
allows us to carried out the proof of Proposition~\ref{prop:PsiHDO.adjoint-manifold} in the holomorphic setting as much so to get: 

\begin{proposition}\label{prop:HolPHDO.adjoint}
    Let $(P_{z})_{z\in \Omega} \subset \pvdo^{*}(M,\cE)$ be a holomorphic family of \psivdos. Then the family 
    $(P_{z}^{*})_{z \in \Omega}\subset \pvdo^{*}(M,\cE^{*})$ is an anti-holomorphic family of \psivdos, in the sense 
    that $(P_{\overline{z}}^{*})_{z \in \Omega}$ is a holomorphic family of \psivdos.
\end{proposition}

\section{Complex powers of hypoelliptic differential operators}
\label{sec.powers1}
In this section we will show that the complex powers of a positive 
hypoelliptic differential operator, \emph{a priori} defined as unbounded operators on $L^{2}(M,\cE)$, give rise 
to a holomorphic family of \psivdos. 

Since we cannot carry out in the Heisenberg setting the standard approach of Seeley~\cite{Se:CPEO} to the complex powers of elliptic operators, we 
will rely on a new approach based on the pseudodifferential representation of the heat kernel of~\cite{BGS:HECRM}, which extended to the Heisenberg 
setting the results of~\cite{Gr:AEHE}. 

This section is divided into two subsections. In the first one we recall the  pseudodifferential representation of the heat kernel of a hypoelliptic 
operator of~\cite{BGS:HECRM}, and in the second one we deal with the complex powers of a positive hypoelliptic operator. 

Throughout this section we let $(M^{d+1},H)$ be a compact Heisenberg manifold endowed with a positive density and let $\cE$ be a Hermitian vector 
bundle over $M$ of rank $r$.

\subsection{Pseudodifferential representation of the heat kernel}
Let $P:C^{\infty}(M,\cE)\rightarrow C^{\infty}(M,\cE)$ be a selfadjoint differential operator of 
even (Heisenberg) order $v$ which is bounded from below and has an invertible principal symbol. In particular,  $P$ is hypoelliptic with loss of 
$\frac{v}{2}$-derivatives by Theorem~\ref{thm:PsiHDO.hypoellipticity}. 

Since $P$ is bounded from below it generates on $L^{2}(M,\cE)$ a heat semigroup $e^{-tP}$, $t\geq 0$.
In fact, for $t>0$ and for any integer $k\geq 1$ the operators 
$P^{k}e^{-tP}$ and $e^{-tP}P^{k}$ coincide and are bounded. Thus $e^{-tP}$ maps continuously 
$L^{2}(M,\cE)$ to $\cap_{k\geq 0} \dom P^{k}$, which is just $C^{\infty}(M,\cE)$ since $P$ is hypoelliptic with loss of $\frac{v}{2}$-derivatives.  

Moreover, as $e^{-tP}$ is selfadjoint it defines by duality a continuous map from $\cD'(M,\cE)$ to $L^{2}(M,\cE)$. Since
$e^{-tP}=e^{-tP/2}e^{-tP/2}$ it follows that $e^{-tP}$ extends to a continuous map from 
$\cD'(M,\cE)$ to $C^{\infty}(M,\cE)$, that is   
 $e^{-tP}$ is smoothing. In particular,  $e^{-tP}$ is given by a  smooth kernel $k_{t}(x,y)$  in 
$C^\infty(M\times M, \cE\boxtimes(\cE \otimes|\Lambda|(M)))$, 
where $|\Lambda|(M)$ denotes the bundle of densities on $M$. 

On the other hand, the heat semi-group allows us to invert the heat equation. Indeed, the operator given by 
\begin{equation}
    Q_{0}f(x,t)=\int_{0}^\infty e^{-sP} f(x,t-s)dt, \qquad f \in C^\infty_{c}(M\times \R, \cE), 
     \label{eq:volterra.inverse-heat-operator}
\end{equation}
maps continuously into $C^{0}(\R, L^{2}(M,\cE)) \subset \cD(M\times\R, \cE)$ and satisfies
\begin{equation}
    (P+\partial_{t})Q_{0}f = Q_{0}(P+\partial_{t})f=f \qquad \forall u \in C^\infty_{c}(M\times\R,\cE).
\end{equation} 

Notice that the operator $Q_{0}$  has the \emph{Volterra property} of~\cite{Pi:COPDTV}, i.e.~it has 
a distribution kernel of the form $K_{Q_{0}}(x,y,t-s)$ with $K_{Q_{0}}(x,y,t)$ supported outside the region $\{t<0\}$. 
Indeed, at the level of distribution kernels the formula~(\ref{eq:volterra.inverse-heat-operator}) implies that we have 
\begin{equation}
    K_{Q_{0}}(x,y,t) = \left\{ 
    \begin{array}{ll}
         k_{t}(x,y) & \quad \text{if $t> 0$},\\
        0 &  \quad \text{if $t<0$}. \label{eq:Powers1.kernel-inverse-heat-operator-heat-kernel}
    \end{array}\right. 
\end{equation}

The above equalities are the main motivation for using pseudodifferential techniques to study the heat kernel $k_{t}(x,y)$. 
The idea is to consider a class of \psivdo's, the Volterra \psivdo's, taking into account:  \smallskip 

(i)  The aforementioned Volterra property;\smallskip 

(ii) The parabolic homogeneity of the heat operator $P+ \partial_{t}$, i.e.~the homogeneity with respect 
             to the dilations,
             \begin{equation}
                  \lambda.(\xi,\tau)=(\lambda.\xi,\lambda^{v}\tau), \qquad (\xi,\tau)\in \R^{d+2}, \ \lambda\neq 0. 
             \end{equation}

 In the sequel for $g\in \cS'(\R^{d+2})$  and $\lambda\neq 0$ we let $g_{\lambda}$ be the tempered distribution 
 defined by   
\begin{equation}
    \acou{g_{\lambda}( \xi,\tau)}{f(\xi,\tau)} =   |\lambda|^{-(d+2+v)}\acou{g(\xi,\tau)} {f(\lambda^{-1}\xi, 
    \lambda^{-v}\tau)}, \qquad f\in \cS(\R^{d+2}).               
\end{equation}

\begin{definition}
    A distribution $ g\in \cS'(\R^{d+2})$ is parabolic 
homogeneous of degree $m$, $m\in \Z$, when we have $g_{\lambda}=\lambda^m g$ for any $\lambda \neq 0$.  
\end{definition}

Let $\C_{-}$ denote the complex halfplane $\{\Im \tau <0\}$ with closure $\overline{\C_{-}}\subset \C$. Then:  
\begin{lemma}[{\cite[Prop.~1.9]{BGS:HECRM}}] \label{lem:volterra.volterra-extension}
Let $q(\xi,\tau)\in C^\infty((\R^{d+1}\times\R)\setminus0)$ be parabolic homogeneous of degree $m$ and such that: \smallskip 

(i) $q$ extends to a smooth function on $(\R^{d+1}\times\overline{\C_{-}})\setminus0$ which restricts on 
$\Rd\times \C_{-}$ to an element of $C^{\infty}(\Rd)\hotimes \Hol(\C_{-})$.
    \smallskip 

\noindent Then  there exists a unique distribution $g\in \cS'(\R^{d+2})$ agreeing with $q$ 
on $\R^{n+1}\setminus 0$ and such that: \smallskip 

(ii) $g$ is parabolic homogeneous of degree $m$; \smallskip 

(iii) The inverse Fourier transform $\check g(x,t)$ vanishes for $t<0$.
\end{lemma}

Let $U$ be an open subset of $\R^{d+1}$ together with a hyperplane bundle $H\subset TU$ and $H$- frame $X_{0},\ldots, 
X_{d}$ of $TU$. We define Volterra symbols and Volterra \psivdos\ as follows. 

\begin{definition}
    $S_{\op{v},m}(U\times\Rd\times\R_{(v)})$, $m\in\Z$,  consists of functions $q(x,\xi,\tau)$ in $C^{\infty}(U\times(\R^{d+2}\setminus 0))$ such 
    that:\smallskip
    
    (i) $q(x,\lambda.\xi,\lambda^{v}\tau)=\lambda^{m}q(x,\xi,\tau)$ for any $(x,\xi,\tau)\in U\times(\R^{d+2}\setminus 0)$ and any $\lambda\neq 0$;\smallskip
    
    (ii) $q(x,\xi,\tau)$ extends to an element of $C^{\infty}(U\times[(\Rd\times\overline{\C_{-}})\setminus 0])$ in such way that its restricts to 
    an element of $C^{\infty}(\URd)\hotimes \Hol(\C_{-})$ on $U\times \Rd\times \C_{-}$.
\end{definition}

\begin{definition}
    $S_{\op{v}}^m(U\times\Rd\times\R_{(v)})$, $m\in\Z$,  consists of functions $q(x,\xi,\tau)$ in 
    $C^{\infty}(U\times\Rd\times\R)$ with an asymptotic expansion  $q \sim \sum_{j\geq 0} q_{m-j}$, $q_{l}\in S_{\op{v},m}(U\times\Rd\times\R_{(v)})$, 
   in the sense that, for any integer $N$ and any compact $K\subset U$,  we have  
            \begin{equation}
                |\partial^{\alpha}_{x}\partial^{\beta}_{\xi} \partial^k_{\tau}(q-\sum_{j< N} 
            q_{m-j})(x,\xi,\tau) | 
                \leq C_{NK\alpha\beta k} (\|\xi\|+|\tau|^{1/v})^{m-N-\brak\beta-v k},
                          \label{eq:volterra.asymptotic-symbols}
            \end{equation}
            for  $x\in K$ and $|\xi|+|\tau|^{\frac{1}{v}}>1$.
\end{definition}

\begin{definition}
  Let $q(x,\xi,\tau)\in S_{\op{v},m}(U\times\Rd\times\R_{(v)})$ and let $g\in 
  C^{\infty}(U)\hotimes \cS'(\Rd)$ be the unique homogeneous extension of $q$ provided by Lemma~\ref{lem:volterra.volterra-extension}. Then we let 
  $\check{q_{m}}(x,y,t)$ denote the inverse Fourier transform of $g(x,\xi,\tau)$ with respect to the variables $(\xi,\tau)$.
\end{definition}

\begin{remark}
    The above definition makes sense since it follows from the proof of Lemma~\ref{lem:volterra.volterra-extension} 
    in \cite{BGS:HECRM} that the extension process of Lemma~\ref{lem:volterra.volterra-extension}  applied to every symbol $q(x,.,.)$, $x \in U$, 
    is smooth with respect to $x$, so really gives rise to an element of $C^\infty(U)\hotimes\cS'(\R^{d+2})$. 
 \end{remark}

\begin{definition}\label{def:volterra.PsiDO}
    $\pvhdo^m(U\times\R_{(v)})$, $m\in\Z$,  consists of continuous operators 
    $Q:C_{c}^\infty(U_{x}\times\R_{t}) \rightarrow C^\infty(U_{x}\times\R_{t})$ such that $Q$ has the Volterra property and can be put into the form
\begin{equation}
    Q=q(x,-iX,D_{t})+R
\end{equation}  
   with $q$ in $S^m_{\op{v}}(U\times\Rd\times\R_{(v)})$ and $R$ in $\psinf(U\times \R)$. 
\end{definition}

It is immediate to extend the properties of \psivdos\ on $U$ alluded to in Section~\ref{sec:PsiHDO} to Volterra \psivdos\ on $U\times \R$ except for 
the asymptotic completeness as in Lemma~\ref{lem:PsiHDO.asymptotic-completeness}, 
which is crucial for constructing. The problem is that the cut-off arguments of the classical proof 
the asymptotic completeness of standard \psidos\ cannot be carried through in Volterra setting because we require analyticity with respect to the time 
covariable. A proof of the asymptotic completeness of Volterra \psidos\ is given in~\cite{Pi:COPDTV}, but simpler proofs which can be carried out 
\emph{verbatim} for Volterra \psivdos\ can be found in~\cite{Po:JAM1}.  

On the other hand, thanks to the Volterra property the kernels of \psivdos\ can be characterized as follows.

\begin{definition}
$\cK_{\op{v}, m}(U\times\Rd\times\R_{(v)})$, $m\in \Z$, consists of distributions $K(x,y,t)$ in 
 $C^{\infty}(U)\otimes\cS'_{\reg}(\R^{d+2})$ such that:\smallskip 
     
    \indent (i) The support of $K(x,y,t)$ is contained in $U\times\Rd\times\R_{+}$;\smallskip 
    
    \indent (ii) $K(x,y,t)$ is parabolic homogeneous of degree $m$ with respect to the variables $(y,t)$. 
\end{definition}

\begin{definition}
 $\cK_{\op{v}}^{m}(U\times\Rd\times\R_{(v)})$, $m\in \Z$, is the space of distributions $K(x,y,t)$ in $\cD'(U\times\R^{d+2})$ which  
 admit an asymptotic expansion $K \sim \sum_{j \geq 0} K_{m+j}$, $K_{m+j}\in \cK_{\op{v}, m+j}(U\times\R^{d+2})$, where
 $\sim$ is taken in the sense of~(\ref{eq:PsiHDO.asymptotics-kernel}).  
\end{definition}

In the sequel, for $x\in U$ we let $\psi_{x}$ and $\varepsilon_{x}$ respectively denote the changes of variable to the privileged coordinates and 
the Heisenberg coordinates at $x$. Then,  along the same lines as that of the proofs of 
Proposition~\ref{prop:PsiVDO.characterisation-kernel1} and Proposition~\ref{prop:PsiVDO.characterisation-kernel2}, 
we obtain the following characterization of Volterra \psivdos.

\begin{proposition}\label{prop:volterra.kernel-charaterization} 
    Let $Q:C_{c}^\infty(U_{x}\times\R_{t}) \rightarrow C^\infty(U_{x}\times\R_{t})$ be a continuous operator with distribution kernel $k_{Q}(x,t;y,s)$. Then 
    the following are equivalent:\smallskip 
    
    (i) The operator $Q$ belongs to $\pvhdo^{m}(U\times \R)$;\smallskip
    
    (ii) The kernel of $Q$ can be put into the form,
    \begin{equation}
        k_{Q}(x,t;y,s)=|\psi_{x}'|K(x,-\psi_{x}(y),t-s)+R(x,y,t-s),
     \label{eq:Powers1.characterization-Volterra-PsiHDO's-psix}
\end{equation}
    for some $K(x,y,t)\in \cK_{\op{v}}^{\hat m}(U\times\Rd\times\R_{(v)})$, $\hat{m}=-(m+d+2+v)$  and $R\in C^{\infty}(U\times 
    \Rd\times\R)$.

(iii) The kernel of $Q$ can be put into the form,
    \begin{equation}
        k_{Q}(x,t;y,s)=|\varepsilon_{x}'|K_{Q}(x,-\varepsilon_{x}(y),t-s)+R(x,y,t-s),
     \label{eq:Powers1.characterization-Volterra-PsiHDO's}
\end{equation}
    for some $K_{Q}(x,y,t)\in \cK_{\op{v}}^{\hat m}(U\times\Rd\times\R_{(v)})$, $\hat{m}=-(m+d+2+v)$  and $R\in C^{\infty}(U\times 
    \Rd\times\R)$.
\end{proposition}

An interesting consequence of Proposition~\ref{prop:volterra.kernel-charaterization} 
is the following small time asymptotics for the kernel of a Volterra \psivdo. 

\begin{proposition}[{\cite[Thm.~4.5]{BGS:HECRM}}]\label{prop:Volterra.asymptotics-kQ}
Let $Q\in \pvhdo^{m}(U\times\R_{(v)})$ have symbol $q \sim \sum_{j \geq 0} q_{m-j}$ and kernel $k_{Q}(x,y,t-s)$. 
Then as $t\rightarrow 0^{+}$ the following asymptotics 
holds in  $C^\infty(U)$, 
\begin{equation}
    k_{Q}(x,x,t) \sim t^{-\frac{2[\frac{m}2]+d+4}{v}} \sum_{j\geq 0} t^j 
    |\varepsilon_{x}'|(q_{2[\frac{m}2]-2j})^{\vee}_{(\xi,\tau)\rightarrow (y,t)}(x,0,1). 
    \label{eq:volterra.asymptotics-Q}
\end{equation}
\end{proposition}

On the other hand, using Proposition~\ref{prop:volterra.kernel-charaterization} and arguing along similar lines that of the proof 
of Proposition~\ref{prop:PsiHDO.invariance} in~\cite{Po:BSM1} allows us to prove:

 \begin{proposition}\label{prop:Powers1.invariance}
     Let $U$ (resp.~$\tilde{U}$) be an open subset of $\Rd$ together with a hyperplane bundle $H\subset TU$ (resp.~$\tilde{H}\subset T\tilde{U}$) and a 
    $H$-frame of $TU$ 
    (resp.~a $\tilde{H}$-frame of $T\tilde{U}$). Let $\phi:(U,H)\rightarrow (\tilde{U},\tilde{H})$ be a Heisenberg diffeomorphism and let $\tilde{Q}\in 
    \Psi_{\tilde{H},\op{v}}^{m}(\tilde{U}\times\R_{(v)})$.\smallskip 
    
   1) The operator $Q=(\phi\oplus 1_{\R})^{*}\tilde{Q}:C^{\infty}_{c}(\URd\times \R)\rightarrow C^{\infty}(\URd\times \R)$ 
   belongs to $\pvhdo^{m}(U\times\R_{(v)})$.\smallskip 
    
   2) If the distribution kernel of $\tilde{P}$ is of the form~(\ref{eq:Powers1.characterization-Volterra-PsiHDO's}) 
   with $K_{\tilde{Q}}(\tilde{x},\tilde{y},t)\in \cK_{\op{v}}^{\hat m}(\tilde{U}\times\Rd\times\R_{(v)})$ then the distribution kernel of 
   $P$ can be written in the form~(\ref{eq:Powers1.characterization-Volterra-PsiHDO's}) with $K_{Q}(x,y,t)\in \cK_{\op{v}}^{\hat m}(U\times\Rd\times\R_{(v)})$ such that 
   \begin{equation}
       K_{Q}(x,y,t) \sim \sum_{\brak\beta\geq \frac{3}{2}\brak\alpha} \frac{1}{\alpha!\beta!} 
       a_{\alpha\beta}(x)y^{\beta}(\partial_{\tilde{y}}^{\beta}K_{\tilde{Q}})(\phi(x),\phi_{H}'(x)y,t),
        \label{eq:Powers1.asymptotic-expansion-KQ}
   \end{equation}
   where the functions $a_{\alpha\beta}(x)$ are as in~(\ref{eq:PsiHDO.asymptotic-expansion-KP}).
 \end{proposition}

This allows us to define Volterra \psivdos\ on the manifold $M\times \R$ and acting on the sections of the 
bundle $\cE$ (or rather on the sections of the pullback of $\cE$ by the projection $M\times \R\rightarrow M$, again denoted $\cE$).

\begin{definition}
  $\pvhdo^{m}(M\times\R_{(v)},\cE)$, $m\in \Z$, consists of continuous operators $Q:C^{\infty}_{c}(M\times \R,\cE)\rightarrow C^{\infty}(M\times \R,\cE)$ such 
  that:\smallskip 
  
  (i) $Q$ has the Volterra property;\smallskip
  
  (ii) The distribution kernel of $Q$ is smooth off the diagonal of $(M\times \R)\times(M\times \R)$;\smallskip
  
  (iii) For any trivialization $\tau:\cE_{|_{U}}\rightarrow U\times \C^{r}$ of $\cE$ over a 
     local Heisenberg chart $\kappa:U \rightarrow V\subset \Rd$ the operator $(\kappa\otimes\op{id})_{*}\tau_{*}(Q_{|_{U\times \R}})$ belongs to
     $\pvhdo^{m}(V\times\R_{(v)}, \C^{r}):=\pvhdo^{m}(V\times \R_{(v)})\otimes \End \C^{r}$. 
\end{definition}

Using Proposition~\ref{prop:Powers1.invariance} 
we can define the global principal symbol of a Volterra \psivdo\ as follows. Let $\fg^{*}M$ denote the dual bundle of the Lie algebra 
bundle $\fg M$ of $M$ and consider the canonical projection 
$\pi:\fg^{*}M \times \R\rightarrow M$. 

\begin{definition}
 $S_{\op{v},m}(\fg^{*}M\times \R_{(v)},\cE)$, $m\in \Z$, is the subspace of $C^{\infty}((\fg^{*}M\times \R)\setminus 0, \pi^{*}\End 
 \cE)$ consisting of sections $q(x,\xi,\tau)$  such that:\smallskip
 
 (i) $q(x,\lambda.\xi,\lambda^{v}\tau)=\lambda^{m}q(x,\xi,\tau)$ for any 
 $(x,\xi,\tau)\in (\fg^{*}M\times \R)\setminus 0$ and any $\lambda \in \R\setminus 0$;\smallskip 
 
 (ii) $q(x,\xi,\tau)$ extends to a section in $C^{\infty}((\fg^{*}M\times \bar\C_{-})\setminus 0, \pi^{*}\End  \cE)$ which restricts to an element of 
 $C^{\infty}(\fg^{*}M, \pi^{*}\End  \cE)\hotimes \Hol(\C_{-})$ on $\fg^{*}M\times \C_{-}$.
\end{definition}

Using~(\ref{eq:Powers1.asymptotic-expansion-KQ}) and arguing as in the proof of Proposition~\ref{prop:PsiHDO.principal-symbol} in~\cite{Po:BSM1} we get:

\begin{proposition}\label{prop:Powers1.principal-symbol}
    For any $Q \in \pvhdo^{m}(M\times\times\R_{(v)},\cE)$ there is a unique symbol $\sigma_{m}(Q)\in S_{\op{v},m}(\fg^{*}M \times \times\R_{(v)}, \cE)$ 
    such that if in a local trivializing Heisenberg 
    chart $U\subset \Rd$ we let $K_{Q,\hat{m}}(x,y,t)\in \cK_{\hat{m}}(\URd)$ be the leading kernel for the kernel $K_{Q}(x,y,t)$  in the 
    form~(\ref{eq:Powers1.characterization-Volterra-PsiHDO's})   
    for $Q$,  then we have
    \begin{equation}
       \sigma_{m}(Q)(x,\xi,\tau)=[K_{Q,\hat{m}}]^{\wedge}_{(y,t)\rightarrow (\xi,\tau)}(x,\xi,\tau), \qquad (x,\xi)\in U\times \Rdo. 
    \end{equation}
    
    Equivalently, on any trivializing Heisenberg coordinates   centered at $a\in M$ the symbol $\sigma_{m}(Q)(a,.,.)$ coincides 
    with the (local) principal symbol of $Q$ at $x=0$. 
\end{proposition}

\begin{definition}\label{def:Powers1.principal-symbol}
   For $Q \in \pvhdo^{m}(M\times\R_{(v)},\cE)$ the symbol $\sigma_{m}(Q)\in S_{\op{v},m}(\fg^{*}M \times \R, \cE)$  provided by 
   Proposition~\ref{prop:Powers1.principal-symbol} is called the (global) principal symbol of $Q$. 
\end{definition}

Granted this we can define the model operator of a Volterra \psivdo\ as follows. 
\begin{definition}
 Let $Q \in \pvhdo^{m}(M\times\R_{(v)},\cE)$ have principal symbol $\sigma_{m}(Q)$.   Then the model operator of $Q$ at $a \in M$ is the left-convolution 
 operator by  $\sigma_{m}(Q)^{\vee}_{(\xi,\tau)\rightarrow (y,t)}(a,.,.)$, i.e.~the continuous endomorphism $Q^{a}:\cS(G_{a}M\times \R,\cE_{a})\rightarrow 
    \cS(G_{a}M\times \R,\cE_{a})$ given by 
    \begin{equation}
        Q^{a}u(y,t)=\acou{\sigma_{m}(Q)^{\vee}_{(\xi,\tau)\rightarrow (y,t)}(x,z,t)}{u(y.z^{-1},t-s)}, \qquad u \in \cS(G_{a}M,\cE_{a}). 
    \end{equation}
\end{definition}
\begin{remark}
    The model operator $Q^{a}$ can be defined as an endomorphism of $\cS(G_{a}M,\cE_{a})$, not just as an endomorphism of $\cS_{0}(G_{a}M,\cE_{a})$ as in 
    Definition~\ref{def:PsiHDO.model-operator}, because $\sigma_{m}(Q)^{\vee}_{(\xi,\tau)\rightarrow (y,t)}(a,.,.)$ makes sense as an element of 
    $\cS'(G_{a}M,\cE_{a})$. 
\end{remark}

\begin{proposition}
  The group laws on the fibers of $GM\times \R$ give rise to a convolution product,
    \begin{equation}
        *: S_{\op{v},m_{1}}(\fg^{*}M\times\R_{(v)},\cE) \times S_{\op{v},m_{2}}(\fg^{*}M\times\R_{(v)},\cE) 
        \longrightarrow S_{\op{v},m_{1}+m_{2}}(\fg^{*}M\times\R_{(v)},\cE),
         \label{eq:Powers1.volterra-Heisenberg-symbol-product}
    \end{equation}
such that for any symbols $q_{m_{j}}\in S_{\op{v},m_{j}}(\fg^{*}M\times\R_{(v)},\cE)$, $j=1,2$, we have
    \begin{gather}
        q_{m_{1}}*q_{m_{2}}(x,\xi,\tau)=[q_{m_{1}}(x,.,.)*^{x}q_{m_{2}}(x,.,.)](\xi,\tau), \qquad (x,\xi,\tau)\in (\fg^{*}M\times \R)\setminus 0,
    \end{gather}
 where $*^{x}$ denote the convolution product for symbols on $G_{x}M\times \R$.
\end{proposition}
 
In a local trivializing Heisenberg chart the symbolic calculus for Volterra \psivdos\ reduces the existence of a Volterra \psivdo\ parametrix to the 
invertibility of the local and global principal symbols. Therefore, we obtain:

\begin{proposition}
    Let $Q\in \pvhdo^{m}(M\times\R_{(v)},\cE)$, $m \in \Z$. Then we have equivalence:\smallskip 
    
    (i) The principal symbol of $Q$ is invertible with respect to the product~(\ref{eq:Powers1.volterra-Heisenberg-symbol-product}) 
    of Volterra-Heisenberg symbols;\smallskip 
    
    (ii) The operator $Q$ admits a parametrix in $\pvhdo^{-m}(M\times\R_{(v)},\cE)$.
\end{proposition}

In the case of the heat operator $P+\partial_{t}$, comparing a parametrix with the inverse~(\ref{eq:volterra.inverse-heat-operator}) and 
using~(\ref{eq:Powers1.kernel-inverse-heat-operator-heat-kernel}) allows us to obtain the pseudodifferential representation of the heat kernel of $P$ 
below.

\begin{proposition}[{\cite[pp.~362--363]{BGS:HECRM}}]\label{thm:volterra.inverse} 
Suppose that the principal symbol of $P+\partial_{t}$ is an invertible Volterra-Heisenberg symbols. Then:\smallskip  

1) The heat operator  $P+\partial_{t}$ has an inverse  $(P+\partial_{t})^{-1}$ in $\pvhdo^{-v}(M\times\R_{(v)},\cE)$.\smallskip

2) Let $K_{(P+\partial_{t})^{-1}}(x,y,t-s)$ denote the kernel of $(P+\partial_{t})^{-1}$. Then the heat kernel $k_{t}(x,y)$ of $P$ satisfies
\begin{equation}
    k_{t}(x,y)=K_{(P+\partial_{t})^{-1}}(x,y,t) \quad \text{for $t>0$.}
     \label{eq:Powers1.kernel-inverse-heat-operator-heat-kernel2}
\end{equation}
\end{proposition}
Combining this with Proposition~\ref{prop:Volterra.asymptotics-kQ} then gives the heat kernel asymptotics for $P$ in the form below.

\begin{proposition}[{\cite[Thm.~5.6]{BGS:HECRM}}]\label{thm:Powers1.heat-kernel-asymptotics}
  If the principal symbol of $P+\partial_{t}$ is an invertible Volterra-Heisenberg symbol, then as $t\rightarrow 0^{+}$  the following asymptotics holds in 
  $C^{\infty}(M,(\End \cE)\otimes|\Lambda|(M))$, 
    \begin{equation}
     k_{t}(x,x) \sim t^{-\frac{d+2}{v}} \sum t^{\frac{2j}{v}} a_{j}(P)(x), \qquad  
    a_{j}(P)(x) =|\varepsilon_{x}'|(q_{-v-2j})^{\vee}_{\xitauyt}(x,0,1), 
    \end{equation}
where the equality on the right shows how to compute $a_{j}(P)(x)$ in a local trivializing Heisenberg chart by means of the symbol 
$q_{-v-2j}(x,\xi,\tau)$ of degree $-v-2j$ of any parametrix of $P+\partial_{t}$ in $\pvhdo^{-v}(M\times\R_{(v)},\cE)$. 
\end{proposition}
We will give in Section~\ref{sec:Rockland-heat} criterions for the invertibility  principal symbol of $P+\partial_{t}$ to be invertible. Nevertheless, 
in the case of a sublaplacian we have: 

\begin{proposition}[{\cite[Thm.~5.22]{BGS:HECRM}}]\label{thm:Powers1.heat-sublaplacians}
 Let $\Delta: C^{\infty}(M,\cE) \rightarrow C^{\infty}(M,\cE)$ be a selfadjoint sublaplacian which is bounded from below  and assume that 
 the condition~(\ref{eq:Rockland.sublaplacian'}) is satisfied at every point of $M$. Then the 
 principal symbol of $\Delta+\partial_{t}$ is an invertible Volterra-Heisenberg symbol, hence Proposition~\ref{thm:volterra.inverse} 
    and Proposition~\ref{thm:Powers1.heat-kernel-asymptotics} hold for $\Delta$. 
 \end{proposition}

\begin{remark}
  As shown by Proposition~\ref{thm:PsiDO.Rockland-sublaplacian} the condition~(\ref{eq:Rockland.sublaplacian'}) at a point of $M$ is equivalent to 
  the Rockland condition for $\Delta$ when the Levi form is degenerate, but is a stronger  
  condition when the Levi form is nondegenerate. It will be shown in Section~\ref{sec:Rockland-heat} that when the Levi form is everywhere 
  nondegenerate the Rockland condition for $\Delta$ is enough to insure us the 
  invertibility of the principal symbol of $\Delta+\partial_{t}$ even when the Levi form is nondegenerate at a point. In addition, it will be 
   also shown in~\cite{Po:JFA3} that any selfadjoint sublaplacian with an invertible principal symbol  is bounded from below, so the assumption on the 
   boundedness from below of $\Delta$ is in fact superfluous.
\end{remark}

\begin{example}
  The above result is true for the following sublaplacians: \smallskip 
  
 (a) A selfadjoint sum of squares $\Delta=X_{1}X_{1}^{*}+\ldots+X_{m}X_{m}^{*}$ where $X_{1},\ldots,X_{d}$ span $H$; \smallskip 
  
 (b) The Kohn Laplacian on a CR manifold and acting on $(p,q)$-forms under condition $Y(q)$;\smallskip 
  
  (c) The horizontal sublaplacian on a Heisenberg manifold acting on horizontal forms of degree $k$ under condition $X(k)$;\smallskip 
   
   (d) The conformal powers of the horizontal sublaplacian acting on functions on a strictly pseudoconvex CR manifold.
\end{example}
 
\subsection{Complex powers}
Let $P:C^{\infty}(M,\cE)\rightarrow C^{\infty}(M,\cE)$ be a selfadjoint differential operator of 
even (Heisenberg) order $v$ such that $P$ has an invertible principal symbol and is positive, i.e.,~we have $ \acou{Pu}{u}\geq 0$ for any $u\in C^{\infty}(M,\cE)$.
Let $\Pi_{0}(P)$ be the orthogonal projection onto $\ker P$. Then the operator $P_{0}:=(1-\Pi_{0}(P))P+\Pi_{0}(P)$ is 
selfadjoint with spectrum contained in $[c,\infty)$ for some $c>0$. Thus by standard functional 
calculus, for any $s\in \C$, the power 
$P_{0}^{s}$ is a well defined unbounded operator on $L^{2}(M,\cE)$. We then define the power $P^{s}$, $s \in \C$, by letting 
\begin{equation}
    P^{s}=(1-\Pi_{0}(P))P_{0}^{s}=P_{0}^{s}-\Pi_{0}(P),
     \label{eq:Powers1.definition}
\end{equation}
so that  $P^{s}$ coincides with $P_{0}^{s}$ on $(\ker P)^{\perp}$ and is zero on $\ker 
P$.  In particular, we have $P^{0}=1-\Pi_{0}(P)$ and $P^{-1}$ is the partial inverse of $P$.

The key result of this section is the following. 

\begin{theorem}\label{thm:Powers1.main} 
Suppose that the principal symbol of $P+\partial_{t}$ is an invertible Volterra-Heisenberg symbol. Then:\smallskip 
    
    (i) For any $s \in \C$ the operator $P^{s}$ defined by~(\ref{eq:Powers1.definition}) is a \psivdo\ of order $vs$;\smallskip 
    
    (ii) The family $(P^{s})_{s\in \C}$ forms a holomorphic 1-parameter group of \psivdos. 
\end{theorem}
\begin{proof}
Let us first assume that $\cE$ is the trivial line bundle over $M$, so that $P$ is a scalar operator.  For  $\Re s >0$ the function $x\rightarrow x^{-s}$ is 
bounded on $[0,\infty)$, so the operators $P_{0}^{-s}$ and 
$P^{-s}$ are bounded. Moreover, by the Mellin formula we have
\begin{equation}
     P^{-s}= (1-\Pi_{0}(P))P_{0}^{s}=\frac{1}{\Gamma(s)} \int_{0}^{\infty} t^{s}(1-\Pi_{0}(P))e^{-tP}\frac{dt}{t}. 
     \label{eq:Powers1.}
\end{equation}
This leads us to define  
\begin{equation}
    A_{s}=\int_{0}^{1}t^{s-1}e^{-tP}dt,  \qquad \Re s>0.
     \label{eq:Powers1.Ds}
\end{equation}
Then we have 
\begin{equation}
    \begin{split}
      \Gamma(s)P^{-s}-A_{s} & 
      = \int_{0}^{1} t^{s-1}\Pi_{0}(P)e^{-tP}dt + \int_{1}^{\infty} t^{s-1} (1-\Pi_{0}(P)) e^{-tP} dt,  \\  
      &=\frac{1}{2}\Pi_{0}(P)+ e^{-P/2} (\int_{0}^{\infty} (1+t)^{s-1}e^{-tP}dt)e^{-P/2}.
     \end{split}
\end{equation}
Since $\Pi_{0}(P)$ and $e^{-P/2}$ are smoothing operators and  $(\int_{0}^{\infty} (1+t)^{s-1}e^{-tP}dt)_{\Re s>0}$ is a 
holomorphic family of bounded operators on $L^{2}(M)$, we get
\begin{equation}
     (\Gamma(s)P^{-s}-A_{s})_{\Re s>0} \in \Hol(\Re s>0, \psinf(M)).
    \label{eq:Powers1.relation-As}
\end{equation}

Let us now show that $(A_{s})_{\Re s>0}$ defined by~(\ref{eq:Powers1.Ds}) is a holomorphic family of \psivdos\ such that $\ord A_{s}=-vs$. To this  end 
observe that, in terms of distribution kernels, the formula~(\ref{eq:Powers1.Ds}) means that $A_{s}$ has distribution kernel 

\begin{equation}
    k_{A_{s}}(x,y)= \int_{0}^{1}t^{s-1}k_{t}(x,y)dt.
\end{equation}
where $k_{t}(x,y)$ denotes the heat kernel of $P$.

On the other hand, since $P$ is bounded from below and the principal symbol of $P+\partial_{t}$ is invertible, Theorem~\ref{thm:volterra.inverse} 
tells us that $P+\partial_{t}$ has an inverse $Q_{0}:=(P+\partial_{t})^{-1}$ in
$\pvhdo^{-v}(M\times\R_{(v)},\cE)$ and that the distribution kernel $K_{Q_{0}}(x,y,t-s)$ of $Q_{0}$ is related to the heat kernel of $P$ by 
\begin{equation}
    K_{Q_{0}}(x,y,t)=k_{t}(x,y) \qquad \text{for $t>0$}.
\end{equation}
Therefore, for $\Re s>0$ we have
\begin{equation}
      k_{A_{s}}(x,y)= \int_{0}^{1}t^{s-1}K_{Q_{0}}(x,y,t)dt..
      \label{eq:Powers1.kAs}
\end{equation}

Let $\varphi$ and $\psi$ be smooth functions on $M$ 
with disjoint supports. Then using~(\ref{eq:Powers1.kAs}) we see that $\varphi A_{s}\psi$ has distribution kernel 
\begin{equation}
    k_{\varphi A_{s}\psi}(x,y)= \int_{0}^{1}t^{s-1}\varphi(x)K_{Q_{0}}(x,y,t)\psi(y)dt.
     \label{eq:Powers1.phi-kDs-psi}
\end{equation}
Since the distribution kernel of a Volterra-\psivdo\ is smooth off the diagonal of $(M\times \R)\times (M\times \R)$ the distribution $K_{Q_{0}}(x,y,t)$ is smooth on 
the region $\{x\neq y\}\times \R$, so~(\ref{eq:Powers1.phi-kDs-psi}) defines a holomorphic family of smooth kernels. Thus, 
\begin{equation}
    (\varphi A_{s}\psi)_{\Re s>0} \in \Hol(\Re s>0,\Psi^{-\infty}(M)).
    \label{eq:Powers1.smoothing}
\end{equation}
 
Next, the following holds. 

\begin{lemma}\label{lem:Powers1.analyticity-Bs}
 Let $V\subset \Rd$ be a Heisenberg chart, let $Q\in \pvhdo^{-v}(V\times\R_{(v)})$ have distribution kernel 
$K_{Q}(x,y,t-s)$ and for $\Re s>0$ let $B_{s}:C_{c}(V)\rightarrow C(V)$  be given by the distribution kernel,
\begin{equation}
    k_{B_{s}}(x,y)= \int_{0}^{1}t^{s-1}K_{Q}(x,y,t)dt, \quad \Re s>0.
\end{equation} 
Then $(B_{s})_{\Re s >0}$ is a holomorphic family of \psivdos\ such that $\ord B_{s}=-vs$
\end{lemma}
\begin{proof}[Proof of the lemma]
Let $\varepsilon_{x}$ denote the change to the Heisenberg coordinates at $x$. By Proposition~\ref{prop:volterra.kernel-charaterization} 
on $V\times V\times \R$ the distribution $K_{Q}(x,y,t)$ is of the form
\begin{equation}
    K_{Q}(x,y,t)=|\epsilon_{x}'|K(x,-\varepsilon_{x}(y),t)+R(x,y,t),
\end{equation}
for some $K\in \cK_{\op{v}}^{-(d+2)}(V\times \Rd\times\R_{(v)})$ and some $R\in C^{\infty}(U\times U\times \R)$. Let us write 
$K\sim \sum_{j\geq 0} K_{j-(d+2)}$ 
with $K_{l}\in \cK_{\op{v},l}(V\times \Rd\times\R_{(v)})$. Thus, for any integer $N$, as soon as $J$ large enough we have
\begin{equation}
    K(x,y,t) =\sum_{j\leq J} K_{j-(d+2)}(x,y,t) + R_{NJ}(x,y,t), \quad R_{NJ}\in C^{N}(U\times\R^{d+2}).
\end{equation}
In particular, on $V\times V$ we have 
\begin{equation}
    k_{B_{s}}(x,y)=|\epsilon_{x}'|K_{s}(x,\varepsilon_{x}(y))+R_{s}(x,y), \quad 
    K_{s}(x,y)=\sum_{j\leq J}K_{j,s}(x,y)+R_{NJ,s}(x,y), 
    \label{eq:Powers1.kernelBs}
\end{equation}
where we have let 
\begin{equation}
    K_{s}(x,y)=\int_{0}^{1} t^{s-1} K (x,y,t)dt, \qquad K_{j,s}(x,y)=\int_{0}^{1} t^{s-1} K_{j-(d+2)} 
    (x,y,t)dt, \quad  j\geq 0,
\end{equation}
and $(R_{s})_{\Re s>0}$ and $(R_{NJ,s})_{\Re s>0}$ are in $\Hol(\Re s>0, C^{\infty}(V \times V))$ and 
$\Hol(\Re s>0, C^{N}(V \times V))$ respectively. 
 
Notice that $K_{j-(d+2)}(x,y,t)$ is in $C^{\infty}(V)\hotimes \cD_{\reg}'(\Rd\times \R)$ and is parabolic homogeneous of 
degree $j-(d+2)\geq -(d+2)$. Thus the family $(K_{j,s})_{\Re s>0}$ belongs to 
$\Hol(\Re s>0,C^{\infty}(U)\hotimes \cD_{\reg}'(\Rd))$. 
Moreover, for any $\lambda>0$, the difference $K_{j,s}(x, \lambda.y)-\lambda^{vs+j-(d+2)}K_{j,s}(x,y)$ is equal to 
\begin{equation}
  \int_{1}^{\lambda^{2}}t^{s-1} K(x,y,t) dt \in 
    \Hol(\Re s>0,C^{\infty}(V\times \R^{d+2})).
\end{equation}
Hence $(K_{j,s})_{\Re s>0}$ is a holomorphic family of almost homogeneous distributions of degree $vs-(d+2)+j$. Combining this with~(\ref{eq:Powers1.kernelBs}) 
then shows that $(K_{s})_{\Re 
s>0}$ belongs to $\Hol(\Re s>0, \cK^{*}(V\times \Rd))$ and has order $vs-(d+2)$. Therefore,
using~(\ref{eq:Powers1.kernelBs}) and Proposition~\ref{prop:HolPHDO.kernel.characterization}
we see that $(B_{s})_{\Re s>0}$ is a holomorphic family of \psivdos\ such that $\ord B_{s}=-(\ord K_{s}+d+2)=-vs$.
\end{proof}

It follows from Lemma~\ref{lem:Powers1.analyticity-Bs} that for any  local Heisenberg 
chart $\kappa:U\rightarrow V$ the family $(\kappa_{*}A_{s|_{U}})_{\Re s>0}$ is a 
holomorphic family of \psivdos\ on $V$ of order~$-vs$. Combining this with~(\ref{eq:Powers1.smoothing}) and~(\ref{eq:Powers1.relation-As}) then 
shows that $(A_{s})_{\Re s>0}$  and $(P^{s})_{\Re s<0}$ are
holomorphic families of \psivdos\ of orders $-vs$ and $vs$ respectively.

Now,  let $s \in \C$ and let $k$ be a positive  integer such that $k> \Re s$. Then we have 
\begin{equation}
    P^{s}u=P^{s-k}P^{k}u \quad \text{for any $u \in C^{\infty}(M,\cE)$}.
\end{equation}
As $P^{k}$ is a differential operator and $P^{s-k}$ is a \psivdo\ of order $m(s-k)$ this proves that 
$P^{s}$ is a \psivdo\ of order $ms$. In fact, as by Proposition~\ref{prop:HolPHDO.composition} the product of \psivdos\ is 
analytic this actually shows that $(P^{s})_{s \in \C}$ is a holomorphic family of \psivdos\ such that $\ord P^{s}=vs$ for every $s\in \C$.

Finally, when $\cE$ is a general vector bundle we can similarly prove that the complex powers $P^{s}$, $s \in \C$,  forms a holomorphic family of \psivdos\ 
such that $\ord P^{s}= vs$ for any $s\in \C$. The proof is thus complete. 
\end{proof}

\begin{example}\label{ex:Powers1.sublaplacian}
Theorem~\ref{thm:Powers1.main} holds for the following sublaplacians:\smallskip 

    (a) A selfdajoint sums of squares $  \Delta=\nabla_{X_{1}}^{*}\nabla_{X_{1}}+\ldots+\nabla_{X_{m}}^{*}\nabla_{X_{m}}$, 
    where  $X_{1},\ldots,X_{m}$ span $H$ and $\nabla$ is a connection on $\cE$, under the condition that the Levi form is nonvanishing.\smallskip 

    (b) The Kohn Laplacian  $\Box_{b}$ on a compact CR manifold acting on $(p,q)$-forms when the condition~$Y(q)$ holds everywhere.\smallskip 
    
   (c) The  horizontal sublaplacian $\Delta_{b}$ on a compact Heisenberg manifold $(M^{d+1}, H)$ acting on sections of $\Lambda^{k}_{\C}H^{*}$ 
   when the condition~$X(k)$ holds everywhere.
\end{example}

In Section~\ref{sec:Rockland-heat} we will actually make use of Theorem~\ref{thm:Powers1.main} to show that when the Levi form has constant rank the 
Rockland condition is enough to insure us that the  the principal symbol of $P+\partial_{t}$ is an invertible Volterra-Heisenberg symbol (see 
Theorem~\ref{thm:Heat1.main}). Therefore, we obtain: 

\begin{theorem}\label{thm:Powers1.main2}
Assume that the Levi form of $(M,H)$ has constant rank and let $P:C^{\infty}(M,\cE)\rightarrow C^{\infty}(M,\cE)$ be a positive differential operator of 
even (Heisenberg) order $v$ such that $P$ satisfies the Rockland condition at every point of $M$. Then:\smallskip 
    
    (i) For any $s \in \C$ the operator $P^{s}$ defined by~(\ref{eq:Powers1.definition}) is a \psivdo\ of order $vs$;\smallskip 
    
    (ii) The family $(P^{s})_{s\in \C}$ forms a holomorphic 1-parameter group of \psivdos. 
\end{theorem}
\begin{example}
   Theorem~\ref{thm:Powers1.main2} is valid for the contact Laplacian $\Delta_{R}$ on a compact orientable contact manifold $(M^{2n+1},\theta)$. 
   In this case $\Delta_{R}^{s}$ has order $2s$ 
   on degree $k=0,\ldots, 2n$ with $k\neq n$ and has order $4s$ on degree $n$.
\end{example}

\begin{remark}
    The above example allows us to fill a technical gap in the proof in~\cite{JK:OKTGSU} of the Baum-Connes conjecture for $SU(n,1)$ (see~\cite{Po:Crelle1}). 
\end{remark}

\section{Weighted Sobolev Spaces}
\label{sec.Sobolev}
Let $(M^{d+1}, H)$ be a compact Heisenberg manifold endowed with a smooth positive density and assume that the Levi form of $(M,H)$ is non-vanishing. We shall  
now construct weighted Sobolev spaces 
$W_{H}^{s}(M)$, $s \in \R$, which extend to any real 
parameter the weighted Sobolev spaces $S_{k}^{2}(M)$, $k\in \N$, of Folland-Stein~\cite{FS:EDdbarbCAHG} (see also~\cite{RS:HDONG}). 
As a consequence these Sobolev spaces will provide us with sharp regularity estimates for hypoelliptic \psivdos. 

Let $X_{1},\ldots,X_{m}$ be real vector fields spanning $H$ and consider the 
positive sum of squares, 
\begin{equation}
    \Delta_{X}= X_{1}^{*}X_{1}+\ldots+X_{m}^{*}X_{m}.
\end{equation}
As mentioned in Example~\ref{ex:Powers1.sublaplacian} (a), 
since the Levi form of $(M,H)$ is non-vanishing Theorem~\ref{thm:Powers1.main} is valid for $1+\Delta_{X}$. Thus, 
the complex powers $(1+\Delta_{X})^{s}$, $s \in \C$,  
gives rise to a holomorphic 1-parameter group  of \psivdos\ such that $\ord (1+\Delta_{X})^{s}=2s$ for any $s\in \C$. 

\begin{definition}
  $W_{H}^{s}(M)$, $s\in \R$, consists of all distributions $u\in \cD'(M)$ such that 
  $(1+\Delta_{X})^{\frac{s}{2}}u$ is in $L^{2}(M)$. It is endowed with the Hilbert norm given by
  \begin{equation}
      \|u\|_{W_{H}^{s}}=\|(1+\Delta_{X})^{\frac{s}{2}}u\|_{L^{2}}, \qquad u \in W_{H}^{s}(M). 
  \end{equation}
 \end{definition}

\begin{proposition}\label{prop:Sobolev.embeddings}
1) Neither $W_{H}^{s}(M)$, nor its topology, depend on the choice of the vector fields $X_{1},\ldots,X_{m}$.\smallskip 
    
    2) We have the following continuous embeddings: 
    \begin{equation}
   \begin{array}{rccccl}
       L^{2}_{s}(M)  & \hookrightarrow & W_{H}^{s}(M)&\hookrightarrow & L^{2}_{s/2}(M) & \qquad \text{if $s\geq 0$},\\
       L^{2}_{s/2}(M)  & \hookrightarrow & W_{H}^{s}(M) & \hookrightarrow & L^{2}_{s}(M) & \qquad \text{if $s< 0$}.
   \end{array}
         \label{eq:Sobolev.embeddings}
    \end{equation}
\end{proposition}
\begin{proof}
    1) Let $Y_{1},\ldots,Y_{p}$ be other vector fields spanning $H$. The operator 
    $(1+\Delta_{Y})^{s}(1+\Delta_{X})^{-s}$ is a \psivdo\ 
    of order $0$, so is bounded on $L^{2}(M)$ by Proposition~\ref{prop:PsiHDO.Sobolev-regularity}. Therefore,  we get the
    estimates,
    \begin{multline}
        \|(1+\Delta_{Y})^{s}u\|_{L^{2}}=\|(1+\Delta_{Y})^{s}(1+\Delta_{X})^{-s}(1+\Delta_{X})^{s}u\|_{L^{2}} \leq C_{XYs} 
        \|(1+\Delta_{X})^{s}u\|_{L^{2}},
    \end{multline}
     which hold for any $u \in C^{\infty}(M)$. Interchanging the roles of the $X_{j}$'s and of the $Y_{k}$'s also gives the estimates
     \begin{equation}
         \|(1+\Delta_{X})^{s}u\|_{L^{2}} \leq C_{YXs}\|(1+\Delta_{Y})^{s}u\|_{L^{2}}, \qquad u \in C^{\infty}(M). 
     \end{equation}
     Therefore, wether we use the $X_{j}$'s or the $Y_{k}$'s to define the $W_{H}^{s}(M)$ changes neither the space, nor its topology. \smallskip 
     
     2) Let $ s\in [0,\infty)$. Since $(1+\Delta_{X})^{{\frac{s}{2}}}$ is a \psivdo\ of order $s$, 
     Proposition~\ref{prop:PsiHDO.Sobolev-regularity}   
     tells us that 
     it is bounded from $L^{2}_{s}(M)$ to $L^{2}(M)$. Thus,
     \begin{equation}
         \|u\|_{W_{H}^{s}}=\|(1+\Delta_{X})^{{\frac{s}{2}}}u\|_{L^{2}}\leq C_{s} \|u\|_{L^{2}_{s}}, \qquad u \in W_{H}^{s}(M), 
     \end{equation}
     which shows that  $L^{2}_{s}(M)$  embeds continuously into $W_{H}^{s}(M)$.
     
    On the other hand, as $(1+\Delta_{X})^{-{\frac{s}{2}}}$ has order $-s$ Proposition~\ref{prop:PsiHDO.Sobolev-regularity}   
     also tells us that $(1+\Delta_{X})^{-{\frac{s}{2}}}$ is bounded from $L^{2}(M)$ to 
     $L^{2}_{s/2}(M)$. Therefore, for we get the estimates
     \begin{equation}
         \|u\|_{L^{2}_{s/2}}=\|(1+\Delta_{X})^{-{\frac{s}{2}}}(1+\Delta_{X})^{{\frac{s}{2}}}u\|_{L^{2}_{s/2}} \leq 
         C_{s}\|(1+\Delta_{X})^{-{\frac{s}{2}}}\|_{L^{2}}=C_{s}\|u\|_{W_{H}^{s}},
     \end{equation}
     which hold for any $u \in W_{H}^{s}(M)$ and show that $W_{H}^{s}(M)$ embeds continuously into $L^{2}_{s/2}(M)$.  
     
     Finally, when  $s<0$ we can similarly show that we have continuous embeddings 
     $L^{2}_{s/2}(M)  \hookrightarrow W_{H}^{s}(M)$ and $ W_{H}^{s}(M)\hookrightarrow L^{2}_{s}(M)$. 
\end{proof}

As an immediate consequence of Proposition~\ref{prop:Sobolev.embeddings} we obtain: 

\begin{proposition}\label{prop:Sobolev.Cinfty-cD'}
   The following equalities between topological spaces hold:
   \begin{equation}
       C^{\infty}(M)=\cap_{s \in \R} W_{H}^{s}(M) \qquad \text{and} \qquad \cD'(M)= \cup_{s \in \R}  W_{H}^{s}(M).
   \end{equation}
\end{proposition}

Let us now compare the Weighted Sobolev spaces $W_{H}^{s}(M)$ to the weighted Sobolev spaces $S_{k}^{2}(M)$, $k=1,2,\ldots$, 
of Folland-Stein~\cite{FS:EDdbarbCAHG}. 

In we sequel we let $\N_{m}=\{1,\ldots,m\}$ and for any $I=(i_{1}, \ldots, i_{k})$ in $\N_{m}^{k}$ we set 
\begin{equation}
    X_{I}=X_{i_{1}}\ldots X_{i_{l}}.
\end{equation}

\begin{definition}[\cite{FS:EDdbarbCAHG}]
The Hilbert space $S_{k}^{2}(M)$, $k\in\N$, consists of functions $u\in L^{2}(M)$ such that 
  $(X_{I})u \in L^{2}(M)$ for any $I \in \cup_{j=1}^{k}\N_{m}^{j}$. It is endowed with the Hilbertian norm given by
  \begin{equation}
      \|u\|_{S_{k}^{2}}^{2}=\|u\|_{L^{2}}^{2} + \sum_{1\leq j\leq k} \sum_{I \in \N_{m}^{j}} 
      \|X_{I}u\|_{L^{2}}^{2}, \qquad u \in S_{k}^{2}(M). 
  \end{equation}
\end{definition}

\begin{proposition}\label{prop:Sobolev.equality-FS}
  For $k=1,2,\ldots$ the weighted Sobolev spaces $W_{H}^{k}(M)$ and $S_{k}^{2}(M)$ agree as spaces and bear the same topology.

\end{proposition}
\begin{proof}
    First, if $I \in \cup_{j=1}^{k}\N_{m}^{j}$ then by Proposition~\ref{prop:Sobolev.regularity-PsiHDOs}
    the differential operator $X_{I}$ is bounded from $W_{H}^{k}(M)$ to $L^{2}(M)$, so we get the get the estimate, 
    \begin{equation}
        \|u\|_{S^{k}}^{2}=\|u\|_{L^{2}}^{2} + \sum_{1\leq j\leq k} \sum_{I \in \N_{m}^{j}} 
      \|X_{I}u\|_{L^{2}}^{2} \leq C_{k}^{2}\|u\|_{W_{H}^{k}} \qquad u \in 
        C^{\infty}(M).
         \label{eq:Sobolev.SkWHk}
    \end{equation}
    
    On the other hand,  for any $l\in \N$  the differential operators $X_{j}$ and 
    $X_{j}^{*}$, $j=1,\ldots,m$, are bounded linear maps from $S_{l+1}^{2}(M)$ to $S_{l}^{2}(M)$. Thus,  for any integer $p$,  
    the operator 
    $(1+\Delta_{X})^{p}$ is bounded from $S^{2}_{l+2p}(M)$ to $S^{2}_{l}(M)$. It follows that, when  $k$ is even, we have 
    \begin{equation}
        \|u\|_{W_{H}^{k}}=\|(1+\Delta_{X})^{\frac{k}{2}}u\|_{L^{2}} \leq C_{k}\|u\|_{S^{k}}, \qquad u \in C^{\infty}(M).
    \end{equation}
    Moreover, for any $u \in C^{\infty}(M)$ we have 
    \begin{equation}
       \|(1+\Delta_{X})^{\frac{1}{2}}u\|_{L^{2}}^{2}=\acou{(1+\sum_{1\leq j\leq m}X^{*}_{j}X_{j})u}{u}_{L^{2}}= 
        \|u\|_{L^{2}}^{2}+\sum_{1\leq j\leq m}\|X_{j}^{2}\|_{L^{2}}^{2} = \|u\|_{S^{2}_{1}}^{2}. 
         \label{eq:Sobolev.WHkSk}
    \end{equation}
    Thus $(1+\Delta_{X})^{\frac{1}{2}}$ is bounded from $S^{2}_{1}(M)$ to $L^{2}(M)$. Therefore,  if $k$ is odd, say 
    $k=2p+1$, then the operator 
    $(1+\Delta_{X})^{\frac{k}{2}}=(1+\Delta_{X})^{p}(1+\Delta_{X})^{\frac{1}{2}}$ 
    is bounded from $S_{k}^{2}(M)$ to $L^{2}(M)$. Hence~(\ref{eq:Sobolev.WHkSk}) is valid in the odd case as well. Together 
    with~(\ref{eq:Sobolev.SkWHk}) this implies that  $W_{H}^{k}(M)$ and $S_{k}^{2}(M)$ agree as spaces and bear the same topology. 
\end{proof}

Now, let $\cE$ be a Hermitian vector bundle over $M$. We can also define weighted Sobolev spaces of sections of $\cE$ as follows. 
Let $\nabla: C^{\infty}(M,\cE)\rightarrow C^{\infty}(M, T^{*}M\times \cE)$  be a connection on $\cE$ and define 
\begin{equation}
    \Delta_{\nabla,X}=\nabla_{X_{1}}^{*}\nabla_{X_{1}}+\ldots+\nabla_{X_{m}}^{*}\nabla_{X_{m}}.
\end{equation}
As in the scalar case,  
the complex powers $(1+\Delta_{\nabla,X})^{s}$, $s\in \C$, form
an analytic 1-parameter group of \psivdos\ such that $\ord (1+\Delta_{\nabla,X})^{s}=2s$ for any $s \in \C$.  

\begin{definition}
  $W_{H}^{s}(M,\cE)$, $s\in \R$, consists of all distributional sections $u\in \cD'(M,\cE)$ such that 
  $(1+\Delta_{\nabla,X})^{\frac{s}{2}}u\in L^{2}(M,\cE)$. It is endowed with the Hilbertian norm given by
  \begin{equation}
      \|u\|_{W_{H}^{s}}=\|(1+\Delta_{\nabla,X})^{\frac{s}{2}}u\|_{L^{2}}, \qquad u \in W_{H}^{s}(M,\cE). 
  \end{equation}
 \end{definition}

Along similar lines as that of the proof of Proposition~\ref{prop:Sobolev.embeddings} we can prove: 

\begin{proposition}
   1) As a topological space  $W_{H}^{s}(M,\cE)$, $s\in \R$, is independent of the choice of the connection $\nabla$ and 
    of the vector fields $X_{1},\ldots,X_{m}$.\smallskip 
    
    2)  We have the following continuous embeddings: 
    \begin{equation}
   \begin{array}{rccccl}
       L^{2}_{s}(M,\cE)  & \hookrightarrow & W_{H}^{s}(M,\cE) &\hookrightarrow& L^{2}_{s/2}(M,\cE) & \qquad \text{if $s\geq 0$},\\
       L^{2}_{s/2}(M,\cE)  & \hookrightarrow & W_{H}^{s}(M,\cE)& \hookrightarrow &L^{2}_{s}(M,\cE) & \qquad \text{if $s< 0$}.
   \end{array}
         \label{eq:Sobolev.embeddings2}
    \end{equation}
\end{proposition}

As a consequence we obtain the equalities of topological spaces,
\begin{equation}
       C^{\infty}(M,\cE)=\cap_{s \in \R} W_{H}^{s}(M,\cE) \qquad \text{and} \qquad \cD'(M,\cE)= \cup_{s \in \R}  W_{H}^{s}(M,\cE).
        \label{eq:Sobolev.smooth-distributions-WHs}
\end{equation}

Notice that we can also define Folland-Stein spaces $S_{k}^{2}(M,\cE)$, $k=1,2,\ldots$, as in the scalar case, by using the 
differential operators $\nabla_{X_{I}}=\nabla_{X_{i_{1}}}\ldots \nabla_{X_{i_{k}}}$, $I \in 
\cup_{j=1}^{k}\N_{m}^{j}$. Then, by arguing as in the proof  of Proposition~\ref{prop:Sobolev.equality-FS}, we can show that the spaces $W_{H}^{k}(M,\cE)$ 
and $S_{k}^{2}(M,\cE)$ agree and bear the same topology.  

Now, the Sobolev spaces $W_{H}^{s}(M,\cE)$ yield sharp regularity results for \psivdos. 

\begin{proposition}\label{prop:Sobolev.regularity-PsiHDOs}
    Let $P\in \pvdo^{m}(M,\cE)$, set $k=\Re m$ and let $s\in \R$.\smallskip 
    
    1) The operator $P$ extends to a continuous linear mapping from $W_{H}^{s+k}(M,\cE)$ to $W_{H}^{s}(M,\cE)$. 
    \smallskip 
    
    2) Assume that the principal symbol of $P$ is invertible. Then for any $u \in \cD'(M,\cE)$ we have
    \begin{equation}
        Pu \in W_{H}^{s}(M,\cE) \Longrightarrow u \in W_{H}^{s+k}(M,\cE).
         \label{eq:Sobolev.hypoellipricity-WHs}
    \end{equation}
   In fact, for any $s'\in \R$ we have the estimate,
    \begin{equation}
        \|u\|_{W^{s+k}_{H}} \leq C_{ss'}(\|Pu\|_{W^{s}_{H}}+\|u\|_{W^{s'}_{H}}), \qquad u \in W^{s+k}_{H}(M,\cE).
    \end{equation}
\end{proposition}
\begin{proof}
    1) As $P_{s}=(1+\Delta_{\nabla,X})^{\frac{s}{2}}P(1+\Delta_{\nabla,X})^{-\frac{(s+k)}{2}}$ is a \psivdo\ with purely imaginary order, 
    Proposition~\ref{prop:PsiHDO.Sobolev-regularity} tells us that it gives rise to a bounded operator on $L^{2}(M,\cE)$. Therefore, we have 
    \begin{equation}
        \|Pu\|_{W^{s}_{H}}=\|P_{s}(1+\Delta_{\nabla,X})^{s+k}u\|_{L^{2}} \leq C_{s}\|u\|_{W^{s+k}_{H}}, \qquad u \in C^{\infty}(M,\cE), 
    \end{equation}
    It then follows that $P$ extends to a continuous linear mapping from $W_{H}^{s+k}(M,\cE)$ to $W_{H}^{s}(M,\cE)$. \smallskip 
    
    2) Since the principal symbol of $P$ is invertible 
    by Proposition~\ref{thm:PsiHDO.hypoellipticity} there exist $Q$ in $\pvdo^{-m}(M,\cE)$ and $R$ in $\psinf(M,\cE)$ such that $QP=1-R$. 
    Therefore, for any $u \in \cD'(M,\cE)$ we have $ u=QPu+Ru$.  
     Thanks to the first part we know that  $Q$ maps 
    $W_{H}^{s}(M,\cE)$ to $W_{H}^{s+k}(m,\cE)$. Since $R$ is smoothing, and so $Ru$ always is smooth, it follows
    that if $Pu$ is in $W_{H}(M,\cE)$ then $u$ 
    must be in $W_{H}(M,\cE)$.
    
    In fact, as $Q$ is actually bounded from 
    $W_{H}^{s}$ to $W_{H}^{s+k}$ and~(\ref{eq:Sobolev.smooth-distributions-WHs}) 
    implies that $R$ is bounded from any space $W_{H}^{s'}(M,\cE)$ to $W_{H}^{s+k}(M,\cE)$ we have 
    estimates, 
    \begin{equation}
        \|u\|_{W^{s+k}_{H}}\leq \|QPu\|_{W^{s+k}_{H}} + \|Ru\|_{W^{s+k}_{H}} \leq 
        C_{ss'}(\|Pu\|_{W^{s}_{H}}+\|u\|_{W^{s'}_{H}}) ,\qquad u \in W^{s+k}_{H}(M,\cE). 
    \end{equation}
    The proof is thus complete.
\end{proof}
\begin{remark}
    Combining Proposition~\ref{prop:Sobolev.regularity-PsiHDOs} and the embeddings~(\ref{eq:Sobolev.embeddings}) 
    and~(\ref{eq:Sobolev.embeddings2}) we recover the Sobolev regularity of \psivdos\ as in Proposition~\ref{prop:PsiHDO.Sobolev-regularity} 
    as well as the hypoelliptic estimates~(\ref{eq:PsiHDO.subellipticity.subelliptic-estimates}).
\end{remark}

On the other hand, the spaces $W_{H}^{s}(M,\cE)$ can be localized as follows. 

\begin{definition}\label{def:Sobolev.loaclization-WHs}
    We say that $u \in \cD'(M,\cE)$ is $W_{H}^{s}$ near a point $x_{0}\in M$ whenever there exists $ \varphi\in C^{\infty}(M)$ such that 
    $\varphi(x_{0})\neq 0$  and $\varphi u$ is in $W_{H}^{s}(M,\cE)$. 
\end{definition}
 
This definition depends only on the germ of $u$ at $x_{0}$ because we have:
\begin{lemma}\label{lem:Sobolev.localization-WHs}
    Let $u \in \cD'(M,\cE)$ be  $W_{H}^{s}$ near $x_{0}$. Then for any $\varphi\in C^{\infty}(M)$ with a small enough support 
    about $x_{0}$ the distribution $\varphi u$ lies in $W_{H}^{s}(M,\cE)$.
\end{lemma}
\begin{proof}
    Let $\varphi\in C^{\infty}(M)$ be such that $\varphi(x_{0})\neq 0$ and $\varphi u$ is in $W_{H}^{s}(M,\cE)$ and let $\psi\in C^{\infty}(M)$ be 
    so that $\psi(x_{0})\neq 0$ and $\supp \psi \cap \varphi^{-1}(0)\neq \emptyset $. 
    Then $\chi:=\frac{\psi}{\varphi}$ is in $C^{\infty}(M)$ and we have
    \begin{equation}
       (1+\Delta_{\nabla,X})^{\frac{s}{2}}\psi u=  (1+\Delta_{\nabla,X})^{\frac{s}{2}}\chi \varphi u = 
       (1+\Delta_{\nabla,X})^{\frac{s}{2}}\chi  (1+\Delta_{\nabla,X})^{-\frac{s}{2}}.  (1+\Delta_{\nabla,X})^{\frac{s}{2}}\varphi u.
    \end{equation}
    Since $(1+\Delta_{\nabla,X})^{\frac{s}{2}}\varphi u$ is in $L^{2}(M,\cE)$ and $ (1+\Delta_{\nabla,X})^{\frac{s}{2}}\chi  
    (1+\Delta_{\nabla,X})^{-\frac{s}{2}}$ is a zero'th order \psivdo, so maps $L^{2}(M,\cE)$ to itself, it follows that $\psi u$ lies in 
    $W_{H}^{s}(M,\cE)$. Hence the lemma.
\end{proof}

We can now get a localized version of~(\ref{eq:Sobolev.hypoellipricity-WHs}).
\begin{proposition}\label{prop:Sobolev.hypoellipricity-WHs-localized}
  Let $P \in \pvdo^{m}(M,\cE)$ have an invertible principal symbol, set $k=\Re m$ and let $s\in \R$. Then for any $u \in \cD'(M,\cE)$ we have 
  \begin{equation}
      \text{$Pu$ is $W_{H}^{s}$ near $x_{0}$}\ \Longrightarrow \ \text{$u$ is $W_{H}^{s}$ near $x_{0}$}.
  \end{equation}
\end{proposition}
\begin{proof}
    Assume that $Pu$ is $W_{H}^{s}$ near $x_{0}$ and let $\varphi\in C^{\infty}(M)$ be such that $\varphi(x_{0})\neq 0$ and $\varphi u$ is in 
    $W_{H}^{s}(M,\cE)$. Thanks to Lemma~\ref{lem:Sobolev.localization-WHs} we may assume that 
 $\varphi=1$ near $x_{0}$. Let $\psi \in C^{\infty}(M)$ be such that $\psi(x_{0})\neq 0$ and $\varphi=1$ near $\supp \psi$. 
    Since the principal symbol of $P$ is invertible there exist $Q$ in $\pvdo^{-m}(M,\cE)$ and $R$ in $\psinf(M,\cE)$ such that $QP=1-R$. 
    Thus for any $u \in \cD'(M,\cE)$ we have 
    \begin{equation}
        \psi u=\psi QPu+\psi Ru= \psi Q\varphi Pu+ \psi Q(1-\varphi) Pu +\psi Ru.
    \end{equation}
   
    In the above equality $Ru$ is a smooth function since $R$ is a smoothing operator. Similarly, as  $\psi$ and $1-\varphi$ have disjoint supports, the 
    operator $\psi Q(1-\varphi) P$ is a smoothing operator and so $\psi Q(1-\varphi) Pu$ is a smooth function.  
    In addition, since $\varphi Pu$ is in $W_{H}^{s}(M,\cE)$ and $ \psi Q$ is a \psivdo\ of order $-m$, and so
    maps $W_{H}^{s}(M,\cE)$ to $W_{H}^{s+k}(M,\cE)$ , it follows 
    that $\psi u$ is in $W_{H}^{s+k}$. Hence $u$ is $W_{H}^{s+k}$ near $x_{0}$. 
\end{proof}

Next, the first part of Proposition~\ref{prop:Sobolev.regularity-PsiHDOs} admits generalization to holomorphic families of \psivdos. 

\begin{proposition}\label{prop:Sobolev.regularity-PsiHDOs-families}
    Let $\Omega \subset \C$ be open and let $(P_{z})_{z\in \Omega}\in \Hol(\Omega, \Psi^{*}_{H}(M,\cE))$. Assume there exists a real $m$ 
    such that $\Re\ord P_{z}\leq m<\infty$. Then for any $s\in \R$ the family $(P_{z})_{z\in \Omega}$ defines 
    a holomorphic family with values in $\cL(W_{H}^{s+m}(M,\cE),W_{H}^{s}(M,\cE))$. 
\end{proposition}
\begin{proof}
    First, let $V\subset \Rd$ be a Heisenberg chart with $H$-frame $Y_{0},Y_{1},\ldots,Y_{d}$ and let $(Q_{z})_{z\in 
    \Omega}\in \Hol(\Omega, \pvdo^{*}(V))$ be such that $\Re \ord Q_{z}\leq 0$. 
    Then  we can write
    \begin{equation}
        Q_{z}=p_{z}(x,-iY)+R_{z},
    \end{equation}
    for some families $(p_{z})_{z\in \Omega}\in \Hol(\Omega, S^{*}(V\times \Rd))$ and $(R_{z})_{z\in \Omega}\in \Hol(\Omega, \Psi^{\infty}(V))$.
    Since $\Re \ord p_{z}=\Re \ord Q_{z}\leq 0$ we see that $(p_{z})_{z\in \Omega}$ belongs to $\Hol(\Omega,S^{0}_{\|}(\URd))$. 
    
    Next, for $j=1,\ldots,d$ let $\sigma_{j}$ denote the standard symbol of 
    $-iY_{j}$ and set $\sigma=(\sigma_{0},\ldots,\sigma_{d})$. 
    Then it follows from the proof of Proposition~\ref{prop:Complexpowers.operators.properties} that the family $(p_{z}(x,\sigma(x,\xi)))_{z\in 
    \Omega}$ lies in $\Hol(\Omega, S_{\frac{1}{2},\frac{1}{2}}(V\times \Rd))$. 
  Since by~\cite{Hw:L2BPO} the 
    quantization map $q\rightarrow q(x,D_{x})$ is continuous from $S_{\frac{1}{2}\frac{1}{2}}(\URd)$ to 
    $\cL(L^{2}_{\op{loc}}(U), L^{2}(U))$, we deduce that $(p_{z}(x,-iX))_{z\in \Omega}$ and $(Q_{z})_{z\in \Omega}$ are 
    holomorphic families with values in $\cL(L^{2}_{\op{loc}}(U), L^{2}(U))$. 
    
    It follows from the above result that any family $(Q_{z})_{z\in \Omega}\in \Hol(\Omega, \pvdo*(M,\cE))$ 
    such that $\Re \ord Q_{z}\leq 0$ gives rise to a holomorphic family with values in $\cL(L^{2}(M,\cE))$. 
    This applies in particular to the family 
    \begin{equation}
        Q_{z}^{(s)}=(1+\Delta_{\nabla,X})^{\frac{s}{2}}Q_{z}(1+\Delta_{\nabla,X})^{-\frac{m+s}{2}}, \quad z \in \Omega. 
    \end{equation}
    As $Q_{z}=(1+\Delta_{\nabla,X})^{-\frac{s}{2}}Q_{z}^{(s)}(1+\Delta_{\nabla,X})^{\frac{m+s}{2}}$ it follows that $(Q_{z})_{z\in 
    \Omega}$ gives rise to a holomorphic family with values in $\cL(W_{H}^{s+m}(M,\cE),W_{H}^{s}(M,\cE))$.
\end{proof}

Combining Proposition~\ref{prop:Sobolev.regularity-PsiHDOs-families} with Theorem~\ref{thm:Powers1.main2} 
we immediately get:

\begin{proposition}
   Let $P:C^{\infty}(M,\cE)\rightarrow C^{\infty}(M,\cE)$ be a positive differential operator of 
even (Heisenberg) order $w$ with an invertible principal symbol. 
Then, for any reals $m$ and $s$, the family of the complex powers of $P$ satisfies
\begin{equation}
    (P^{z})_{\Re z<m} \in \Hol(\Re z<m, \cL(W_{H}^{s+mw}(M,\cE),W_{H}^{s}(M,\cE))).
\end{equation}
\end{proposition}

\section{Rockland condition and the heat equation}
 \label{sec:Rockland-heat}
Let $(M^{d+1},H)$ be a compact Heisenberg manifold equipped with a smooth positive density and let $\cE$ be a Hermitian vector bundle over $M$. 
 In this section we shall show that for a selfadjoint differential operator $P:C^{\infty}(M,\cE)\rightarrow C^{\infty}(M,\cE)$ which is bounded from 
 below, the invertibility 
 of the principal symbol of $P$ is enough to insure us the invertibility of the principal symbol of $P+\partial_{t}$. This is needed for the 
 pseudodifferential representation of the heat kernel and the invertibility of the heat operator $P+\partial_{t}$ in~\cite{BGS:HECRM}, as well as in 
 Theorem~\ref{thm:Powers1.main} for the complex powers of $P$ to give rise to a holomorphic family of \psivdos. 

 Let us first look at symbols that are positive in the sense below. 
 
\begin{definition}\label{def:Heat1.positive-symbol}
   A symbol $p\in S_{m}(\fg^{*}M,\cE)$ is said to be positive when it can be put into the form $p=\bar{q}*q$ for some symbol 
   $q\in S_{\frac{m}{2}}(\fg^{*}M,\cE)$. 
\end{definition}

\begin{lemma}\label{lem:Rockland-heat.positivity}
   Let $P\in \pvdo^{m}(M,\cE)$ have principal symbol $\sigma_{m}(P)\in S^{m}(\fg^{*}M,\cE)$. Then we have equivalence:\smallskip
   
   (i) The symbol $\sigma_{m}(P)$ is positive;\smallskip

   (ii) There exist $Q\in \pvdo^{\frac{m}{2}}(M,\cE)$ and $R\in \pvdo^{m-1}(M,\cE)$ so that $P=Q^{*}Q+R$. 
\end{lemma}
\begin{proof}
  Assume that $\sigma_{m}(P)$ is a positive symbol, so that there exists $q_{\frac{m}{2}}\in S_{\frac{m}{2}}(\fg^{*}M,\cE)$ such that  
  $p=\overline{q_{\frac{m}{2}}}*q_{\frac{m}{2}}$. By Proposition~\ref{prop:PsiHDO.surjectivity-principal-symbol-map}   the principal 
 symbol map $\sigma_{\frac{m}{2}}:\pvdo^{\frac{m}{2}}(M,\cE)\rightarrow S^{\frac{m}{2}}(\fg^{*}M,\cE)$ is surjective, so there  
 exists $Q\in \pvdo^{\frac{m}{2}}(M,\cE)$ such that $\sigma_{\frac{m}{2}}(Q)=q_{\frac{m}{2}}$. Then by Proposition~\ref{prop:PsiHDO.composition2} 
 and~Proposition~\ref{prop:PsiHDO.adjoint-manifold}   the operator 
$Q^{*}Q$ has principal symbol 
 $\overline{q_{\frac{m}{2}}}*q_{\frac{m}{2}}=\sigma_{m}(P)$, hence coincides with $P$ modulo $\pvdo^{m-1}(M,\cE)$. 
  
 Conversely, if $P$ is of the form $P=Q^{*}Q+R$ with $Q\in \pvdo^{\frac{m}{2}}(M,\cE)$ and $R\in \pvdo^{m-1}(M,\cE)$ then $P$ and $Q^{*}Q$ have same 
 principal symbol. Then $\sigma_{m}(P)=\overline{\sigma_{\frac{m}{2}}(Q)}*\sigma_{\frac{m}{2}}(Q)$, that is,  
$\sigma_{m}(P)$ is a positive symbol.
\end{proof}

From now on we assume that the Levi form of $(M,H)$ has constant rank and we let $P:C^{\infty}(M,\cE)\rightarrow C^{\infty}(M,\cE)$ be a 
differential operator of even Heisenberg order $m$. 

 \begin{lemma}\label{lem:Rockland.invertibility}
      Assume that  $P$ satisfies the Rockland condition at every point and has a positive principal symbol. 
      Then the principal symbol of the heat operator $P+\partial_{t}$ is an invertible Volterra-Heisenberg symbol.
 \end{lemma}
 \begin{proof}
      Let us first assume that $\cE$ is the trivial line bundle. Since proving the invertibility of the principal symbol of $P+\partial_{t}$ is a local problem with respect to 
      the space covariable, we may as well work in a local Heisenberg chart.  
      
    Moreover, as the Levi form has constant rank it follows from~\cite{Po:Pacific1} that $GM$ is a fiber bundle of Lie group with fiber 
   $G=\bH^{2n+1}\times \R^{d-2n}$, where $2n$ is the rank of the Levi form. Therefore, by considering a trivialization of this fiber bundle by means 
   of a suitable local $H$-frame (see~\cite{Po:Pacific1}) we may further assume that $GU$ is the trivial bundle $U\times G$. In particular,  
   the families of model operators $(P^{x})_{x \in U}$ can be seen as a smooth family of self-adjoint convolution operators on 
   $\cS_{0}(G)$.  
   
   Next, as the global principal symbol $\sigma_{m}(P)$ of $P$ is positive there exists $\tilde{p}_{\frac{m}{2}}\in S_{\frac{m}{2}}(U\times G)$ such that 
   $\sigma_{m}(P)=\overline{\tilde{p}_{\frac{m}{2}}}*\tilde{p}_{\frac{m}{2}}$, so if we let $\tilde{P}^{x}$ be the left-convolution operator with 
   symbol $\tilde{p}_{\frac{m}{2}}(x,.)$ then by Proposition~\ref{prop:PsiHDO.composition2} and~Proposition~\ref{prop:PsiHDO.adjoint-manifold} 
   we have $P^{x}=(\tilde{P}^{x})^{*}\tilde{P}^{x}$. Thus, by 
   Proposition~\ref{PsiHDO.properties-symbol-representation} for 
   every non-trivial irreducible representation $\pi$ of $G$ we have $\overline{\pi_{P^{x}}}=(\overline{\pi_{\tilde{P}^{x}}})^{*}\overline{\pi_{\tilde{P}^{x}}}$, 
   which shows that $\overline{\pi_{P}^{x}}$ is positive. 
   
   Now, as $P^{x}$ is positive and satisfies the Rockland condition it then follows from~\cite{FS:HSHG} that:\smallskip
   
   - The operator $P^{x}+ \partial_{t}$ satisfies the Rockland condition on $G\times \R$ (see~\cite[Lem.~4.21]{FS:HSHG}); 
   
   - There exists $K^{x}(y,t)\in C^{\infty}((G\times \R)\setminus 0)$ homogeneous of degree $-(d+2)$ such that $K^{x}(y,t)=0$ for $t<0$ and 
   $K^{x}(y,t)$ is a fundamental solution of $P^{x}+ \partial_{t}$, that is, 
   \begin{equation}
       (P^{x}+\partial_{t})K^{x}(y,t)=\delta(y,t),
       \label{eq:Heat1.fundamental-solution-model-operator}
   \end{equation}
    where $\delta(y,t)$ denotes the Dirac distribution at the origin on $G\times \R$ (see~\cite[Lem.~4.24]{FS:HSHG}
    in the case $m>\frac{d+2}{2}$ and \cite[pp.~136--137]{FS:HSHG} for the general case).

       Let $Q^{x}$ be the left-convolution operator on $G\times \R$ with $K^{x}(y,t)$. As $K^{x}$ belongs to $ \cK_{\op{v},-(d+2)}(G\times \R_{(m)})$ this is 
       a left-invariant Volterra \psivdo\ with symbol $q^{x}_{-m}(\xi,\tau)=(K^{x})^{\wedge}(\xi,\tau)$. Moreover, thanks 
       to~(\ref{eq:Heat1.fundamental-solution-model-operator}) for any $f \in \cS(G\times \R)$ we have 
       \begin{equation}
           (P^{x}+\partial_{t})Q^{x}f= (P^{x}+\partial_{t})(K*f)=[(P^{x}+\partial_{t})K]*f=\delta*f=f.
       \end{equation}
        Hence  $(P^{x}+\partial_{t})Q^{x}=1$, which gives $(\sigma_{m}(P)(x,.)+i\tau)*q^{x}_{-m}=1$. 
        
        Next, the distribution $K^{x}(y,-t)$ belongs to $ \cK_{\op{v},-(d+2)}(G\times \R_{(m)})$ and is a fundamental solution of 
        $P^{x}-\partial_{t}=(P^{x}+\partial_{t})^{t}$, so the same arguments as above show that the left-convolution operator $\tilde{Q}^{x}$ with $K^{x}(y,-t)$ 
        is a left-invariant Volterra \psivdo\ such that $(P^{x}+\partial_{t})^{t}\tilde{Q}^{x}=1$. Taking transposes we obtain 
        $(\tilde{Q}^{x})^{t}(P^{x}+\partial_{t})$, so 
         if we let $\tilde{q}_{-m}^{x}\in S_{\op{v},-m}(G\times \R_{(m)})$ be the symbol of $(\tilde{Q}^{x})^{t}$ then we get 
        $\tilde{q}_{-m}^{x}*(\sigma_{m}(P)(x,.)+i\tau)=1$. Henceforth $\sigma_{m}(P)(x,.)+i\tau$ is a two-sided invertible 
        Volterra-Heisenberg symbol on $G\times \R$ with inverse $q_{-m}^{x}=\tilde{q}_{-m}^{x}$. 

       On the other hand, since for every $x \in U$ the operator $P^{x}+\partial_{t}$ satisfies the Rockland condition on $G\times \R$,  it follows 
        from~\cite[Thm.~5(d)]{CGGP:POGD} that there exists $\tilde{K}(x,y,t)\in C^{\infty}(U\times[(G\times \R)\setminus 0])$ homogeneous of degree 
        $-(d+2)$ with respect to $(y,t)$ such that, for every $x\in U$, the left-convolution operator by $\tilde{K}(x,.,.)$ is a right-inverse for 
        $P^{x}+\partial_{t}$. Necessarily, this operator agrees with the two-sided 
        inverse $Q^{x}$, so $K^{x}$ agrees with $\tilde{K}(x,.,.)$. 
        
        This shows that $K^{x}$ depends smoothly only on $x$, that is, 
        $K(x,y,t)=K^{x}(x,t)$ is an element of $\cK_{\op{v},-(d+2)}(U\times G\times 
        \R)$. It then follows that $q_{-m}(x,\xi,\tau)=q_{-m}^{x}(xi,\tau)=K^{\wedge}_{(y,t)\rightarrow (\xi,\tau)}(x,\xi,\tau)$ belongs to 
        class $S_{\op{v},-m}(U\times G\times \R_{(m)})$ and satisfies
        \begin{equation}
             q_{-m}*(\sigma_{m}(P)+i\tau)=(\sigma_{m}(P)+i\tau)*q_{-m}=1.
        \end{equation}
        Hence $\sigma_{m}(P)(x,\xi)+i\tau$ is an invertible Volterra-Heisenberg symbol with inverse $q_{-m}(x,\xi,\tau)$. 
        
        Finally, suppose that $\cE$ is not the trivial line bundle. Then, as the aforementioned results of~\cite{FS:HSHG} 
        and~\cite{CGGP:POGD} remain valid for systems, by working in a local trivializing Heisenberg chart and by arguing as above we can show that the principal 
        symbol $P+\partial_{t}$ is an invertible Volterra-Heisenberg symbol.  
    \end{proof}

 Granted Proposition~\ref{lem:Rockland.invertibility} we can give the following  criterion for a selfadjoint hypoelliptic 
 differential operator to have a positive principal symbol. 
 
 \begin{proposition}\label{prop:Heat1.positivity-criterion}
Suppose that $P$ satisfies the Rockland condition at every point. Then the following are equivalent:\smallskip
   
   (i) The operator $P$ is bounded from below.\smallskip
   
   (ii) The principal symbol $\sigma_{m}(P)$ of $P$ is positive.
  \end{proposition}
\begin{proof}
      Assume that $P$ is bounded from below. Possibly by replacing $P$ by $P+c$ with $c$ large enough we may assume that 
    $P$ is positive. Observe that $P^{2}$ satisfies the Rockland condition at every point. Indeed, at every point $x\in M$ the model operator of 
    $P^{2}$ is $(P^{2})^{x}=(P^{x})^{2}$, so for any nontrivial irreducible representation $\pi$ of $G_{x}M$ on $C^{\infty}(\pi)$ 
    we have $\overline{\pi_{(P^{2})^{x}}}=\overline{\pi_{(P^{x})^{2}}}=\overline{\pi_{P^{x}}}\overline{\pi_{P^{x}}}$. Since $\overline{\pi_{P^{x}}}$ 
    is injective on $C^{\infty}(\pi)$ we see that so is $\overline{\pi_{(P^{2})^{x}}}$, that is, $P^{2}$ satisfies the Rockland condition 
    at $x$. 
    
    Granted this, we can apply Lemma~\ref{lem:Rockland.invertibility} to deduce that the principal symbol of $P^{2}+\partial_{t}$ is an invertible 
    Volterra-Heisenberg symbol. As $P^{2}$ is positive we then may use Theorem~\ref{thm:Powers1.main} 
    to see that $(P^{2})^{\frac{1}{4}}=P^{\frac{1}{2}}$ is a \psivdo\ of order $\frac{m}{2}$. 
    Since $P=(P^{\frac{1}{2}})^{2}$ it then follows from Lemma~\ref{lem:Rockland-heat.positivity} that the principal symbol of $P$ is a positive symbol. 

    Conversely, suppose that $P$ has a positive principal symbol. Then thanks to Lemma~\ref{lem:Rockland-heat.positivity} the operator $P$ can be written as
    $P=Q^{*}Q+R$ with $Q\in \pvdo^{\frac{m}{2}}(M,\cE)$ and $R\in \pvdo^{m-1}(M,\cE)$.  Let $P_{1}=Q^{*}Q$. 
    Since $P$ and $P_{1}$ have same principal symbol, the operator $P_{1}$ also satisfies the Rockland condition at every point. 
    
      Next, since $P_{1}$ is positive for $\lambda<-1$ we have
    \begin{equation}
        (P_{1}-\lambda)^{-1}=\int_{0}^{\infty}e^{-tP_{1}}e^{t\lambda}dt,
    \end{equation}
    where the integral converges in $\cL(L^{2}(M,\cE))$. Let $\alpha=\frac{m-1}{m}$. Then we have
    \begin{equation}
        R(P_{1}-\lambda)^{-1}=RP^{-\alpha}_{1}R_{(\lambda)}, \qquad R_{(\lambda)}=\int_{0}^{\infty}P^{\alpha}_{1}e^{-tP_{1}}e^{t\lambda}dt.
         \label{eq:Heat1.decomposition-Rlambda}
    \end{equation}

    Since $P_{1}=Q^{*}Q$ has a positive principal symbol and satisfies the Rockland condition, Lemma~\ref{lem:Rockland.invertibility} tells us that 
    the principal symbol of $P_{1}+\partial_{t}$ is an invertible Volterra-Heisenberg symbol. This allows us to apply  
    Theorem~\ref{thm:Powers1.main2} to deduce that $P^{-\alpha}_{1}$ is a \psivdo\ of order $-m\alpha=-(m-1)$. 
    Therefore $RP^{-\alpha}_{1}$ is a \psivdo\ of order~$0$, hence is bounded on $L^{2}(M,\cE)$.  
  
    On the other hand, we have 
  \begin{equation}
      \|R_{(\lambda)}\|_{L^{2}}\leq \int_{0}^{\infty}t^{-\alpha} \|(tP_{1})^{\alpha}e^{-tP_{1}}\|e^{t\lambda}dt.
  \end{equation}
  As $\alpha \in(0,1)$  the function $x\rightarrow x^{\alpha}e^{-x}$ maps $[0,\infty)$ to $[0,1]$. Therefore, we have
  $\|(tP_{1})^{\alpha}e^{-tP_{1}}\|\leq 1$, from which we get
  \begin{equation}
       \|R_{(\lambda)}\|_{L^{2}}\leq \int_{0}^{\infty}t^{-\alpha} e^{t\lambda}dt=|\lambda|^{\alpha-1} \int_{0}^{\infty}u^{-\alpha} e^{-u}du.
       \label{eq:Heat1.Rlambda-estimates}
  \end{equation}
  
 Now, by combining~(\ref{eq:Heat1.decomposition-Rlambda})--(\ref{eq:Heat1.Rlambda-estimates}) we obtain
  \begin{equation}
      \|R(P_{1}-\lambda)^{-1}\|_{L^{2}}\leq  \|RP^{-\alpha}_{1}\|_{L^{2}}\|R_{(\lambda)}\|_{L^{2}}\leq C_{\alpha}|\lambda|^{\alpha-1},
  \end{equation}
  where the constant $C_{\alpha}$ does not depend on $\lambda$. Since $\alpha<1$ it follows that for $\lambda$ negatively large enough we have 
  $\|R(P_{1}-\lambda)^{-1}\|_{L^{2}}\leq \frac{1}{2}$, so that $1+R(P_{1}-\lambda)^{-1}$ is invertible. 
 Since we have 
 \begin{equation}
     P-\lambda=P_{1}+R-\lambda=(P_{1}+R-\lambda)(P_{1}-\lambda)^{-1},
 \end{equation}
it follows that $P-\lambda$ has a right inverse for $\lambda$ negatively large enough. Since we can similarly  show that for $\lambda$ negatively large enough 
$P-\lambda$ is left-invertible, we deduce that 
as soon as $\lambda$ negatively large enough $P-\lambda$ admits a bounded two-sided inverse. This means
 that the spectrum of $P$ is contained in some interval $[c,\infty)$, that is, $P$ is bounded from below.
\end{proof}

Now, by combining Lemma~\ref{lem:Rockland.invertibility} and Proposition~\ref{prop:Heat1.positivity-criterion} 
with Proposition~\ref{thm:volterra.inverse} and 
Proposition~\ref{thm:Powers1.heat-kernel-asymptotics} we get the main result of this section: 

\begin{theorem}\label{thm:Heat1.main}
 Assume that the Levi form of $(M,H)$ has constant rank and let $P:C^{\infty}(M,\cE)\rightarrow C^{\infty}(M,\cE)$ be a selfadjoint differential operator of even 
 Heisenberg order $m$ such that $P$ is bounded from below and satisfies the Rockland condition at every point. Then:\smallskip 
   
   1) The principal symbol $P+\partial_{t}$ is an invertible Volterra-Heisenberg symbol;\smallskip 
   
   2) The heat operator $P+\partial_{t}$ has an inverse in $\pvhdo^{-m}(M\times\R_{(m)},\cE)$;\smallskip  
   
   3) The heat kernel $k_{t}(x,y)$ of $P$ has an asymptotics in  $C^{\infty}(M,(\End \cE)\otimes|\Lambda|(M))$ of the form
    \begin{equation}
     k_{t}(x,x) \sim_{t\rightarrow 0^{+}} t^{-\frac{d+2}{m}} \sum t^{\frac{2j}{m}} a_{j}(P)(x), \qquad  
    a_{j}(P)(x) =|\varepsilon_{x}'|\check{q}_{-m-2j}(x,0,1),
    \label{eq:Rockland-Heat.heat-kernel-asymptotics}
    \end{equation}
where the equality on the right-hand side shows how to compute $a_{j}(P)(x)$ in a local trivializing Heisenberg chart by means 
of the (local) symbol $q_{-m-2j}(x,\xi,\tau)$ of degree $-m-2j$ of any Volterra-\psivdo\ parametrix for $P+\partial_{t}$ in this chart.
  \end{theorem}

  \begin{example}
    Theorem~\ref{thm:Heat1.main} holds for the following examples:\smallskip 
    
    (a) The conformal powers $\boxdot_{\theta}^{(k)}$ of the horizontal sublaplacian acting on functions on a compact strictly pseudoconvex CR manifold;\smallskip 

    (b) The contact Laplacian on a compact orientable contact manifold.

\end{example}

 \section{Spectral asymptotics for hypoelliptic operators}
 \label{sec:Spectral}
In this section we will make use of Theorem~\ref{thm:Heat1.main} to derive spectral asymptotics for hypoelliptic operators on Heisenberg manifolds. 

Let $(M^{d+1},H)$ be a compact Heisenberg manifold equipped with a smooth positive density and let $\cE$ be a Hermitian vector bundle over $M$. 
Let $P:C^{\infty}(M,\cE)\rightarrow C^{\infty}(M,\cE)$ be a selfadjoint differential operator of Heisenberg order $m$ such that $P$ is bounded from below and 
and the principal symbol of $P+\partial_{t}$ is an invertible Volterra-Heisenberg symbol. Recall that by Proposition~\ref{thm:Powers1.heat-sublaplacians} 
and Theorem~\ref{thm:Heat1.main} the latter condition is satisfied in the following cases:\smallskip

- The operator $P$ is a sublaplacian and satisfies the condition~(\ref{eq:Rockland.sublaplacian'}) at every point.\smallskip

- The Levi form of $(M,H)$ has constant rank and $P$ satisfies the Rockland condition at every point. \smallskip

Let $\lambda_{0}(P)\leq \lambda_{1}(P)\leq \ldots$ denote the eigenvalues of $P$ counted with multiplicity and let $N(P;\lambda)$ be the counting 
function of $P$, that is, 
\begin{equation}
     N(P;\lambda)=\#\{k\in \N;\ \lambda_{k}(P)\leq \lambda \}, \qquad \lambda \in \R. 
\end{equation}

\begin{theorem}\label{thm:Heat1.spectral-asymptotics}
Under the above assumptions the following hold.\smallskip
 1) As $t\rightarrow 0^{+}$ we have
     \begin{equation}
         \Tr e^{-tP} \sim t^{-\frac{d+2}{m}} \sum t^{\frac{2j}{m}} A_{j}(P), \qquad  
   A_{j}(P)=\int_{M}\tr_{\cE} a_{j}(P)(x),
          \label{eq:Heat1.heat-trace-asymptotics}
     \end{equation}
     where the density $a_{j}(P)(x)$ is the coefficient of $t^{\frac{j-d+2}{m}}$ in the heat kernel 
     asymptotics~(\ref{eq:Rockland-Heat.heat-kernel-asymptotics}) for $P$.\smallskip
     
     2) We have $A_{0}(P)>0$.\smallskip 
     
     3) As $\lambda\rightarrow \infty$ we have 
\begin{equation}
    N(P;\lambda) \sim \nu_{0}(P)\lambda^{\frac{d+2}{m}}, \qquad 
    \nu_{0}(P)=\Gamma(1+\frac{d+2}{m})^{-1}A_{0}(P).
    \label{eq:Heat1.counting-function-asymptotics}
\end{equation}

  4) As $k\rightarrow \infty$ we have 
     \begin{equation}
         \lambda_{k}(P)\sim \left(\frac{k}{\nu_{0}(P)}\right)^{\frac{m}{d+2}}. 
         \label{eq:Heat1.eigenvalue-asymptotics}
     \end{equation}     
 \end{theorem}
 \begin{proof}
 First,  as $\Tr e^{-tP}=\int_{M}\tr_{\cE}k_{t}(x,x)$ the asymptotics~(\ref{eq:Heat1.heat-trace-asymptotics}) immediately follows 
     from~(\ref{eq:Rockland-Heat.heat-kernel-asymptotics}). 
     
     Second, we  have $A_{0}(P)=\lim_{t\rightarrow 0^{+}}t^{\frac{d+2}{m}}\Tr e^{-tP}\geq 0$, so if we can show that $A_{0}(P)\neq 0$ then we get 
     $A_{0}(P)>0$ and the 
     asymptotics~(\ref{eq:Heat1.counting-function-asymptotics}) will follow from~(\ref{eq:Heat1.heat-trace-asymptotics}) by Karamata's Tauberian theorem 
     (see \cite[Thm.~108]{Ha:DS}). This will also yield~(\ref{eq:Heat1.eigenvalue-asymptotics}), since the latter is equivalent 
     to~(\ref{eq:Heat1.counting-function-asymptotics}) (e.g.~\cite[Sect.~13]{Sh:POST}). Therefore, the bulk of the proof is to show that $A_{0}(P)\neq 0$.
       
  Now, notice that there is at least one integer~$< \frac{m}{d+2}$ such that $A_{j}(P)\neq 0$. Otherwise, by~(\ref{eq:Heat1.heat-trace-asymptotics}) 
    there would exist a constant $C> 0$ such that $\Tr e^{-tP}\leq C$ for $0\leq t\leq 1$. Thus, 
    \begin{equation}
        k e^{-t\lambda_{k}(P)}\leq \sum_{j<k} e^{-t\lambda_{j}(P)}\leq \Tr e^{-tP}\leq C, \qquad 0<t<1.
    \end{equation}
   Therefore, letting $t\rightarrow 0^{+}$ would give  $k\leq C$ for every $k\in \N$, which is not possible.
   
   On the other hand, when $m\geq d+2$ the only integer~$<\frac{m}{d+2}$ is $j=0$, so in this case we must have $A_{0}(P)\neq  0$. 
    
    Next, assume $m<d+2$ and $A_{0}(P)=0$. Let $\mu=\frac{d+2}{m}-j_{0}$ where $j_{0}$ is the smallest integer $j$ such that $A_{j}(P)\neq 0$. Since 
   $1\leq j_{0}< \frac{d+2}{m}$ we have $0<\mu\leq \frac{d+2}{m}-1$. Moreover, the asymptotics~(\ref{eq:Heat1.heat-trace-asymptotics}) becomes
   \begin{equation}
       \Tr e^{-tP}=A_{j_{0}}(P)t^{-\mu}+\op{O}(t^{1-\mu}) \qquad \text{as $t\rightarrow 0^{+}$}. 
   \end{equation}
   This implies that we have $A_{j_{0}}(P)\geq 0$. Since we have $A_{j_{0}}(P)\neq 0$ by definition of $j_{0}$, we get 
   $A_{j_{0}}(P)>0$. Therefore, as alluded to above, it follows from Karamata's Tauberian theorem that as $k\rightarrow \infty$ we have
   \begin{equation}
       \lambda_{k}(P)\sim  \left(\frac{k}{\beta}\right)^{\frac{1}{\mu}}, \qquad \beta=\Gamma(1+\mu)^{-1}A_{j_{0}}(P). 
   \end{equation}
It follows that $\lambda_{k}(P^{- \frac{d+2}{2m}})=\op{O}(k^{-\frac{1}{2}-\delta})$, with 
   $\delta=\frac{1}{2}(\frac{1}{\mu}\frac{d+2}{m}-1)>0$. In particular, we have $\sum_{k\geq 0}\lambda_{k}(P^{- \frac{d+2}{2m}})^{2}<\infty$, that is,
 $P^{- \frac{d+2}{2m}}$ is a Hilbert-Schmidt operator on $L^{2}(M,\cE)$. 
   
   In addition, observe that $P^{- \frac{d+2}{2m}}$ is a \psivdo\ of order $-\frac{(d+2)}{2}$ and that any $Q \in \pvdo^{-\frac{(d+2)}{2}}(M,\cE)$ can be written as 
   \begin{equation}
       Q=\Pi_{0}(P)Q +(1-\Pi_{0}(P))P^{-\frac{d+2}{2m}}P^{\frac{d+2}{2m}}Q,
        \label{eq:Spectral.Hilbert-Schmidt}
   \end{equation}
   where $\Pi_{0}(P)$ denotes the orthogonal projection onto $\ker P$. Recall that the space of Hilbert-Schmidt operators is a two-sided ideal of 
   $\cL(L^{2}(M<\cE)$. Observe that in~(\ref{eq:Spectral.Hilbert-Schmidt}) the projection $\Pi_{0}(P)$ is a smoothing operator, so is a Hilbert-Schmidt 
   operator and  $P^{\frac{d+2}{2m}}Q$ is a bounded operator on $L^{2}(M,\cE)$ because this 
   a zero'th order \psivdo. It thus follows that any $Q\in \pvdo^{-\frac{(d+2)}{2}}(M,\cE)$ is a Hilbert-Schmidt operator on $L^{2}(M,\cE)$.

    Now, we get a contradiction as follows. Let $\kappa: U \rightarrow V$ be a Heisenberg chart over which there is a trivialization 
    $\tau:\cE_{_{U}}\rightarrow U\times \C^{r}$ of $\cE$ 
    and such that the open $V\subset \Rd$ is bounded. 
    Let $\varphi \in C^{\infty}_{c}(\Rd)$ have non-empty support $L$, let $\chi \in C^{\infty}_{c}(V\times V)$ be such that $\chi(x,y)=1$ near $L\times L$, and
   let $Q:C_{c}^{\infty}(V)\rightarrow C^{\infty}(V)$ be given by the kernel,
    \begin{equation}
        k_{Q}(x,y)=|\varepsilon_{x}'|\varphi(x)\|\varepsilon_{x}(y)\|^{-\frac{d+2}{2}}\chi(x,y).
    \end{equation}
    
    The kernel $k_{Q}(x,y)$ has a compact support contained in $V\times V$ and, as $\varphi(x)(1-\chi(x,y))$ vanishes near the diagonal of $V\times V$, we have 
\begin{equation}
         k_{Q}(x,y)=|\varepsilon_{x}'|\varphi(x)\|-\varepsilon_{x}(y)\|^{-\frac{d+2}{2}} \bmod C^{\infty}(V\times V). 
\end{equation}
     Since $\varphi(x)\|y\|^{-\frac{d+2}{2}}$ belongs to $\cK^{-\frac{d+2}{2}}(V\times \Rd)$,
     it follows from Proposition~\ref{prop:PsiVDO.characterisation-kernel2} 
    that $Q$ is a \psivdo\ of order $-\frac{d+2}{2}$. 
    
    Let $Q_{0}=\tau^{*}\kappa^{*}(Q\otimes 1)$. Then $Q_{0}$ belongs to $\pvdo^{-\frac{d+2}{2}}(M,\cE)$, hence is a Hilbert-Schmidt operator on 
    $L^{2}(M,\cE)$. This implies that  $Q$ is a Hilbert-Schmidt operator on $L^{2}(V)$, so by~\cite[p.~109]{GK:ITLNSO} its kernel lies in $L^{2}(V\times V)$. 
    This cannot be true, however, because we have 
   \begin{multline}
       \int_{V\times V}|k_{Q}(x,y)|^{2}dxdy\geq \int_{L\times L} |\varepsilon_{x}'|^{2}|\varphi(x)|^{2}\|\varepsilon_{x}(y)\|^{-(d+2)}dxdy \\  
       = \int_{L}|\varepsilon_{x}'||\varphi(x)|^{2}(\int_{\varepsilon_{x}(L)}\|y\|^{-(d+2)}dy)dx=\infty.
   \end{multline}
   We have thus obtained a contradiction, so we must have $A_{0}(P)\neq 0$. The proof is now complete. 
\end{proof}
 
\begin{example}
    Theorem~\ref{thm:Heat1.spectral-asymptotics} is valid for the following operators:\smallskip
  
   (a) Real selfadjoint sublaplacian $ \Delta=\nabla_{X_{1}}^{*}\nabla_{X_{1}}+\ldots+\nabla_{X_{m}}^{*}\nabla_{X_{m}}^{2}+\mu(x)$, 
   where $\nabla$ is connection on $\cE$, the vector fields $X_{1},\ldots,X_{m}$ span $H$ and $\mu(x)$ is a selfadjoint section of $\End \cE$, 
   provided that the Levi form of $(M,H)$ is non-vanishing; \smallskip 
  
   (b) The Kohn Laplacian on a compact CR manifold and acting on $(p,q)$-forms when the condition $Y(q)$ holds everywhere;\smallskip 
  
   (c) The horizontal sublaplacian on a compact Heisenberg manifold acting on horizontal forms of degree $k$ when the condition $X(k)$ holds 
   everywhere;\smallskip 
   
   (d) The contact Laplacian on a compact orientable contact manifold.\smallskip

   (e) The conformal powers of the horizontal sublaplacian acting on functions on a compact strictly pseudoconvex CR manifold.\smallskip 
\end{example}
 
Several authors have obtained Weyl asymptotics closely related to~(\ref{eq:Intro1.counting-function-asymptotics}) for bicharacteristic hypoelliptic 
 operators~(\cite{II:PDPEAABSFSP}, \cite{Me:HOCVC2WE}, \cite{MS:ECH2}), including sublaplacians on Heisenberg manifolds, and for more general hypoelliptic 
 operators~(\cite{Mo:ESOHCM1}, \cite{Mo:ESOHCM2}) using different approaches involving other pseudodifferential calculi. 
 
 While these authors deal in a setting more general than the Heisenberg setting, as far as the Heisenberg setting is concerned, our approach using the 
 Volterra-Heisenberg calculus presents the following advantages:\smallskip
 
 (i) The pseudodifferential analysis is somewhat simpler, since the Volterra-Heisenberg calculus yields for free 
 a heat kernel asymptotics, once the principal symbol of the heat operator is shown to be invertible, for which it is enough to use the Rockland condition 
 in many cases;\smallskip 
 
 (ii) As the Volterra-Heisenberg calculus take fully into account the underlying Heisenberg geometry of the manifold and is invariant by change of 
 Heisenberg coordinates, we can explicitly compute the coefficient $\nu_{0}(P)$ in~(\ref{eq:Heat1.counting-function-asymptotics}) 
 for operators admitting a normal form. 
 This is illustrated  below by Proposition~\ref{prop:Spectral.normal-form-nu0P}, which will be used in the next two section to give geometric 
 expressions for the Weyl asymptotics~(\ref{eq:Heat1.counting-function-asymptotics}) for geometric operators on CR and contact manifolds.\smallskip

Next, prior to dealing with operators admitting a normal form we need the following.
 \begin{lemma}\label{lem:Spectral.a0(P)}
    Let $P:C^{\infty}(M,\cE)\rightarrow C^{\infty}(M,\cE)$ be a selfadjoint differential operator of Heisenberg order $m$ such that $P$ is bounded 
    from below and the principal symbol of $P+\partial_{t}$ is an invertible Volterra-Heisenberg symbol. Then, for any $a\in M$, the following 
    hold.\smallskip
    
    1) The model operator $P^{a}+\partial_{t}$ admits a unique fundamental solution $K^{a}(x,t)$ which, in addition, belongs  to the class 
    $\cK_{\op{v},-(d+2)}(G_{a}M\times\R_{(m)})$.\smallskip
    
    2) In Heisenberg coordinates centered at $a$ we have
    \begin{equation}
        a_{0}(P)(0)=K^{a}(0,1).
    \end{equation}
\end{lemma}
\begin{proof}
    Since it is enough to prove the result in a trivializing Heisenberg chart near $a$, we may as well assume that $\cE$ is the trivial line bundle, 
    since the proof in the general case follows along similar lines. 
    
    Let $\sigma_{-m}(Q)\in S_{\op{v},-m}(\fg^{*}U\times 
    \R_{(m)})$ be the inverse of $\sigma_{m}(P)+i\tau$ and let $Q^{a}$ be the left-invariant Volterra-\psivdo\ on $G_{a}U\times \R$ with symbol 
    $\sigma_{-m}(Q)(a,.,.)$. In particular, the operator $Q^{a}$ is the inverse of $P^{a}+\partial_{t}$ on $\cS(G_{a}U\times \R)$ and it agrees 
    with the left-convolution by $K^{a}=[\sigma_{-m}(a,.,t)]_{(\xi,\tau)\rightarrow (y,t)}^{\vee}$. 
   Therefore,  the left-convolution operator with $[(P^{a}+\partial_{t})K^{a}]$ agrees with $(P^{a}+\partial_{t})Q^{a}=1$. Thus, 
    \begin{equation}
        (P^{a}+\partial_{t})K^{a}(y,t)=\delta(y,t),
    \end{equation}
    that is, $K^{a}(y,t)$ is a fundamental solution for $P^{a}+\partial_{t}$. 
    
    Let $K\in \cS'(G_{a}U\times \R)$ be another fundamental solution for $P^{a}+\partial_{t}$.
    As it follows from the proof of Lemma~\ref{lem:Rockland.invertibility}, the left-convolution operator $Q$ by $K$ 
   is a right-inverse for $P^{a}+\partial_{t}$, so agrees with $Q^{a}$. Hence $K=K^{a}$, which shows that $K^{a}$ is the unique fundamental solution
  of  $P^{a}+\partial_{t}$.
  
  In addition, $K^{a}$ is in the class $\cK_{\op{v},-(d+2)}(G_{a}M\times\R_{(m)})$ because this the inverse 
  Fourier transform of a symbol in $S_{\op{v},-m}(\fg^{*}U\times \R_{(m)})$. 

   Finally, let $q_{-m}\in S_{\op{v},-m}(\URd\times\R_{(m)})$ be the local principal symbol for a parametrix $Q$ for $P+\partial_{t}$ on $U\times \R$. 
  As $Q$ has global principal symbol $\sigma_{-m}(Q)$ it follows from Proposition~\ref{prop:Powers1.principal-symbol} 
  that $q_{-m}(0,.,.)=\sigma_{-m}(Q)(a,.,.)$. Moreover, 
   since we are in the Heisenberg coordinates centered at $a$, the map $\varepsilon_{0}$ to the Heisenberg coordinates centered at $0$ is just the 
   identity. Combining this with~(\ref{eq:Heat1.counting-function-asymptotics}) then gives
   \begin{equation}
       a_{0}(P)(0)=|\varepsilon_{0}'|\check{q}_{m}(0,0,1)=[\sigma_{-m}(Q)]_{(\xi,\tau)\rightarrow (y,t)}(a,0,1)=K^{a}(0,1).
   \end{equation}
   The proof is thus achieved.
 \end{proof}

Now, assume that the Levi form of $(M,H)$ has constant rank $2n$. Therefore, by Proposition~\ref{prop:Bundle.intrinsic.fiber-structure} 
the tangent Lie group bundle $GM$ is a fiber bundle with typical fiber $G=\bH^{2n+1}\times \R^{d-2n}$. 

Let $P:C^{\infty}(M,\cE)\rightarrow C^{\infty}(M,\cE)$ be a selfadjoint differential operator of Heisenberg order $m$ such that $P$ is bounded 
    from below and satisfies the Rockland condition at every point. 
    
    We further assume that the density $d\rho(x)$ of $M$ and the model operator of $P$ 
    at a point admit normal forms. By this it is meant that there exist $\rho_{0}>0$ and a differential operator $P_{0}:C^{\infty}(G,\C^{r})\rightarrow 
    C^{\infty}(G,\C^{r})$  such that near every point $x_{0}$, there exist trivializing Heisenberg coordinates, herewith called 
    normal trivializing Heisenberg coordinates centered at $x_{0}$, in which the following hold:\smallskip 
    
    (i) We have $d\rho(x)_{|_{x=0}}=[\rho_{0}dx]_{|_{x=0}}$;\smallskip 
    
    (ii) At $x=0$ the tangent group is $G$ and the model operator of $P$ is $P_{0}$.\smallskip
    
In the sequel, we let $\op{vol}_{\rho}M$ denote the volume of $M$ with respect to $\rho$, that is, 
\begin{equation}
    \op{vol}_{\rho}M=\int_{M}d\rho(x).
\end{equation}
    As it turns out the assumptions (i) and (ii) allows us to relate the Weyl asymptotics~(\ref{eq:Heat1.counting-function-asymptotics}) 
    for $P$ to $\op{vol}_{\rho}M$ as follows.

\begin{proposition} \label{prop:Spectral.normal-form-nu0P}
Under the assumptions (i) and (ii) above, as $\lambda \rightarrow \infty$ we have
\begin{equation}
  N(P;\lambda)\sim \rho_{0}^{-1}\nu_{0}(P_{0})(\op{vol}_{\rho}M)\lambda^{\frac{d+2}{m}}, 
  \qquad \nu_{0}(P_{0})=\Gamma(1+\frac{d+2}{m})^{-1}\tr_{\C^{r}}K_{0}(0,1),
     \label{eq:Spectral.normal-form-nu0P}
\end{equation}
where $K_{0}(x,t)$ denotes the fundamental solution of $P_{0}+\partial_{t}$.
\end{proposition}
\begin{proof}
   By Proposition~\ref{thm:Heat1.spectral-asymptotics} as $\lambda \rightarrow \infty$ we have 
    \begin{equation}
        N(P;\lambda) \sim \nu_{0}(P)\lambda^{\frac{d+2}{m}}, \qquad \nu_{0}(P)=\Gamma(1+\frac{d+2}{m})^{-1}\int_{M}\tr_{\cE}a_{0}(P)(x).
        \label{eq:Heat1.counting-function-asymptotics2}
    \end{equation}

    On the other hand, by Lemma~\ref{lem:Spectral.a0(P)} in normal trivializing Heisenberg coordinates centered at point $a\in M$ we have
    \begin{equation}
      [ \tr_{\C^{r}} a_{0}(P)(x)dx]_{|_{x=0}}=[\tr_{\C^{r}}K_{0}(0,1)dx]_{|_{x=0}}=\tr_{\C^{r}}K_{0}(0,1)\rho_{0}^{-1}[d\rho(x)]_{|_{x=0}}.
    \end{equation}
  Hence $\tr_{\cE}a_{0}(P)(x)=\rho_{0}^{-1}\tr_{\C^{r}}K_{0}(0,1)d\rho(x)$. Thus, 
  \begin{equation}
   \nu_{0}(P)=\Gamma(1+\frac{d+2}{m})^{-1}\int_{M}\rho_{0}^{-1}\tr_{\C^{r}}K_{0}(0,1)d\rho(x)=\rho_{0}^{-1}\nu_{0}(P_{0})\op{vol}_{\rho}M,
  \end{equation}
  with $\nu_{0}(P_{0})=\Gamma(1+\frac{d+2}{m})^{-1}\tr_{\C^{r}}K_{0}(0,1)$. Combining this~(\ref{eq:Heat1.counting-function-asymptotics2}) 
  then proves the lemma.
\end{proof}

\section{Weyl asymptotics and CR geometry}
\label{sec:Spectral-CR}
The aim of this section is to express in geometric terms the Weyl asymptotics~(\ref{eq:Heat1.counting-function-asymptotics}) 
for the Kohn Laplacian and the horizontal sublaplacian on a nondegenerate CR manifold with a Levi metric.  

Let $(M^{2n+1},\theta)$ be a  $\kappa$-strictly pseudoconvex pseudohermitian manifold with $0\leq \kappa \leq \frac{n}{2}$. 
Let $T_{1,0}\subset T_{\C}M$ be the CR structure of $M$ and define $T_{0,1}=\overline{T_{0,1}}$ and $H=\Re (T_{1,0}\otimes T_{0,1})$. 
The Levi form $L_{\theta}$ on $T_{1,0}$ is given by 
\begin{equation}
     L_{\theta}(Z,W)=-id\theta (Z,\bar W)=i\theta([Z,\bar W]), \qquad Z,W \in T_{1,0}.
     \label{eq:Heat1.Levi-form}
\end{equation}

Let $X_{0}$ be the Reeb field associated to $\theta$, so that $\imath_{X_{0}} \theta=1$ and $\imath_{X_{0}} d\theta=0$. Then we have
  \begin{equation}
    [Z,\bar W]= -iL_{\theta}(Z,W)X_{0} \quad \bmod T_{1,0}\oplus T_{0,1} \qquad \forall Z, W \in T_{1,0}M.
     \label{eq:Heat1.commutators-metric}
\end{equation}

We endow $M$ with a Levi metric as follows (see also~\cite{FS:EDdbarbCAHG}). 
Let $\tilde{h}$ be a (positive definite) Hermitian metric on $T_{1,0}$. Then there exists a Hermitian-valued section $A$ of $\End_{\C}T_{1,0}$ such that 
\begin{equation}
    L_{\theta}(Z,W) =\tilde{h}(AZ,W), \qquad Z,W\in T_{1,0}.
\end{equation}

Let $A=U|A|$ be the polar decomposition of $A$. Since by assumption the Hermitian form $L_{\theta}$ has signature $(n-\kappa,\kappa,0)$, hence is 
nondegenerate, the section $A$ is invertible and the section $U$ is an orthogonal matrix with only the eigenvalues $ 1$ and $-1$ with multiplicities $n-\kappa$ 
and $\kappa$ respectively. Therefore, we have the splitting,
\begin{equation}
    T_{1,0}=T_{1,0}^{+}\oplus T_{1,0}^{-}, \qquad T_{1,0}^{\pm}=\ker (U\mp 1), 
    \label{eq:Heat1.positive-negative-splitting}
\end{equation}
where the restriction of $L_{\theta}$ to $T_{1,0}^{+}$ (resp.~$T_{1,0}^{-}$) is positive definite (resp.~negative definite). We then get a Levi 
metric on $T_{1,0}$ by letting
\begin{equation}
    h(Z,W)=\tilde{h}(|A|Z,W), \qquad Z,W\in T_{1,0}.
\end{equation}
In particular, on the direct summands $T_{1,0}^{\pm}$ of the splitting~(\ref{eq:Heat1.positive-negative-splitting}) we have
\begin{equation}
    L_{\theta}(Z,W)=h(UZ,W)=\pm h(Z,W), \qquad Z,W\in T_{1,0}^{\pm}.
     \label{eq:Heat1.Levi-form.Levi-metric}
\end{equation}

We now extend $h$ into a Hermitian metric $h$ on $T_{\C}M$ by making the following requirements:
\begin{gather}
  h(X_{0},X_{0})=1, \qquad   h(Z,W)=\overline{h(\bar Z,\bar W)} \quad \forall Z, W\in T_{0,1},\\ 
  \text{The splitting $T_{\C}M=\C X_{0}\oplus T_{1,0}\oplus T_{0,1}$ is orthogonal with respect to $h$}.
\end{gather}

This allows us to express the Levi form $\cL_{\C}:(H\otimes \C)\times (H\otimes \C)\rightarrow T_{\C}M/(H\otimes \C)$ as follows. 
Since $\theta(X_{0})=1$  we have 
\begin{equation}
    \cL_{\C}(X,Y)=\theta([X,Y])X_{0}=-d\theta(X,Y)X_{0}, \qquad X,Y\in H\otimes \C. 
\end{equation}
Therefore, if follows from~(\ref{eq:Heat1.Levi-form.Levi-metric}) that we have 
\begin{equation}
    \cL(Z,\bar W)=-iL_{\theta}(Z,W)=h(Z,iUW), \qquad Z,W\in T_{1,0}.
\end{equation}
Since $\cL$ is antisymmetric and the integrability condition $[T_{1,0},T_{1,0}]\subset T_{1,0}$ implies that $\cL_{\C}$ vanishes on $T_{1,0} \times 
T_{1,0}$ and on $T_{0,1} \times T_{0,1}$, we get
\begin{equation}
       \cL_{\C}(X,Y)=h(X,LY), \qquad X,Y\in  H\otimes \C,
\end{equation}
where $L$ is the antilinear antisymmetric section of $\End_{\R} (H\otimes \C)$ such that
\begin{equation}
    L(Z+\bar{W})=iUW-i(\overline{UZ}), \qquad Z,W\in T_{1,0}.
     \label{eq:Spectral.Levi-form-Levi-metric2}
\end{equation}
In particular, since $U^{*}=U$ and $U^{2}=1$ we have $|L|=1$. 
Let $Z_{1},\ldots,Z_{n}$ be a local orthonormal frame for $T_{1,0}$ (with respect to $h$) and such that $Z_{1},\ldots,Z_{n-\kappa}$ 
span $T_{1,0}^{+}$ and $Z_{n-\kappa+1},\ldots,Z_{n}$ span $T_{1,0}^{-}$. Then $\{X_{0},Z_{j},\bar Z_{j}\}$ is an orthonormal frame. In the sequel we 
will call such a frame an \emph{admissible orthonormal frame of} $T_{\C}M$. Then from~(\ref{eq:Heat1.Levi-form.Levi-metric}) we get: 
\begin{equation}
 L_{\theta}(Z_{j},\bar Z_{k})=\epsilon_{j}h(Z_{j},Z_{k})=\epsilon_{j}\delta_{jk},
\label{eq:Heat1.dtheta.admissible-frame}
\end{equation}
where $\epsilon_{j}=1$ for $j=1,\ldots,n-\kappa$ and $\epsilon_{j}=-1$ for $j=n-\kappa+1,\ldots,n$.

Let $\{\theta,\theta^{j},\theta^{\bar j}\}$ be the dual coframe of $T^{*}_{\C}$ associated to $\{X_{0},Z_{j},\bar Z_{j}\}$. Then the volume form of 
$h$ is locally given by 
\begin{equation}
    \sqrt{h(x)} dx=i^{n}\theta\wedge \theta^{1}\wedge \theta^{\bar 1}\wedge \cdots \wedge \theta^{n}\wedge \theta^{\bar n}.
\end{equation}
On the other hand, because of~(\ref{eq:Heat1.dtheta.admissible-frame}) we have 
$d\theta =i \sum_{j=1}^{n} \epsilon_{j}\theta^{j}\wedge \theta^{\bar j} \bmod \theta \wedge T^{*}M$, so the form $\theta \wedge d\theta^{n}$ is equal to
\begin{multline}
    n! i^{n} \epsilon_{1}\ldots \epsilon_{n}\theta\wedge \theta^{1}\wedge \theta^{\bar 1}\wedge \cdots \wedge 
    \theta^{n}\wedge \theta^{\bar n}= n! i^{n}(-1)^{\kappa}\theta\wedge \theta^{1}\wedge \theta^{\bar 1}\wedge \cdots \wedge 
    \theta^{n}\wedge \theta^{\bar n}. 
\end{multline}
Therefore, we get the following global formula for the volume form,
\begin{equation}
     \sqrt{h(x)} dx=\frac{(-1)^{\kappa}}{n!}  \theta \wedge d\theta^{n}. 
\end{equation}
In particular, we see that the volume form depends only on $\theta$ and not on the choice of the Levi metric. 

\begin{definition}
    The pseudohermitian volume of $(M,\theta)$ is 
    \begin{equation}
        \op{vol}_{\theta}M=\frac{(-1)^{\kappa}}{n!} \int_{M} \theta \wedge d\theta^{n}.
    \end{equation}
\end{definition}
We shall now relate the asymptotics~(\ref{eq:Heat1.counting-function-asymptotics}) for the Kohn Laplacian and the Tanaka sublaplacian to the volume 
$\op{vol}_{\theta}M$. To this end consider the Heisenberg group $\bH^{2n+1}=\R \times \R^{2n}$ equipped with the standard left-invariant basis of $\fh^{2n+1}$ 
given by
\begin{equation}
    X_{0}=\frac{\partial}{\partial x_{0}}, \quad X_{j}=\frac{\partial}{\partial x_{j}}+x_{n+j}\frac{\partial}{\partial 
    x_{0}}, \quad X_{n+j}=\frac{\partial}{\partial x_{n+j}}-x_{j}\frac{\partial}{\partial 
    x_{0}}, \quad 1\leq j\leq n.
     \label{eq:Heat1.Heisenberg-vector-fields}
\end{equation}

For $\mu \in \C$ let $\cL_{\mu}$ denote  the Folland-Stein sublaplacian, 
\begin{equation}
    \cL_{\mu}=-\frac{1}{2}(X_{1}^{2}+\ldots+X_{2n}^{2})-i\mu X_{0}. 
\end{equation}
The Rockland condition for $\cL_{\mu}$ reduces to $\mu\neq \pm n, \pm(n+1),\ldots$ (see~\cite{FS:EDdbarbCAHG}, \cite{BG:CHM}, \cite{Ta:NCMA}), 
so in this case the operator $\cL_{\mu}$ is hypoelliptic. Moreover, when $\mu\neq \pm n, \pm(n+1),\ldots$ the heat operator $\cL_{\mu}+\partial_{t}$ is also 
hypoelliptic 
and admits a unique fundamental solution since its symbol is an invertible Volterra-Heisenberg symbol (see~\cite{BGS:HECRM}). 

In the sequel it will be convenient to use the variable $x'=(x_{1},\ldots,x_{2n})$. 

\begin{lemma}
 For $|\Re \mu|<n$ the fundamental solution of $\cL_{\mu}+\partial_{t}$ is given by 
    \begin{equation}
       k_{\mu}(x_{0},x',t) = \chi(t) (2\pi t)^{-(n+1)} \int_{-\infty}^{\infty}e^{i t^{-1}x_{0}\xi_{0}-\mu\xi_{0}}(\frac{\xi_{0}}{\sinh \xi_{0}})^{n}
        \exp [-\frac{1}{2t} \frac{\xi_{0}}{\tanh \xi_{0}}|x'|^{2}]d\xi_{0}.
        \label{eq:Heat1.fundamental-solution-cLlambda}
    \end{equation}
    where $\chi(t)$ denotes the characteristic function of $(0,\infty)$. 
\end{lemma}
\begin{proof}
  The fundamental solution $k_{\mu}(x_{0},x',t)$ is solution to the equation, 
   \begin{equation}
       (\cL_{\mu}+\partial_{t})k(x_{0},x',t)=\delta(x_{0})\otimes \delta(x')\otimes \delta(t).
        \label{eq:Heat1.heat-cLmu}
   \end{equation}
   Let us make a Fourier transform with respect to $x_{0}$. Thus, letting $\hat{k}(\xi_{0},x',t)=\hat{k}_{x_{0}\rightarrow 
   \xi_{0}}(\xi_{0},x',t)$, the equation~(\ref{eq:Heat1.heat-cLmu}) becomes
    \begin{gather}
        (\hat{\cL}_{0}+\mu \xi_{0}+\partial_{t})\hat{k}_{\mu}=\delta(x')\otimes \delta(t), 
        \label{eq:Heat1.heat-cLmu-Fourier}\\ 
        \hat{\cL}_{0}= -\frac{1}{2}\sum_{j=1}^{n}(\frac{\partial}{\partial x_{j}}-ix_{n+j}\xi_{0})^{2}- 
        \frac{1}{2}\sum_{j=1}^{n}(\frac{\partial}{\partial x_{n+j}}+ix_{j}\xi_{0})^{2}.
    \end{gather}
    
    Notice that if we fix $\xi_{0}$ we have $\hat{\cL}_{0}=\frac{1}{2}H_{A(\xi_{0})}$, where 
    $H_{A(\xi_{0})}=-\sum_{j=1}^{2n}(\partial_{j}-\sum_{k=1}^{2n}A(\xi_{0})_{jk}x_{k})^{2}$ is the harmonic oscillator associated to the 
    real antisymmetric matrix, 
    \begin{equation}
        A(\xi_{0})=\xi_{0}J, \qquad 
        J=\left(  
        \begin{array}{cc}
            0 & -I_{n}  \\
            I_{n} & 0
        \end{array} \right).
    \end{equation}

    Therefore, the fundamental solution of $\hat{\cL}_{0}+\mu \xi_{0}+\partial_{t}$ is given by a version of the Melher formula 
    (see~\cite{GJ:QPFIPV}, \cite[p.~225]{Po:CMP1}), that is, 
    \begin{multline}
        \hat{k}_{0}(\xi_{0},x',t)= \chi(t) (2\pi t)^{-n} \det{}^{\frac{1}{2}}( \frac{t A(\xi_{0})}{\sinh (tA(\xi_{0}))}) 
        \exp [ -\frac{1}{2t} \acou{ \frac{tA(\xi_{0})}{\tanh (tA(\xi_{0})}x'}{x'}],\\
        =\chi(t) (2\pi t)^{-n} (\frac{t\xi_{0}}{\sinh t\xi_{0}})^{n} \exp [-\frac{1}{2t} \frac{t\xi_{0}}{\tanh t\xi_{0}}|x'|^{2}]. 
    \end{multline}

 A solution of~(\ref{eq:Heat1.heat-cLmu-Fourier}) is now given by $\hat{k}_{\mu}(\xi_{0},x',t)=e^{-\mu\xi_{0}t}\hat{k}_{0}(\xi_{0},x',t)$. 
Moreover, as we have 
\begin{equation}
    |h_{0}(\xi_{0},x',t)|\leq \pi^{-n}|\xi_{0}|^{n}e^{-tn|\xi_{0}|}, 
\end{equation}
we further see that for  $|\Re \mu| <n$ the function $e^{-t\mu 
    \xi_{0}}\hat{k}_{0}$ is integrable with respect to $\xi_{0}$. Since $k_{\mu}(x_{0},x',t)$ is  the inverse Fourier transform with respect to $\xi_{0}$ of 
    $\hat{k}_{\mu}(\xi_{0},x',t)$ it follows that we have
     \begin{equation}
        k_{\mu}(x_{0},x',t) = \chi(t) (2\pi t)^{-(n+1)} \int_{-\infty}^{\infty}e^{i t^{-1}x_{0}\xi_{0}-\mu\xi_{0}}(\frac{\xi_{0}}{\sinh \xi_{0}})^{n}
        \exp [-\frac{1}{2t} \frac{\xi_{0}}{\tanh \xi_{0}}|x'|^{2}]d\xi_{0}.
    \end{equation}
The lemma is thus proved. 
\end{proof}

Next, for $|\Re \mu|<n$ we let 
\begin{equation}
    \nu(\mu)=\frac{1}{(n+1)!}k_{\mu}(0,0,1)= \frac{(2\pi)^{-(n+1)}}{(n+1)!} \int_{-\infty}^{\infty}e^{-\mu\xi_{0}}(\frac{\xi_{0}}{\sinh \xi_{0}})^{n}d\xi_{0}. 
\end{equation}

\begin{lemma}\label{lem:Heat1.alpha0-volume}
  Let $\Delta:C^{\infty}(M,\cE)\rightarrow C^{\infty}(M,\cE)$ be a selfadjoint sublaplacian  which is bounded from below and assume there exists $\mu 
  \in (-n,n)$ such that near any point of $M$ there is an admissible orthonormal frame $Z_{1},\ldots,Z_{n}$ of $T_{1,0}$ with respect to which 
  $\Delta$ takes the form, 
    \begin{equation}
        \Delta= - \sum_{j=1}^{n}(\overline{Z}_{j}Z_{j}+Z_{j}\overline{Z}_{j}) - i\mu X_{0}+ \text{lower order terms}.
         \label{eq:Heat1.principal-term-Deltamu}
    \end{equation}
Then as $\lambda \rightarrow \infty$ we have
     \begin{equation}
        N(\Delta; \lambda) \sim \nu(\mu) \rk \cE(\op{vol}_{\theta}M)\lambda^{n+1}. 
           \label{eq:Heat1.asymptotics-sublaplacian-mu}
     \end{equation}

\end{lemma}
\begin{proof}
Let $Z_{1},\ldots,Z_{n}$ be a local admissible orthonormal frame of $T_{1,0}$. Then from~(\ref{eq:Heat1.Levi-form}) and~(\ref{eq:Heat1.Levi-form.Levi-metric})  
we obtain 
\begin{equation}
    [Z_{j},\bar Z_{k}]= -L_{\theta}(Z_{j},Z_{k})X_{0}=-i\epsilon_{j}\delta_{jk}X_{0} \quad \bmod T_{1,0}\oplus T_{0,1}. 
     \label{eq:Heat1.commutators-admissible-frame}
\end{equation}
In addition, we let $X_{1},\ldots,X_{2n}$ be the vector fields in $H$ such that
 \begin{equation}
     Z_{j}=\left\{ 
      \label{eq:Spectral.convention-Xj-orientation}
 \begin{array}{ll}
     \frac{1}{2}(X_{j}-iX_{n+j})& \text{for $j=1,\ldots,n-\kappa$},\\
     \frac{1}{2}(X_{n+j}-iX_{j}) & \text{for $j=n-\kappa+1,\ldots,n$}.
 \end{array}\right.
 \end{equation}
Then $X_{1},\ldots,X_{0}$ is a local frame of $H=\Re (T_{1,0}\oplus T_{0,1})$ and from~(\ref{eq:Heat1.commutators-admissible-frame}) we get 
\begin{gather}
    [X_{j},X_{n+k}]=-2\delta_{jk}X_{0} \ \bmod H, 
    \label{eq:Heat1.almost-Heisenberg-relations1}
    \\ [X_{0},X_{j}]=[X_{j},X_{k}]=[X_{n+j},X_{n+k}]=0\ \bmod H.
     \label{eq:Heat1.almost-Heisenberg-relations}
\end{gather}

Moreover, in terms of the vector fields $X_{1},\ldots,X_{2n}$ the formula~(\ref{eq:Heat1.principal-term-Deltamu}) becomes
    \begin{equation}
        \Delta= - \frac{1}{2}(X_{1}^{2}+\ldots+X_{2n}^{2}) + i\mu X_{0}+ \text{lower order terms}.
         \label{eq:Heat1.principal-term-Deltamu2}
    \end{equation}
Combining this with~(\ref{eq:Spectral.Levi-form-Levi-metric2}) shows that the condition~(\ref{eq:Rockland.sublaplacian'}) for $\Delta$, 
is given in terms of the eigenvalues of $|L|=1$ and becomes 
\begin{equation}
    \mu \not \in  \{\pm(n+k);\ k=0,1,2,\ldots\}.
\end{equation}
Since by assumption we have $\mu\in (-n,n)$, the condition~(\ref{eq:Rockland.sublaplacian'}) is fulfilled at every point, 
so by Proposition~\ref{thm:Powers1.heat-sublaplacians} the principal symbol of $\Delta+\partial_{t}$ is an invertible 
Volterra-Heisenberg symbol. This allows us to apply Theorem~\ref{thm:Heat1.spectral-asymptotics} to deduce that as $\lambda \rightarrow \infty$ we have
\begin{equation}
    N(\Delta; \lambda) \sim \nu_{0}(\Delta) \lambda^{n+1},
     \label{eq:Spectral.Weyl-Delta-lemma}
\end{equation}
where $\nu_{0}(P)$ is given by~(\ref{eq:Heat1.counting-function-asymptotics}). 
    
Let us now work in Heisenberg coordinates centered at a point $a\in M$ related to a local $H$-frame $X_{0},X_{1},\ldots,X_{2n}$ as above. 
Because of~(\ref{eq:Heat1.almost-Heisenberg-relations}) the model vector fields of $X_{0},\ldots,X_{2n}$ at $a$ are
\begin{equation}
    X_{0}^{a}=\frac{\partial}{\partial x_{0}}, \quad X_{j}^{a}=\frac{\partial}{\partial x_{j}}+x_{n+j}\frac{\partial}{\partial 
    x_{0}}, \quad X_{n+j}^{a}=\frac{\partial}{\partial x_{n+j}}-x_{j}\frac{\partial}{\partial 
    x_{0}}, \quad 1\leq j\leq n.
\end{equation}
These are the vector fields~(\ref{eq:Heat1.Heisenberg-vector-fields}), 
so in the Heisenberg coordinates $G_{a}M$ agrees with the Heisenberg group $\bH^{2n+1}$. In addition, using~(\ref{eq:Heat1.principal-term-Deltamu2}) 
we see that the model operator of $\Delta$ is 
\begin{equation}
    \Delta^{a}= - \frac{1}{2}((X_{1}^{a})^{2}+\ldots+(X_{2n}^{a})^{2}) - i\mu X_{0}^{0}=\cL_{\mu}.
\end{equation}
Since $\mu\in (-n,n)$ it follows that $\Delta$ satisfies the Rockland at every point. 
On the other hand, since we are in Heisenberg coordinates we have $X_{j}(0)=\frac{\partial}{\partial x_{j}}$. In particular, we have 
$\theta(0)=dx_{0}$. Moreover, as $d\theta(X,Y)=-\theta([X,Y])$ we deduce from~(\ref{eq:Heat1.almost-Heisenberg-relations}) that 
\begin{equation}
    d\theta(X_{j},X_{n+k})=2\delta_{jk} \qquad \text{and} \qquad d\theta(X_{j},X_{k})=0,
\end{equation}
for $j=0,1,\ldots,n$ and $k=1,\ldots,n$. It then follows that $d\theta(0)=2\sum_{j=1}^{n} dx_{j}\wedge dx_{n+j}$. Therefore, the volume form is given 
by
\begin{equation}
  \frac{(-1)^{\kappa}}{n!2^{n}}\theta \wedge d\theta^{n}(0)  
  =(-1)^{\kappa}dx_{0}\wedge dx_{1}\wedge dx_{n+1}\wedge \ldots \wedge dx_{n}\wedge dx_{2n}.
     \label{eq:Spectral.volume-form-Heisenberg-coordinates}
\end{equation}
Note also that because of~(\ref{eq:Spectral.convention-Xj-orientation}) 
the orientation of $M$ is that given by $i^{n}X_{0}\wedge Z_{1}\wedge\bar{Z}_{1}\wedge \ldots \wedge 
Z_{n}\wedge\bar{Z}_{n}=(-1)^{\kappa}X_{0}\wedge X_{1}\wedge X_{n}\wedge \ldots X_{n+1}\wedge X_{2n}$. Therefore, at $x=0$ the volume form of $M$ 
is in the same orientation class as $(-1)^{\kappa}dx_{0}\wedge dx_{1}\wedge dx_{n+1}\wedge \ldots \wedge dx_{n}\wedge dx_{2n}$. In terms of density, 
together with~(\ref{eq:Spectral.volume-form-Heisenberg-coordinates}), this means that we have 
\begin{equation}
    \frac{(-1)^{\kappa}}{n!2^{n}}\theta \wedge d\theta^{n}(0)  =dx|_{x=0}. 
\end{equation}
All this shows that in the above Heisenberg coordinates the volume form and the model operator of $\Delta$ have normal forms in the sense of 
Proposition~\ref{prop:Spectral.normal-form-nu0P} with $\rho_{0}=1$ and $P_{0}=\cL_{\mu}\otimes 1_{\C^{r}}$. Therefore, 
we may apply~(\ref{eq:Spectral.normal-form-nu0P}) to get
\begin{equation}
    \nu_{0}(\Delta)=\nu_{0}(\cL_{\mu}\otimes 1_{\C^{r}})\int_{M}  \frac{(-1)^{\kappa}}{n!2^{n}}\theta \wedge d\theta^{n}=\nu_{0}(\cL_{\mu})\rk\cE 
    \op{vol}_{\theta}M,. 
\end{equation}
where  $\nu_{0}(\cL_{\mu})=\frac{1}{(n+1)!}k_{\mu}(0,0,1)=\nu(\mu)$. The lemma is thus proved.

\end{proof}
We are now ready to relate the asymptotics~(\ref{eq:Heat1.counting-function-asymptotics}) for the Kohn Laplacian to the pseudohermitian volume of $M$. 
\begin{theorem}
   Let $\Box_{b}:C^{\infty}(M,\Lambda^{*,*})\rightarrow C^{\infty}(M,\Lambda^{*,*})$ be the Kohn Laplacian on $M$ associated to a  
   Levi metric on $M$. Then for $p,q=1,\ldots,n$,  
   with  $q\neq \kappa$ and $q\neq n-\kappa$, as $\lambda \rightarrow \infty$ we have 
\begin{gather}
    N(\Box_{b|_{\Lambda^{p,q}}};\lambda) \sim \alpha_{n\kappa pq}(\op{vol}_{\theta}M)\lambda^{n+1},\\
     \alpha_{n\kappa pq}= \frac{1}{2^{n+1}} \binom{n}{p} \sum_{\max(0,q-\kappa)\leq  k\leq \min(q,n-\kappa)} \binom{n-\kappa}{k}\binom{\kappa}{q-k} 
    \nu(n-2(\kappa-q+2k)).
     \label{eq:Heat1.alphapq}  
\end{gather}    
\end{theorem}
\begin{proof}
As $M$ is $\kappa$-strictly pseudoconvex the $Y(q)$-condition reduces to $q\not\in \{\kappa, n-\kappa\}$, so in this case $\Box_{b}^{p,q}$ has an 
invertible principal symbol. Since $\Box_{b}$ is a positive it follows from Theorem~\ref{thm:Heat1.spectral-asymptotics} 
that as $\lambda \rightarrow \infty$ we have 
\begin{equation}
    N(\Box_{b}^{p,q};\lambda)  \sim \nu_{0}(\Box_{b}^{p,q}) \lambda^{n+1}.
     \label{eq:Heat1.Weyl-asymptotic-Boxb1}
\end{equation}
It then remains to show that $\nu_{0}(\Box_{b})=\alpha_{n\kappa pq}\op{vol}_{\theta}M$ with $\alpha_{n \kappa pq}$ given by~(\ref{eq:Heat1.alphapq}). 

Let $\{X_{0},Z_{j},\bar Z_{j}\}$ be a local admissible orthonormal frame of $T_{\C}M$ and let $\{\theta,\theta^{j}, \theta^{\bar j}\}$ be the dual 
coframe of $T^{*}_{\C}M$. For ordered subsets $J=\{j_{1},\ldots,j_{p}\}$ and $K=\{k_{1},\ldots,k_{q}\}$ of $\{1,\ldots,n\}$ with $j_{1}<\ldots<j_{p}$ 
and $k_{1}<\ldots<k_{q}$ we let 
\begin{equation}
    \theta^{J,\bar K}=\theta^{j_{1}}\wedge \ldots \wedge \theta^{j_{p}}\wedge \theta^{\bar k_{1}}\wedge \ldots \wedge \theta^{\bar k_{q}}.
\end{equation}
Then $\{\theta^{J,\bar K}\}$ form an orthonormal frame of $\Lambda^{*,*}$ and, as shown in~\cite[Sect.~20]{BG:CHM}, with respect to this frame 
on $(p,q)$-forms $\Box_{b}$ takes the form, 
\begin{gather}
    \Box_{b|_{\Lambda^{p,q}}}=\op{diag}\{ \Box_{J\bar K}\} +  \text{lower order terms},
    \label{eq:Heat1.Boxb-local-form1}\\
    \Box_{J\bar K}=-\frac{1}{2}\sum_{1\leq j \leq n} (Z_{j}\bar Z_{j}+\bar Z_{j} Z_{j}) +\frac{1}{2}\sum_{j \in K}[Z_{j},\bar Z_{j}] - 
    \frac{1}{2}\sum_{j \not \in K}[Z_{j},\bar Z_{j}].
      \label{eq:Heat1.Boxb-local-form2}
\end{gather}

Moreover, using~(\ref{eq:Heat1.commutators-admissible-frame}) 
we see that the leading part of $\Box_{J \bar K}$ is equal to
\begin{equation}
    -\frac{1}{2}\sum_{1\leq j \leq n} (Z_{j}\bar Z_{j}+\bar Z_{j} Z_{j}) -\frac{i}{2}\mu_{K}X_{0}, \qquad \mu_{K}=\sum_{j \in K} \epsilon_{j} 
    -\sum_{j\not\in K}\epsilon_{j}.
    \label{eq:Heat1.principal-term-BoxJK}
\end{equation}
Notice that since $\epsilon_{j}=1$ for $j=1,\ldots,n-\kappa$ and $\epsilon_{j}=-1$ for $j=n-\kappa+1,\ldots,n$ we have $\mu_{K}\in [-n,n]$ and 
$\mu_{K}=\pm n$ if, and only if, we have $K=\{1,\ldots,n-\kappa\}$ or $K=\{n-\kappa+1,\ldots,n\}$. Thus if $|K|=q$ with $q \not\in\{\kappa,n-\kappa\}$ 
then $\Box_{J \bar K}$ is two times a sublaplacian of the form~(\ref{eq:Heat1.principal-term-Deltamu}) with $\mu=\mu_{K}$ in $(-n,n)$. 

On the other hand, complex conjugation is an isometry, so from the orthogonal splitting~(\ref{eq:Heat1.positive-negative-splitting}) 
we get the orthogonal splitting $T_{0,1}=T_{0,1}^{+}\oplus T_{0,1}^{-}$ with $T_{0,1}^{\pm}=\overline{T_{1,0}^{\pm}}$. 
By duality these splittings give rise to the orthogonal decompositions,
\begin{gather}
    \Lambda^{1,0}=\Lambda^{1,0}_{+}\oplus \Lambda^{1,0}_{-}, \qquad \Lambda^{0,1}=\Lambda^{0,1}_{+}\oplus \Lambda^{0,1}_{-},\\
    \Lambda^{p,q}=  \bigoplus_{\max(0,q-\kappa)\leq  k\leq \min(q,n-\kappa)} \Lambda^{p;q,k}, \qquad 
    \Lambda^{p;q,k}=\Lambda^{p,0}\wedge (\Lambda^{0,1}_{+})^{k}\wedge (\Lambda^{0,1}_{-})^{q-k}.
      \label{eq:Heat1.splitting-pqk}
\end{gather}
Since~(\ref{eq:Heat1.Boxb-local-form1}) and~(\ref{eq:Heat1.Boxb-local-form2}) show that $\Box_{b}$ is scalar modulo lower order terms, on 
$\Lambda^{p,q}$ we can write 
\begin{equation}
    \Box_{b}= \sum_{\max(0,q-\kappa)\leq  k\leq \min(q,n-\kappa)} \Box_{p;qk} +  \text{lower order terms}, 
    \qquad \Box_{pqk}=\Pi_{pqk}\Box_{b} \Pi_{pqk},
\end{equation}
where $\Pi_{p;qk}$ denotes the orthogonal projection onto $\Lambda^{p;q,k}$. In particular,  
 as $\nu_{0}(\Box_{b}^{p,q})$ depends only on the principal 
symbol of $\Box_{b}^{p,q}$ we get
\begin{equation}
    \nu_{0}(\Box_{b|_{\Lambda^{p,q}}})=\sum_{\max(0,q-\kappa)\leq  k\leq \min(q,n-\kappa)}  \nu_{0}(\Box_{p;qk}).
     \label{eq:Heat1-splitting-alpha0}
\end{equation}
We are thus reduced to express each coefficient $\nu_{0}(\Box_{p;qk})$ in terms of $\op{vol}_{\theta}M$. 

Next, if $\theta^{J\bar K}$ is a section of $\Lambda^{p,q,k}$ then we have
\begin{gather}
    \# K\cap \{1,\ldots,n-\kappa\}=k, \qquad  \#K\cap \{n-\kappa+1,\ldots,n\}=q-k,\\
   \#K^{c}\cap \{1,\ldots,n-\kappa\}=n-\kappa-k, \qquad  \#K^{c}\cap \{n-\kappa+1,\ldots,n\}=\kappa-q+k,
\end{gather}
from which we get $ \mu_{K}= k-(q-k)-[(n-\kappa-k)-(\kappa-q+k)]=n+2q-2\kappa-4k.$
Combining this with~(\ref{eq:Heat1.Boxb-local-form1})--(\ref{eq:Heat1.principal-term-BoxJK}) then gives
\begin{equation}
    \Box_{p;qk}=-\frac{1}{2}\sum_{1\leq j \leq n} (Z_{j}\bar Z_{j}+\bar Z_{j} Z_{j}) -\frac{i}{2}(n+2q-2\kappa -4k)X_{0}+  \text{lower order terms}.
     \label{eq:Heat1.Boxpqk-local}
\end{equation}

We are now in position to apply Lemma~\ref{lem:Heat1.alpha0-volume} to $2 \Box_{p;qk}$ to get 
\begin{multline}
 \nu_{0}( \Box_{p;qk})= 2^{-(n+1)}\nu_{0}( 2\Box_{p;qk})=
 2^{-(n+1)}(\rk \Lambda^{p;q,k}) \nu(n+2q-2\kappa-4k).\op{vol}_{\theta}M \\
    = \frac{1}{2^{n+1}} \binom{n}{p} \binom{n-\kappa}{k}\binom{\kappa}{q-k}\nu(n+2q-2\kappa-4k).\op{vol}_{\theta}M.
\end{multline}
Combining this with~(\ref{eq:Heat1.Weyl-asymptotic-Boxb1}) and~(\ref{eq:Heat1-splitting-alpha0}) then shows that as $\lambda \rightarrow \infty$ we have 
\begin{equation}
    N(\Box_{b|_{\Lambda^{p,q}}};\lambda) \sim \alpha_{n\kappa pq}(\op{vol}_{\theta}M)\lambda^{n+1},  
\end{equation}    
with $\alpha_{n\kappa pq}$ given by~(\ref{eq:Heat1.alphapq}). The proof is thus achieved.
\end{proof}

We now turn to the horizontal sublaplacian on $(p,q)$-forms.

\begin{theorem}\label{thm:Heat1.Deltab}
  Let $\Delta_{b}:C^{\infty}(M,\Lambda^{*,*})\rightarrow C^{\infty}(M,\Lambda^{*,*})$ be the horizontal  sublaplacian associated to the   
   Levi metric on $M$. Then for $p,q=1,\ldots,n$, with $(p,q)\neq (\kappa,n-\kappa)$ and $(p,q)\neq (n-\kappa,\kappa)$, as $\lambda  \rightarrow \infty$ we have 
    \begin{gather}
         N(\Delta_{b|_{\Lambda^{p,q}}};\lambda) \sim \beta_{npq}(\op{vol}_{\theta}M) \lambda^{n+1},\\
    \beta_{n\kappa pq}= \!  \!  \!  \! \sum_{\substack{\max(0,q-\kappa)\leq  k\leq \min(q,n-\kappa)\\ \max(0,p-\kappa)\leq l\leq \min(p,n-\kappa)}}  \!  \!  \!  \!
    \binom{n-\kappa}{l}\binom{\kappa}{p-l} \binom{n-\kappa}{k}\binom{\kappa}{q-k} 
    \nu(2(q-p)+4(l-k)).
     \label{eq:Heat1.betapq}
    \end{gather}
\end{theorem}
\begin{proof}
 Thanks to~(\ref{eq:Operators.Tanaka-Kohn}) we know that we have 
 \begin{equation}
     \Delta_{b}=\Box_{b}+ \overline{\Box}_{b},
 \end{equation}
 where $\overline{\Box}_{b}$ denotes the conjugate operator of $\Box_{b}$, that is $\overline{\Box}_{b}\alpha=\overline{\Box_{b}\overline{\alpha}}$, 
 or equivalently the Laplacian of the $\partial_{b}$-complex. 
 
As in~(\ref{eq:Heat1.splitting-pqk}) we have orthogonal splittings, 
\begin{gather}
    \Lambda^{p,q}=    \bigoplus_{\max(0,p-\kappa)\leq  l\leq \min(p,n-\kappa)} \Lambda^{p,l;q}  
    = \bigoplus_{\substack{\max(0,q-\kappa)\leq  k\leq \min(q,n-\kappa)\\ \max(0,p-\kappa)\leq  l\leq \min(p,n-\kappa)}} \Lambda^{p,l;q,k},\\
    \Lambda^{p,l;q}= (\Lambda^{1,0}_{+})^{l}\wedge (\Lambda^{1,0}_{-})^{p-l} \wedge \Lambda^{0,q}, \qquad 
     \Lambda^{p,l;q,k}= \Lambda^{p,l;0}\wedge \Lambda^{p,l;0}.
\end{gather}
Since $\Box_{b}$ is a scalar operator modulo lower order terms, the same is true for $\overline{\Box}_{b}$ and $\Delta_{b}$. Therefore,
on $\Lambda^{p,q}$ we can write
\begin{gather}
    \overline{\Box}_{b}= \sum_{\max(0,p-\kappa)\leq  l\leq \min(p,n-\kappa)} \overline{\Box}_{p,l;q} +  \text{lower order terms}, \qquad  \overline{\Box}_{p,l;q}= 
    \Pi_{p,l;q} \overline{\Box}_{b}\Pi_{p,l;q},\\ 
    \Delta_{b}= \sum_{\substack{\max(0,q-\kappa)\leq  k\leq \min(q,n-\kappa)\\ \max(0,p-\kappa)\leq l\leq \min(p,n-\kappa)}} \Delta_{p,l;q,k} 
    +  \text{lower order terms},  \qquad    \Delta_{p,l;q,k}=  \Pi_{p,l;q,k}\Delta_{b}\Pi_{p,l;q,k},
\end{gather}
where $\Pi_{p,l;q}$ and $ \Pi_{p,l;q,k}$ denote the orthogonal projections onto $ \Lambda^{p,l;q}$ and $\Lambda^{p,l;q,k}$ respectively.
In particular, as in~(\ref{eq:Heat1-splitting-alpha02}) we have 
\begin{equation}
    \nu_{0}( \Delta_{b|_{\Lambda^{p,q}}})=  \sum_{\substack{\max(0,q-\kappa)\leq  k\leq \min(q,n-\kappa)\\ \max(0,p-\kappa)\leq l\leq \min(p,n-\kappa)}} 
    \nu_{0}(\Delta_{p,l;q,k}).
     \label{eq:Heat1-splitting-alpha02}
\end{equation}

Next, let $\{X_{0},Z_{j},\bar Z_{j}\}$ be a local admissible orthonormal frame for $T_{\C}M$. Since 
$\overline{\Box}_{p,l;q}=\overline{\Box_{q;p,l}}$, using~(\ref{eq:Heat1.Boxpqk-local}) we see that on $\Lambda^{p,l;q}$ we have 
\begin{equation}
    \overline{\Box}_{p,l;q}= - \frac{1}{2}\sum_{1\leq j \leq n} (Z_{j}\bar Z_{j}+\bar Z_{j} Z_{j}) +\frac{i}{2}(n+2p-2\kappa -4l)X_{0}+  
    \text{lower order terms}.
\end{equation}
Therefore, on $\Lambda^{p,l;q,k}$ we can write
\begin{multline}
      \Delta_{p,l;q,k}= \Box_{p;q,k}+\overline{\Box}_{p,l;q}\\ = 
      \sum_{1\leq j \leq n} (Z_{j}\bar Z_{j}+\bar Z_{j} Z_{j}) -i(2(q-p)+4(l-k))X_{0} +  \text{lower order terms}.
     \label{eq:Heat1.Deltaplqk}
\end{multline}

Thanks to~(\ref{eq:Heat1.Deltaplqk}) we can apply Lemma~\ref{lem:Heat1.alpha0-volume} to get
\begin{multline}
    \nu_{0}(\Delta_{p,l;q,k})= \rk \Lambda^{p,l;q,k}\nu(2(q-p)+4(l-k)) \op{vol}_{\theta}M\\ 
    = \binom{n-\kappa}{l}\binom{\kappa}{p-l} \binom{n-\kappa}{k}\binom{\kappa}{q-k}  \nu(2(q-p)+4(l-k)) \op{vol}_{\theta}M.
\end{multline}
Combining this with~(\ref{eq:Heat1-splitting-alpha02}) then shows that as $\lambda\rightarrow \infty$ we have 
   \begin{equation}
         N(\Delta_{b};\lambda) \sim \beta_{npq}(\op{vol}_{\theta}M) \lambda^{n+1}, 
    \end{equation}
with $\beta_{n\kappa pq}$ given by~(\ref{eq:Heat1.betapq}). The theorem is thus proved.
\end{proof}

Finally, we deal with with the conformal powers of the horizontal sublaplacian.
\begin{theorem}
Assume that $M$ is strictly pseudoconvex (i.e.~$\kappa=0$) and for $k=1,\ldots,n+1$ let $\boxdot_{\theta}^{(k)}:C^{\infty}(M)\rightarrow C^{\infty}(M)$ be a 
$k$'th conformal power of $\Delta_{b}$ acting on functions. Then  as $\lambda 
    \rightarrow \infty$ we have 
    \begin{equation}
         N(\boxdot_{\theta}^{(k)};\lambda) \sim  \nu(0)(\op{vol}_{\theta}M) \lambda^{\frac{n+1}{k}}. 
    \end{equation}

 \end{theorem}
 \begin{proof}
     Since $\nu_{0}(\boxdot_{\theta}^{(k)})$ depends only on the principal symbol of 
     $\boxdot_{\theta}^{(k)}$, which is the same as that of $\Delta_{b}^{k}$, we have
     $\nu_{0}(\boxdot_{\theta}^{(k)})=\nu_{0}(\Delta_{b}^{k})$. Therefore, using Theorem~\ref{thm:Heat1.Deltab} 
     we see that as $\lambda\rightarrow \infty$ we have 
    \begin{equation}
        N(\boxdot_{\theta}^{(k)};\lambda)\sim N(\Delta_{b}^{k};\lambda)=N(\Delta_{b};\lambda)^{\frac{1}{k}} \sim  \beta_{n000}(\op{vol}_{\theta}M) 
        \lambda^{\frac{n+1}{k}}.  
    \end{equation}
    Since $\beta_{n000}=2^{n} \nu(0)$ the result follows.
\end{proof}

\section{Weyl asymptotics and contact geometry}
\label{sec:Spectral-contact}
In this section we express in more geometric terms the Weyl asymptotics for the horizontal sublaplacian and the contact Laplacian on a compact 
orientable contact manifold $(M^{2n+1},\theta)$. 

We let $H=\ker \theta$ and let $X_{0}$ be the Reeb fields associated to $\theta$, so that $\imath_{X_{0}}d\theta=0$ and $\imath_{X_{0}}\theta=1$.  
Since $M$ is orientable $H$ admits a a calibrated almost complex structure $J\in C^{\infty}(M,\End H)$, $J^{2}=-1$, so that for any nonzero section $X$ of $H$ 
we have $d\theta(X,JX)=-d\theta(JX,X)>0$. 
We then endow $M$ with the orientation  defined by $\theta$ and the almost complex structure $J$, so that we have $\theta\wedge d\theta^{n}>0$, and with  
the Riemannian metric, 
\begin{equation}
    g_{\theta}=d\theta(.,J.)+\theta^{2}.
\end{equation}

A local orthonormal frame $X_{1},\ldots, X_{2n}$ of $H$ will be called admissible when we have $X_{n+j}=JX_{j}$ for $j=1,\ldots,n$. If 
$\theta^{1},\ldots,\theta^{2n}$ denotes the dual frame then we have $d\theta=\sum_{j=1}^{n}\theta^{j}\wedge \theta^{n+j}$, so the  
volume form of $g_{\theta}$ is equal to 
\begin{equation}
    \theta^{1}\wedge \theta^{n+1}\wedge \ldots\wedge \theta^{n}\wedge \theta^{2n}\wedge \theta=\frac{1}{n!}d\theta^{n}\wedge \theta.
\end{equation}
In particular, the volume form is independent of the choice of the almost complex structure and depends only on the contact form.

\begin{definition}
    The contact volume of $(M^{2n+1},\theta)$ is given by 
    \begin{equation}
        \op{vol}_{\theta}M=\frac{1}{n!}\int_{M}d\theta^{n}\wedge \theta.
    \end{equation}
\end{definition}

\begin{lemma}\label{lem:Heat1.alpha0-volume-contact}
  Let $\Delta:C^{\infty}(M,\cE)\rightarrow C^{\infty}(M,\cE)$ be selfadjoint sublaplacian  such that $\Delta$ is bounded from below. We further assume 
  that there exists $\mu 
  \in (-n,n)$ so that near any point of $M$ there is an admissible orthonormal frame $X_{1},\ldots,X_{2n}$ of $H$ with respect to which 
  $\Delta$ takes the form, 
    \begin{equation}
        \Delta= - (X_{1}^{2}+\ldots+X_{2n}^{2}) - i\mu X_{0}+ \text{lower order terms}.
         \label{eq:Heat1.principal-term-Deltamu-contact}
    \end{equation}
Then as $\lambda \rightarrow \infty$ we have
     \begin{equation}
        N(\Delta; \lambda) \sim 2^{-n}\nu(\mu) \rk \cE (\op{vol}_{\theta}M)\lambda^{n+1}. 
           \label{eq:Heat1.asymptotics-sublaplacian-mu-contact}
     \end{equation}

\end{lemma}
\begin{proof}
  Let $X_{1},\ldots,X_{2n}$ be a admissible local orthonormal frame and for $j=1,\ldots,2n$ let $\tilde{X}_{j}=\sqrt{2}X_{j}$. Then 
  $X_{0},\tilde{X}_{1},\ldots,\tilde{X}_{2n}$ is a local $H$-frame of $TM$ with respect to which $\Delta$ takes the form,
  \begin{equation}
      \Delta=-\frac{1}{2}(\tilde{X}_{1}^{2}+\ldots+\tilde{X}_{2n}^{2})-i\mu(x)X_{0}+\text{lower order terms}.
       \label{eq:Spectral.local-form-admissible}
  \end{equation}
  
 Moreover, for $j,k=1,\ldots,2n$ we have 
 \begin{equation}
     \theta([\tilde{X}_{j},J\tilde{X}_{k}])=-2d\theta(X_{j},JX_{k})=-2g_{\theta}(X_{j},X_{k})=-2\delta_{jk}.
 \end{equation}
  Therefore, for $j,k=1,\ldots,n$ we get:
\begin{gather}
    [X_{j},X_{n+k}]=-2\delta_{jk}X_{0} \ \bmod H,\\ 
    [X_{0},X_{j}]=[X_{j},X_{k}]=[X_{n+j},X_{n+k}]=0\ \bmod H.
     \label{eq:Heat1.almost-Heisenberg-relations-contact}
\end{gather}
The equalities~(\ref{eq:Spectral.local-form-admissible})--(\ref{eq:Heat1.almost-Heisenberg-relations-contact}) 
are the same as~(\ref{eq:Heat1.almost-Heisenberg-relations1})--(\ref{eq:Heat1.principal-term-Deltamu2}) in the case $\kappa=0$. Therefore, 
along the same lines as that of the proof  Lemma~\ref{lem:Heat1.alpha0-volume} we get
\begin{equation}
    \nu_{0}(\Delta)=\frac{\nu(\mu)}{n!2^{n}}\int_{M}d\theta^{n}\wedge \theta=2^{-n}\nu(\mu)\op{vol}_{\theta}M.
\end{equation}
Combining this with Theorem~\ref{thm:Heat1.spectral-asymptotics} then proves the asymptotics~(\ref{eq:Heat1.asymptotics-sublaplacian-mu-contact}).
\end{proof}

\begin{theorem}
    Let $\Delta_{b}:C^{\infty}(M,\Lambda^{*}_{\C}H^{*})\rightarrow C^{\infty}(M,\Lambda^{*}_{\C}H^{*})$ be the horizontal sublaplacian on 
    $M$ associated to the metric $g_{\theta}$ above and assume $k\neq n$. Then as $\lambda 
    \rightarrow \infty$ we have
    \begin{gather}
         N(\Delta_{b|_{\Lambda^{k}_{\C}H^{*}}};\lambda) \sim \gamma_{nk} (\op{vol}_{\theta}M) \lambda^{n+1},\qquad 
          \label{eq:Spectral.Weyl-Deltab-contact}
 \gamma_{nk}=2^{-n}\sum_{p+q=k}\binom{n}{p} \binom{n}{q}\nu(p-q).         
    \end{gather}
\end{theorem}
\begin{proof}
   As explained in Section~\ref{sec:Operators} the almost complex structure of $H$ gives rise to an orthogonal decomposition 
   $\Lambda^{k}_{\C}H^{*}=\bigoplus_{p+q=k}\Lambda^{p,q}$.   If $X_{1},\ldots,X_{2n}$ is a local admissible orthonormal frame of $H$ then, as shown by 
   Rumin~\cite[Prop.~2]{Ru:FDVC}, on $\Lambda^{p,q}$ the operator $\Delta_{b}$ takes the form,
   \begin{equation}
       \Delta_{b}=-(X_{1}^{2}+\ldots+X_{2n}^{2})+i(p-q)X_{0}+\text{lower order terms}.
        \label{eq:Spectral.contact-local-form-Deltapq}
   \end{equation}
   where the lower order part is not scalar. Therefore, modulo lower order terms, $\Delta_{b}$ preserves the bidegree. We thus may write
   \begin{equation}
       \Delta_{b|_{\Lambda^{k}_{\C}H^{*}}}=\sum_{p+q=k} \Delta_{p,q} +\text{lower order terms}, \qquad \Delta_{p,q}=\Pi_{p,q}\Delta_{b}\Pi_{p,q},
   \end{equation}
   where $\Pi_{p,q}$ denotes the orthogonal projection of $\Lambda^{*}_{\C}H^{*}$ onto $\Lambda^{p,q}$. In particular, we have 
   \begin{equation}
       \nu_{0}(\Delta_{b|_{\Lambda^{k}_{\C}H^{*}}})=\sum_{p+q=k}\nu_{0}(\Delta_{p,q}).
        \label{eq:Spectral.decomposition-nuDeltab-contact}
   \end{equation}
   
   Moreover, since $\Delta_{p,q}$ takes the form~(\ref{eq:Spectral.contact-local-form-Deltapq}) 
   with respect to  any admissible orthonormal frame of $H$, we may apply 
   Lemma~\ref{lem:Heat1.alpha0-volume-contact} to get
   \begin{equation}
       \nu_{0}(\Delta_{p,q})=2^{-n}\sum_{p+q=k}\binom{n}{p} \binom{n}{q}\nu(p-q). 
   \end{equation}
   Combining this with~(\ref{eq:Spectral.decomposition-nuDeltab-contact}) then gives the asymptotics~(\ref{eq:Spectral.Weyl-Deltab-contact}). 
\end{proof}
Finally, in the case of the contact Laplacian we can prove:

\begin{theorem}
    Let $\Delta_{R}:C^{\infty}(M,\Lambda^{*}\oplus \Lambda^{n}_{*})\rightarrow C^{\infty}(M,\Lambda^{*}\oplus \Lambda^{n}_{*})$ 
    be the contact Laplacian on $M$.\smallskip 
    
    1)  For $k=0,\ldots,2n$ with $k\neq n$ there exists a  universal constant $\nu_{nk}>0$ depending only on $n$ and $k$ 
    such that as $\lambda \rightarrow \infty$ we have
     \begin{equation}
           N(\Delta_{R|_{\Lambda^{k}}})\sim \nu_{nk} (\op{vol}_{\theta}M)\lambda^{n+1}. 
                \label{eq:Spectral.Weyl-contact-Laplacian1}
      \end{equation}
    
    2) For $j=1,2$ there exists a  universal constant $\nu_{n}^{(j)}>0$ depending only on $n$ and $j$  such that as $\lambda \rightarrow \infty$ we have 
    \begin{equation}
      N(\Delta_{R|_{\Lambda^{n}_{j}}})\sim \nu_{n}^{(j)} (\op{vol}_{\theta}M)\lambda^{\frac{n+1}{2}}. 
         \label{eq:Spectral.Weyl-contact-Laplacian2}
    \end{equation}
\end{theorem}
\begin{proof}
Let $a\in M$ and consider a chart around $a$ together with an admissible orthonormal frame $X_{1},\ldots,X_{2n}$ of $H$. Since 
$X_{0},X_{1},\ldots,X_{2n}$ form a $H$-frame this chart is a Heisenberg chart. Moreover, as shown in the proofs of Lemma~\ref{lem:Heat1.alpha0-volume} 
and Lemma~\ref{lem:Heat1.alpha0-volume-contact} the following hold:\smallskip

(i) We have $X_{j}^{a}$, where $X_{0}^{0},\ldots,X_{2n}^{0}$ are the left-invariant vector fields~(\ref{eq:Heat1.Heisenberg-vector-fields}) 
on $\bH^{2n+1}$. In particular, we have 
$G_{a}M=\bH^{2n+1}$ and $H_{a}=H^{0}_{0}$, where $H^{0}_{0}$ denotes the left-invariant Heisenberg hyperplane bundle of $\bH^{2n+1}$.\smallskip 

(ii) We have $\theta(0)=dx_{0}=\theta^{0}(0)$ and $d\theta(0)=2\sum_{j=1}^{n}dx_{j}\wedge dx_{n+j}=d\theta^{0}(0)$, where 
$\theta^{0}=dx_{0}+\sum_{j=1}^{n}(x_{j}dx_{n+j}-x_{n+j}dx_{j})$ is the standard left-invariant contact form of $\bH^{2n+1}$.\smallskip 

(iii) The density on $M$ given by the contact volume form $\frac{1}{n!}\theta\wedge d\theta^{n}$ agrees at $x=0$ with the density $dx$ on 
$\Rd$.\smallskip 

On the other hand, for $k=0,1,\ldots,2n$ the fiber at $a$ of the bundle $\Lambda^{k}_{*}$ depends only on $H_{a}$ and on the values of $\theta$ and 
$d\theta$ at $a$. Therefore, it follows from the statements (i) and (ii) that  in the Heisenberg coordinates centered at $a$  
the fibers at $x=0$ of the bundles $\Lambda^{*}\oplus \Lambda^{n}_{*}$ of $M$ and $\bH^{2n+1}$ agree. 

Next, let $\Delta_{R}^{0}$ denote the contact Laplacian on $\cH^{2n+1}$. Then we have:
\begin{lemma}
       In the Heisenberg coordinates centered at $a$ the model operator $\Delta_{R}^{a}$ agrees with $\Delta_{R}^{0}$.
\end{lemma}
\begin{proof}[Proof of the lemma]
First, note that in view of the formulas~(\ref{eq:Operators.equalities-DeltaR-Deltab.1})--(\ref{eq:Operators.equalities-DeltaR-Deltab.5}) 
for $\Delta_{R}$ and of Proposition~\ref{prop:PsiHDO.composition2} and Proposition~\ref{prop:PsiHDO.transpose-global}, we only have to show that in the 
Heisenberg coordinates centered at $a$ the   
model operators $d_{R}^{a}$ and $D^{a}_{R}$ agree with the operators $d_{R}^{0}$ and $D^{0}_{R}$ on $\bH^{2n+1}$.

Let $\theta^{1},\ldots,\theta^{2n}$ be the coframe of $H^{*}$ dual to $X_{1},\ldots,X_{2n}$. This coframe gives rise to a trivialization of 
$\Lambda^{*}_{\C}H^{*}$ in the chart with respect to which we have $d_{b}=\sum_{j=1}^{2n}\varepsilon(\theta^{j})X_{j}$. Furthermore, 
since $X_{j}(0)=\frac{d}{dx_{j}}$ we have $\theta_{j}(0)=dx_{j}$, so the model operator of $d_{b}$ is 
$d_{b}^{a}=\sum_{j=1}^{2n}\varepsilon(dx_{j})X_{j}^{0}=d_{b}^{0}$, where $d_{b}^{0}$ is the $d_{b}$-operator on $\bH^{2n+1}$.  
In particular, as for $k=n+1,\ldots,2n$ we have $d_{R}=d_{b}$ on $\Lambda^{k}$, we get $d_{R}^{a}=d_{b}^{a}=d_{b}^{0}=d_{R}^{0}$. 

On the other hand, as shown in~\cite{Ru:FDVC} for $k=0,\ldots,n-1$ we have $d_{R}^{*}=d_{b}^{*}$ on $\Lambda^{k}$. Thus 
$(d_{R}^{a})^{*}=(d_{R}^{*})^{a}=(d_{b}^{*})^{a}=(d_{b}^{a})^{*}=(d_{b}^{0})^{*}=(d_{R}^{0})^{*}$, 
which by taking adjoints gives $d_{R}^{a}=d_{R}^{0}$. 

Finally, it is proved in~\cite{Ru:FDVC} that $D_{R}=X_{0}+d_{b}\varepsilon(d\theta)^{*}d_{b}$ on $\Lambda^{n}_{1}$. Therefore, we get
$D^{a}_{R}=X_{0}^{a}+d_{b}^{a}(\varepsilon(d\theta)^{a})^{*}d_{b}^{a}= X_{0}^{0}+d_{b}^{0}\varepsilon(d\theta^{0})^{*}d_{b}^{0}=D^{0}_{R}$. 
The proof  that $\Delta_{R}^{a}=\Delta_{R}^{0}$ is thus complete. 
\end{proof}

Thanks to the statements (i) and (iii) and the claim above we may apply Proposition~\ref{prop:Spectral.normal-form-nu0P}. Letting $K_{0}(x,t)$ be 
the fundamental solution of $\Delta_{R}^{0}+\partial_{t}$ we then obtain:\smallskip 

    - For $k=0,\ldots,2n$, with $k\neq n$, as $\lambda \rightarrow \infty$ we have
     \begin{equation}
           N(\Delta_{R|_{\Lambda^{k}}})\sim \nu_{nk} (\op{vol}_{\theta}M)\lambda^{n+1}, \qquad  
           \nu_{nk}=\frac{2^{-n}}{(n+1)!}\tr_{\Lambda^{k,0}}K_{0|_{\Lambda^{k,0}}}(0,1).
      \end{equation}
 
      - For $j=1,2$ as $\lambda \rightarrow \infty$, we have 
    \begin{equation}
      N(\Delta_{R|_{\Lambda^{n}_{j}}})\sim \nu_{n,j} (\op{vol}_{\theta}M)\lambda^{\frac{n+1}{2}}, \qquad
     \nu_{n}^{(j)}=2^{-n}\Gamma(1+\frac{n+1}{2})^{-1}\tr_{\Lambda^{n,0}_{j}}K_{0|_{\Lambda^{n,0}_{j}}}(0,1).
    \end{equation}
In particular, the constant $\nu_{nk}$ (resp.~$\nu_{n}^{(j)}$) depends only on $n$ and $k$ (resp.~$n$ and $j$), hence is a universal constant. Furthermore, 
Theorem~\ref{thm:Heat1.spectral-asymptotics} implies that $\nu_{nk}$ and $\nu_{n}^{(j)}$ are positive numbers.
\end{proof}

\end{document}